# KOBAYASHI-HITCHIN CORRISPONDENCE FOR TWISTED VECTOR BUNDLES

ARVID PEREGO

ABSTRACT. We prove the Kobayashi-Hitchin correspondence and the approximate Kobayashi-Hitchin correspondence for twisted holomorphic vector bundles on compact Kähler manifolds. More precisely, if $X$ is a compact manifold and $g$ is a Gauduchon metric on $X$, a twisted holomorphic vector bundle on $X$ is $g-$polystable if and only if it is $g-$Hermite-Einstein, and if $X$ is a compact Kähler manifold and $g$ is a Kähler metric on $X$, then a twisted holomorphic vector bundle on $X$ is $g-$semistable if and only if it is approximate $g-$Hermite-Einstein.

## Contents

1. Introduction

The Kobayashi-Hitchin correspondence for vector bundles is a nowadays well-established result in complex geometry, saying that a holomorphic vector bundle on a compact complex manifold $X$ is polystable if and only if it admits a Hermite-Einstein metric. Here a holomorphic vector bundle is polystable if it is the direct sum of stable holomorphic vector bundles (where stability is the slope-stability, or Mumford-Takemoto stability) with the same slope, and a Hermite-Einstein metric is a Hermitian metric whose mean curvature is a constant multiple of the identity.

This result was proved in an increasing order of generalization by several authors. First, in 1980 Kobayashi introduced in [15] the notion of Hermite-Einstein metric on a holomorphic vector bundle over a complex manifold. In [16] he showed that an irreducible Hermite-Einstein vector bundle on a compact Kähler manifold is polystable with respect to a Kähler metric. A different proof of this was given by Lübke in [20].

Shortly after [16], Donaldson proved in [6] that on a Riemann surface even the opposite is true, i. e. that a stable holomorphic vector bundle carries a Hermite-Einstein metric. This gave a new proof of the Narasimhan-Seshadri theorem (see [23]) saying that a holomorphic vector bundle on a Riemann surface $X$ is stable if and only if it has an irreducible projective unitary representation of the fundamental group of $X$.

Donaldson's result motivated Kobayashi and Hitchin, indipendently, to conjecture that this result holds for every holomorphic vector bundle on a compact Kähler manifold: it is this correspondence which is usually referred to as *Kobayashi-Hitchin correspondence*. The first proof of this correspondence for higher dimensional manifolds was given by Donaldson in [7] for algebraic surfaces, and then in [8] for algebraic manifolds.

Uhlenbeck and Yau proved in [27] and [28] that the Kobayashi-Hitchin correspondence holds on arbitrary compact Kähler manifolds, completing the proof of the original conjecture of Kobayashi and Hitchin. A few years after that, Buchdahl proved in [1] that the Kobayashi-Hitchin correspondence holds on any compact complex surface, and in [18] Li and Yau proved that it holds on every compact complex manifold, i. e. that if $X$ is a compact complex manifold and $g$ is a Gauduchon metric on $X$, then a holomorphic vector bundle $E$ is $g-$Hermite-Einstein if and only if it is $g-$polystable.

Instead of Hermite-Einstein holomorphic vector bundle, i. e. a holomorphic vector bundle admitting a Hermite-Einstein metric, one can consider the weaker notion of approximate Hermite-Einstein holomorphic vector bundle, i. e. a holomorphic vector bundle $E$ such that for each $\epsilon > 0$ there is a Hermitian metric $h_\epsilon$ on $E$ whose mean curvature $K_g(E, h_\epsilon)$ verifies

$$\max_{x \in X} |Tr((K_g(E, h_\epsilon) - c \cdot id_E)^2)| < \epsilon,$$

where $c \in \mathbb{R}$ depends only on $X$, $g$, $c_1(E)$ and the rank of $E$.

It was shown by Kobayashi in [17] that an approximate Hermite-Einstein holomorphic vector bundle on a compact Kähler manifold $X$ is semistable (with respect to a Kähler metric $g$), and that if $X$ is a projective manifold, then even the converse holds. This equivalence is often referred to as *approximate Kobayashi-Hitchin correspondence*. In [13] Jacob proved that the approximate Kobayashi-Hitchin correspondence holds for every holomorphic vector bundle on a compact Kähler manifold.

A natural generalization of holomorphic vector bundles on a compact complex manifold $X$ is given by twisted holomorphic vector bundles, where the twist is a $2-$cocyle representing an element in the Brauer group $Br(X)$ of $X$ (i. e. the torsion of $H^2(X, \mathscr{O}_X^*)$). Twisted sheaves were introduced by Giraud in [10], and can be defined in several equivalent ways: as family of sheaves on an open covering of $X$ together with a twisted gluing, as sheaves of modules over an Azumaya algebra on $X$ (see [3]), as sheaves on a gerb on $X$ (see [10], [11], [5]), as sheaves on a $\mathscr{O}_X^*-$gerbe (see [19]) or as sheaves on a projective bundle over $X$ (see [30]).

Stability for coherent twisted sheaves was introduced first by Lieblich in [19] in the language of sheaves on $\mathscr{O}_X^*-$gerbes, and by Yoshioka in [30] in the language of sheaves on a projective bundle over $X$. In [24] stability of coherent twisted sheaves is discussed in the language of twisted gluing of coherent sheaves and in the language of modules over an Azumaya algebra. In all of these categories in order to define stability one needs a definition of Chern classes of coherent twisted sheaves.

The notion of connection on a twisted holomorphic vector bundles appears in [11], [5] (see even [22]), and was used in [29] to prove the Kobayashi-Hitchin correspondence for twisted holomorphic vector bundles on a compact Kähler manifold. All the definitions there are in the category of twisted holomorphic vector bundles as holomorphic bundles over a gerb.

In the present paper we will consider twisted vector bundles following Căldăraru's point of view, i. e. local vector bundles on a open covering together with a twisted gluing. Connections on such vector bundles may be found in [14], and in the present paper we present a definition of Hermitian metric. The twist $\alpha$ will be a 2−cocycle (given once an open covering $\mathscr{U} = \{U_i\}_{i \in I}$ of $X$ is fixed) whose cohomology class lies in $Br(X)$, and which is associated a $B$−field, i. e. a family of closed $(1,1)$−forms $B_i \in A^{1,1}(U_i)$ such that $B_i - B_j$ is an exact form.

The aim of the paper is to prove the following:

**Theorem 1.1.** *Let $X$ be a compact Kähler manifold with a Kähler metric $g$, $\alpha$ the twist associated to a $B$−field and $E$ an $\alpha$−twisted holomorphic vector bundle on $X$.*

(1) *$E$ is $g$−polystable if and only if it is $g$−Hermite-Einstein.*
(2) *$E$ is $g$−semistable if and only if it is approximate $g$−Hermite-Einstein .*

As already mentioned, the first item of the statement of Theorem 1.1 was already proved by Wang in [29]. Here we present a proof in the language of local vector bundles with twisted gluing, and we provide a generalization of Wang's result to compact complex manifolds with a Gauduchon metric $g$ (see Theorem 5.17). The second item is the twisted version of the approximate Kobayashi-Hitchin correspondence, and is new. The proof is an adaptation to twisted holomorphic vector bundles of the original proofs of Kobayashi and Jacob.

The structure of the paper is the following: section 2 contains the main basic facts about twisted vector bundles and sheaves, and about connections and metrics on them. As this was not available in the literature, we gave a global overview of this.

We prove in particular that every twisted holomorphic vector bundle $E$ over which we fix a Hermitian metric $h$, carries a unique connection which is compatible with the holomorphic structure of $E$ and the metric $h$. In analogy with the untwisted case, we will call this connection the *Chern connection* of the pair $(E, h)$.

The $B$−field fixed at the beginning will allow us to define the curvature of every connection, which is a global 2−form with values in the (untwisted) vector bundle $\underline{End}(E)$. The curvature of the Chern connection will be called *Chern curvature*, and it will be a global $(1,1)$−form with values in $\underline{End}(E)$. It is the fact that $\underline{End}(E)$ is a true (i. e. untwisted) holomorphic vector bundle on $X$ that will allow us to provide a proof of Theorem 1.1.

We will moreover discuss how connections, curvatures and Hermitian metrics behave under various operations of twisted bundles, like direct sum, tensor product, dual, pull-back, sub-bundle and quotient bundles, proving in particular a twisted version of the Gauss-Codazzi formulas.

Section 3 is devoted to introduce the notion of $g$−Hermite-Einstein and approximate $g$−Hermite-Einstein twisted bundles. To do so, we need to introduce Chern forms and Chern classes for twisted sheaves: as already done by Wang in [29], we define Chern forms and Chern classes by means of the curvature of a connection, similarly to what happens for holomorphic vector bundles.

Once the Chern forms and classes are introduced, we define the mean curvature of a pair $(E, h)$ of a twisted holomorphic vector bundle $E$ and a Hermitian metric $h$: exactly as in the untwisted case, this will be a smooth endomorphism of $E$ (which is Hermitian with respect to $h$). The notion of (weak) $g-$Hermite-Einstein metric is then as in the untwisted case, and we prove that in the conformal class of a weak $g-$Hermite-Einstein metric there is always a $g-$Hermite-Einstein metric.

We define $g-$Hermite-Einstein and approximate $g-$Hermite-Einstein vector bundles as in the untwisted case, and we will provide several properties of (approximate) $g-$Hermite-Einstein bundles following closely the analogous properties for untwisted bundles.

Section 4 is devoted to the notion of $g-$semistable and $g-$stable twisted holomorphic vector bundles, proving several properties of these bundles, and in particular that $g-$Hermite-Einstein bundles are $g-$polystable, and that approximate $g-$Hermite-Einstein bundles are $g-$semistable. The proof is essentially the same as in the untwisted case, and we follow closely Lübke's argument in [20]. This proves half of Theorem 1.1.

In section 5 we prove that a $g-$stable twisted holomorphic vector bundle is $g-$Hermite-Einstein: this is the content of Theorem 5.1, which completes the proof of point 1 of Theorem 1.1. The proof we present is identical to the one given by Uhlenbeck and Yau in [27], and its adaptation to Gauduchon metrics on compact complex manifold as presented in [18] and in section 3 of [21].

Wang's approach in [29] was to adapt to twisted bundles the original argument of Donaldson, adapted by Simpson in [25]. Since [18] and [21] work more generally if $g$ is a Gauduchon metric on a compact complex manifold $X$, we will finally prove that the Kobayashi-Hitchin correspondence for twisted holomorphic vector bundles holds on every compact complex manifold (with respect to a chosen Gauduchon metric on it), thus generalizing [29].

The remaining part of the paper is devoted to the proof of the approximate Kobayashi-Hitchin correspondence, namely that each $g-$semistable twisted holomorphic vector bundle is approximate $g-$Hermite-Einstein: this is the content of Theorem 6.1, which completes the proof of point 2 of Theorem 1.1.

The proof will follow closely the original argument in the untwisted case as presented in [17] and [13]. As in [17] we first define the Donaldson Lagrangian for Hermitian metrics on a twisted holomorphic vector bundle, and prove that if it is bounded below, then the bundle is approximate $g-$Hermite-Einstein. We will then prove that the Donaldson Lagrangian of a $g-$semistable twisted holomorphic vector bundle is bounded below. The proof is based on the one proposed by Jacob in [13] for untwisted vector bundles, and we adapt it to the twisted case.

Notation

All along the paper we consider a $C^\infty-$differentiable manifold $M$. If $J$ is a complex structure on $M$, we will let $X = (M, J)$ be the induced complex manifold. We moreover fix a sufficiently fine open covering $\mathscr{U} = \{U_i\}_{i \in I}$ of $M$, where $I$ is a set of indexes. We write $U_{ij} := U_i \cap U_j$ and $U_{ijk} := U_i \cap U_j \cap U_k$.

Once the open covering $\mathscr{U}$ is fixed, we choose a $B-$field $B$ on $X$ with respect to $\mathscr{U}$: for further reference and more details on $B-$fields, see [11], [5], [29] and [4].

A $B-$**field** on $X$ with respect to $\mathscr{U}$ is a family $B = \{B_i\}_{i \in I}$ where $B_i$ is a 2$-$form on $U_i$, such that there are 1$-$forms $\omega_{ij}$ on $U_{ij}$ with the property that
$$B_i - B_j = d\omega_{ij}.$$
We notice that
$$dB_i - dB_j = d^2\omega_{ij} = 0,$$
hence the 3$-$forms $dB_i$ glue together to give a closed 3$-$form $dB$ on $X$, whose cohomology class is an element in $H^3(X, \mathbb{Z})$.

Notice that
$$d(\omega_{ij} + \omega_{jk} + \omega_{ki}) = 0,$$
so if the covering $\mathscr{U}$ is sufficiently fine we may find $U(1)-$valued functions $\alpha_{ijk}$ on $U_{ijk}$ such that
$$\omega_{ij} + \omega_{jk} + \omega_{ki} = -\alpha_{ijk}^{-1} d\alpha_{ijk}.$$
Then $\alpha_B = \{\alpha_{ijk}\}$ is a 2$-$cocycle whose cohomology class lies in $H^2(X, \mathscr{O}_X^*)$.

The 2$-$cocyle $\alpha_B$ will be called **twist induced by** $B$, and all along the paper we will use the notation $\alpha$ instead of $\alpha_B$. The natural morphism $H^2(X, \mathscr{O}_X^*) \longrightarrow H^2(X, \mathbb{Z})$ from the exponential sequence sends the cohomology class $[\alpha]$ to the cohomology class $[dB]$ of $dB$.

If $[dB]$ is torsion in $H^3(X, \mathbb{Z})$, then $[\alpha]$ is torsion in $H^2(X, \mathscr{O}_X^*)$, i. e. it corresponds to an element in the Brauer group $Br(X)$ of $X$. In this case $[dB]$ is trivial in $H^3(X, \mathbb{R})$, and we may and will choose the forms $B_i$ to be $d-$closed for every $i \in I$. We will moreover ask that $B_i$ is a purely imaginary $(1,1)-$form on $U_i$, and that $\omega_{ij}$ is a $(1,0)-$form on $U_{ij}$.

If $H^3(X, \mathbb{Z})$ is free, then every element in $Br(X)$ may be represented by a twist induced by some $B-$field.

## 2. Connections and metrics

In this section we introduce the definitions of twisted vector bundle and of twisted coherent sheaf that we will use all along the paper. After having reviewed all the basic operations we will use on twisted sheaves, we will introduce the notion of connection, of curvature and of Hermitian metric on a twisted vector bundle, showing that once a holomorphic twisted vector bundle and a Hermitian metric on it are given, then there is a unique connection on it which is compatible with the metric and with the holomorphic structure. In analogy to the untwisted case, this connection will be called the *Chern connection* of the twisted bundle, whose curvature will be the most important tool in the paper, exactly as in the untwisted case.

2.1. **Twisted vector bundles.** Let $M$ be a $C^\infty-$differentiable manifold over which we have a complex structure, and let $X$ be the induced complex manifold. We first recall the definition of twisted vector bundle on $X$.

**Definition 2.1.** *An $\alpha-$**twisted complex $C^\infty$ vector bundle** on $X$ is a family $E = \{E_i, \phi_{ij}\}_{i,j \in I}$ where*
  (1) *for each $i \in I$, $E_i$ is a complex $C^\infty$ vector bundle on $U_i$,*
  (2) *for each $i, j \in I$, $\phi_{ij} : E_{i|U_{ij}} \longrightarrow E_{j|U_{ij}}$ is an isomorphism of complex $C^\infty$ vector bundles on $U_{ij}$,*
  (3) *we have $\phi_{ii} = id_{E_i}$, $\phi_{ij}^{-1} = \phi_{ji}$ and $\phi_{ki} \circ \phi_{jk} \circ \phi_{ij} = \alpha_{ijk} \cdot id_{E_{i|U_{ijk}}}$ for every $i, j, k \in I$.*

Morphisms of twisted bundles are defined in a natural way:

**Definition 2.2.** *Let $E = \{E_i, \phi_{ij}\}$ and $F = \{F_i, \psi_{ij}\}$ be two $\alpha-$twisted complex $C^\infty$ vector bundles on $X$. A **morphism of $\alpha-$twisted complex $C^\infty$ vector bundles** $f : E \longrightarrow F$ is a family $f = \{f_i\}_{i \in I}$ where*
  (1) *for each $i \in I$, we have that $f_i : E_i \longrightarrow F_i$ is a morphism of complex $C^\infty$ vector bundles on $U_i$,*
  (2) *for each $i, j \in I$, we have $\psi_{ij} \circ f_i = f_j \circ \phi_{ij}$.*

The category of $\alpha-$twisted complex $C^\infty$ vector bundles on $X$ will be denoted $Bun_{C^\infty}(X, \alpha)$. The objects that will be under investigation in this paper will anyway more precisely be twisted *holomorphic* vector bundles, defined as follows:

**Definition 2.3.** *An $\alpha-$twisted complex $C^\infty$ vector bundle $E = \{E_i, \phi_{ij}\}_{i,j \in I}$ on $X$ will be called $\alpha-$**twisted holomorphic vector bundle** on $X$ if $E_i$ is a holomorphic vector bundle on $U_i$ and $\phi_{ij} : E_{i|U_{ij}} \longrightarrow E_{j|U_{ij}}$ is an isomorphism of holomorphic vector bundles.*

Morphisms among twisted holomorphic vector bundles are then defined as follows:

**Definition 2.4.** *If $E$ and $F$ are two $\alpha-$twisted holomorphic vector bundles on $X$, a **morphism of $\alpha-$twisted holomorphic vector bundles** from $E$ to $F$ is a morphism $f = \{f_i\}_{i \in I} : E \longrightarrow F$ of $\alpha-$twisted complex $C^\infty$ bundles such that for every $i \in I$ we have that $f_i : E_i \longrightarrow E_j$ is a morphism of holomorphic vector bundles on $U_i$.*

The category of $\alpha-$twisted holomorphic vector bundles on $X$ will be denoted $Bun(X, \alpha)$. If instead of looking at vector bundles we are willing to look at sheaves, we will talk about twisted sheaves as follows.

**Definition 2.5.** *An $\alpha-$**twisted sheaf** on $X$ is a family $\mathscr{E} = \{\mathscr{E}_i, \phi_{ij}\}_{i,j \in I}$ where*

(1) *for each $i \in I$, $\mathscr{E}_i$ is sheaf of Abelian groups on $U_i$,*
(2) *for each $i, j \in I$, $\phi_{ij} : \mathscr{E}_{i|U_{ij}} \longrightarrow \mathscr{E}_{j|U_{ij}}$ is an isomorphism of sheaves of Abelian groups on $U_{ij}$,*
(3) *we have $\phi_{ii} = id_{E_i}$, $\phi_{ij}^{-1} = \phi_{ji}$ and $\phi_{ki} \circ \phi_{ji} \circ \phi_{ij} = \alpha_{ijk} \cdot id_{E_{i|U_{ijk}}}$ for every $i, j, k \in I$.*

*If $\mathscr{F}_i$ is a sheaf of $\mathcal{O}_{U_i}-$modules for every $i \in I$, then we say that $\mathscr{F}$ is an $\alpha-$**twisted sheaf of** $\mathcal{O}_X-$**modules**. If moreover $\mathscr{F}_i$ is coherent (resp. quasi-coherent), we will say that $\mathscr{F}$ is an $\alpha-$**twisted coherent sheaf** (resp. an $\alpha-$**twisted quasi-coherent sheaf**).*

Moreover, we have the notion of morphism between $\alpha-$twisted sheaves.

**Definition 2.6.** *If $\mathscr{E} = \{\mathscr{E}_i, \phi_{ij}\}_{i,j \in I}$ and $\mathscr{F} = \{\mathscr{F}_i, \psi_{ij}\}_{i,j \in I}$ are two $\alpha-$twisted sheaves (of $\mathcal{O}_X-$modules), a **morphism of $\alpha-$twisted sheaves (of $\mathcal{O}_X-$modules)** from $\mathscr{E}$ to $\mathscr{F}$ is a family $f = \{f_i\}_{i \in I}$ where*

(1) *for each $i \in I$, $f_i : \mathscr{E}_i \longrightarrow \mathscr{F}_i$ is morphism of sheaves (of $\mathcal{O}_{U_i}-$modules),*
(2) *we have $\psi_{ij} \circ f_i = f_j \circ \phi_{ij}$ for every $i, j \in I$.*

We therefore have the categories $Sh(X, \alpha)$ (resp. $Sh_{\mathcal{O}_X}(X, \alpha)$) of $\alpha-$twisted sheaves (resp. of $\alpha-$twisted sheaves of $\mathcal{O}_X-$modules). The full subcategory of $Sh_{\mathcal{O}_X}(X, \alpha)$ whose objects are $\alpha-$twisted (quasi-)coherent sheaves are denoted $Coh(X, \alpha)$ (resp. $QCoh(X, \alpha)$). Exactly as for untwisted sheaves, we have the notion of locally free twisted sheaf.

**Definition 2.7.** *An $\alpha-$twisted sheaf $\mathscr{E} = \{\mathscr{E}_i, \phi_{ij}\}$ of $\mathcal{O}_X-$modules is said to be **locally free** (of rank $r$) if for each $i \in I$ we have that $\mathscr{E}_i$ is locally free (of rank $r$).*

The full-subcategory of $Sh_{\mathcal{O}_X}(X, \alpha)$ whose objects are $\alpha-$twisted locally free sheaves is denoted $Lf(X, \alpha)$. It is easy to prove that there is an equivalence of categories between $Bun(X, \alpha)$ and $Lf(X, \alpha)$.

**Remark 2.8.** If $\mathscr{U}'$ is a refinement of $\mathscr{U}$, by restriction we see that a $B-$field relative to $\mathscr{U}$ gives a $B-$field relative to $\mathscr{U}'$, whose associated twist is a Čech 2$-$cocycle $\alpha'$ relative to $\mathscr{U}'$. We moreover get a canonical equivalence between $Bun_{C^\infty}(X, \alpha)$ and $Bun_{C^\infty}(X, \alpha')$ (and similarly for the other categories we mentioned before).

If $E = \{E_i, \phi_{ij}\}$ is an $\alpha-$twisted complex $C^\infty$ vector bundle relative to $\mathscr{U}$, we may refine $\mathscr{U}$ so that the twisted vector bundle corresponding to $E$ will be $\{E'_i, \phi'_{ij}\}_{i \in I'}$ where $E'_i$ is the trivial vector bundle.

2.2. **Operations with twisted bundles.** The usual operations between vector bundles ($C^\infty$ or holomorphic) and sheaves can be defined as well in the twisted setting. We will only consider the case of $\alpha-$twisted complex $C^\infty$ vector bundles, but the same definitions work for $\alpha-$twisted holomorpic vector bundles and for $\alpha-$twisted sheaves (of $\mathcal{O}_X-$modules, coherent or quasi-coherent). We refer the reader to [3] for further details.

**Dual bundle**. Let $E = \{E_i, \phi_{ij}\}$ be an $\alpha-$twisted complex $C^\infty$ vector bundle. The **dual** of $E$ is $E^* = \{E_i^*, \phi_{ij}^*\}$ where $E_i^*$ is the dual vector bundle of $E_i$ on $U_i$, and $\phi_{ij}^* : E_i^* \longrightarrow E_j^*$ is the dual of $\phi_{ij}$: more precisely, if $\eta$ is a local section of $E_i^*$, then $\phi_{ij}^*(\eta)$ is the local section of $E_j^*$ mapping a local section $\xi$ of $E_j$ to $\eta(\phi_{ij}^{-1}(\xi))$. It is easy to see that $E^*$ is a complex $C^\infty$ vector bundle twisted by $\alpha^{-1} = \{\alpha_{ijk}^{-1}\}$.

**Conjugate bundle**. Let $E = \{E_i, \phi_{ij}\}$ be an $\alpha-$twisted complex $C^\infty$ vector bundle. The **conjugate** of $E$ is $\overline{E} = \{\overline{E}_i, \overline{\phi}_{ij}\}$ where $\overline{E}_i$ is the conjugate vector bundle of $E_i$ on $U_i$, and $\overline{\phi}_{ij} : \overline{E}_i^* \longrightarrow \overline{E}_j^*$ is the conjugate of $\overline{\phi}_{ij}$. Then $\overline{E}$ is a complex $C^\infty$ vector bundle twisted by $\overline{\alpha} = \{\overline{\alpha_{ijk}}\}$. If $E$ is holomorphic, then $\overline{E}$ is holomorphic over $\overline{X}$ (the complex manifold obtained by putting on $M$ the conjugate complex structure $\overline{J}$).

**Direct sum**. If $E = \{E_i, \phi_{ij}\}$ and $F = \{F_i, \psi_{ij}\}$ are two $\alpha-$twisted complex $C^\infty$ vector bundles, their **direct sum** is $E \oplus F := \{E_i \oplus F_i, \phi_{ij} \oplus \psi_{ij}\}$, which is an $\alpha-$twisted complex $C^\infty$ vector bundle as well.

**Tensor product**. Consider $B$ and $B'$ two $B-$fields with respect to $\mathscr{U}$, and let $\alpha$ and $\alpha'$ be the respective twists. If $E = \{E_i, \phi_{ij}\}$ is an $\alpha-$twisted complex $C^\infty$ vector bundle and $F = \{F_i, \psi_{ij}\}$ an $\alpha'-$twisted complex $C^\infty$ vector bundle, their **tensor product** is $E \otimes F := \{E_i \otimes F_i, \phi_{ij} \otimes \psi_{ij}\}$, which is an $\alpha\alpha'-$twisted complex $C^\infty$ vector bundle.

In particular $E \otimes E^*$ and $\overline{E}^* \otimes E^*$ are untwisted complex $C^\infty$ vector bundles: indeed $E^*$ is $\alpha^{-1}-$twisted, $\overline{E}^*$ is $\overline{\alpha}^{-1}$ twisted, hence $E \otimes E^*$ is twisted by $\alpha\alpha^{-1} = 1$, and $\overline{E}^* \otimes E^*$ is twisted by $\overline{\alpha}^{-1}\alpha^{-1} = |\alpha|^{-2} = 1$ (since $\alpha_{ijk}$ is a function taking values in $U(1)$).

**Wedge product**. If $E$ is an $\alpha-$twisted complex $C^\infty$ vector bundle on $X$, for every $p \geq 0$ we may consider the $p-$**th wedge product** $\bigwedge^p E = \{\bigwedge^p E_i, \wedge^p \phi_{ij}\}$. This is a direct summand of $E^{\otimes p}$, hence it is an $\alpha^p-$twisted complex $C^\infty$ vector bundle. In particular, if $E$ has rank $r$, then we have that $\det(E) := \bigwedge^r E$ is an $\alpha^r-$twisted complex $C^\infty$ vector bundle.

**Bundle of morphisms**. If $E = \{E_i, \phi_{ij}\}$ are is an $\alpha-$twisted complex $C^\infty$ vector bundle and $F = \{F_i, \psi_{ij}\}$ is an $\alpha'-$twisted complex $C^\infty$ vector bundle on $M$, we let $\underline{Hom}(E, F) := \{\underline{Hom}(E_i, F_i), \Phi_{ij}\}$, where

$$\Phi_{ij} : \underline{Hom}(E_i, F_i) \longrightarrow \underline{Hom}(E_j, F_j), \quad \Phi_{ij}(f) := \psi_{ij} \circ f \circ \phi_{ij}^{-1}.$$

It is easy to show that $\underline{Hom}(E, F)$ is an $\alpha^{-1}\alpha'-$twisted complex $C^\infty$ vector bundle.

The canonical isomorphism $E_i^* \otimes F_i \simeq \underline{Hom}(E_i, F_i)$ (coming from the universal property of tensor product) induces a canonical isomorphism between $E^* \otimes F$ and $\underline{Hom}(E, F)$. In particular we see that $\underline{End}(E)$ and $\underline{Hom}(\overline{E}, E^*)$ are untwisted complex $C^\infty$ vector bundles.

Notice that if $E$ and $F$ are $\alpha-$twisted holomorphic vector bundles, then $\underline{End}(E)$ and $\underline{Hom}(E, F)$ are holomorphic vector bundles. As such we may consider their global sections as complex $C^\infty$ vector bundles, or as holomorphic vector bundles.

In the first case, the global sections of $\underline{End}(E)$ (resp. of $\underline{Hom}(E, F)$) are the smooth endomorphisms of $E$ (resp. the smooth morphisms from $E$ to $F$), and will be denoted $A^0(End(E))$ or simply $End(E)$ (resp. $A^0(Hom(E, F))$, $Hom(E, F)$). In the second case, we will use the notation $H^0(End(E))$ and $H^0(Hom(E, F))$.

**Pull-back**. Let $X$ and $Y$ be two $C^\infty$ differentiable manifolds and $f : X \longrightarrow Y$ be a smooth map between them. If $\mathscr{U} = \{U_i\}_{i \in I}$ is an open covering of $Y$, then $f^*\mathscr{U} := \{f^{-1}(U_i)\}_{i \in I}$ is an open covering of $X$. If $E = \{E_i, \phi_{ij}\}$ is an $\alpha-$twisted complex $C^\infty$ vector bundle on $Y$, the **pull-back of $E$ under $f$** is $f^*E = \{f^*E_i, f^*\phi_{ij}\}$, which is a $f^*\alpha-$twisted complex $C^\infty$ vector bundle on $X$.

2.2.1. *Properties of morphisms.* If $f : E \longrightarrow E$ is an endomorphism of an $\alpha-$twisted complex $C^\infty$ (resp. holomorphic) vector bundle, then for every $i, j \in I$ we have
$$Tr(f_i) = Tr(\phi_{ij}^{-1} f_j \phi_{ij}) = Tr(f_j),$$
so that we may glue together the traces of the endomorphisms $f_i$'s to get a global smooth (resp. holomorphic) function $Tr(f)$, called **trace of $f$**.

Similarily we have
$$\det(f_i) = \det(\phi_{ij}^{-1} f_j \phi_{ij}) = \det(f_j),$$
so we get a global smooth (resp. holomorphic) function $\det(f)$, called **determinant of $f$**.

If $\mathscr{F} = \{\mathscr{F}_i, \phi_{ij}\}_{i,j \in I}$ and $\mathscr{G} = \{\mathscr{G}_i, \psi_{ij}\}_{i,j \in I}$ are two $\alpha-$twisted coherent sheaves on $X$ and $f = \{f_i\}_{i \in I} : \mathscr{F} \longrightarrow \mathscr{G}$ is a morphism, then for every $x \in U_i$ we have that $f_{i,x} : \mathscr{F}_{i,x} \longrightarrow \mathscr{G}_{i,x}$ is a morphism of $\mathscr{O}_{X,x}-$modules.

If $x \in U_{ij}$, as $\phi_{ij}$ and $\psi_{ij}$ are isomorphisms of vector bundles we have
$$rk_x(f_i) = rk_x(\psi_{ij} \circ f_i) = rk_x(f_j \circ \phi_{ij}) = rk_x(f_j).$$

It follows that $rk_x(f_i)$ does not depend on the choice of $i \in I$: we will write it $rk_x(f)$ and call it the **rank of $f$ at $x$**.

We now need to make some remarks about eigenvalues of endomorphisms of twisted bundles.

**Remark 2.9.** If $E$ is an $\alpha-$twisted holomorphic vector bundle and $f = \{f_i\}$ is a smooth endomorphism of $E$, then it makes sense to consider the eigenvalues of $f$ (which are smooth functions on $X$).

Indeed, suppose that $\lambda_i$ is an eigenvalue of $f_i$, i. e. $\lambda_i$ is a smooth function on $U_i$ for which there is a nowhere vanishing smooth section $s$ of $E_i$ with $f_i(s) = \lambda_i s$. Then $\lambda_i$ is an eigenvalue for $f_j$ over $U_{ij}$: indeed $\phi_{ij}(s)$ is a nowhere vanishing smooth section of $E_j$ over $U_{ij}$, and we have

$$f_j(\phi_{ij}(s)) = \phi_{ij}(f_i(s)) = \phi_{ij}(\lambda_i s) = \lambda_i \phi_{ij}(s).$$

Hence the eigenvalues of the $f_i$'s glue together to give global smooth functions on $X$ that on each $U_i$ rectrict to the eigenvalues of $f_i$, and that will be referred to as eigenvalues of $f$.

As for morphisms of untwisted sheaves, the trace of a morphism of twisted sheaves is the sum of the eigenvalues, and its determinant is their product. The previous Remark 2.9 shows moreover that $f_i$ is diagonalizable if and only if $f_j$ is, hence it makes sense to talk about diagonalizable endomorphisms of $\alpha$-twisted (holomorphic) vector bundles.

If $f$ is a diagonalizable endomorphism of $E$ whose eigenvalues are $\lambda_1, \cdots, \lambda_r$, consider a smooth function $\varphi : \mathbb{R} \longrightarrow \mathbb{R}$ and suppose that the images of $\lambda_1, \cdots, \lambda_r$ are all contained in the definition domain of $\varphi$. In particular we see that $\varphi \circ \lambda_i$ is a smooth function on $X$.

This allows us to perform the following general construction: for every $i \in I$ consider a local frame $\sigma_i$ of $E_i$ which diagonalizes $f_i$. With respect to $\sigma_i$ we then have that $f_i$ is represented by a diagonal matrix $F_i$ whose diagonal entries are the eigenvalues of $f_i$ (each one appearing with its respective multiplicity). We then let $\varphi(F_i)$ be the diagonal matrix whose diagonal entries are the $\varphi \circ \lambda_i$ (each one with the respective multiplicity), and consider the endomorphism $\varphi(f_i)$ of $E_i$ corresponding to $\varphi(F_i)$.

**Lemma 2.10.** *The family $\varphi(f) := \{\varphi(f_i)\}_{i \in I}$ is a diagonalizable endomorphism of $E$.*

*Proof.* The endomorphism $\varphi(f_i)$ is diagonalizable, so $\varphi(f)$ is diagonalizable. Moreover, if $\sigma_i$ is a local frame of $E_i$ diagonalizing $\varphi(f_i)$, then $\phi_{ij}(\sigma_i)$ is a local frame of $E_j$ diagonalizing $\varphi(f_j)$. It then follows that $\phi_{ij} \circ \varphi(f_i) = \varphi(f_j) \circ \phi_{ij}$, and we are done. $\square$

Particular cases are $\exp(f)$, the **exponential** of $f$ (which may be defined for every endomorphism of $E$), $\log(f)$, the **logarithm** of $f$ (which may be defined for positive definite endomorphisms) and $f^\sigma$ for every $\sigma \in (0, 1]$ (which may be defined for positive semidefinite endomorphisms).

2.3. **Connections and curvatures.** We now define connections and curvatures on twisted vector bundles. Before doing this, we recall some very basic facts about connections on vector bundles: we refer the reader to [17] for further details. If $V$ is a complex $C^\infty$ vector bundle on $X$ of rank $r$, we use the notation $A^p(V)$ for the space of $p$-forms on $X$ with values in $V$, and $A^p(X)$ for the space of $p$-forms on $X$.

A connection on $V$ is a $\mathbb{C}$−linear map $D : A^0(V) \longrightarrow A^1(V)$ such that for every $f \in A^0(X)$ and every $s \in A^0(E)$ we have

$$D(s \cdot f) = D(s) \cdot f + s \cdot df.$$

If $s = \{s_1, \cdots, s_r\}$ is a local frame of $V$, then the connection form of $D$ relative to $s$ is a matrix $\Gamma$ of 1−forms on $X$ such that

$$D(s_i) = \sum_{j=1}^r s_j \cdot \gamma_i^j,$$

where $\Gamma = [\gamma_i^j]$. If $\xi \in A^0(X)$, write

$$\xi = \sum_{j=1}^r \xi^j s_j,$$

so that

$$D(\xi) = \sum_{j=1}^r s_j \left( d\xi^j + \sum_{k=1}^r \gamma_k^j \xi^k \right),$$

that we write simply $D(\xi) = d\xi + \Gamma\xi$.

It is easy to see that to give a connection on $V$ is equivalent to give an open covering $\mathscr{U} = \{U_i\}_{i \in I}$ of $X$ such that $V$ is trivialized over $U_i$, a local frame $s_i$ of $V$ over $U_i$ and a $r \times r$−matrix $\Gamma_i$ of 1−forms on $U_i$ such that if $\phi_{ij} : V_i \longrightarrow V_j$ is the transition function, then we have

$$\Gamma_i = a_{ij}^{-1} \Gamma_j a_{ij} + a_{ij}^{-1} da_{ij},$$

where $a_{ij}$ is the $r \times r$−matrix of smooth functions on $V$ representing $\phi_{ij}$ with respect to $s_i$ and $s_j$.

We introduce now the notion of connection on a twisted bundle (see [11], [5] and [29] for connections on gerbs, and [14] for connections on twisted vector bundles). Let $E = \{E_i, \phi_{ij}\}_{i,j \in I}$ be an $\alpha$−twisted complex $C^\infty$ vector bundle of rank $r$.

**Definition 2.11.** *A **connection** on $E$ is a family $D = \{D_i\}_{i \in I}$ where*
  (1) *for each $i \in I$, $D_i$ is a connection on $E_i$,*
  (2) *for every $i, j \in I$ if $\Gamma_i$ is a connection form of $D_i$ with respect to a local frame of $E_i$, and if $a_{ij}$ is the matrix of smooth functions representing $\phi_{ij}$ with respect to the chosen local frames, we have*

$$\Gamma_i = a_{ij}^{-1} \Gamma_j a_{ij} + a_{ij}^{-1} da_{ij} + \omega_{ij} \cdot I_r.$$

**Remark 2.12.** The motivation of the previous definition comes from the following remark: if $V$ is a complex $C^\infty$ vector bundle and $D$ is a connection on it, take a family of connections form $\Gamma_i$ associated to local frames on an open covering of $X$, and let $a_{ij}$ be the matrix representing the transition function with respect to the given local frames. We then have

$$\Gamma_i = a_{ij}^{-1} \Gamma_j a_{ij} + a_{ij}^{-1} da_{ij} =$$

$$= a_{ij}^{-1}(a_{jk}^{-1}\Gamma_k a_{jk} + a_{jk}^{-1}da_{jk})a_{ij} + a_{ij}^{-1}da_{ij} =$$
$$= a_{ij}^{-1}a_{jk}^{-1}\Gamma_k a_{jk}a_{ij} + a_{ij}^{-1}a_{jk}^{-1}da_{jk}a_{ij} + a_{ij}^{-1}da_{ij} =$$
$$= a_{ij}^{-1}a_{jk}^{-1}(a_{ki}^{-1}\Gamma_i a_{ki} + a_{ki}^{-1}da_{ki})a_{jk}a_{ij} + a_{ij}^{-1}a_{jk}^{-1}da_{jk}a_{ij} + a_{ij}^{-1}da_{ij} =$$
$$= (a_{ki}a_{jk}a_{ij})^{-1}\Gamma_i(a_{ki}a_{jk}a_{ki}) + a_{ki}^{-1}da_{ki})a_{jk}a_{ij} + a_{ij}^{-1}a_{jk}^{-1}da_{jk}a_{ij} + a_{ij}^{-1}da_{ij} =$$
$$= (a_{ki}a_{jk}a_{ij})^{-1}\Gamma_i(a_{ki}a_{jk}a_{ki}) + (a_{ki}a_{jk}a_{ij})^{-1}d(a_{ki}a_{jk}a_{ki}).$$

If $a_{ki}a_{jk}a_{ij} = I_r$, this last line is $\Gamma_i$. But in the twisted case we have that $a_{ki}a_{jk}a_{ij} = \alpha_{ijk} \cdot I_r$, hence we get $\Gamma_i = \Gamma_i + \alpha_{ijk}^{-1}d\alpha_{ijk}$. In order to avoid this discrepancy we need to add $\omega_{ij} \cdot I_r$ in the relation between $\Gamma_i$ and $\Gamma_j$.

The existence of a connection on any $\alpha$-twisted complex $C^\infty$ vector bundle on $M$ is granted by Example 7.2 of [14]. We present here a more general construction that will be used in what follows.

**Proposition 2.13.** *Let $E$ be an $\alpha$-twisted complex $C^\infty$ vector bundle on $M$. Then $E$ admits a connection.*

*Proof.* Write $E = \{E_i, \phi_{ij}\}$, and let $p = \{p_i\}_{i \in I}$ be a partition of the unity with respect to $\mathscr{U}$. Choose a connection $D_i$ on $E_i$, and let $\Gamma_i$ be the connection form of $D_i$ with respect to a chosen local frame of $E_i$. We write $a_{ij}$ for the matrix of smooth functions representing $\phi_{ij}$ with respect to these local frames.

We consider
$$\Gamma'_i := \sum_{j \in I} p_j a_{ij}^{-1}\Gamma_j a_{ij}, \quad \Phi_i := \sum_{j \in I} p_j a_{ij}^{-1}da_{ij},$$

which are two matrices of 1-forms on $U_i$. Notice that
$$a_{ij}^{-1}\Gamma'_j a_{ij} = a_{ij}^{-1}\left(\sum_{k \in I} p_k a_{jk}^{-1}\Gamma_k a_{jk}\right)a_{ij} = \sum_{k \in I} p_k a_{ij}^{-1}a_{jk}^{-1}\Gamma_k a_{jk}a_{ij}.$$

Now recall that $a_{jk}a_{ij} = \alpha_{ijk}a_{ik}$, so
$$a_{ij}^{-1}\Gamma'_j a_{ij} = \sum_{k \in I} p_k(\alpha_{ijk}a_{ik})^{-1}\Gamma_k(\alpha_{ijk}a_{ik}) = \Gamma'_i.$$

Now, recall that $a_{ki}a_{jk}a_{ij} = \alpha_{ijk} \cdot I_r$, so that $a_{ik} = \alpha_{ijk}^{-1}a_{jk}a_{ij}$. It follows that
$$a_{ik}^{-1}da_{ik} = (\alpha_{ijk}^{-1}a_{jk}a_{ij})^{-1}d(\alpha_{ijk}^{-1}a_{jk}a_{ij}) =$$
$$= \alpha_{ijk}a_{ij}^{-1}a_{jk}^{-1}((d\alpha_{ijk}^{-1})a_{jk}a_{ij} + \alpha_{ijk}^{-1}(da_{jk})a_{ij} + \alpha_{ijk}^{-1}a_{jk}(da_{ij})) =$$
$$= \alpha_{ijk}d\alpha_{ijk}^{-1} \cdot I_r + a_{ij}^{-1}(a_{jk}^{-1}da_{jk})a_{ij} + a_{ij}^{-1}da_{ij} =$$
$$= -\alpha_{ijk}^{-1}d\alpha_{ijk} \cdot I_r + a_{ij}^{-1}(a_{jk}^{-1}da_{jk})a_{ij} + a_{ij}^{-1}da_{ij}.$$

But then
$$\Phi_i = \sum_{k \in I} p_k a_{ik}^{-1}da_{ik} = a_{ij}^{-1}\left(\sum_{k \in I} p_k a_{jk}^{-1}da_{jk}\right)a_{ij} + a_{ij}^{-1}da_{ij} - \sum_{k \in I} p_k \alpha_{ijk}^{-1}d\alpha_{ijk} =$$

$$= a_{ij}^{-1} \Phi_j a_{ij} + a_{ij}^{-1} da_{ij} + \omega_{ij} \cdot I_r.$$

We now let $\widetilde{\Gamma}_i := \Gamma'_i + \Phi_i$ for every $i \in I$, so that

$$\widetilde{\Gamma}_i = \Gamma'_i + \Phi_i = a_{ij}^{-1} \Gamma'_j a_{ij} + a_{ij}^{-1} \Phi_j a_{ij} + a_{ij}^{-1} da_{ij} + \omega_{ij} \cdot I_r =$$
$$= a_{ij}^{-1} \widetilde{\Gamma}_j a_{ij} + a_{ij}^{-1} da_{ij} + \omega_{ij} \cdot I_r.$$

Consider now the family $\widetilde{D} = \{\widetilde{D}_i\}_{i \in I}$ where $\widetilde{D}_i$ is the connection whose connection form is $\widetilde{\Gamma}_i$ with respect to the given local frame: we then see that $\widetilde{D}$ is a connection on $E$. $\square$

**Remark 2.14.** The set of connections on an $\alpha$-twisted complex $C^\infty$ vector bundle $E$ is an affine space over the vector space $A^1(\underline{End}(E))$. Indeed, if $D = \{D_i\}$ and $D' = \{D'_i\}$ are two connections on $E$, we have that $D_i - D'_i \in A^1(\underline{End}(E_i))$ and

$$D_i - D'_i = \phi_{ij}^{-1}(D_j - D'_j)\phi_{ij},$$

so the $D_i - D'_i$'s glue together to a global $D - D' \in A^1(\underline{End}(E))$. As a consequence of this, any affine linear combination of connections on an $\alpha$-twisted complex $C^\infty$ vector bundle $E$ is again a connection on $E$.

Using $D_i : A^0(E_i) \longrightarrow A^1(E_i)$ we define a $\mathbb{C}$-linear map

$$D_i^1 : A^1(E_i) \longrightarrow A^2(E_i), \quad D_i^1(s \cdot \varphi) := D_i(s) \wedge \varphi + s \cdot d\varphi,$$

and more generally we define $D_i^p : A^p(E_i) \longrightarrow A^{p+1}(E_i)$ using the previous formula in a recursive way (letting $D_i^0 := D_i$).

If now $E$ is an $\alpha$-twisted complex $C^\infty$ vector bundle on $X$ and $D = \{D_i\}$ is a connection on $E$, for every $i \in I$ we consider $R_i := D_i^1 \circ D_i$, which is the curvature of the connection $D_i$ on $E_i$, so we know that $R_i \in A^2(\underline{End}(E_i))$, i. e. $R_i$ is a 2-form with values in the complex $C^\infty$ vector bundle $\underline{End}(E_i)$ of endomorphisms of $E_i$ (over $U_i$). Let

$$\widetilde{R}_i := R_i - B_i \cdot id_{E_i} \in A^2(\underline{End}(E_i)).$$

**Lemma 2.15.** *There is a unique $R_D \in A^2(\underline{End}(E))$ such that $R_{D|U_i} = \widetilde{R}_i$ for every $i \in I$.*

*Proof.* We have $\widetilde{R}_i = R_i - B_i \cdot id_{E_i}$. With respect to a local frame of $E_i$ we represent $R_i$ by a matrix $\Omega_i$ of 2-forms. For each $i \in I$ recall that $\Omega_i = d\Gamma_i + \Gamma_i \wedge \Gamma_i$ (see as instance section 1 in Chapter I of [17]). But then

$$R_i = \phi_{ij}^{-1} R_j \phi_{ij} + d\omega_{ij} \cdot id_{E_i},$$

we then get

$$\widetilde{R}_i = \phi_{ij}^{-1} R_j \phi_{ij} + (d\omega_{ij} - B_i) \cdot id_{E_i}.$$

But as $d\omega_{ij} - B_i = -B_j$, we get

$$\widetilde{R}_i = \phi_{ij}^{-1} R_j \phi_{ij} - B_j \cdot id_{E_i} = \phi_{ij}^{-1} R_j \phi_{ij} - \phi_{ij}^{-1}(B_j \cdot id_{E_j})\phi_{ij} = \phi_{ij}^{-1} \widetilde{R}_{ij} \phi_{ij},$$

so we glue together the $\widetilde{R}_i$'s to produce the 2-form $R_D$. $\square$

**Definition 2.16.** *The 2-form $R_D \in A^2(\underline{End}(E))$ is the **curvature of** $D$.*

2.4. **Connections and holomorphic structures.** Let us now fix a holomorphic structure, getting a complex manifold $X$. Let $E = \{E_i, \phi_{ij}\}$ be an $\alpha$-twisted complex $C^\infty$ vector bundle on $X$ and $D = \{D_i\}$ a connection on $E$. The holomorphic structure on $X$ gives a direct sum decomposition
$$A^0(E_i) = A^{1,0}(E_i) \oplus A^{0,1}(E_i).$$
Composing $D_i$ with the two projections we get
$$D_i^{1,0} : A^0(E_i) \longrightarrow A^{1,0}(E_i), \quad D_i^{0,1} : A^0(E_i) \longrightarrow A^{0,1}(E_i)$$
and we clearly have $D_i = D_i^{1,0} + D_i^{0,1}$. We will let
$$D^{1,0} := \{D_i^{1,0}\}, \quad D^{0,1} := \{D_i^{0,1}\}.$$

Similarly, we have
$$A^2(\underline{End}(E)) = A^{2,0}(\underline{End}(E)) \oplus A^{1,1}(\underline{End}(E)) \oplus A^{0,2}(\underline{End}(E)),$$
hence if $R_D$ is the connection of $D$, we have three components $R_D^{2,0}$, $R_D^{1,1}$ and $R_D^{0,2}$ such that
$$R_D = R_D^{2,0} + R_D^{1,1} + R_D^{0,2}.$$

If $\Gamma_i$ is a connection form for $D_i$ with respect to a given frame, we have a natural decomposition $\Gamma_i = \Gamma_i^{1,0} + \Gamma_i^{0,1}$ since $\Gamma_i$ is a matrix of 1-forms. It follows from the definition of connection and the fact that $\omega_{i,j}$ is a $(1,0)$-form that
$$\Gamma_i^{1,0} = a_{ij}^{-1} \Gamma_j^{1,0} a_{ij} + a_{ij}^{-1} \partial a_{ij} + \omega_{ij}^{1,0} \cdot id_{E_i},$$
and
$$\Gamma_i^{0,1} = a_{ij}^{-1} \Gamma_j^{0,1} a_{ij} + a_{ij}^{-1} \overline{\partial} a_{ij}.$$

Suppose now furthermore that $E$ is an $\alpha$-twisted holomorphic vector bundle, i. e. $E_i$ has a holomorphic structure and $\phi_{ij}$ is an isomorphism of holomorphic vector bundles. We then represent the connection $D_i$ by a matrix $\Gamma_i$ of 1-forms with respect to a holomorphic local frame, and $\phi_{ij}$ by a matrix $a_{ij}$ whose entries are holomorphic functions. In this case we get
$$\Gamma_i^{1,0} = a_{ij}^{-1} \Gamma_j^{1,0} a_{ij} + a_{ij}^{-1} \partial a_{ij} + \omega_{ij}^{1,0} \cdot id_{E_i}, \quad \Gamma_i^{0,1} = a_{ij}^{-1} \Gamma_j^{0,1} a_{ij}.$$

The holomorphic structure of $E_i$ corresponds to a semi-connection
$$\overline{\partial}_i : A^0(E_i) \longrightarrow A^{0,1}(E_i),$$
i. e. such that $\overline{\partial}(f \cdot s) = \overline{\partial}(f) \cdot s + f \cdot \overline{\partial}_i(s)$ for every $f \in A^0(U_i)$ and $s \in A^0(E_i)$.

**Definition 2.17.** *A connection $D = \{D_i\}_{i \in I}$ on $E$ is* **compatible with the holomorphic structure of** *$E$ if for every $i \in I$ we have $D_i^{0,1} = \overline{\partial}_i$.*

This is equivalent to asking that $\Gamma_i^{0,1} = 0$ for every $i$, or even that for every holomorphic section $\xi$ of $E_i$ we have $D_i(\xi) = D_i^{1,0}(\xi)$ (see Proposition 3.9 in Chapter I of [17]).

The following shows that each twisted holomorphic vector bundle carries a connection compatible with its holomorphic structure.

**Lemma 2.18.** *Let $E$ be an $\alpha-$twisted holomorphic vector bundle on $X$. Then $E$ admits a connection $D$ compatible with its holomorphic structure, and if the $B-$field $B$ is such that $B_i^{0,2} = 0$ for every $i \in I$, then $R_D^{0,2} = 0$.*

*Proof.* Write $E = \{E_i, \phi_{ij}\}$, where $E_i$ is a holomorphic vector bundle and $\phi_{ij}$ is holomorphic for every $i, j \in I$. We know that $E_i$ admits a connection $D_i$ compatible with its holomorphic structure (see Proposition 3.5 in Chapter I of [17]). We let $\Gamma_i$ be its connection form with respect to a holomorphic local frame of $E_i$. Consider moreover a partition of the unity $p = \{p_i\}$ relative to $\mathscr{U}$.

The proof of Proposition 2.13 tells us that if we let $\widetilde{D} = \{\widetilde{D}_i\}$ whose connection form, with respect to the given holomorphic local frames, is

$$\widetilde{\Gamma}_i = \sum_{j \in I} p_j (a_{ij}^{-1} \Gamma_j a_{ij} + a_{ij}^{-1} da_{ij}),$$

then $\widetilde{D}$ is a connection on $E$. Notice that

$$\widetilde{\Gamma}_i^{0,1} = \sum_{j \in I} p_j a_{ij}^{-1} \Gamma_j^{0,1} a_{ij} + \sum_{j \in I} p_j a_{ij}^{-1} \overline{\partial} a_{ij} = 0,$$

since $\Gamma_j^{0,1} = 0$ by the fact that $D$ is compatible with the holomorphic structure of $E$, and $\overline{\partial} a_{ij} = 0$ since $a_{ij}$ is a matrix of holomorphic functions. It follows that $\widetilde{D}$ is compatible with the holomorphic structure of $E$.

To conclude, notice that $R_{\widetilde{D}}^{0,2} = 0$ if and only if its restriction to $U_i$ is 0 for every $i \in I$, i. e. if and only if $(R_i - B_i \cdot id_{E_i})^{0,2} = 0$. But since $B_i^{0,2} = 0$, we get

$$R_{\widetilde{D}|U_i}^{0,2} = (R_i - B_i \cdot id_{E_i})^{0,2} = R_i^{0,2} - B_i^{0,2} = R_i^{0,2}.$$

Now, recall that as $D_i^{0,1} = \overline{\partial}_i$ we have $R_i^{0,2} = 0$ (see Proposition 3.5 in Chapter I of [17]), and we are done. $\square$

A converse of the previous Lemma holds too.

**Lemma 2.19.** *Let $E$ be an $\alpha-$twisted complex $C^\infty$ vector bundle on $X$ and $D$ a connection on $E$. Suppose that the $B-$field $B = \{B_i\}$ is such that $B_i^{0,2} = 0$ for every $i \in I$, and that $R_D^{0,2} = 0$. Then there is a unique holomorphic structure on $E$ with which $D$ is compatible.*

*Proof.* As in the proof of Lemma 2.18, the fact that $R_D^{0,2} = 0$ and that $B_i^{0,2} = 0$ imply that $R_i^{0,2} = 0$ for every $i \in I$. Proposition 3.7 in Chapter I of [17] then implies the existence of a unique holomorphic structure on $E_i$ with which $D_i$ is compatible.

We now need to prove that $\phi_{ij}$ is holomorphic with respect to the holomorphic structures of $E_i$ and $E_j$. Let $\Gamma_i$ be the connection form of $D_i$ with respect to a holomorphic local frame of $E_i$, and $a_{ij}$ the matrix of smooth functions representing $\phi_{ij}$ with respect to the chosen local frames of $E_i$ and $E_j$. Then $\Gamma_i^{0,1} = 0$, and since
$$\Gamma_i = a_{ij}^{-1}\Gamma_j a_{ij} + a_{ij}^{-1}da_{ij} + \omega_{ij} \cdot id_{E_i},$$
by multiplying by $a_{ij}$ on both sides we then get
$$da_{ij} = a_{ij}\Gamma_i - \Gamma_j a_{ij} - \omega_{ij} a_{ij}.$$
But as $\omega_{ij}$ is a $(1,0)$−form, we see that $da_{ij}$ is a matrix of $(1,0)$−forms: hence $a_{ij}$ is a matrix of holomorphic functions, and $\phi_{ij}$ is holomorphic. □

2.5. **Hermitian metrics and connections.** We now introduce the notion of Hermitian metric on a twisted bundle. We recall that if $V$ is a complex $C^\infty$ vector bundle on $X$, a Hermitian metric on $V$ is a $C^\infty$ field of positive definite Hermitian products on the fibers of $V$.

**Definition 2.20.** *Let $E = \{E_i, \phi_{ij}\}$ be an $\alpha$−twisted complex $C^\infty$ vector bundle on $X$. A **Hermitian metric** on $E$ is a collection $h = \{h_i\}_{i \in I}$ where*
  (1) *for every $i \in I$, $h_i$ is a Hermitian metric on $E_i$,*
  (2) *for every $i, j \in I$ we have $h_i = {}^T\phi_{ij}h_j\overline{\phi}_{ij}$, i. e. for every sections $\xi$ and $\eta$ of $E_i$ we have*
$$h_i(\xi, \eta) = h_j(\phi_{ij}(\xi), \phi_{ij}(\eta)).$$

**Remark 2.21.** As $\alpha_{ijk}$ is $U(1)$−values, we have
$$h_i = {}^T\phi_{ij}h_j\overline{\phi}_{ij} = {}^T\phi_{ij}{}^T\phi_{jk}h_k\overline{\phi}_{jk}\overline{\phi}_{ij} = {}^T\phi_{ij}{}^T\phi_{jk}{}^T\phi_{ki}h_i\overline{\phi}_{ki}\overline{\phi}_{jk}\overline{\phi}_{ij} =$$
$$= {}^T(\phi_{ki}\phi_{jk}\phi_{ij})h_i\overline{\phi_{ki}\phi_{jk}\phi_{ij}} = |\alpha_{ijk}|^2 h_i = h_i.$$
It follows that there is no discrepancy on $U_{ijk}$, and the definition makes sense.

**Remark 2.22.** If $H_i$ is the matrix of smooth functions representing $h_i$ with respect to a given local frame of $E_i$, and $a_{ij}$ is the matrix of smooth functions representing $\phi_{ij}$ with respect to the chosen local frames of $E_i$ and $E_j$, then $H_i$ is a Hermitian matrix and $H_i = {}^Ta_{ij}H_j\overline{a}_{ij}$.

We first show that Hermitian metrics exist on every $\alpha$−twisted $C^\infty$ vector bundle:

**Lemma 2.23.** *Let $E$ be an $\alpha$−twisted $C^\infty$ vector bundle. Then $E$ admits a Hermitian metric.*

*Proof.* Let $h_i$ be a Hermitian metric on $E_i$, and $p = \{p_i\}_{i \in I}$ a partition of the unity with respect to $\mathscr{U}$. Let
$$\widetilde{h}_i := \sum_{j \in I} p_j {}^T\phi_{ij}h_j\overline{\phi}_{ij}.$$
Hence we have
$${}^T\phi_{ij}\widetilde{h}_j\overline{\phi}_{ij} = \sum_k p_k {}^T\phi_{ij}{}^T\phi_{jk}h_k\overline{\phi}_{jk}\overline{\phi}_{ij} = \sum_k p_k {}^T(\phi_{jk} \circ \phi_{ij})h_k\overline{(\phi_{jk} \circ \phi_{ij})}.$$

But since $\phi_{jk} \circ \phi_{ij} = \alpha_{ijk} \cdot \phi_{ik}$ and $|\alpha_{ijk}| = 1$, we see that

$$^T\phi_{ij}\widetilde{h}_j\overline{\phi}_{ij} = \sum_k p_k |\alpha_{ijk}|\,^T\phi_{ik} h_k \overline{\phi}_{ik} = \sum_k p_k\,^T\phi_{ik} h_k \overline{\phi}_{ik} = \widetilde{h}_i,$$

so that $\widetilde{h} = \{\widetilde{h}_i\}$ is a Hermitian metric on $E$. $\square$

As for untwisted vector bundles, we look for relations between Hermitian metrics and connections. More precisely, let $E = \{E_i, \phi_{ij}\}$ be an $\alpha$−twisted complex $C^\infty$ vector bundle, $D = \{D_i\}$ a connection on $E$ and $h = \{h_i\}$ a Hermitian metric on $E$.

**Definition 2.24.** *We say that $D$ and $h$ are **compatible** (or that $D$ is a $h$−**connection**) if for every $i \in I$ we have that $D_i$ is a $h_i$−connection, i. e. for every sections $\xi$ and $\eta$ of $E_i$ we have*

$$d(h_i(\xi, \eta)) = h_i(D_i(\xi), \eta) + h_i(\xi, D_i(\eta)).$$

Representing $h_i$ and $D_i$ by matrices $H_i$ (of smooth functions) and $\Gamma_i$ (of $1-$forms) with respect to a chosen local frame of $E_i$, this reads as

$$dH_i = {}^T\Gamma_i H_i + H_i \overline{\Gamma}_i.$$

Using this we prove the following:

**Lemma 2.25.** *Let $E$ be an $\alpha$−twisted holomorphic vector bundle and $h$ a Hermitian metric on $E$.*

1. *There is a unique connection $D$ on $E$ which is compatible both with the holomorphic structure of $E$ and the Hermitian metric $h$.*
2. *If the $B$−field $B$ is such that $B_i$ is a $(1,1)$−form for every $i \in I$, then $R_D \in A^{1,1}(End(E))$.*

*Proof.* As $E_i$ is a holomorphic vector bundle on $U_i$ and $h_i$ is a Hermitian connection on $E_i$, there is a unique connection $D_i$ on $E_i$ which is compatible with the holomorphic structure of $E_i$ and with the Hermitian metric $h_i$. If $\Gamma_i$ is a connection form of $D_i$ and $H_i$ is a matrix representing $h_i$ with respect to a given local frame of $E_i$, we know that

$$dH_i = {}^T\Gamma_i \cdot H_i + H_i \cdot \overline{\Gamma}_i.$$

We will moreover let $a_{ij}$ be the matrix of smooth functions representing $\phi_{ij}$ with respect to the chose local frames.

Let $p = \{p_i\}$ be a partition of the unity with respect to $\mathscr{U}$. The proof of Lemma 2.18 tell us that if we let $\widetilde{D}$ be the connection on $E$ whose connection form (with respect to the local frame given above) is

$$\widetilde{\Gamma}_i := \sum_{j \in I} p_j(a_{ij}^{-1}\Gamma_j a_{ij} + a_{ij}^{-1}da_{ij}),$$

then $\widetilde{D}$ is compatible with the holomorphic structure of $E$.

It only remains to show that $\widetilde{D}$ and $h$ are compatible. Notice that

$$^T\widetilde{\Gamma}_i \cdot H_i + H_i \cdot \overline{\widetilde{\Gamma}}_i = \sum_{j\in I} p_j\, ^T a_{ij}\, ^T\Gamma_j\, ^T a_{ij}^{-1} H_i + \sum_{j\in I} p_j\, ^T da_{ij}\, ^T a_{ij}^{-1} H_i +$$

$$+ \sum_{j\in I} p_j H_i \overline{a}_{ij}^{-1} \overline{\Gamma}_j \overline{a}_{ij} + \sum_{j\in I} p_j H_i \overline{a}_{ij}^{-1} \overline{da}_{ij}.$$

Using the fact that $H_i = {}^T a_{ij} H_j \overline{a}_{ij}$ we see that ${}^T a_{ij}^{-1} H_i = H_j \overline{a}_{ij}$ and that $H_i \overline{a}_{ij}^{-1} = {}^T a_{ij} H_j$, hence we find that

$$^T\widetilde{\Gamma}_i \cdot H_i + H_i \cdot \overline{\widetilde{\Gamma}}_i = \sum_{j\in I} p_j\, ^T a_{ij}\, ^T\Gamma_j H_j \overline{a}_{ij} + \sum_{j\in I} p_j\, ^T da_{ij} H_j \overline{a}_{ij} +$$

$$+ \sum_{j\in I} p_j\, ^T a_{ij} H_j \overline{\Gamma}_j \overline{a}_{ij} + \sum_{j\in I} p_j\, ^T a_{ij} H_j \overline{da}_{ij} =$$

$$= \sum_{j\in I} p_j\, ^T a_{ij} ({}^T\Gamma_j H_j + H_j \overline{\Gamma}_j) \overline{a}_{ij} + \sum_{j\in I} p_j ({}^T da_{ij} H_j \overline{a}_{ij} + {}^T a_{ij} H_j \overline{da}_{ij}).$$

But as ${}^T\Gamma_j H_j + H_j \overline{\Gamma}_j = dH_j$, we find

$$^T\widetilde{\Gamma}_i \cdot H_i + H_i \cdot \overline{\widetilde{\Gamma}}_i = \sum_{j\in I} p_j ({}^T a_{ij} dH_j \overline{a}_{ij} + {}^T da_{ij} H_j \overline{a}_{ij} + {}^T a_{ij} H_j \overline{da}_{ij}) =$$

$$= \sum_{j\in I} p_j d({}^T a_{ij} H_j \overline{a}_{ij}) = \sum_{j\in I} p_j dH_i = dH_i,$$

which proves that $\widetilde{D}$ is compatible with Hermitian metric $h$.

Let now $R_{\widetilde{D}}$ be the curvature of $\widetilde{D}$: as $\widetilde{D}$ is compatible with the holomorphic structure of $E$, we know from Lemma 2.18 that $R_{\widetilde{D}}^{0,2} = 0$. Moreover, for every $i \in I$ we have

$$R_{\widetilde{D}|U_i} = R_i - B_i \cdot id_{E_i},$$

where $R_i$ is the curvature of $\widetilde{D}_i$. As $\widetilde{D}_i$ is compatible with $h_i$ we know that also $R_i^{2,0} = 0$ (see section 4 in Chapter I of [17]). But since $B_i$ is a $(1,1)$-form by hypothesis, it follows that $R_{\widetilde{D}} \in A^{1,1}(End(E))$. □

Now, if $E$ is an $\alpha$-twisted holomorphic vector bundle and $h$ is a Hermitian metric on $E$, the previous Lemma allows us to give the following:

**Definition 2.26.** *The unique connection $D$ on $E$ which is compatible with $h$ and with the holomorphic structure of $E$ is called the **Chern connection** of the pair $(E, h)$, and its curvature will be called **Chern curvature** of the pair $(E, h)$. We will sometimes use the notation $D_h$ and $R_h$ for them.*

Notice that by definition we have that $D = \{D_i\}$ is the Chern connection of $(E, h)$ if and only if $D_i$ is the Chern connection of $(E_i, h_i)$. As an immediate consequence of Lemma 2.19 we get the following converse of Lemma 2.25:

**Lemma 2.27.** *Let $E$ be an $\alpha$−twisted complex $C^\infty$ vector bundle on $X$, $h$ a Hermitian metric on $E$ and $D$ a connection on $E$ compatible with $h$. If the $B$−field $B$ is such that $B_i$ is a $(1,1)$−form for every $i \in I$, then there is a unique holomorphic structure on $E$ so that $D$ is the Chern connection of $(E, h)$.*

2.6. **Connections and metrics on associated bundles.** We resume here the basic facts about how a connection (or a Hermitian metric) on a twisted vector bundle $E$ induces a connection (or a Hermitian metric) on twisted vector bundles that may be constructed from $E$.

**Dual bundle**. Let $E = \{E_i, \phi_{ij}\}$ be an $\alpha$−twisted complex $C^\infty$ vector bundle, and $D = \{D_i\}$ a connection on $E$. In particular $D_i$ is a connection on $E_i$, so we may use it to produce a connection $D_i^*$ on $E_i^*$: for every local section $\xi$ of $E_i^*$, we need to define a 1−form $D_i^*(\xi)$ with coefficients in $E_i^*$. We then define $D_i^*(\xi)$ by expressing the 1−form $(D_i^*(\xi), \eta)$ obtained by evaluating the coefficients of $D_i^*(\xi)$ on $\eta$. We then let

$$(D_i^*(\xi), \eta) := d(\xi, \eta) - (\xi, D_i(\eta)).$$

If $\Gamma_i$ is the connection form of $D_i$ with respect to a local frame of $E_i$, the connection form of $D_i^*$ is $-\Gamma_i$ with respect to the dual local frame (see section 5 in Chapter I of [17]).

**Lemma 2.28.** *The family $D^* = \{D_i^*\}_{i \in I}$ is a connection on $E^*$, and $R_{D^*} = -R_D$ as elements of $A^{1,1}(\underline{End}(E)) \simeq A^{1,1}(\underline{End}(E^*))$.*

*Proof.* For every local sections $\xi$ (of $E_i^*$) and $\eta$ (of $E_i$) we have

$$(D_i^*(\xi), \eta) = d(\xi, \eta) - (\xi, D_i(\eta)) =$$
$$= d(\xi, \eta) - (\xi, d\eta) - (\xi, a_{ij}^{-1}\Gamma_j a_{ij}\eta) - (\xi, a_{ij}^{-1}da_{ij}\eta) - (\xi, \omega_{ij}\eta) =$$
$$= (d\xi, \eta) - (a_{ij}^*\xi, \Gamma_j a_{ij}\eta) - (a_{ij}^*\xi, da_{ij}\eta) - (\omega_{ij}\xi, \eta) =$$
$$= (d\xi, \eta) + (\Gamma_j^* a_{ij}^*\xi, a_{ij}\eta) - d(a_{ij}^*\xi, a_{ij}\eta) + (da_{ij}^*\xi, a_{ij}\eta) - (\omega_{ij}\xi, \eta) =$$
$$= (d\xi, \eta) + ((a_{ij}^*)^{-1}\Gamma_j^* a_{ij}^*\xi, \eta) + ((a_{ij}^*)^{-1}da_{ij}^*\xi, \eta) - (\omega_{ij}\xi, \eta) =$$
$$= (d\xi + ((a_{ij}^*)^{-1}\Gamma_j^* a_{ij}^* + (a_{ij}^*)^{-1}da_{ij}^* - \omega_{ij})\xi, \eta),$$

where we let $a_{ij}^*$ be the matrix representing $\phi_{ij}^*$, and $\Gamma_i^*$ the connection form of $D_i^*$ with respect to the dual local frame. In conclusion, we get

$$\Gamma_i^* = (a_{ij}^*)^{-1}\Gamma_j^* a_{ij}^* + (a_{ij}^*)^{-1}da_{ij}^* - \omega_{ij} \cdot I_r,$$

so that $D^* = \{D_i^*\}$ is a connection on $E^*$.

If $B = \{B_i\}$ is a $B$−field inducing the twist $\alpha$, the family $-B = \{-B_i\}$ is a $B$−field inducing the twist $\alpha^{-1}$, and we have

$$R_{D^*|U_i} = R_i^* - (-B_i) = -R_i + B_i = -R_{D|U_i},$$

where $R_i^*$ is the curvature of the dual connection $D_i^*$ (here we use the fact that $R_i^* = -R_i$, see section 5 in Chapter I of [17]). □

The connection $D^*$ is called the **dual connection** of $D$, or equivalently **connection induced by $D$ on $E^*$**.

Let now $h = \{h_i\}$ be a Hermitian metric on $E$. Then $h_i$ is a Hermitian metric on $E_i$, i. e. an isomorphism $h_i : \overline{E}_i \longrightarrow E_i^*$ of complex $C^\infty$ vector bundles: the dual of this gives $h_i^* : \overline{E}_i^* \longrightarrow E_i^{**}$, which is then a Hermitian metric on $E_i^*$ represented by the matrix $H_i^* = \overline{H}_i^{-1}$ (since $H_i$ is Hermitian).

**Lemma 2.29.** *Let $E$ be an $\alpha-$twisted holomorphic vector bundle on $X$ and $h$ a Hermitian metric on $E$. The family $h^* = \{h_i^*\}$ is a Hermitian metric on $E^*$, and $D_{h^*} = D_h^*$.*

*Proof.* We notice that as $H_i = {}^T a_{ij} H_j \overline{a}_{ij}$, then $H_i^* = {}^T a_{ij}^* H_j^* \overline{a}_{ij}^*$, so that $h^* = \{h_i^*\}_{i \in I}$ is a Hermitian metric on $E^*$. If $D = \{D_i\}$ is the Chern connection of $(E, h)$, then $D_i$ is the Chern connection of $(E_i, h_i)$. The dual of $D$ is the $D^* = \{D_i^*\}$, where $D_i^*$ is the dual connection of $D_i$. But this implies that $D_i^*$ is the Chern connection of $(E_i^*, h_i^*)$, so that $D^*$ is the Chern connection of $(E^*, h^*)$ (see Lemma 2.25). $\square$

The Hermitian metric $h^*$ is called the **dual Hermitian metric** of $h$, or equivalently **Hermitian metric induced by $h$ on $E^*$**.

**Direct sum**. Let $E = \{E_i, \phi_{ij}\}$ and $F = \{F_i, \psi_{ij}\}$ be two $\alpha-$twisted complex $C^\infty$ vector bundles, and consider a connection $D = \{D_i\}$ on $E$ and a connection $D' = \{D_i'\}$ on $F$. In particular $D_i$ is a connection on $E_i$ and $D_i'$ is a connection on $F_i$, so we may use them to produce a connection $\overline{D}_i \oplus D_i'$ on $\overline{E}_i \oplus F_i$: a local section of $E_i \oplus F_i$ is of the form $\xi \oplus \xi'$ for a local section $\xi$ of $E_i$ and a local section $\xi'$ of $F_i$, so we let

$$(\overline{D}_i \oplus D_i')(\xi \oplus \xi') := D_i(\xi) \oplus D_i'(\xi') \in A^1(E_i \oplus F_i).$$

If $\Gamma_i$ is the connection form of $D_i$ with respect to a local frame of $E_i$ and $\Gamma_i'$ is the connection form of $D_i'$ with respect to a local frame of $F_i$, then

$$\Gamma_i \oplus \Gamma_i' := \begin{bmatrix} \Gamma_i & 0 \\ 0 & \Gamma_i' \end{bmatrix}$$

is the connection form of $D_i \oplus D_i'$ with respect to the corresponding local frame of $E_i \oplus F_i$.

**Lemma 2.30.** *The family $D \oplus D' = \{D_i \oplus D_i'\}_{i \in I}$ is a connection on $E \oplus F$, and we have $R_{D \oplus D'} = R_D \oplus R_{D'}$ as elements of $A^{1,1}(\underline{End}(E \oplus F))$.*

*Proof.* It is easy to see that $D \oplus D' = \{D_i \oplus D_i'\}$ is a connection on $E \oplus F$. Moreover, we have

$$R_{D \oplus D'|U_i} = R_{D_i \oplus D_i'} - B_i \cdot id_{E_i \oplus F_i} = R_{D_i} \oplus R_{D_i'} - B_i \cdot id_{E_i \oplus F_i} =$$

$$= (R_i - B_i \cdot id_{E_i}) \oplus (R_i' - B_i id_{F_i}) = R_{D|U_i} \oplus R_{D'|U_i} = (R_D \oplus R_{D'})_{|U_i}$$

(where we used the fact that $R_{D_i \oplus D_i'} = R_{D_i} \oplus R_{D_i'}$, see section 5 in Chapter I of [17]). We then get $R_{D \oplus D'} = R_D \oplus R_{D'}$, and we are done. $\square$

The connection $D \oplus D'$ is called the **direct sum connection** of $D$ and $D'$, or equivalently **connection induced by $D$ and $D'$ on $E \oplus F$**.

Let now $h = \{h_i\}$ be a Hermitian metric on $E$ and $h' = \{h'_i\}$ a Hermitian metric on $F$. Then $h_i$ is a Hermitian metric on $E_i$ and $h'_i$ is a Hermitian metric on $F_i$, and we define the sum Hermitian metric $h_i \oplus h'_i$ on $E_i \oplus F_i$ as

$$(h_i \oplus h'_i)(\xi \oplus \xi', \eta \oplus \eta') := h_i(\xi, \eta) + h'_i(\xi', \eta')$$

for every local sections $\xi, \xi'$ of $E_i$ and $\eta, \eta'$ of $F_i$.

**Lemma 2.31.** *Let $E$ and $F$ be two $\alpha-$twisted holomorphic vector bundles on $X$, $h$ a Hermitian metric on $E$ and $h'$ a Hermitian metric on $F$. The family $h \oplus h' = \{h_i \oplus h'_i\}_{i \in I}$ is a Hermitian metric on $E \oplus F$, and $D_{h \oplus h'} = D_h \oplus D_{h'}$.*

*Proof.* If $H_i$ and $H'_i$ represent $h_i$ and $h'_i$ with respect to local frames of $E_i$ and $F_i$, then

$$H_i \oplus H'_i = \begin{bmatrix} H_i & 0 \\ 0 & H'_i \end{bmatrix}$$

represents $h_i \oplus h'_i$ with respect to the corresponding local frame. It is then easy to see that $h \oplus h' = \{h_i \oplus h'_i\}$ is a Hermitian metric on $E \oplus F$.

If $D = \{D_i\}$ is the Chern connection of $(E, h)$, then $D_i$ is the Chern connection of $(E_i, h_i)$. Similarily, if $D' = \{D'_i\}$ is the Chern connection of $(F, h')$, then $D'_i$ is the Chern connection of $F_i$. The sum of $D$ and $D'$ is $D \oplus D' = \{D_i \oplus D'_i\}$, and we know that $D_i \oplus D'_i$ is the Chern connection of $(E_i \oplus F_i, h_i \oplus h'_i)$, and we are done. □

The Hermitian metric $h \oplus h'$ is called the **direct sum Hermitian metric** of $h$ and $h'$, or equivalently **Hermitian metric induced by $h$ and $h'$ on $E \oplus F$**.

**Tensor product**. Let $E = \{E_i, \phi_{ij}\}$ be an $\alpha-$twisted complex $C^\infty$ vector bundle and $F = \{F_i, \psi_{ij}\}$ be an $\alpha'-$twisted complex $C^\infty$ vector bundle. Consider a connection $D = \{D_i\}$ on $E$ and a connection $D' = \{D'_i\}$ on $F$. In particular $D_i$ is a connection on $E_i$ and $D'_i$ is a connection on $F_i$, so we may use them to produce a connection $D_i \otimes D'_i$ on $E_i \otimes F_i$: we let

$$D_i \otimes D'_i = D_i \otimes id_{F_i} + id_{E_i} \otimes D'_i.$$

If $\Gamma_i$ is the connection form of $D_i$ with respect to a local frame of $E_i$ and $\Gamma'_i$ is the connection form of $D'_i$ with respect to a local frame of $F_i$, then $\Gamma_i \otimes I_s + I_r \otimes \Gamma'_i$ is the connection form of $D_i \otimes D'_i$ with respect to the corresponding local frame of $E_i \otimes F_i$ (where $A \otimes B$ is the Kronecker product).

**Lemma 2.32.** *The family $D \otimes D' = \{D_i \otimes D'_i\}$ is a connection on $E \otimes F$, and*

$$R_{D \otimes D'} = R_D \otimes id_{F_i} + id_{E_i} \otimes R_{D'}$$

*as elements of $A^{1,1}(End(E \otimes F))$.*

*Proof.* We know that $D_i \otimes D'_i$ is a connection on $E_i \otimes F_i$, and an easy calculation shows that

$$\Gamma_i \otimes I_s + I_r \otimes \Gamma'_i = (a_{ij} \otimes b_{ij})^{-1}(\Gamma_j \otimes I_s + I_r \otimes \Gamma'_j)(a_{ij} \otimes b_{ij})+$$

$$+(a_{ij} \otimes b_{ij})^{-1}d(a_{ij} \otimes b_{ij}) + (\omega_{ij} + \omega'_{ij}) \cdot I_{rs},$$

where $a_{ij}$ and $b_{ij}$ are matrices of smooth functions representing $\phi_{ij}$ and $\psi_{ij}$ respectively with respect to local frames of $E_i$ and $F_i$. Hence $D \otimes D'$ is a connection on $E \otimes F$.

Now, let $B$ a $B$−field inducing the twist $\alpha$ and $B'$ a $B$−field inducing the twist $\alpha'$. Then $B + B' = \{B_i + B'_i\}_{i \in I}$ is a $B$−field inducing the twist $\alpha\alpha'$, and we have

$$R_{D \otimes D'|U_i} = R_{D_i \otimes D'_i} - (B_i + B'_i) \cdot id_{E_i \otimes F_i} =$$

$$= R_i \otimes id_{F_i} + id_{E_i} \otimes R'_i - B_i \cdot id_{E_i \otimes F_i} - B'_i \cdot id_{E_i \otimes F_i} =$$

$$= (R_i - B_i \cdot id_{E_i}) \otimes id_{F_i} + id_{E_i} \otimes (R'_i - B'_i id_{F_i})$$

(where we use the fact that $R_{D_i} \otimes R_{D'_i} = R_i \otimes id_{F_i} + id_{E_i} \otimes R'_i$, see section 5 in Chapter I of [17]), which implies the statement. $\square$

The connection $D \otimes D'$ is called the **tensor product connection** of $D$ and $D'$, or equivalently **connection induced by $D$ and $D'$ on $E \otimes F$**.

Let now $h = \{h_i\}$ be a Hermitian metric on $E$ and $h' = \{h'_i\}$ a Hermitian metric on $F$. Then $h_i$ is a Hermitian metric on $E_i$ and $h'_i$ is a Hermitian metric on $F_i$, and we define the product Hermitian metric $h_i \otimes h'_i$ on $E_i \otimes F_i$ as

$$(h_i \otimes h'_i)(\xi \otimes \xi', \eta \otimes \eta') := h_i(\xi, \eta) \cdot h'_i(\xi', \eta')$$

for every local sections $\xi, \xi'$ of $E_i$ and $\eta, \eta'$ of $F_i$, and then extending this by linearity. If $H_i$ and $H'_i$ represent $h_i$ and $h'_i$ with respect to local frames, then $H_i \otimes H'_i$ represents $h_i \otimes h'_i$.

**Lemma 2.33.** *Let $E$ be an $\alpha$−twisted holomorphic vector bundle on $X$, $F$ an $\alpha'$−twisted holomorphic vector bundle on $X$, $h$ a Hermitian metric on $E$ and $h'$ a Hermitian metric on $F$. The family $h \otimes h' = \{h_i \otimes h'_i\}_{i \in I}$ is a Hermitian metric on $E \otimes F$, and $D_{h \otimes h'} = D_h \otimes D_{h'}$.*

*Proof.* By the very definition we have

$$H_i \otimes H'_i = (^T a_{ij} H_j \overline{a}_{ij}) \otimes (^T b_{ij} H_j \overline{b}_{ij}) = {}^T(a_{ij} \otimes b_{ij}) H_i \otimes H_j \overline{a_{ij} \otimes b_{ij}},$$

so $h \otimes h'$ is a Hermitian metric on $E$. The remaining part of the proof is straightforward. $\square$

The Hermitian metric $h \otimes h'$ is called the **tensor product Hermitian metric** of $h$ and $h'$, or equivalently **Hermitian metric induced by $h$ and $h'$ on $E \otimes F$**.

As a particular case, if $E$ is an $\alpha-$twisted complex $C^\infty$ vector bundle $E$ and two integers $p, q \geq 0$, we let

$$E^{p,q} := \underbrace{E \otimes \cdots \otimes E}_{p} \otimes \underbrace{E^* \otimes \cdots \otimes E^*}_{q}$$

which is then an $\alpha^{p-q}-$twisted complex $C^\infty$ vector bundle. If $D$ is a connection on $E$, it induces a connection $D^{p,q}$ on $E^{p,q}$, and if $h$ is a Hermitian metric on $E$, it induces a Hermitian metric $h^{p,q}$ on $E^{p,q}$.

The most important case to consider is when $p = q = 1$, in which case $E^{1,1} = E \otimes E^* = \underline{End}(E)$: this is a usual complex $C^\infty$ vector bundle. If $D$ is a connection on $E$, the connection $D^{1,1}$ is a usual connection on a vector bundle, and if $h$ is a Hermitian metric on $E$, then $h^{1,1}$ is a usual Hermitian metric on a vector bundle. The Chern connection of $(E^{1,1}, h^{1,1})$ is the Chern connection of $(\underline{End}(E), h^{1,1})$.

**Wedge product**. Let $E$ be an $\alpha-$twisted complex $C^\infty$ vector bundle on $X$ and $p$ a strictly positive integer. Consider a connection $D = \{D_i\}$ on $E$. The product connection $D^{\otimes p}$ on $E^{\otimes p}$ is easily seen to verify the following: if $\xi$ is a section of $\bigwedge^p E_i$, then $D_i^{\otimes p}(\xi) \in A^1(\wedge^p E_i)$. It follows that $D^{\otimes p}_{|\wedge^p E}$ is a connection on $\bigwedge^p E$, denoted $D^p$ and called **wedge connection** on $E$, or equivalently **connection induced by** $D$ **on** $\wedge^p E$.

Let now $h = \{h_i\}$ be a Hermitian metric on $E$. The induced Hermitian metric $h^{\otimes p}$ restricted to $\wedge^p$ induces a Hermitian metric, denoted $h^p$, and called **wedge Hermitian metric** on $E$, or equivalently **Hermitian metric induced by** $h$ **on** $\wedge^p E$. The following is immediate:

**Lemma 2.34.** *Let $E$ be an $\alpha-$twisted holomorphic vector bundle on $X$ and $h$ a Hermitian metric on $E$. Then $D_{h^p} = D_h^p$.*

Particular case is when $p$ is the rank $r$ of $E$, in which case we have

$$\det(E) = \bigwedge^r E,$$

which is an $\alpha^r-$twisted complex $C^\infty$ line bundle on $X$. If $D$ is a connection on $E$, it induces a connection $\det(D)$ on $\det(E)$, called **determinant connection**, and if $h$ is a Hermitian metric on $E$, it induces a Hermitian metric $\det(h)$ on $\det(E)$, called **determinant Hermitian metric**. We moreover have $D_{\det(h)} = \det(D_h)$.

**Pull-back**. Let now $X$ and $Y$ be two complex manifolds and $f : X \longrightarrow Y$ be a holomorphic map between them. If $\mathscr{U} = \{U_i\}_{i \in I}$ is an open covering of $Y$, then $f^*\mathscr{U} := \{f^{-1}(U_i)\}_{i \in I}$ is an open covering of $X$. Let $E$ be an $\alpha-$twisted complex $C^\infty$ vector bundle on $Y$ and consider a connection $D = \{D_i\}$ on $E$.

In particular $D_i$ is a connection on $E_i$, so we may use it to produce a connection $f^*D_i$ on $f^*E_i$: to define it, we notice that if $\xi$ is a local section of $f^*E_i$, then there is a unique local section $\xi'$ of $E_i$ such that $\xi' \circ f = f' \circ \xi$

(where $f' : f^*E_i \longrightarrow E_i$ is the natural morphism induced by $f$). Hence we define
$$f^*D_i(\xi) := f^*(D_i(\xi')),$$
where on the right we have the pull-back under $f$ of the $1-$form $D_i(\xi')$ with coefficients in $E_i$. If $\Gamma_i$ is the connection form of $D_i$ with respect to a local frame of $E_i$, then $f^*\Gamma_i$ is the connection form of $f^*D_i$ with respect to the pull-back local frame.

**Lemma 2.35.** *The family $f^*D = \{f^*D_i\}_{i \in I}$ is a connection on $f^*E$, and we have $R_{f^*D} = f^*R_D$ as elements of $A^{1,1}(\underline{End}(f^*E))$.*

*Proof.* We know that $f^*D_i$ is a connection on $f^*E_i$, it is easy to show that $f^*D = \{f^*D_i\}$ is a connection on $f^*E$. If $B = \{B_i\}_{i \in I}$ is a $B-$field inducing the twist $\alpha$, then $f^*B = \{f^*B_i\}$ is a $B-$field inducing the twist $f^*\alpha$, and
$$R_{f^*D|f^{-1}(U_i)} = R_{f^*D_i} - f^*B_i \cdot id_{f^*E_i} = f^*R_i - f^*(B_i id_{E_i}) = f^*(\widetilde{R_i})$$
(where we use the fact that $R_{f^*D_i} = f^*R_i$, see section 5 in Chapter I of [17]). □

The connection $f^*D$ is called the **pull-back connection** of $D$, or equivalently **connection induced by** $D$ **on** $f^*E$.

Let now $h = \{h_i\}$ be a Hermitian metric on $E$. Then $h_i$ is a Hermitian metric on $E_i$, i. e. an isomorphism of complex $C^\infty$ vector bundles $h_i : \overline{E}_i \longrightarrow E_i^*$. Then $f^*h_i : f^*\overline{E}_i \longrightarrow f^*E_i^*$ is an isomorphism of complex $C^\infty$ vector bundles, getting a Hermitian metric $f^*h_i$. The following is immediate.

**Lemma 2.36.** *Let $E$ be an $\alpha-$twisted holomorphic vector bundle on $Y$ and $h$ a Hermitian metric on $E$. The family $f^*h = \{f^*h_i\}_{i \in I}$ is a Hermitian metric on $f^*E$, and $D_{f^*h} = f^*D_h$.*

The Hermitian metric $f^*h = \{f^*h_i\}$ is called **pull-back Hermitian metric** or **Hermitian metric induced by** $h$ **on** $f^*E$.

2.7. **Subbundles and quotients.** Let $E = \{E_i, \phi_{ij}\}$ be an $\alpha-$twisted holomorphic vector bundle on a complex manifold $X$, and let $r$ be its rank.

**Definition 2.37.** *A **twisted holomorphic subbundle of** $E$ is an $\alpha-$twisted holomorphic vector bundle $S = \{S_i, \psi_{ij}\}$ on $X$ such that for every $i \in I$ we have an injective morphism of $\alpha-$twisted holomorphic vector bundles $f : S \longrightarrow E$, i. e. a morphism $f = \{f_i\}$ such that $f_i : S_i \longrightarrow E_i$ is an injective morphism of holomorphic vector bundles on $U_i$ for every $i \in I$. The morphism $f$ is called **inclusion** of $S$ in $E$.*

Let now $S$ be a twisted holomorphic sub-bundle of $E$, and let $f : S \longrightarrow E$ be the inclusion. For every $i \in I$ we then may consider the quotient vector bundle $Q_i := E_i/S_i$, and we let
$$\varphi_{ij} : Q_{i|U_{ij}} \longrightarrow Q_{j|U_{ij}}, \quad \varphi_{ij}(x) := [\phi_{ij}(\overline{x})],$$
where $\overline{x}$ is a point in $E_i$ such that $[\overline{x}] = x$.

It is easy to verify that $\varphi_{ij}$ is a well-defined isomorphism of holomorphic vector bundles on $U_{ij}$, and that $Q = \{Q_i, \varphi_{ij}\}$ is an $\alpha-$twisted holomorphic vector bundle, called **quotient of $E$ by $S$**.

Moreover, for every $i \in I$ we have a natural projection $p_i : E_i \longrightarrow Q_i$, and we have $\varphi_{ij} \circ p_i = p_j \circ \phi_{ij}$. The family $p = \{p_i\}_{i \in I}$ is then a morphism $p : E \longrightarrow Q$ of $\alpha-$twisted holomorphic vector bundles, called **projection**.

We notice that we have an exact sequence of $\alpha-$twisted holomorphic vector bundles
$$0 \longrightarrow S \xrightarrow{f} E \xrightarrow{p} Q \longrightarrow 0.$$

2.7.1. *Hermitian metrics and orthogonals.* Let now $h = \{h_i\}$ be a Hermitian metric on $E$. As $f_i : S_i \longrightarrow E_i$, is a morphism of holomorphic vector bundles, we may define
$$h_i^S(\xi, \eta) := h_i(f_i(\xi), f_i(\eta))$$
for every sections $\xi$ and $\eta$ of $S_i$.

**Lemma 2.38.** *The family $h^S := \{h_i^S\}_{i \in I}$ is a Hermitian metric on $S$.*

*Proof.* As $f_i$ is injective and $h_i$ is a Hermitian metric on $S_i$, it is easy to see that $h_i^S$ is a Hermitian metric on $S_i$. For every sections $\xi, \eta$ of $S_i$ we have
$$h_i^S(\xi, \eta) = h_i(f_i(\xi), f_i(\eta)) = h_j(\phi_{ij}(f_i(\xi)), \phi_{ij}(f_i(\eta))) =$$
$$h_j(f_j(\psi_{ij}(\xi)), f_j(\psi_{ij}(\eta))) = h_j^S(\psi_{ij}(\xi), \psi_{ij}(\eta))$$
where we used the fact that $f$ is a morphism of $\alpha-$twisted holomorphic vector bundles. As a consequence, the family $h^S$ is a Hermitian metric on $S$. $\square$

For every $i \in I$ and every $x \in U_i$ we define
$$S_{i,x}^\perp := \{s \in E_{i,x} \,|\, h_{i,x}(s, f_{i,x}(t)) = 0 \text{ for every } t \in S_{i,x}\}.$$

We then get a complex $C^\infty$ vector bundle $S_i^\perp$ on $U_i$, called $h_i-$**orthogonal complement of $S_i$ in $E_i$**. Let now
$$\psi_{ij}^\perp := \phi_{ij|S_i^\perp} : S_{i|U_{ij}}^\perp \longrightarrow E_{j|U_{ij}}.$$

**Lemma 2.39.** *The family $S^\perp := \{S_i^\perp, \psi_{ij}^\perp\}_{i,j \in I}$ is an $\alpha-$twisted complex $C^\infty$ sub-bundle of $E$.*

*Proof.* First we need to prove that the image of $\psi_{ij}^\perp$ is contained in $S_j^\perp$. To do so, let $x \in U_{ij}$ and $s \in S_{i,x}^\perp$, we prove that $\psi_{ij,x}^\perp(s) \in S_{j,x}^\perp$, i. e. that $\phi_{ij}(s) \in S_{j,x}^\perp$. This means that $h_{j,x}(\phi_{ij,x}(s), f_{j,x}(t')) = 0$ for every $t \in S_{j,x}$.

As $\psi_{ij} : S_{i|U_{ij}} \longrightarrow S_{j|U_{ij}}$ is an isomorphism, there is $t \in S_{i,x}$ such that $t' = \psi_{ij}(t)$, hence we need to show that $h_{j,x}(\phi_{ij,x}(s), f_{j,x}(\psi_{ij,x}(t))) = 0$. Recall moreover that $f_j \circ \psi_{ij} = \phi_{ij} \circ f_i$, so we need to show that
$$h_{j,x}(\phi_{ij,x}(s), \phi_{ij,x}(f_{i,x}(t))) = 0.$$

As $h$ is a Hermitian metric on $E$, this holds if and only if $h_{i,x}(s, f_{i,x}(t)) = 0$, and this last holds as $s \in S_{i,x}^\perp$.

As a consequence we see that $\psi_{ij}^\perp : S_{i|U_{ij}}^\perp \longrightarrow S_{j|U_{ij}}^\perp$. By definition, this map is the restriction of a biholomorphism to a complex $C^\infty$ sub-bundle, hence $\psi_{ij}^\perp$ is injective and $C^\infty$. We need to show that it is surjective.

To do so, let $x \in U_{ij}$ and choose $s' \in S_{j,x}^\perp$. As $S_j^\perp \subseteq E_j$ and $\phi_{ij} : E_i \longrightarrow E_j$ is surjective, it follows that there is $s \in E_{i,x}$ such that $s' = \phi_{ij,x}(s)$. As $s' \in S_{j,x}^\perp$, for every $t' \in S_{j,x}$ we have $h_{j,x}(s', f_{j,x}(t')) = 0$.

Now, consider $t \in S_{i,x}$. We have

$$h_{i,x}(s, f_{i,x}(t)) = h_{j,x}(\phi_{ij,x}(s), \phi_{ij,x}(f_{i,x}(t))) = h_{j,x}(s', f_{j,x}(\psi_{ij,x}(t))) = 0,$$

so that $s \in S_{i,x}^\perp$. It follows that $\psi_{ij}^\perp : S_{i|U_{ij}}^\perp \longrightarrow S_{j|U_{ij}}^\perp$ is an isomorphism of complex $C^\infty$ vector bundles on $U_{ij}$.

As $\psi_{ij}^\perp = \phi_{ij|S_i^\perp}$, it is now easy to verify that $S^\perp = \{S_i^\perp, \psi_{ij}^\perp\}$ is an $\alpha$−twisted complex $C^\infty$ vector bundle on $X$, and that the natural inclusion $\iota_i : S_i^\perp \longrightarrow E_i$ makes it a twisted complex $C^\infty$ subbundle of $E$. $\square$

We will call $S^\perp$ the $h$−**orthogonal** of $S$ in $E$. A priori there is no reason why $S^\perp$ is holomorphic.

In any case, as for every $i \in I$ we have $E_i = S_i \oplus S_i^\perp$, it follows that $E = S \oplus S^\perp$ as $\alpha$−twisted complex $C^\infty$ vector bundles. Hence the exact sequence

$$0 \longrightarrow S \xrightarrow{f} E \xrightarrow{p} Q \longrightarrow 0$$

splits as an exact sequence of $\alpha$−twisted complex $C^\infty$ vector bundles, i. e. we have an isomorphism of $\alpha$−twisted complex $C^\infty$ vector bundles between $Q$ and $S^\perp$. We then have an exact sequence

$$0 \longrightarrow Q \xrightarrow{\varphi} E \xrightarrow{\pi} S \longrightarrow 0$$

of $\alpha$−twisted complex $C^\infty$ vector bundles, where $p \circ \varphi = id_Q$ and $\pi \circ f = id_S$.

The morphism $\varphi$ is then an injective morphism of $\alpha$−twisted complex $C^\infty$ vector bundles, i. e. $Q$ is an $\alpha$−twisted complex $C^\infty$ subbundle of $E$. We then may use $\pi$ to define a Hermitian metric $h^Q$ on $Q$. The Hermitian metrics $h^S$ and $h^Q$ are called **Hermitian metrics induced by $h$ on $S$ and $Q$**.

2.7.2. *Connections and orthogonals.* Let now $E$, $S$ and $Q$ as before, consider a Hermitian metric $h$ on $E$, and $h^S$ and $h^Q$ the induced Hermitian metrics. Let $D$ be the Chern connection of $(E, h)$, and consider the exact sequence

$$0 \longrightarrow S \xrightarrow{f} E \xrightarrow{p} Q \longrightarrow 0$$

of $\alpha$−twisted holomorphic vector bundles, and the exact sequence

$$0 \longleftarrow S \xleftarrow{\pi} E \xleftarrow{\varphi} Q \longleftarrow 0$$

of $\alpha$−twisted complex $C^\infty$ vector bundles. We write $\pi = \{\pi_i\}$ and $\varphi = \{\varphi_i\}$.

As $E_i = S_i \oplus Q_i$, for every section $\xi$ of $S_i$ we have

$$D_i(\xi) \in A^1(E_i) = A^1(S_i) \oplus A^1(Q_i).$$

We then write $D_i(\xi) = D_i^S(\xi) + A_i(\xi)$ where $D_i^S(\xi) \in A^1(S_i)$ and $A_i(\xi) \in A^1(Q_i)$. More precisely, we have

$$D_i^S(\xi) = \pi_i(D_i(f_i(\xi))), \quad A_i(\xi) = p_i(D_i(f_i(\xi))).$$

We then get two maps

$$D_i^S : A^0(S_i) \longrightarrow A^1(S_i), \quad A_i : A^0(S_i) \longrightarrow A^1(Q_i).$$

**Lemma 2.40.** *The family $D^S = \{D_i^S\}_{i \in I}$ is the Chern connection of $(S, h^S)$, and the maps $A_i$ glue together to form an element $A \in A^{1,0}(\underline{Hom}(S,Q))$.*

*Proof.* It is known that $D_i^S$ is the Chern connection of $(S_i, h_i^S)$ (see Proposition 6.4 in Chapter I of [17]). Choose now a local frame $t_i = \{t_1, \cdots, t_s\}$ of $S_i$, so that $f_i(t_i) = \{f_i(t_1), \cdots, f_i(t_s)\}$ may be completed to a local frame $t_i'$ of $E_i$. Let $\Gamma_i$ be the connection form of $D_i$ with respect to $t_i'$, and let $F_i$ be the matrix representing $f_i$ with respect to $t_i$ and $t_i'$, and $\Pi_i$ be the matrix representing $\pi_i$ with respect to $t_i'$ and $t_i$.

If $\Gamma_i^S$ is the connection form of $D_i^S$ with respect to $t_i$, we then have $\Gamma_i^S = \Pi_i \Gamma_i F_i$. It follows that

$$\Gamma_i^S = \Pi_i \Gamma_i F_i = \Pi_i a_{ij}^{-1} \Gamma_j a_{ij} F_i + \Pi_i a_{ij}^{-1} da_{ij} F_i + \Pi_i \omega_{ij} \cdot I_r F_i =$$

$$= b_{ij}^{-1} \Pi_j \Gamma_j F_j b_{ij} + b_{ij}^{-1} \Pi_j F_j db_{ij} + \omega_{ij} \Pi_i F_i,$$

where $a_{ij}$ is the matrix representing $\phi_{ij}$ with respect to $t_i'$ and $t_j'$, and $b_{ij}$ is the matrix representing $\psi_{ij}$ with respect to $t_i$ and $t_j$. But as $\pi_i \circ f_i = id_{S_i}$ we see that

$$\Gamma_i^S = b_{ij}^{-1} \Gamma_j^S b_{ij} + b_{ij}^{-1} db_{ij} + \omega_{ij} \cdot I_s,$$

showing that $D^S$ is a connection on $S$. It then follows that $D^S$ is the Chern connection of $(S, h^S)$.

Moreover, we know that $A_i \in A^{1,0}(\underline{Hom}(S_i, Q_i))$ (see Proposition 6.4 in Chapter I of [17]). Complete $t_i'$ as $t_i' = \{t_1, \cdots, t_s, t_{s+1}, \cdots, t_r\}$, so that $t_i'' = \{p_i(t_{s+1}), \cdots, p_i(t_r)\}$ is a local frame of $Q_i$. Let $P_i$ be the matrix representing $p_i$ with respect to $t_i'$ and $t_i''$. If we represent $A_i$ by a matrix $A_i$ with respect to $t_i$ and $t_i''$, we have $A_i = P_i \Gamma_i F_i$. Then

$$A_i = P_i \Gamma_i F_i = P_i a_{ij}^{-1} \Gamma_j a_{ij} F_i + P_i a_{ij}^{-1} da_{ij} F_i + P_i \omega_{ij} \cdot I_r F_i =$$

$$= c_{ij}^{-1} P_j \Gamma_j F_j c_{ij} + c_{ij}^{-1} P_j F_j dc_{ij} + \omega_{ij} P_i F_i,$$

where $c_{ij}$ is the matrix representing $\varphi_{ij}$ with respect to $t_i''$ and $t_j''$.

But as $p_i \circ f_i = 0$ we then get $A_i = c_{ij}^{-1} A_j c_{ij}$. It follows that the $A_i$'s glue together to give a global $A \in A^{1,0}(\underline{Hom}(S,Q))$. $\square$

In a similar way one produces a connection $D^Q$ on $Q$ as $D^Q = \{D_i^Q\}$, which turns out to be the Chern connection of $(Q, h^Q)$, and an element $C \in A^{0,1}(\underline{Hom}(Q,S))$. The form $A$ is called **second fundamental form of $S$**, and the form $Q$ is called **second fundamental form of $Q$**.

Let now $E$ be an $\alpha-$twisted holomorphic vector bundle on $X$ and $D$ a connection on it.

**Definition 2.41.** *An $\alpha-$twisted complex $C^\infty$ subbundle $E'$ of $E$ is called $D-$**invariant** if for every $i \in I$ and every section $\xi$ of $E'_i$ we have that $D_i(\xi) \in A^1(E'_i)$.*

The following will be used in the proof of the Kobayashi-Hitchin correspondence.

**Lemma 2.42.** *Le $E$ be an $\alpha-$twisted holomorphic vector bundle on $X$, $h$ a Hermitian metric on $E$ and $D$ the Chern connection of $(E,h)$. Suppose that $E'$ is a $D-$invariant subbundle of $E$, and let $E'' := (E')^\perp$. Then $E'$ and $E''$ are both $\alpha-$twisted holomorphic subbundles of $E$, and the direct sum decomposition $E = E' \oplus E''$ is holomorphic.*

*Proof.* This is an immediate consequence of the untwisted analogue, see Proposition 4.18 in Chapter I of [17]. □

An immediate corollary of this is the following:

**Corollary 2.43.** *Let $E$ be an $\alpha-$twisted holomorphic vector bundle on $X$, $h$ a Hermitian metric on $E$ and $D$ the Chern connection of $(E,h)$. Suppose that $S$ is an $\alpha-$twisted holomorphic subbundle, and let $A \in A^{1,0}(\underline{Hom}(S,S^\perp))$ be as before. If $A = 0$, then $S^\perp$ is an $\alpha-$twisted holomorphic subbundle which is isomorphic, as an $\alpha-$twisted holomorphic bundle, to $Q$.*

We end this section with the Gauss-Codazzi equations in the twisted setting. Let $E$ be an $\alpha-$twisted holomorphic vector bundle, $h$ a Hermitian metric on $E$, $D$ the Chern connection of $(E,h)$, $S$ a twisted holomorphic sub-bundle of $E$, $Q$ the quotient of $E$ by $S$ and $A \in A^{1,0}(\underline{Hom}(S,Q))$, $C \in A^{0,1}(\underline{Hom}(Q,S))$ as before.

For every $i \in I$ we let $R_i$ be the curvature of $D_i$, $R_i^S$ the curvature of $D_i^S$ and $R_i^Q$ the curvature of $D_i^Q$ of $(Q_i, h_i^Q)$. By Lemma 2.40, the Gauss-Codazzi equations in the untwisted setting give

$$R_i = \begin{bmatrix} R_i^S - C_i \wedge C_i^* & D_i^{1,0} C_i \\ -D_i^{0,1} C_i^* & R_i^Q - C_i^* \wedge C_i \end{bmatrix} = \begin{bmatrix} R_i^S - A_i^* \wedge A_i & D_i^{1,0} A_i^* \\ -D_i^{0,1} A_i & R_i^Q - A_i \wedge A_i^* \end{bmatrix},$$

where we let $D_i^{1,0}$ and $D_i^{0,1}$ for the $(1,0)-$part and the $(0,1)-$part of the connection induced by $D_i$ on $\underline{Hom}(S_i, Q_i)$ and $\underline{Hom}(Q_i, S_i)$ (see section 6 in Chapter I of [17]).

We notice that the $Hom(S_i, Q_i)$'s glue together to give a holomorphic vector bundle $\underline{Hom}(S,Q)$, and $D$ induces the Chern connection on it: hence the $D_i^{1,0}$'s glue together to give the $(1,0)-$part of the Chern connection on this bundle (and similarly for the $D_i^{0,1}$'s). If $B = \{B_i\}$ is the $B-$field, then

$$R_i - B_i \cdot id_{E_i} = \begin{bmatrix} R_i^S - B_i \cdot id_{S_i} - C_i \wedge C_i^* & D_i^{1,0} C_i \\ -D^{0,1} C_i^* & R_i^Q - B_i \cdot id_{Q_i} - C_i^* \wedge C_i \end{bmatrix} =$$

$$= \begin{bmatrix} R_i^S - B_i \cdot id_{S_i} - A_i^* \wedge A_i & D_i^{1,0} A_i^* \\ -D^{0,1} A_i & R_i^Q - B_i \cdot id_{Q_i} - A_i \wedge A_i^* \end{bmatrix}.$$

Now all the forms in the formula glue together to give

$$R_D = \begin{bmatrix} R_S - C \wedge C^* & D^{1,0} C \\ -D^{0,1} C^* & R_Q - C^* \wedge C \end{bmatrix} = \begin{bmatrix} R_S - A^* \wedge A & D^{1,0} A^* \\ -D^{0,1} A & R_Q - A \wedge A^* \end{bmatrix},$$

where $R_S := R_{D^S}$ and $R_Q = R_{D^Q}$. These are called the **twisted Gauss-Codazzi equations**.

2.8. **Hermitian forms and Hermitian endomorphisms.** We now define the notion of Hermitian form on an $\alpha-$twisted complex $C^\infty$ vector bundle $E$.

**Definition 2.44.** *A **Hermitian form** on $E$ is a family $u = \{u_i\}_{i \in I}$ where*
  (1) *for every $i \in I$, $u_i$ is a Hermitian form on $E_i$, i. e. a $C^\infty$ field of Hermitian products $u_{i,x}$ on the fibers of $E_i$ over $x \in U_i$,*
  (2) *for every $i, j \in I$ we have $u_i = {}^T\phi_{ij} u_j \overline{\phi}_{ij}$.*

Hermitian metrics on $E$ are clearly examples of Hermitian forms. Another useful example is the following.

**Example 2.45.** Take $A \subseteq \mathbb{R}$ an interval and consider a differentiable family $\hbar = \{h_t\}_{t \in A}$ of Hermitian metrics on $E$, i. e. a family of Hermitian metrics such that if we write $h_t = \{h_{t,i}\}_{i \in I}$ and represent $h_{t,i}$ by a matrix $H_{t,i}$ with respect to a local frame, then the entries of $H_{t,i}$ are differentiable in $t$.

For every $t \in A$ let
$$V_{t,i} := \partial_t H_{i,t},$$
i. e. the matrix whose entries are the derivatives of the entries of $H_{i,t}$ with respect to $t$. We then let $v_{t,i}$ be the form on $E_i$ represented by the matrix $V_i$ with respect to the given local frame, i. e. $v_{t,i} = \partial_t h_{i,t}$.

As $h_{t,i}$ is a Hermitian metric, it is easy to see that $v_{t,i}$ is a Hermitian form on $E_i$. Moreover, on $U_{ij}$ we have
$$v_{t,i} = \partial_t h_{i,t} = \partial_t({}^T\phi_{ij} h_{j,t} \overline{\phi}_{ij}) = {}^T\phi_{ij}(\partial_t h_{j,t})\overline{\phi}_{ij} = {}^T\phi_{ij} v_{t,j} \overline{\phi}_{ij},$$
so the family $v_t = \{v_{t,i}\}$ is a Hermitian form on $E$.

We will use the notation $v = \partial_t \hbar$, $v_t = \hbar'(t)$, or $v_t = h'_t$. The family $v = \{v_t\}_{t \in A}$ is a family of Hermitian forms on $E$, called **derivation of $\hbar$**.

2.8.1. *Endomorphism from a Hermitian metric.* Let $h$ be a Hermitian metric in $E$, and $v$ a Hermitian form on $E$. For every $i \in I$ let us define an endomorphism $f_i^{h,v} : E_i \longrightarrow E_i$ as follows: for every $x \in U_i$ and every $a \in E_{i,x}$, let $f_{i,x}^{h,v}(a)$ to be the unique element of $E_{i,x}$ such that
$$h_{i,x}(b, f_{i,x}^{h,v}(a)) = v_{i,x}(b, a).$$
The fact that $h_i$ is a Hermitian metric and $v_i$ is a Hermitian form imply that
$$f_i^{h,v} : E_i \longrightarrow E_i$$

is a well-defined smooth endomorphism of $E_i$.

For every $x \in U_{ij}$ and every $a, b \in E_{i,x}$ we have
$$v_{i,x}(b, a) = h_{i,x}(b, f^{h,v}_{i,x}(a)) = h_{j,x}(\phi_{ij,x}(b), \phi_{ij,x}(f^{h,v}_{i,x}(a))),$$
by using the fact that $h$ is a Hermitian metric, and
$$v_{i,x}(b, a) = v_{j,x}(\phi_{ij,x}(b), \phi_{ij,x}(a)) = h_{j,x}(\phi_{ij,x}(b), f^{h,v}_{j,x}(\phi_{ij,x}(a))),$$
hence for every $x \in U_{ij}$ we get
$$h_{j,x}(\phi_{ij,x}(b), \phi_{ij,x}(f^{h,v}_{i,x}(a))) = h_{j,x}(\phi_{ij,x}(b), f^{h,v}_{j,x}(\phi_{ij,x}(a))).$$

As $\phi_{ij}$ is an isomorphism, this implies that for every $a \in E_{i,x}$ and every $b' \in E_{j,x}$ we have
$$h_{j,x}(b', \phi_{ij,x}(f^{h,v}_{i,x}(a))) = h_{j,x}(b', f^{h,v}_{j,x}(\phi_{ij,x}(a))),$$
and as $h_{j,x}$ is non-degenerate we get
$$\phi_{ij,x}(f^{h,v}_{i,x}(a)) = f^{h,v}_{j,x}(\phi_{ij,x}(a)).$$

As this holds for every $x \in U_{ij}$ and for every $a \in E_{i,x}$ we then finally get
$$f^{h,v}_i = \phi_{ij}^{-1} f^{h,v}_j \phi_{ij},$$
i. e. the family $f^{h,v} = \{f^{h,v}_i\}$ is a smooth endomorphism of $E$.

**Definition 2.46.** *The endomorphism $f^{h,v}$ will be called **endomorphism associated to $h$ and $v$**.*

Choosing a local frame on $U_i$ and representing $v_i$ by a matrix $V_i$ and $h_i$ by a matrix $H_i$, we see that $f^{h,v}$ is represented by the matrix $H_i^{-1} V_i$.

**Remark 2.47.** If $h$ and $k$ are two Hermitian metrics, the endomorphism $f^{h,k}$ is an automorphism whose inverse is $f^{k,h}$. Indeed, for every $i \in I$ and sections $\xi, \eta$ of $E_i$, we have
$$k_i(\xi, \eta) = h_i(f^{h,k}_i(\xi), \eta), \quad h_i(\xi, \eta) = k_i(f^{k,h}_i(\xi), \eta).$$
Then
$$k_i(\xi, \eta) = h_i(f^{h,k}_i(\xi), \eta) = k_i(f^{k,h}_i(f^{h,k}_i(\xi)), \eta).$$
Since $k_i$ is a Hermitian metric, this implies that $f^{k,h}_i(f^{h,k}_i(\xi)) = \xi$, and hence that $f^{k,h} \circ f^{h,k} = id_E$.

**Remark 2.48.** If $v$ is a Hermitian form on $E$ and $h, k$ are two Hermitian metrics on $E$, we have $f^{k,h} \circ f^{h,v} = f^{k,v}$. Indeed, for every $i \in I$ consider two sections $\xi$ and $\eta$ of $E_i$. We then have
$$v_i(\xi, \eta) = h_i(f^{h,v}_i(\xi), \eta) = k_i(f^{k,h}_i(f^{h,v}_i(\xi)), \eta).$$
Since $v_i(\xi, \eta) = k_i(f^{k,v}_i(\xi), \eta)$ by definition, we then have $f^{k,h}_i \circ f^{h,v}_i = f^{k,v}_i$, and hence $f^{k,h} \circ f^{h,v} = f^{k,v}$.

**Remark 2.49.** If $v_1, v_2$ are Hermitian forms on $E$, $\lambda_1, \lambda_2 \in \mathbb{R}$ and $h$ is a Hermitian metric on $E$, then $\lambda_1 v_1 + \lambda_2 v_2$ is a Hermitian form on $E$ and we have $f^{h,\lambda_1 v_1 + \lambda_2 v_2} = \lambda_1 f^{h,v_1} + \lambda_2 f^{h,v_2}$.

Indeed, for every $i \in I$ and every two sections $\xi, \eta$ of $E_i$ we have

$$(\lambda_1 v_{1,i} + \lambda_2 v_{2,i})(\xi, \eta) = \lambda_1 v_{1,i}(\xi, \eta) + \lambda_2 v_{2,i}(\xi, \eta) =$$

$$= \lambda_1 h_i(f_i^{h,v_1}(\xi), \eta) + \lambda_2 h_i(f_i^{h,v_2}(\xi), \eta) = h_i((\lambda_1 f_i^{h,v_1} + \lambda_2 f_i^{h,v_2})(\xi), \eta).$$

But as we have

$$(\lambda_1 v_{1,i} + \lambda_2 v_{2,i})(\xi, \eta) = h_i(f_i^{h,\lambda_1 v_1 + \lambda_2 v_2}(\xi), \eta),$$

we conclude.

**Remark 2.50.** If $h$ is a Hermitian metric on $E$ and $v$ is a Hermitian form on $E$, then $f^{h,v}$ is diagonalizable and has the same signature of $v$. Indeed, consider a local frame $\sigma$ of $E_i$ with respect to which we represent $h_i$ and $v_i$ by Hermitian matrices $H_i$ and $V_i$, so that $f_i^{h,k}$ is represented by the matrix $F_i^{h,v} = H_i^{-1} V_i$ with respect to $\sigma$.

As $H_i$ is positive definite and $V_i$ is diagonalizable, by classical linear algebra (see as instance Theorem 7.6.3 of [12]) their product $F_i^{h,v}$ is diagonalisable, has real eigenvalues and same signature of $V_i$, proving the claim.

In particular, if $k$ is a Hermitian metric, then $f^{h,k}$ is diagonalizable and its eigenvalues are all strictly positive smooth function. It then makes sense to consider $\log(f^{h,v})$ and $(f^{h,v})^\sigma$ for every $\sigma \in (0,1]$.

A particular example of this construction is obtained by taking a differentiable family $\hbar = \{h_t\}_{t \in A}$. By Example 2.45 we know that $h'_t$ is a Hermitian form on $E$ for every $t \in A$, and hence we may consider the endomorphism $f^{h_t, h'_t}$ of $E$: we then get a function

$$f^\hbar : A \longrightarrow End(E), \quad f^\hbar(t) := f^{h_t, h'_t}.$$

If the family $\hbar'$ of Hermitian forms is differentiable as well, then $f^\hbar$ is differentiable.

2.8.2. *Hermitian endomorphisms.* A converse of the previous construction is also possible. Before, recall the following definition:

**Definition 2.51.** Let $E$ be an $\alpha-$twisted complex $C^\infty$ vector bundle on $X$, $h$ a Hermitian metric on $E$ and $f : E \longrightarrow E$ an endomorphism. We say that $f$ is $h-$**Hermitian** if for every $i \in I$ and every sections $\xi, \eta$ of $E_i$ we have

$$h_i(f_i(\xi), \eta) = h_i(\xi, f_i(\eta)).$$

We let $End_h(E)$ be the set of $h-$Hermitian endomorphisms of $E$, which is easily seen to be a real vector space. We first provide an easy example of $h-$Hermitian endomorphism:

**Example 2.52.** Let $h$ be a Hermitian metric on $E$ and $v$ a Hermitian form on $E$. The endomorphism $f^{h,v}$ is $h-$Hermitian: indeed, by definition for every $i \in I$ and every sections $\xi$ and $\eta$ of $E_i$ we have

$$h_i(f_i^{h,v}(\xi),\eta) = v_i(\xi,\eta) = \overline{v_i(\eta,\xi)} = \overline{h_i(f_i^{h,v}(\eta),\xi)} = h_i(\xi, f_i^{h,v}(\eta)).$$

If $h$ and $k$ are both Hermitian metrics on $E$, then $f^{h,k}$ is moreover $k-$Hermitian. Indeed, for every $i \in I$ and $\xi, \eta$ two sections of $E_i$ we have

$$k_i(f_i^{h,k}(\xi),\eta) = h_i(f_i^{h,k}(f_i^{h,k}(\xi)),\eta) = h_i(f_i^{h,k}(\xi), f_i^{h,k}(\eta)) = k_i(\xi, f_i^{h,k}(\eta)).$$

If $E$ is an $\alpha-$twisted complex $C^\infty$ vector bundle, $h$ is a Hermitian metric on $E$ and $f \in End_h(E)$, then we define a Hermitian form $\widehat{f}_h$ on $E$ as follows: for every sections $\xi, \eta$ of $E_i$ we let

$$\widehat{f}_{h,i}(\xi,\eta) := h_i(\xi, f_i(\eta)).$$

It is easy to see that as $f_i$ is $h_i-$Hermitian, then $\widehat{f}_{h,i}$ is a Hermitian form on $E_i$. Moreover we have

$$\widehat{f}_{h,i}(\xi,\eta) = h_i(\xi, f_i(\eta)) = h_j(\phi_{ij}(\xi), \phi_{ij} f_i(\eta)) =$$

$$= h_j(\phi_{ij}(\xi), f_j \phi_{ij}(\eta)) = \widehat{f}_{h,j}(\phi_{ij}(\xi), \phi_{ij}(\eta)),$$

so that $\widehat{f}_h := \{\widehat{f}_{h,i}\}_{i \in I}$ is a Hermitian form on $E$.

**Definition 2.53.** *The Hermitian form $\widehat{f}_h$ is called **Hermitian form associated to $h$ and $f$**.*

Let us first give an example of this construction.

**Example 2.54.** Let $h$ be a Hermitian metric on $E$ and $v$ a Hermitian form on $E$. By Example 2.52 the endormorphism $f^{h,v}$ is $h-$Hermitian. The Hermitian form $\widehat{f^{h,v}}_h$ associated to $h$ and $f^{h,v}$ is $v$. Indeed we have

$$\widehat{f^{h,v}}_{h,i}(\xi,\eta) = h_i(\xi, f_i^{h,v}(\eta)) = h_i(f_i^{h,v}(\xi),\eta) = v_i(\xi,\eta),$$

for every $i \in I$ and every sections $\xi, \eta$ of $E_i$.

Conversely, if $f$ a $h-$Hermitian endomorphism of $E$ we have $f = f^{h, \widehat{f}_h}$. Indeed, for every $i \in I$ and every sections $\xi, \eta$ of $E_i$ we have

$$\widehat{f}_{h,i}(\xi,\eta) = h_i(f_i(\xi),\eta), \quad \widehat{f}_{h,i}(\xi,\eta) = h_i(f_i^{h,\widehat{f}_h}(\xi),\eta).$$

A useful remark is about the eigenvalues of Hermitian endomorphisms:

**Remark 2.55.** The eigenvalues of $f \in End_h(E)$ are all real functions. Indeed if $\lambda$ is an eigenvalue of $f$, and if $s$ is an eigenvector of eigenvalue $\lambda$ over $U_i$, then

$$\lambda h_i(s,s) = h_i(\lambda_i s, s) = h_i(f_i(s), s) = h_i(s, f_i(s)) = h_i(s, \lambda s) = \overline{\lambda} h_i(s,s).$$

But as $s$ is nowhere vanishing we have $h_i(s,s) > 0$, so $\lambda = \overline{\lambda}$.

As a consequence, it makes sense to consider the subset $End_h^+(E)$ of $End_h(E)$ given by $h-$Hermitian endomorphisms of $E$ whose eigenvalues are all strictly positive: it is easy to see that it is a convex domain in $End_h(E)$.

Finally, we remark the following:

**Lemma 2.56.** *Let $f$ be a $h-$Hermitian diagonalizable endomorphism of $E$ whose eigenvalues are $\lambda_1, \cdots, \lambda_r$, and $\varphi : \mathbb{R} \longrightarrow \mathbb{R}$ is a smooth function such that the image of $\lambda_1, \cdots, \lambda_r$ lies in the definition domain of $\varphi$. Then $\varphi(f)$ is a $h-$Hermitian diagonalizable endomorphism.*

*Proof.* By Lemma 2.10 we have that $\varphi(f)$ is a diagonalizable endomorphism. Let $\sigma_i = \{\xi_1, \cdots, \xi_r\}$ be a local frame of $E_i$ diagonalizing $f_i$, and $\lambda_j$ be the eigenvalue corresponding to $\xi_j$. If we take $\xi, \eta$ two sections of $E_i$, write

$$\xi = \sum_{j=1}^r \alpha_j \xi_j, \quad \eta = \sum_{j=1}^r \beta_j \xi_j,$$

so that we have

$$h_i(\varphi(f_i)(\xi), \eta) = \sum_{j,k=1}^r \alpha_j \beta_k \varphi(\lambda_j) h_i(\xi_j, \xi_k),$$

$$h_i(\xi, \varphi(f_i)(\eta)) = \sum_{j,k=1}^r \alpha_j \beta_k \varphi(\lambda_k) h_i(\xi_j, \xi_k).$$

Notice that as $f$ is $h-$Hermitian, by Remark 2.55 we get

$$\lambda_j h_i(\xi_j, \xi_k) = h_i(\lambda_j \xi_j, \xi_k) = h_i(f_i(\xi_j), \xi_k) = h_i(\xi_j, f_i(\xi_k)) =$$

$$= h_i(\xi_j, \lambda_k \xi_k) = \lambda_k h_i(\xi_j, \xi_k).$$

If $h_i(\xi_j, \xi_k) \neq 0$, we then get $\lambda_j = \lambda_k$, so $\varphi(\lambda_j) = \varphi(\lambda_k)$ and hence $\varphi(f_i)$ is $h_i-$Hermitian, concluding the proof. □

**2.9. Space of Hermitian metrics.** Let now $E = \{E_i, \phi_{ij}\}$ be an $\alpha-$twisted complex $C^\infty$ vector bundle of rank $r$ on $X$. We let $Herm(E)$ be the set of Hermitian forms on $E$, and $Herm^+(E)$ be the set of Hermitian metrics on $E$. If $v, w \in Herm(E)$ and $\lambda \in \mathbb{R}$ we define

$$v + w := \{v_i + w_i\}_{i \in I}, \quad \lambda v := \{\lambda v_i\}_{i \in I}.$$

It is easy to see that $v+w$ and $\lambda v$ are Hermitian forms on $E$, and that under these two operations $Herm(E)$ is a real vector space (of dimension $r^2$).

Moreover, the subset $Herm^+(E)$ of $Herm(E)$ is a convex domain in $Herm(E)$, since the same holds for Hermitian metrics on vector bundles. As a consequence, we see that if $h \in Herm^+(E)$, then we may view $Herm(E)$ as the tangent space of $Herm^+(E)$ at $h$, i. e.

$$T_h(Herm^+(E)) = Herm(E).$$

We have an action of the group of automorphisms of $E$ on $Herm(E)$. More precisely, we let $GL(E)$ be the group of automorphisms of $E$, that we will call **complex gauge group of** $E$, and define the action

$$GL(E) \times Herm(E) \longrightarrow Herm(E) \quad (f, v) \mapsto fv := \{{}^T f_i v_i \overline{f}_i\}.$$

We notice that ${}^T f_i v_i \overline{f}_i$ is a Hermitian form on $E_i$, and that

$$ {}^T f_i v_i \overline{f}_i = {}^T f_i {}^T \phi_{ij} v_j \overline{\phi}_{ij} \overline{f}_i = {}^T(\phi_{ij} f_i) v_j \overline{\phi_{ij} f_i} = {}^T \phi_{ij} ({}^T f_j v_j \overline{f}_j) \overline{\phi}_{ij},$$

so that $fv := \{{}^T f_i v_i \overline{f}_i\}$ is a Hermitian form on $E$.

We will moreover let $\mathfrak{gl}(E)$ be the Lie algebra of $GL(E)$, i. e. the Lie algebra of global sections of $\underline{End}(E)$. If $k \in Herm^+(E)$, we let

$$U_k := \{f \in GL(E) \,|\, fk = k\},$$

i. e. the stabilizer of $k$ under the action of $GL(E)$: this will be called **complex gauge group of the pair** $(E, k)$, and we let

$$\mathfrak{u}_k(E) := \{f \in \mathfrak{gl}(E) \,|\, f \text{ is } k-\text{Hermitian}\}.$$

**Lemma 2.57.** *The action of $GL(E)$ on $Herm^+(E)$ is transitive.*

*Proof.* If $h, k \in Herm^+(E)$, then $h_i, k_i \in Herm^+(E_i)$, and on $Herm^+(E_i)$ the group $GL(E_i)$ acts transitively, i. e. there is a unique $f_i \in GL(E_i)$ such that $h_i = f_i k_i = {}^T f_i k_i \overline{f}_i$. Notice that

$$h_i = {}^T \phi_{ij} h_j \overline{\phi}_{ij} = {}^T \phi_{ij} {}^T f_j k_j \overline{f}_j \overline{\phi}_{ij} = {}^T(f_j \phi_{ij}) k_j \overline{f_j \phi}_{ij}.$$

and that

$$h_i = {}^T f_i k_i \overline{f}_i = {}^T f_i {}^T \phi_{ij} k_j \overline{\phi}_{ij} \overline{f}_i = {}^T(\phi_{ij} f_i) k_j \overline{\phi}_{ij} \overline{f}_i.$$

This implies that for every two local sections $\xi$ and $\eta$ of $E_i$ we have

$$k_j(f_j \phi_{ij} \xi, f_j \phi_{ij} \eta) = k_j(\phi_{ij} f_i \xi, \phi_{ij} f_i \eta),$$

so that we finally get $\phi_{ij} f_i = f_j \phi_{ij}$, i. e. $f = \{f_i\} \in GL(E)$ and $h = fk$. □

The previous Lemma allows us to identify $Herm^+(E)$ with the quotient $GL(E)/U_k(E)$, so we may consider $Herm^+(E)$ as a symmetric space with respect to the involution mapping $h = fk$ to $\overline{f}^{-1} k {}^T f^{-1}$.

**Remark 2.58.** We notice that if $h$ is a Hermitian metric on $E$, the map

$$\lambda : Herm(E) \longrightarrow End_h(E), \quad \lambda(v) := f^{h,v}$$

is an isomorphism of real vector spaces: it is linear by Remark 2.49, and by Example 2.54 its inverse is

$$\mu : End_h(E) \longrightarrow Herm(E), \quad \mu(f) := \widehat{f}_h.$$

If $k$ is a Hermitian metric, then by Remark 2.50 we know that $f^{h,k}$ is positive definite, i. e. we have

$$\lambda : Herm^+(E) \longrightarrow End_h^+(E).$$

Conversely, the Spectral Theorem implies that if $f$ is a positive definite $h-$Hermitian endomorphism, we have
$$\widehat{f}_{h,i}(\xi, \xi) = h_i(\xi, f_i(\xi)) > 0$$
if $\xi \neq 0$. It follows that $\widehat{f}_h$ is a Hermitian metric, and that $\lambda$ induces an identification between $Herm^+(E)$ and $End_h^+(E)$.

The previous Remark 2.58 allows us to provide an action of $GL(E)$ on $End_h(E)$ and $End_h^+(E)$, which is easily described by the following:

**Remark 2.59.** If $h \in Herm^+(E)$, $v \in Herm(E)$ and $a \in GL(E)$, then we have $ah \in Herm^+(E)$ and $av \in Herm(E)$. The associated endomorphism $f^{ah,av}$ is $a^{-1} \circ f^{h,v} \circ a$. Indeed, for every $i \in I$ and every two sections $\xi, \eta$ of $E_i$ we have
$$(av)_i(\xi, \eta) = (ah)_i(f_i^{ah,av}(\xi), \eta) = h_i(a_i f_i^{ah,av}(\xi), a_i(\eta)),$$
and
$$(av)_i(\xi, \eta) = v_i(a_i(\xi), a_i(\eta)) = h_i(f_i^{h,v} a_i(\xi), a_i(\eta)).$$
We then get $a_i \circ f_i^{ah,av} = f_i^{h,v} \circ a_i$.

2.9.1. *Riemannian metric and geodesics.* Now, let $h \in Herm^+(E)$, and take $v, w \in Herm(E)$. As already noticed, we have $Herm(E) = T_h(Herm^+(E))$, i. e. we may view $v, w$ as tangent vectors to $Herm^+(E)$ at $h$. We now want to define a metric on $T_h(Herm^+(E))$, so to have a metric on the space $Herm^+(E)$.

To do so, suppose that $X$ is a compact complex manifold of dimension $n$ and that $g$ is a Kähler metric on $X$, whose associated $(1,1)-$form is denoted $\sigma_g$. As $f^{h,v}, f^{h,w} \in End(E)$, we have that $f^{h,v} \circ f^{h,w} \in End(E)$ and that $Tr(f^{h,v} \circ f^{h,w})$ is a smooth function on $X$. We define
$$(v,w)_h := \int_X Tr(f^{h,v} \circ f^{h,w}) \sigma_g^n.$$

**Lemma 2.60.** *For every $h \in Herm^+(E)$ the map*
$$(\cdot, \cdot)_h : Herm(E) \times Herm(E) \longrightarrow \mathbb{R}$$
*is a positive definite bilinear symmetric product, which depends smoothly on $h$ and which is $GL(E)-$invariant.*

*Proof.* Linearity follows from Remark 2.49, while symmetry is a consequence of the properties of the trace. As $f^{h,v}$ is $h-$Hermitian, we see that $Tr(f^{h,v} \circ f^{h,w})$ is the Hilbert-Schmidt inner product, and hence it is positive definite. The $GL(E)-$invariance is an immediate consequence of Remark 2.59, while the smooth dependence on $h$ is obvious. □

As a consequence of Lemma 2.60 we get a $GL(E)-$invariant Riemannian metric on $Herm^+(E)$. For every $h \in Herm^+(E)$ and for every $v \in T_h Herm^+(E)$ we let
$$||v||_h := \sqrt{(v,v)_h}.$$

Now, consider $a < b$ two real numbers and $h, k \in Herm^+(E)$. We let $\Omega_{h,k}^{a,b}(E)$ be the space of piecewise differentiable fonctions
$$\hbar : [a,b] \longrightarrow Herm^+(E)$$
such that $\hbar(a) = h$ and $\hbar(b) = k$. We will write $h_t := \hbar(t)$, so that $h = h_a$ and $k = h_b$.

If $\hbar \in \Omega_{h,k}^{a,b}(E)$ is differentiable, by Example 2.45 we have a function
$$\hbar' : [a,b] \longrightarrow Herm(E), \quad \hbar'(s) = \frac{d\hbar}{dt}(s),$$
and a function
$$f^{\hbar} : [a,b] \longrightarrow End(E), \quad f^{\hbar}(t) = f^{h_t, h_t'},$$
where $h_t' = \hbar'(t)$. We let
$$\mathscr{E} : \Omega_{h,k}^{a,b}(E) \longrightarrow \mathbb{R}, \quad \mathscr{E}(\hbar) := \int_a^b \left( \int_X Tr(f^{h_t, h_t'} \circ f^{h_t, h_t'}) \sigma_g^n \right) dt.$$

We notice that
$$\int_X Tr(f^{h_t, h_t'} \circ f^{h_t, h_t'}) \sigma_g^n = \|h_t'\|_{h_t}^2,$$
hence $\mathscr{E}(\hbar)$ is the lenght of the curve parametrized by $\hbar$ in $Herm^+(E)$. The critical points of this functional correspond then to the geodetics in $Herm^+(E)$.

To describe the geodetics, consider $\hbar \in \Omega_{h,k}^{a,b}(E)$ and
$$v : [a,b] \longrightarrow Herm(E)$$
a differentiable function such that $v(a) = v(b) = 0$. For $s \in \mathbb{R}$ consider
$$\hbar + sv : [a,b] \longrightarrow Herm(E),$$
which is a piecewise differentiable function for every $s$. Moreover, if $s \ll 1$ for every $t \in [a,b]$ we have that $h_t + sv_t := \hbar(t) + sv(t) \in Herm^+(E)$. As $h_0 + sv_0 = h$ and $h_1 + sv_1 = k$, we have $\hbar + sv \in \Omega_{h,k}^{a,b}(E)$: this element is a small deformation of $\hbar$ in the direction of $v$.

We then have
$$\mathscr{E}(\hbar + sv) = \int_a^b \left( \int_X Tr(f^{(h+sv)_t, (h+sv)_t'}, f^{(h+sv)_t, (h+sv)_t'}) \sigma_g^n \right) dt,$$
so an easy calculation gives
$$\frac{d}{ds} \mathscr{E}(\hbar + sv) \bigg|_{s=0} = -2 \int_a^b \left( \int_X Tr(f^{h_t, v_t} \circ (f^{h_t, h_t''} - f^{h_t, h_t'} \circ f^{h_t, h_t'})) \sigma_g^n \right) dt.$$

As $\hbar$ is a critical point for the functional $\mathscr{E}$ if and only if
$$\frac{d}{ds} \mathscr{E}(\hbar + sv) \bigg|_{s=0} = 0$$

for every $v$, the previous calculation shows that $\hbar$ is a critical point of $\mathscr{E}$ if and only if
$$f^{h_t,h''_t} - f^{h_t,h'_t} \circ f^{h_t,h'_t} = 0.$$

If we represent $h_{t,i}$ by the matrix $H_{t,i}$ with respect to a chosen local frame of $E_i$, we see that $f_i^{h_t,h'_t}$ is represented by the matrix $H_{t,i}^{-1}H'_{t,i}$ with respect to the same local frame: here we let $H'_{t,i}$ be the matrix whose entries are the derivatives in $t$ of the entries of $H_{t,i}$. Hence we have that
$$\frac{d}{dt}(H_{t,i}^{-1}H'_{t,i}) = -H_{t,i}^{-1}H'_{t,i}H_{t,i}^{-1}H'_{t,i} + H_{t,i}^{-1}H''_{t,i}.$$

It follows that
$$\frac{d}{dt}f^{\hbar} = \frac{d}{dt}f^{h_t,h'_t} = -f^{h_t,h'_t} \circ f^{h_t,h'_t} + f^{h_t,h''_t}.$$

But then we see that $\hbar$ is a critical point of the functional $\mathscr{E}$, i. e. a geodetic in $Herm^+(E)$ from $h$ to $k$, if and only if
$$\frac{d}{dt}f^{h_t,h'_t} = 0,$$

i. e. if and only if $f^{h_t,h'_t}$ is an endomorphism of $E$ independent of $t$.

2.9.2. *The vector bundle of Hermitian forms.* Let now $X$ be a complex manifold of dimension $n$ and $E$ an $\alpha$−twisted holomorphic vector bundle on $X$ of rank $r$. For every open subset $U$ of $X$ we let $\mathscr{H}_E(U)$ be the set of Hermitian forms on $E_{|U}$. Since the sum of Hermitian forms is a Hermitian form, and the multiplication of a Hermitian form by a smooth real function is again a Hermitian form, it is easy to remark that the functor associating to any open subset $U$ of $X$ the set $\mathscr{H}_E(U)$, and to every $V \subseteq U$ the restriction map, is a sheaf of $C_X^\infty$−modules.

**Lemma 2.61.** *The sheaf $\mathscr{H}_E$ is locally free of rank $r^2$.*

*Proof.* For every $i \in I$ we have a natural morphism $\eta_i : \mathscr{H}_{E|U_i} \longrightarrow \mathscr{H}_{E_i}$ of sheaves of $C_{U_i}^\infty$−modules, the last being the sheaf of Hermitian forms on the (untwisted) vector bundle $E_i$. The sheaf $\mathscr{H}_{E_i}$ is known to be locally free of rank $r^2$, and we now show that $\eta_i$ is an isomorphism: this will imply the statement.

For every $x \in U_i$ the stalk $\mathscr{H}_{E_i,x}$ is the real vector space $Herm(E_{i,x})$ of Hermitian forms on $E_{i,x}$. The morphism
$$\eta_{i,x} : \mathscr{H}_{E,x} \longrightarrow Herm(E_{i,x})$$

maps $\alpha \in \mathscr{H}_{E,x}$ to $h(\alpha)_{i,x}$, where $h(\alpha)$ is a Hermitian form on $E_{|U}$ (for some open subset $U$ of $U_i$ containing $x$).

It is immediate to verify that $\eta_{i,x}$ is injective. For the surjectivity, let $\beta \in Herm(E_{i,x})$, and take an open subset $U$ of $U_i$ containing $x$ and a Hermitian form $h_i$ on $U$ such that $h_{i,x} = \beta$. For every $j \in I$ we then let $h_j := {}^T\phi_{ij}^{-1}h_i\overline{\phi}_{ij}^{-1}$, which is a Hermitian form on $E_{j|U}$, and we see that $h := \{h_j\}_{j \in I}$ is a Hermitian metric on $E_{|U}$ whose germ at $x$ has image $\beta$ under $\eta_{i,x}$. □

As $\mathscr{H}_E$ is a locally free sheaf of $C_X^\infty$−modules of rank $r^2$ on $X$, there is a $C^\infty$ real vector bundle $H_E$ of rank $r^2$ corresponding to it, whose space of global sections is $Herm(E)$.

2.9.3. *Norms*. We now introduce various norms on the space of $p$−forms with values in a vector bundle: we refer the reader to [17] and to [21], chapter 7, for more details.

Let $V$ be a real vector bundle of rank $s$ on a differentiable compact manifold $X$ of dimension $d$. Let $g$ be a Riemannian metric on $X$ and $h$ a fiber metric on $V$. Recall that on $A^*(V)$ we have a pointwise inner product associated to $g$ and $h$, and defined as follows: for every $\alpha \in A^p(X)$, $\beta \in A^q(X)$ and $s, t \in A^0(V)$ we let

$$(\alpha \cdot s, \beta \cdot t) := g(\alpha, \beta) \cdot h(s, t),$$

which is a smooth function on $M$. In particular, to every $\xi \in A^*(V)$ we associate a smooth function

$$|\xi| := (\xi, \xi) : X \longrightarrow \mathbb{R}.$$

We then have a $L^p$−norm on $A^*(V)$ defined as

$$||\xi||_{L^p} := \sqrt[p]{\int_X |\xi|^p \sigma_g^d}.$$

**Remark 2.62.** As a particular case, if $E$ is an $\alpha$−twisted vector bundle on $X$ and $\xi \in A^p(\underline{End}(E))$, if $h$ is a Hermitian metric on $E$ and $g$ is a Hermitian metric on $X$, we associate to $\xi$ a smooth function $|\xi|$ on $X$, which has the property that

$$|\xi|^2 \sigma_g^d = Tr(\xi \wedge *\xi^*),$$

where $\xi^*$ is the adjoint of $\xi$, and $*\xi^*$ is the $*$−Hodge of $\xi^*$ (see section 3.2).

Let now $\nabla$ be the Levi-Civita connection induced by $g$ on the tangent bundle $T_X$ of $X$, and consider a connection $D$ on $V$ which is compatible with $h$. For every $k \in \mathbb{N}$, the covariant derivative of the connection induced on $Sym^k(\Omega_X) \otimes V$ (where $\Omega_X$ is the cotangent bundle of $X$) gives a linear map

$$A^0(Sym^k(\Omega_X) \otimes V) \longrightarrow A^1(Sym^k(\Omega_X) \otimes V) = A^0(\Omega_X \otimes Sym^k(\Omega_X) \otimes V).$$

The natural map $\Omega_X \otimes Sym^k(\Omega_X) \longrightarrow Sym^{k+1}(\Omega_X)$ induces then a linear map

$$A^0(\Omega_X \otimes Sym^k(\Omega_X) \otimes V) \longrightarrow A^0(Sym^{k+1}(\Omega_X) \otimes V).$$

The composition of these maps gives then a linear map

$$d_{\nabla,D}^k : A^0(V) \longrightarrow A^0(Sym^k(\Omega_X) \otimes V).$$

Using the metric on $Sym^k(\Omega_X)$ induced by $h$ and $g$ we then may define for every $\xi \in A^0(V)$, and for every $p, q \in \mathbb{N}$ a $L_q^p$−norm as follows:

$$||\xi||_{L_q^p} := \sqrt[p]{\int_X \sum_{j=1}^q |d_{\nabla,D}^j(\xi)|^p \sigma_g^d}.$$

We will let $L_q^p(V)$ be the completion of $A^0(V)$ with respect to this norm.

A similar construction may be done in a relative context. More precisely, let $a > 0$ be a real number, and consider the projection $\pi_a : X \times [0, a] \longrightarrow X$: the pull-back $\pi_a^* V$ is a real smooth vector bundle of rank $s$ on $X \times [0, a]$, whose global sections are smooth curves in the space of global sections of $V$, defined over $[0, a]$.

Take a smooth family of fiber metrics $\{h_t\}_{t \in [0,a]}$ on $V$, and let $D_t$ be a connection on $V$ compatible with $h_t$, and we suppose that the family $\{D_t\}_{t \in [0,a]}$ is smooth. Consider a global section $f$ of $\pi_a^* V$, and notice that

$$f^{(j)} := \frac{d^j}{d^j t} f$$

is again a global section of $\pi_a^* V$.

We then define

$$||f||_{L_q^p} := \sqrt[p]{\sum_{i+2j \leq q} \int_{X \times [0,a]} |d_{\nabla,D_t}^i f^{(j)}|^p \sigma_g^d dt}.$$

This is a norm on the space of global sections of $\pi^* V$, and we let

- $C^\infty(a, V)$: the space of smooth sections of $\pi_a^* V$,
- $L_q^p(a, V)$: the completion of $C^\infty(a, V)$ with respect to the $L_q^p$−norm,
- $C_0^\infty(a, V)$: the space of smooth sections $f$ of $\pi_a^* V$ such that for every $i + 2j \leq q$ we have $d_{\nabla,D_t}^i f^{(j)} = 0$ on $X \times \{0\}$,
- $L_{q,0}^p(a, V)$: the completion of $C_0^\infty(a, V)$ with respect to the norm $L_q^p$.
- For every $b \in \mathbb{N}$ and every $\alpha \in (0, 1]$ we let $C^{b+\alpha}(a, V)$ be the space of sections $f$ of $\pi_a^* V$ of class $C^b$ such that for every $i, j \in \mathbb{N}$ such that $i + 2j = b$ we have that $d_{\nabla,D_t}^i f^{(j)}$ is Hölder continuous with exponent $\alpha$.

On the space $C^{b+\alpha}(a, V)$ define a norm

$$||f||_{b+\alpha} := \sum_{i+2j \leq b} \sup_{x \in X} |d_{\nabla,D_t}^i f^{(j)}| + \sum_{i+2j=b} \sup_{z \neq w \in X_a} \frac{|d_{\nabla,D_t}^i f^{(j)}(z) - d_{\nabla,D_t}^i f^{(j)}(w)|}{d(z,w)^\alpha},$$

where $d(z, w)$ is the distance between $z$ and $w$, and $X_a = X \times [0, a]$.

All these definitions will be applied to the real vector bundle $H_E$ (of rank $r^2$) on a compact Kähler manifold $X$ of (real) dimension $2n$. The metric $h$ on $H_E$ is the natural Riemannian metric we defined on $Herm(E)$.

## 3. Hermite-Einstein condition

In this section we introduce the notion of Hermite-Einstein and approximate Hermite-Einstein twisted vector bundles: both notions will be identical to the corresponding notions for untwisted bundles.

To do so, we will first need to introduce Chern forms and the Chern classes of twisted holomorphic vector bundles, which will be defined starting from the choice of a connection on a twisted bundle. We start with some preliminary notation.

Let $V$ be a complex vector space and $k \in \mathbb{N}$. We let $f_k : V^k \longrightarrow \mathbb{C}$ be a symmetric multilinear form of degree $k$ on $V$, and define

$$F_k : V \longrightarrow \mathbb{C}, \qquad F_k(v) := f_k(v, \cdots, v)$$

the associated degree $k$ homogeneous polynomial.

If $G$ is a linear group acting freely on $V$, we say that $f_k$ is $G-$**invariant** if for every $g \in G$ and for every $v_1, \cdots, v_k \in V$ we have

$$f_k(g \cdot v_1, \cdots, g \cdot v_k) = f_k(v_1, \cdots, v_k),$$

in which case we have $F_k(g \cdot v) = F_k(v)$.

A particular case is when $V = \mathfrak{gl}_r(\mathbb{C})$, the Lie algebra of $G = GL_r(\mathbb{C})$, on which $G$ acts by conjugation, i. e.

$$G \times V \longrightarrow V, \quad (A, X) \mapsto AXA^{-1}.$$

The form $f_k$ is $G-$invariant if and only if for every $v_1, \cdots, v_k, w \in V$ we have

$$\sum_{i=1}^{k} f_k(v_1, \cdots, v_{i-1}, [w, v_i], v_{i+1}, \cdots, v_k) = 0$$

(see section 2 in Chapter II of [17]).

We now define $r$ homogeneous polynomials $F_1, \cdots, F_r$ on $V$ by letting $F_k(X)$ be the homogeneous part of degree $k$ of

$$\det\left(I_r - \frac{1}{2\pi i}X\right),$$

i. e. we have

$$\det\left(I_r - \frac{1}{2\pi i}X\right) = 1 + F_1(X) + \cdots + F_r(X).$$

We notice that if $A \in G$ we have

$$\det\left(I_r - \frac{1}{2\pi i}AXA^{-1}\right) = \det\left(I_r - \frac{1}{2\pi i}X\right),$$

hence it follows that $F_1, \cdots, F_r$ are $G-$invariants.

3.1. **Chern classes from Chern connection.** Let $E$ be an $\alpha-$twisted holomorphic vector bundle of rank $r$ on a complex manifold $X$, and consider $\omega \in A^2(\underline{End}(E))$. Choose a local frame $s$ of $E_i$ over an open subset $U \subseteq U_i$, with respect to which we represent $\omega_{|U}$ by a matrix $\Omega_U$ of 2−forms. For every $k \in \mathbb{N}$ we then let
$$\gamma_{U,k}(\omega) := F_k(\Omega_U) \in A^{2k}(U).$$
If $s'$ is a local frame of $E_j$ over an open subset $U' \subseteq U_j$, then:
- if $i = j$, then $s$ are $s'$ are two local frames of $E_i$ over $U \cap U'$, and $\Omega_U$ and $\Omega_{U'}$ represent the same 2−form with values in $\underline{End}(E_i)$ with respect to these two local frames. There is then an invertible matrix $M$ of smooth functions such that $\Omega_{U'} = M\Omega_U M^{-1}$ (see section 1 in Chapter I of [17]).
- If $i \neq j$ and $U \cap U' \neq \emptyset$, let $a_{ij}$ be the matrix representing $\phi_{ij}$ with respect to $s$ and $s'$. We then have $\Omega_{U'} = a_{ij}^{-1}\Omega_U a_{ij}$.

In any case, we see that there is an invertible matrix $A$ of smooth functions on $U \cap U'$ such that $\Omega_{U'} = A\Omega_U A^{-1}$. As $F_k$ is invariant under the action of $GL_r(\mathbb{C})$, we see that
$$\gamma_{U,k}(\omega) = \gamma_{U',k}(\omega).$$
It follows that the 2k−forms $\gamma_{U,k}(\omega)$ glue together to give $\gamma_k(\omega) \in A^{2k}(X)$.

If now $E$ is an $\alpha-$twisted holomorphic vector bundle $E$ and $D$ is a connection on it, then $R_D \in A^2(\underline{End}(E))$, so the previous construction gives us $\gamma_k(R_D) \in A^{2k}(X)$ for every $k = 1, \cdots, r$.

**Lemma 3.1.** *For every $k = 1, \cdots, r$ the 2k−form $\gamma_k(R_D)$ is d−closed.*

*Proof.* Let $D = \{D_i\}$, and choose an open covering $\mathscr{U}'$ of $X$ as before, i. e. $\mathscr{U}' = \{U'_j\}_{j \in J}$ is such that for every $j \in J$ there is $i \in I$ such that $U'_j \subseteq U_i$, and there is a local frame $s_j$ of $E_i$ over $U'_j$. Represent $R_{D|U'_j}$ with respect to $s_j$ by a matrix $\widetilde{\Omega}_j$ of 2−forms, and the curvature $R_i$ of $D_i$ by a matrix $\Omega_j$ of 2−forms. We have
$$\widetilde{\Omega}_j = \Omega_j - B_i \cdot I_r, \quad \gamma_k(R_D)_{|U'_j} = F_k(\widetilde{\Omega}_j).$$
As $B_i$ is $d$−closed, we get
$$d\widetilde{\Omega}_j = d\Omega_j - dB_i \cdot id_{E_i} = d\Omega_j.$$
If $\Gamma_j$ is the connection form of $D_i$ with respect to $s_j$, the Bianchi identity for $D_i$ (see section 1 in Chapter I of [17]) gives
$$d\widetilde{\Omega}_j = d\Omega_j = \Omega_j \cdot \Gamma_j - \Gamma_j \cdot \Omega_j = [\Omega_j, \Gamma_j] = [\widetilde{\Omega}_j, \Gamma_j].$$
But since $f_k$ is $GL_r(\mathbb{C})$−invariant, this implies that
$$d\gamma_k(R_D)_{|U'_j} = dF_k(\widetilde{\Omega}_j) = df_k(\widetilde{\Omega}_j, \cdots, \widetilde{\Omega}_j) =$$
$$= \sum_{i=1}^{k} f_k(\widetilde{\Omega}_j, \cdots, \widetilde{\Omega}_j, \underbrace{d\widetilde{\Omega}_j}_{i}, \widetilde{\Omega}_j, \cdots, \widetilde{\Omega}_j) =$$

$$= \sum_{i=1}^{k} f_k(\widetilde{\Omega}_j, \cdots, \widetilde{\Omega}_j, \underbrace{[\widetilde{\Omega}_j, \Gamma_j]}_{i}, \widetilde{\Omega}_j, \cdots, \widetilde{\Omega}_j) = 0.$$

We finally get $d\gamma_k(R_D) = 0$. $\square$

**Definition 3.2.** *The $d-$closed smooth $2k-$form $\gamma_k(R_D)$ is called $k-$**th Chern form of $E$ with respect to** $D$. Its cohomology class $c_k(R_D) \in H^{2k}(X, \mathbb{C})$ is the $k-$**th Chern class of $E$ with respect to** $D$.*

We now analyze how $\gamma_k(R_D)$ varies with $D$.

**Lemma 3.3.** *If $D$ and $D'$ are two connections on $E$, then for every $k \in \mathbb{Z}$ we have that $\gamma_k(R_D) - \gamma_k(R_{D'})$ is an exact $2k-$form on $X$.*

*Proof.* By Remark 2.14, for every $t \in [0,1]$ the affine linear combination $D_t := (1-t)D + tD'$ is a connection on $E$. Let $\mathscr{U}' = \{U'_j\}_{j \in J}$ be an open covering of $X$ as in the proof of Lemma 3.1, and let $\Gamma_j$ be the connection form of $D$, $\Gamma'_j$ the connection form of $D'$ and $\Gamma_{j,t}$ the connection form of $D_t$ with respect to a local frame over $U'_j$. We clearly have $\Gamma_{j,t} := (1-t)\Gamma_j + t\Gamma'_j$.

If $\Delta_j := \Gamma'_j - \Gamma_j$, then $\Gamma_{j,t} = \Gamma_j + t\Delta_j$ and by Remark 2.14 we see that the $\Delta_j$ represents a global $\delta \in A^1(\underline{End}(E))$ with respect to the given local frame. As in the proof of Lemma 3.1 we represent $R_D$ by a matrix $\widetilde{\Omega}_j$, $R_{D'}$ by a matrix $\widetilde{\Omega}'_j$, and the curvature $R_t$ of $D_t$ by a matrix $\widetilde{\Omega}_{j,t}$. Moreover, we represent the curvature of $D_i$ (resp. of $D'_i$, of $D_{i,t}$) by a matrix $\Omega_j$ (resp. $\Omega'_j$, $\Omega_{j,t}$), so that we have

$$\widetilde{\Omega}_j = \Omega_j - B_i \cdot I_r, \quad \widetilde{\Omega}'_j = \Omega'_j - B_i \cdot I_r, \quad \widetilde{\Omega}_{j,t} = \Omega_{j,t} - B_i \cdot I_r.$$

We then have

$$\widetilde{\Omega}_{j,t} = \Omega_{j,t} - B_i \cdot I_r = d\Gamma_{j,t} + \Gamma_{j,t} \wedge \Gamma_{j,t} - B_i \cdot I_r,$$

so that

$$\partial_t \widetilde{\Omega}_{j,t} = \partial_t d\Gamma_{j,t} + \partial_t(\Gamma_{j,t} \wedge \Gamma_{j,t}) =$$
$$= \partial_t d(\Gamma_j + t\Delta_j) + \partial_t((\Gamma_j + t\Delta_j) \wedge (\Gamma_j + t\Delta_j)) =$$
$$= d\Delta_j + \Delta_j \wedge \Gamma_{j,t} + \Gamma_{j,t} \wedge \Delta_j = D_{j,t}\Delta_j.$$

This means that $\partial_t R_t = D_t(\delta)$, where $D_t$ denotes here the connection induced by $D_t$ on $\underline{End}(E)$. This implies that

$$kdf_k(\Delta_j, \widetilde{\Omega}_{j,t}, \cdots, \widetilde{\Omega}_{j,t}) = kD_{j,t}f_k(\Delta_j, \widetilde{\Omega}_{j,t}, \cdots, \widetilde{\Omega}_{j,t}) =$$
$$= kf_k(D_{j,t}\Delta_j, \widetilde{\Omega}_{j,t}, \cdots, \widetilde{\Omega}_{j,t}) = kf_k\left(\partial_t\widetilde{\Omega}_{j,t}, \widetilde{\Omega}_{j,t}, \cdots, \widetilde{\Omega}_{j,t}\right) =$$
$$= \partial_t f_k(\widetilde{\Omega}_{j,t}, \cdots, \widetilde{\Omega}_{j,t}) = \partial_t \gamma_k(R_t)_{|U_j}.$$

Let us now consider the smooth $(2k-1)-$form on $U'_j$ defined by

$$\varphi_{k,j} := k \int_0^1 f_k(\Delta_j, \widetilde{\Omega}_{j,t}, \cdots, \widetilde{\Omega}_{j,t}) dt.$$

As the $\Delta_j$'s glue together to form a global $1-$form, it is easy to prove that if $U'_{j'}$ is another open subset, then $\varphi_{k,j|U'_{jj'}} = \varphi_{k,j'|U_{jj'}}$, so there is a unique smooth $(2k-1)-$form $\varphi_k$ on $X$ such that $\varphi_{k|U_j} = \varphi_{k,j}$.

We now have
$$d\varphi_{k|U_j} = d\varphi_{k,j} = \int_0^1 k df_k(\Delta_j, \widetilde{\Omega}_{j,t}, \cdots, \widetilde{\Omega}_{j,t}) dt = \int_0^1 \partial_t \gamma_k(R_t)_{|U_j} dt =$$
$$= \gamma_k(R_{D'})_{|U_j} - \gamma_k(R_D)_{|U_j}.$$
But this implies that $\gamma_k(R_{D'}) - \gamma_k(R_D) = d\varphi_k$, and we are done. □

As a consequence, the $k-$th Chern class of $E$ does not depend on $D$, so we will write it as $c_k(E)$, and call it $k-$**th Chern class of** $E$.

If now $E$ is an $\alpha-$twisted holomorphic vector bundle and $h$ is a Hermitian metric on $E$, we may choose the Chern connection $D_h$ in order to calculate the Chern forms and the Chern classes. We will use the notation $\gamma_k(E, h)$ instead of $\gamma_k(R_{D_h})$, and call it the $k-$**th Chern form of** $(E, h)$.

As the Chern curvature of $(E, h)$ is a $(1, 1)-$form, it follows that $\gamma_k(E, h)$ is a $(k, k)-$form, and hence the $k-$th Chern class of $E$ is $c_k(E) \in H^{k,k}(X)$.

## 3.2. Mean curvature.

If $M$ is a $C^\infty$ differentiable manifold, a Hermitian metric $g$ on $T_M \otimes \mathbb{C}$ (the complexified tangent bundle of $M$) is called **Hermitian metric on** $M$. A holomorphic structure on $M$ is a holomorphic structure on $T_M$, so that the pair given by the differentiable manifold $M$ and the given complex structure is a complex manifold $X$, whose tangent bundle is denoted $T_X$ (i. e. it is the tangent bundle of $M$ with the given complex structure).

The Hermitian metric $g$ on $X$ induces a Hermitian metric, still denoted $g$, on $\bigwedge^k T_X$ and $\bigwedge^k \Omega_X$ for every $k$, where $\Omega_X$ is the cotangent bundle of $X$. In particular, if $\xi, \eta$ are two $k-$forms on $X$, i. e. two local sections of $\bigwedge^k \Omega_X$, we may calculate $g(\xi, \eta)$, which is a smooth function on $X$.

If $X$ has dimension $n$, then for every smooth $k-$form $\xi$ on $X$ there is a unique smooth $(2n-k)-$form on $X$, denoted $*\xi$, such that for every smooth $k-$form $\eta$ on $X$ we have
$$\eta \wedge *\xi = g(\eta, \xi) \cdot \sigma_g^n,$$
where $\sigma_g$ is the real $(1,1)-$form associated to $g$ and to the complex structure of $X$, i. e. if $z_1, \cdots, z_n$ are local holomorphic coordinate on an open subset $U$ of $X$ and $g_{ij} = g(\partial/\partial z_i, \partial/\partial z_j)$, then
$$\sigma_{g|U} = \sqrt{-1} \sum_{i,j} g_{ij} dz_i \wedge d\bar{z}_j.$$

If $\xi$ is a $(p, q)-$form, then $*\xi$ is a $(n-q, n-p)-$form, so we have the **Hodge $*-$operator**
$$* : A^{p,q}(X) \longrightarrow A^{n-q,n-p}(X),$$

whose conjugate is
$$\bar{*}: A^{p,q}(X) \longrightarrow A^{n-p,n-q}(X).$$

We then have an inner product
$$(\cdot,\cdot)_g: A^{p,q}(X) \times A^{p,q}(X) \longrightarrow \mathbb{R}, \quad (\eta,\xi)_g := \int_X \eta \wedge \bar{*}\xi.$$

Another useful operator is the **Lefschetz operator**
$$L_g: A^{p,q}(X) \longrightarrow A^{p+1,q+1}(X), \quad L_g(\xi) := \xi \wedge \sigma_g,$$
whose adjoint operator
$$\Lambda_g: A^{p+1,q+1}(X) \longrightarrow A^{p,q}(X)$$
is characterized by
$$(L_g(\xi),\eta) = (\xi,\Lambda_g(\eta))$$
for every $\xi \in A^{p,q}(X)$ and $\eta \in A^{p+1,q+1}(X)$. A formula for the adjoint of the Lefschetz operator is
$$\Lambda_g = *^{-1} \circ L_g \circ *.$$
It is easy to extend all these operators to $(p,q)$−forms with coefficients in any vector bundle on $X$.

Let now $E$ be an $\alpha$−twisted holomorphic vector bundle on $X$ and $h$ a Hermitian metric on $E$. Suppose that the $B$−field $B$ is given by $d$−closed, purely imaginary $(1,1)$−forms, and let moreover $R_h$ be the Chern curvature of $(E,h)$. We define
$$K_g(E,h) := i\Lambda_g(R_h).$$
As $R_h \in A^{1,1}(\underline{End}(E))$, we have that $K_g(E,h) \in A^0(\underline{End}(E))$, i. e. it is an endomorphism of the complex $C^\infty$ vector bundle $\underline{End}(E)$.

**Definition 3.4.** *The endomorphism $K_g(E,h)$ is called $g$−**mean curvature** of $(E,h)$.*

By the very definition of $\Lambda_g$ we then see that
$$K_g(E,h) \cdot \sigma_g^n = \sqrt{-1}nR_h \wedge \sigma_g^{n-1}.$$

**Lemma 3.5.** *The mean curvature $K_g(E,h)$ of $(E,h)$ is a $h$−Hermitian endomorphism of $E$.*

*Proof.* By the very definition of $K_g(E,h)$ we have that
$$K_g(E,h)_{|U_i} = K_g(E_i,h_i) - i\Lambda_g(B_i) \cdot id_{E_i},$$
where $K_g(E_i,h_i)$ is the mean curvature of the vector bundle $E_i$ with Hermitian metric $h_i$.

The mean curvature of a holomorphic vector bundle with respect to a Hermitian metric is known to be Hermitian with respect to that metric, so $K_g(E_i,h_i)$ is $h_i$−Hermitian. Moreover, as $B_i$ is a purely imaginary $(1,1)$−form, $i\Lambda_g(B_i)$ is a real smooth function.

Notice that if $f : U_i \longrightarrow \mathbb{R}$ is a real smooth function, for every sections $\xi, \eta$ of $E_i$ we have
$$h_i(f\xi, \eta) = f h_i(\xi, \eta) = \overline{f} h_i(\xi, \eta) = h_i(\xi, f\eta),$$
so $f \cdot id_{E_i}$ is $h_i$−Hermitian. In particular $i\Lambda_g(B_i) \cdot id_{E_i}$ is $h_i$−Hermitian. It follows that $K_g(E,h)_{|U_i}$ is $h_i$−Hermitian, and we are done. $\square$

The Hermitian form associated to the Hermitian metric $h$ and to the $h$−Hermitian endomorphism $K_g(E,h)$ is denoted $\widehat{K}_g(E,h)$ and called **$g$−mean curvature Hermitian form** of $(E,h)$.

Another useful result is the following:

**Lemma 3.6.** *Let $\hbar \in \Omega^{0,1}_{h,k}(E)$ be a differentiable family of Hermitian metrics on $E$. We then have*
$$\partial_t K_g(E, h_t) = i\Lambda_g \overline{\partial} D_t^{1,0} f^{\hbar},$$
*where $f^{\hbar} : [0,1] \longrightarrow End(E)$ maps $t$ to $f^{h_t, h'_t}$.*

*Proof.* Let $h_t := \hbar(t)$ and $D_t$ the Chern connection of $(E, h_t)$. Given a local frame of $E_i$, let $\Gamma_{i,t}$ be the connection form of $D_{i,t}$, and represent $h_{i,t}$ by a matrix $H_{i,t}$. As
$$^T H_{i,t} \cdot \Gamma_{i,t} = \partial\, ^T H_{i,t}$$
(see the proof of Proposition 4.9 in Chapter I of [17]), we get that
$$\partial_t(^T H_{i,t} \cdot \Gamma_{i,t}) = \partial_t \partial\, ^T H_{i,t}.$$
If we let $V_{t,i} := \partial_t H_{i,t}$, we get
$$^T V_{t,i} \cdot \Gamma_{t,i} + {}^T H_{t,i} \cdot \partial_t \Gamma_{t,i} = \partial V_{t,i}.$$
Written in another way we then get
$$^T H_{t,i} \cdot \partial_t \Gamma_{t,i} = \partial V_{t,i} - {}^T V_{t,i} \cdot \Gamma_{t,i} = D^{1,0}_{t,i} V_{t,i},$$
so we get
$$\partial_t \Gamma_{i,t} = D^{1,0}_{t,i}\, {}^T H^{-1}_{t,i} \cdot V_{t,i}.$$
Passing from the matrix notation to the usual notation we then get that
$$\partial_t D_{i,t} = D^{1,0}_{t,i} f_i^{h_t, h'_t},$$
since the matrices $V_{t,i}$'s represent the Hermitian form $h'_t$ with respect to the chosen local frames. Applying $D^{0,1}_i$ to both sides, and using the fact that $D^{0,1}_i = D^{0,1}_{t,i} = \overline{\partial}_i$ (since both are Chern connections, and hence compatible with the holomorphic structure) we get then
$$\partial_t R_{i,t} = \overline{\partial}_i D^{1,0}_{t,i} f_i^{h_t, h'_t},$$
where $R_{i,t}$ is the curvature of the Chern connection $D_{i,t}$ of $(E_i, h_{i,t})$.

As the $B$−field $B = \{B_i\}_{i \in I}$ does not depend on $t$, i. e. we have $\partial_t B_i = 0$. This implies that
$$\partial_t(R_{i,t} - B_i \cdot id_{E_i}) = \overline{\partial}_i D^{1,0}_{t,i} f_i^{h_t, h'_t}.$$

As both sides glue together to give global elements of $A^{1,1}(\underline{End}(E))$, we finally get
$$\partial_t R_t = \overline{\partial} D_t^{1,0} f^{h_t, h'_t}.$$
Applying $i\Lambda_g$ to both sides we conclude. $\square$

### 3.3. Hermite-Einstein metrics.

Let now $E$ be an $\alpha$−twisted holomorphic vector bundle on $X$, and consider a Hermitian metric $h$ on $E$.

**Definition 3.7.** *The pair $(E, h)$ verifies the **weak $g$−Hermite-Einstein condition** if there is a real function $\varphi : X \longrightarrow \mathbb{R}$ such that*
$$K_g(E, h) = \varphi \cdot id_E.$$
*The function $\varphi$ is called **Einstein function** of $(E, h)$ relative to $g$. If $\varphi$ is constant, we will say that $(E, h)$ verifies the $g$−**Hermite-Einstein condition**, and the constant number $\varphi$ will be called **Einstein factor** of $(E, h)$ relative to $g$.*

Let us first see some properties.

**Proposition 3.8.** *If $L$ is an $\alpha$−twisted holomorphic line bundle on $X$ and $h$ is a Hermitian metric on $L$, then $(L, h)$ verifies the weak $g$−Hermite-Einstein condition for every Hermitian metric $g$ on $X$.*

*Proof.* As $K_g(L, h)$ is a smooth section of $\underline{End}(L) \simeq A^0(X)$ (since $L$ is a line bundle), there is a smooth function $\varphi$ such that $K_g(L, h) = \varphi \cdot id_L$. $\square$

**Proposition 3.9.** *If $(E, h)$ verifies the weak $g$−Hermite-Einstein condition with Einstein function $\varphi$, then $(E^*, h^*)$ verifies the weak $g$−Hermite-Einstein condition with Einstein function $-\varphi$.*

*Proof.* By Lemma 2.29 if $D$ is the Chern connection of $(E, h)$, then $D^*$ is the Chern connection of $(E^*, h^*)$, and we have $R_{D^*} = -R_D$. Hence
$$K_g(E^*, h^*) = i\Lambda_g R_{D^*} = -i\Lambda_g R_D = -K_g(E, h),$$
and we are done. $\square$

**Proposition 3.10.** *If $(E_1, h_1)$ verifies the weak $g$−Hermite-Einstein condition with Einstein function $\varphi_1$, and $(E_2, h_2)$ verifies the weak $g$−Hermite-Einstein condition with Einstein function $\varphi_2$, then $(E_1 \otimes E_2, h_1 \otimes h_2)$ verifies the weak $g$−Hermite-Einstein condition with Einstein function $\varphi_1 + \varphi_2$.*

*Proof.* By Lemma 2.33 if $D_i$ is the Chern connection of $(E_i, h_i)$, then $D_1 \otimes D_2$ is the Chern connection of $(E_1 \otimes E_2, h_1 \otimes h_2)$, and we have
$$R_{D_1 \otimes D_2} = R_{D_1} \otimes id_{E_2} + id_{E_1} \otimes R_{D_2}.$$
Hence
$$K_g(E_1 \otimes E_2, h_1 \otimes h_2) = i\Lambda_g(R_{D_1} \otimes id_{E_2} + id_{E_1} \otimes R_{D_2}) =$$
$$= i\Lambda_g(R_{D_1}) \otimes id_{E_2} + id_{E_1} \otimes i\Lambda_g(R_{D_2}) =$$
$$= K_g(E_1, h_1) \otimes id_{E_2} + id_{E_1} \otimes K_g(E_2, h_2).$$
The statement then follows readily. $\square$

**Proposition 3.11.** *Let $E_1, E_2$ be two $\alpha-$twisted holomorphic vector bundles. Then $(E_1, h_1)$ and $(E_2, h_2)$ verify the weak $g-$Hermite-Einstein condition with Einstein function $\varphi$ if and only if $(E_1 \oplus E_2, h_1 \oplus h_2)$ verifies the weak $g-$Hermite-Einstein condition with Einstein function $\varphi$.*

*Proof.* By Lemma 2.31 if $D_i$ is the Chern connection of $(E_i, h_i)$, then $D_1 \oplus D_2$ is the Chern connection of $(E_1 \oplus E_2, h_1 \oplus h_2)$, and we have
$$R_{D_1 \oplus D_2} = R_{D_1} \oplus R_{D_2}.$$
Hence
$$K_g(E_1 \oplus E_2, h_1 \oplus h_2) = i\Lambda_g(R_{D_1}) \oplus i\Lambda_g(R_{D_2}) = K_g(E_1, h_1) \oplus K_g(E_2, h_2).$$
The statement then follows readily. □

The following is an immediate consequence of Propositions 3.9, 3.10 and 3.11:

**Proposition 3.12.** *If $(E, h)$ verifies the weak $g-$Hermite-Einstein condition with Einstein function $\varphi$, then*
  (1) *for every $p, q \in \mathbb{N}$, $(E^{p,q}, h^{p,q})$ verifies the weak $g-$Hermite-Einstein condition with Einstein function $(p - q)\varphi$;*
  (2) *for every $p \in \mathbb{N}$, $(\wedge^p E, \wedge^p h)$ verifies the weak $g-$Hermite-Einstein condition with Einstein function $p\varphi$.*

Finally, we have the following:

**Proposition 3.13.** *Let $f : X \longrightarrow Y$ be a morphism of complex manifolds, $E$ an $\alpha-$twisted holomorphic vector bundle on $Y$ and $h$ a Hermitian metric on $E$. If $(E, h)$ verifies the weak $g-$Hermite-Einstein condition with Einstein function $\varphi$, then $(f^*E, f^*h)$ verifies the weak $f^*g-$Hermite-Einstein condition with Einstein function $\varphi \circ f$.*

*Proof.* By Lemma 2.36 if $D$ is the Chern connection of $(E, h)$, then $f^*D$ is the Chern connection of $(f^*E, f^*h)$, and we have $R_{f^*D} = f^*R_D$. It follows that
$$K_g(f^*E, f^*h) = i\Lambda_g R_{f^*D} = i\Lambda_g f^* R_D = f^*(i\Lambda_g R_D) = f^* K_g(E, h),$$
so the statement follows readily. □

We now give the following definition:

**Definition 3.14.** *An $\alpha-$twisted holomorphic vector bundle on a compact, complex manifold with Hermitian metric $g$ is called $g-$**Hermite-Einstein** if it admits a Hermitian metric $h$ such that $(E, h)$ verifies the $g-$Hermite-Einstein condition.*

The previous Propositions 3.8 to 3.13 tell us that:
- twisted holomorphic line bundles are all $g-$Hermite-Einstein;
- the dual of a $g-$Hermite-Einstein twisted holomorphic vector bundle is $g-$Hermite-Einstein;

- the tensor product of $g-$Hermite-Einstein twisted holomorphic vector bundles is $g-$Hermite-Einstein;
- the direct sum of twisted holomorphic vector bundles is $g-$Hermite-Einstein if and only if the summands are $g-$Hermite-Einstein with the same Einstein factor;
- the pull-back of a $g-$Hermite-Einstein twisted holomorphic vector bundle is Hermite-Einstein (with respect to the pull-back metric).

The following is a key result in the proof of one direction of the Kobayashi-Hitchin correspondence.

**Proposition 3.15.** *Let $E_1$ and $E_2$ be two $\alpha-$twisted holomorphic vector bundles on $X$, and let $h_i$ be a Hermitian metric on $E_i$ for $i = 1, 2$. Let $g$ a Hermitian metric on $X$ and suppose that $(E_i, h_i)$ verifies the weak $g-$Hermite-Einstein condition with Einstein function $\varphi_i$.*

  (1) *If $\varphi_2 < \varphi_1$, then every morphism $f \in \mathscr{H}om(E_1, E_2)$ is 0.*
  (2) *If $\varphi_2 \leq \varphi_1$, then every morphism $f \in \mathscr{H}om(E_1, E_2)$ is a morphism of $\alpha-$twisted holomorphic vector bundles. We moreover have direct sum decompositions $E_1 = \ker(f) \oplus E_1''$ and $E_2 = Im(f) \oplus E_2''$ as $\alpha-$twisted holomorphic vector bundles, and $f$ maps the Chern connection of $(E_1'', h_{1|E_1''})$ to the Chern connection of $(Im(f), h_{2|Im(f)})$.*

*Proof.* Let us first suppose that $\varphi_2 < \varphi_1$, so we prove that if $f : E_1 \longrightarrow E_2$ is a morphism of $\alpha-$twisted sheaves, then $f = 0$.

Recall that $f$ is a global section of $E_1^* \otimes E_2$ (which is an untwisted holomorphic vector bundle). By Lemmas 2.29 and 2.33 the Hermitian metric $h_1$ and $h_2$ induce the Hermitian metric $h_1^* \otimes h_2$ on $E_1^* \otimes E_2$, and if $D_i$ is the Chern connection of $(E_i, h_i)$, then $D_1^* \otimes D_2$ is the Chern connection of $(\underline{Hom}(E_1, E_2), h_1^* \otimes h_2)$.

By Propositions 3.9 and 3.10 it then follows that $(E_1^* \otimes E_2, h_1^* \otimes h_2)$ is a holomorphic vector bundle verifying the weak $g-$Hermite-Einstein condition with Einstein function $\varphi_2 - \varphi_1$. As $\varphi_2 < \varphi_1$ the mean curvature of $E_1^* \otimes E_2$ is negative definite everywhere on $X$: by Theorem 1.9 in Chapter III of [17] it has then no non-zero global sections, so $f = 0$.

Suppose now that $\varphi_2 \leq \varphi_1$ and let $f : E_1 \longrightarrow E_2$ be a morphism of $\alpha-$twisted sheaves. The previous part of the proof tells us that $f$ is a global section of the $g-$Hermite-Einstein holomorphic vector bundle $E_1^* \otimes E_2$, whose Einstein function is $\varphi_2 - \varphi_1$: the mean curvature is then everywhere negative semi-definite, hence by Theorem 1.9 in Chapter III of [17] $f$ has to be parallel with respect to the Chern connection $D_1^* \otimes D_2$.

If $f = \{f_i\}$, this means that $f_i$ is parallel with respect to the Chern connection $D_{1,i}^* \otimes D_{2,i}$, so that the rank of $f_i$ has to be constant. It follows that $f_i$ is a morphism of holomorphic vector bundles, and hence that $f$ is a morphism of $\alpha-$twisted holomorphic vector bundles.

As a consequence, $\ker(f)$ is an $\alpha-$twisted holomorphic subbundle of $E_1$ and $Im(f)$ is an $\alpha-$twisted holomorphic subbundle of $E_2$. As $f$ is parallel, they are both invariant with respect to the Chern connections. We then let $E_1'' := \ker(f)^\perp$ and $E_2'' := Im(f)^\perp$, where the orthogonality is with respect to $h_1$ and $h_2$ respectively. By Lemma 2.42 the statement follows. $\square$

### 3.4. First Chern class and Hermite-Einstein.

Let now $E$ be an $\alpha$−twisted holomorphic vector bundle of rank $r$ on a complex manifold $X$ of dimension $n$, and let us fix a Hermitian metric $h$ on $E$ and a Hermitian metric $g$ on $X$. We let $\sigma_g$ be the real $(1,1)$−form on $X$ associated to $g$ and to the complex structure of $X$.

If $R_h$ is the Chern connection of $(E, h)$, then by its very definition the first Chern form of $(E, h)$ is

$$\gamma_1(E, h) = \frac{i}{2\pi} Tr(R_h) \in A^{1,1}(X).$$

**Lemma 3.16.** *If $X$ is compact and $(E, h)$ verifies the weak $g$−Hermite-Einstein condition with Einstein function $\varphi$, then*

$$\int_X \gamma_1(E, h) \wedge \sigma_g^{n-1} = \frac{r}{2n\pi} \int_X \varphi \sigma_g^n.$$

*Proof.* We have

$$K_g(E, h) \cdot \sigma_g^n = \sqrt{-1} n R_h \wedge \sigma_g^{n-1}.$$

As $(E, h)$ verifies the weak $g$−Hermite-Einstein condition with Einstein function $\varphi$, we have $K_g(E, h) = \varphi \cdot id_E$, hence

$$\sqrt{-1} n R_h \wedge \sigma_g^{n-1} = \varphi \cdot id_E \sigma_g^n$$

as elements of $A^{2n}(\underline{End}(E))$, so that

$$\sqrt{-1} n Tr(R_h) \wedge \sigma_g^{n-1} = r\varphi \sigma_g^n.$$

But since $Tr(R_h) = \frac{2\pi}{i} \gamma_1(E, h)$ we then get

$$\gamma_1(E, h) \wedge \sigma_g^{n-1} = \frac{r}{2n\pi} \varphi \sigma_g^n,$$

and the statement follows. $\square$

If $g$ is a Kähler metric, then $\sigma_g$ is $d$−closed, hence the value

$$\int_X \gamma_1(E, h) \wedge \sigma_g^{n-1}$$

only depends on $c_1(E)$ and the cohomology class $[\sigma_g]$. The same holds more generally if $\sigma_g^{n-1}$ is $\partial\overline{\partial}$−closed.

### 3.5. Hermite-Einstein and weak Hermite-Einstein.

We now show that if a twisted vector bundle verifies the weak $g$−Hermite-Einstein condition, then it is $g$−Hermite-Einstein. Let $E$ be an $\alpha$−twisted holomorphic vector bundle on a complex manifold $X$, $h$ a Hermitian metric on $E$ and $\chi : X \longrightarrow \mathbb{R}$ a positive smooth function.

For every $i \in I$ and for every sections $\xi, \eta$ of $E_i$ we let

$$h_i^\chi(\xi, \eta) := \chi \cdot h_i(\xi, \eta).$$

It is easy to see that $h_i^\chi$ is a Hermitian metric on $E_i$, and since

$$h_i^\chi = \chi \cdot h_i = \chi(^T\phi_{ij} h_j \overline{\phi}_{ij}) = {}^T\phi_{ij}(\chi \cdot h_j)\overline{\phi}_{ij} = {}^T\phi_{ij} h_j^\chi \overline{\phi}_{ij},$$

we see that $h^\chi = \{h_i^\chi\}$ is a Hermitian metric on $E$.

**Definition 3.17.** *The Hermitian metric $h^\chi$ is the **conformal change** of $h$ by $\chi$.*

The conformal change of Hermitian metric affects the curvature, and we have:

**Lemma 3.18.** *If $(E, h)$ verifies the weak $g-$Hermite-Einstein condition with Einstein function $\varphi$, then $(E, h^\chi)$ verifies the weak $g-$Hermite-Einstein condition with Einstein function $\varphi + i\Lambda_g \overline{\partial}\partial \chi$.*

*Proof.* Let $D$ be the Chern connection of $(E, h)$ and $D^\chi$ the Chern connection of $(E, h^\chi)$. Consider an open covering $\mathscr{U}' = \{U'_j\}_{j \in J}$ as in the proof of Lemma 3.1: hence for every $j \in J$ there is $i \in I$ such that $U'_j \subseteq U_i$, and on $U'_j$ there is a local frame $s_j$ of $E_i$.

We let $\Gamma_j$ and $\Gamma_j^\chi$ be the connection forms of $D_i$ and $D_i^\chi$, and we represent $h_i$ and $h_i^\chi$ by the matrices $H_j$ and $H_j^\chi$ respectively over $U'_j$ with respect to the local frame $s_j$. Then

$$^T\Gamma_j = \partial H_j \cdot H_j^{-1}, \quad {}^T\Gamma_j^\chi = \partial H_j^\chi \cdot (H_j^\chi)^{-1}$$

(see the proof of Proposition 4.9 in Chapter I of [17]). As $H_j^\chi = \chi \cdot H_j$ it follows that

$$^T\Gamma_j^\chi = \partial(\chi \cdot H_j) \cdot (\chi \cdot H_j)^{-1} = (\partial(\chi) \cdot H_j + \chi \cdot \partial H_j) \cdot \left(\frac{1}{\chi} H_j^{-1}\right) =$$

$$= \partial \log(\chi) \cdot I_r + \partial H_j \cdot H_j^{-1} = \partial \log(\chi) \cdot I_r + {}^T\Gamma_j.$$

Hence

$$D_i^\chi = D_i + \partial \log(\chi) \cdot id_{E_i},$$

so that

$$R_i^\chi = R_i + \overline{\partial}\partial \log(\chi) \cdot id_{E_i},$$

where $R_i$ (resp. $R_i^\chi$) is the curvature of $D_i$ (resp. of $D_i^\chi$). It follows that

$$R_{D^\chi} = R_D + \overline{\partial}\partial \log(\chi) \cdot id_E,$$

so that

$$K_g(E, h^\chi) = i\Lambda_g R_{D^\chi} = i\Lambda_g(R_D) + i\Lambda_g \overline{\partial}\partial \log(\chi) \cdot id_E =$$

$$= K_g(E, h) + i\Lambda_g(\overline{\partial}\partial \log(\chi)) \cdot id_E = (\varphi + i\Lambda_g(\overline{\partial}\partial \log(\chi))) \cdot id_E,$$

concluding the proof. $\square$

Consequence of this is the following:

**Proposition 3.19.** *Let $X$ be a compact complex manifold and $g$ a Kähler metric on $X$. Let $E$ be an $\alpha-$twisted holomorphic vector bundle on $X$ and $h$ a Hermitian metric on $E$. If $(E, h)$ verifies the weak $g-$Hermite-Einstein condition, then there is a conformal change $h^\chi$ of $h$, which is unique up to homothety, such that $(E, h^\chi)$ verifies the $g-$Hermite-Einstein condition.*

*Proof.* Define
$$c := \frac{\int_X \varphi \sigma_g^n}{\int_X \sigma_g^n},$$
which is well-defined real number since $X$ is compact. Let $\varphi$ be the Einstein function of $(E, h)$: as $g$ is a Kähler metric, we know that there is a $C^\infty$ function $u : X \longrightarrow \mathbb{R}$ such that
$$i\Lambda_g \overline{\partial}\partial(u) = c - \varphi.$$
Let now $\chi := \exp(u)$, which is then a positive $C^\infty$ real function on $X$ such that
$$i\Lambda_g \overline{\partial}\partial \log(\chi) = i\Lambda_g \overline{\partial}\partial(u) = c - \varphi.$$
As $(E, h)$ verifies the weak $g-$Hermite-Einstein condition with Einstein function $\varphi$, by Lemma 3.18 we have that $(E, h^\chi)$ verifies the Hermite-Einstein condition with Einstein function
$$\varphi + i\Lambda_g \overline{\partial}\partial \log(\chi) = \varphi + c - \varphi = c,$$
which is constant, and we are done. $\square$

**Remark 3.20.** A similar proof works without the Kähler assumption on $g$, see Lemma 2.1.5 of [21].

3.6. **Approximate Hermite-Einstein.** Let now $X$ be a compact complex manifold of dimension $n$ and $g$ a Kähler metric on $X$, whose associated Kähler form is $\sigma_g$. Let $E$ be an $\alpha-$twisted holomorphic vector bundle of rank $r$ on $X$, and $h$ a Hermitian metric on $E$.

Recall that $K_g(E, h)$ is a smooth endomorphism of $E$, so its trace
$$T_g(E, h) := Tr(K_g(E, h)),$$
is a smooth function on $X$, called $g-$**scalar curvature** of $(E, h)$. We let
$$c_g(E, h) := \frac{\int_X T_g(E, h) \sigma_g^n}{r \int_X \sigma_g^n},$$
which is the mean value of $\frac{1}{r}T_g(E, h)$ on $X$.

**Lemma 3.21.** *We have*
$$\int_X \gamma_1(E, h) \wedge \sigma_g^{n-1} = \frac{rc_g(E, h)}{2n\pi} \int_X \sigma_g^n.$$

*Proof.* The proof is almost identical to that of Lemma 3.16. We have
$$K_g(E, h) \cdot \sigma_g^n = \sqrt{-1}nR_h \wedge \sigma_g^{n-1}.$$

Taking the trace we then get
$$T_g(E,h)\sigma_g^n = \sqrt{-1}nTr(R_D) \wedge \sigma_g^{n-1} = 2n\pi\gamma_1(E,h) \wedge \sigma_g^{n-1}.$$

Taking the integral over $X$ the statement follows. $\square$

As $g$ is a Kähler metric (or more generally if $\sigma_g^{n-1}$ is $\partial\bar{\partial}$−closed), the value
$$\int_X \gamma_1(E,h) \wedge \sigma_g^{n-1}$$
only depends on $c_1(E)$ and the cohomology class $[\sigma_g]$. Similarily, the value
$$\int_X \sigma_g^n$$
only depends on $[\sigma_g]$. As a consequence, $c_g(E,h)$ only depends on $c_1(E)$ and on $[\sigma_g]$, but not on $h$: we will use the notation $c_g(E)$ for it.

**Remark 3.22.** If $E$ is $g$−Hermite-Einstein, then by Lemmas 3.16 and 3.21 the Einstein factor of $(E,h)$ is $c_g(E)$. In particular, the Einstein factor does not depend on the Hermitian metric.

**Definition 3.23.** *The $g$−**degree** of $E$ is the real number*
$$\deg_g(E) := \int_X \gamma_1(E,h) \wedge \sigma_g^{n-1},$$
*which only depends on $c_1(E)$.*

Lemma 3.21 then reads as
$$\deg_g(E) = \frac{rc_g(E)}{2n\pi} \int_X \sigma_g^n.$$

We now define
$$||K_g(E,h) - c_g(E) \cdot id_E||^2 := Tr((K_g(E,h) - c_g(E) \cdot id_E)^2),$$
which is a smooth function on $X$. This allows us to give the following:

**Definition 3.24.** *We say that an $\alpha$−twisted holomorphic vector bundle $E$ is **approximate $g$−Hermite-Einstein** if for every $\epsilon > 0$ there is a Hermitian metric $h_\epsilon$ on $E$ such that*
$$\max_X ||K_g(E,h_\epsilon) - c_g(E) \cdot id_E||^2 < \epsilon.$$

First, we have the following:

**Proposition 3.25.** *Let $E$ be an $\alpha$−twisted holomorphic vector bundle on $X$. If $E$ is $g$−Hermite-Einstein, then it is approximate $g$−Hermite-Einstein.*

*Proof.* Let $c$ be the Einstein factor of $(E,h)$. By Remark 3.22 we have $c = c_g(E)$. Choose then $h_\epsilon = h$ for every $\epsilon > 0$, so that $K_g(E,h_\epsilon) = K_g(E,h) = c \cdot id_E$, and we are done. $\square$

As for Hermite-Einstein vector bundles, we have some properties about the behavior of approximate $g-$Hermite-Einstein vector bundles with respect to the usual operations, the proof of which is exactly as the one for untwisted vector bundle (we just provide the proof of one of the following result to show the analogies).

**Proposition 3.26.** *Let $E$ be an $\alpha-$twisted holomorphic vector bundle on $X$. If $E$ is approximate $g-$Hermite-Einstein, then $E^*$ is approximate $g-$Hermite-Einstein.*

*Proof.* Let $h$ be a Hermitian metric on $E$. We proved in Proposition 3.9 that $K_g(E^*, h^*) = -K_g(E, h)^*$, so
$$T_g(E^*, h^*) = -T_g(E, h), \quad c_g(E^*) = -c_g(E).$$
But then
$$||K_g(E^*, h^*) - c_g(E^*) \cdot id_E^*||^2 = Tr((K_g(E^*, h^*) - c_g(E^*) \cdot id_E^*)^2) =$$
$$= Tr((-K_g(E, h) + c_g(E))^2) = Tr((K_g(E, h) - c_g(E))^2) =$$
$$= ||K_g(E, h) - c_g(E) \cdot id_E||^2.$$

If $E$ is approximate $g-$Hermite-Einstein, for every $\epsilon > 0$ there is a Hermitian metric $h_\epsilon$ on $E$ such that $\max ||K_g(E, h_\epsilon) - c_g(E) \cdot id_E||^2 < \epsilon$. But then
$$\max_{x \in X} ||K_g(E^*, h_\epsilon^*) - c_g(E^*) \cdot id_{E^*}||^2 = \max_{x \in X} ||K_g(E, h_\epsilon) - c_g(E) \cdot id_E||^2 < \epsilon,$$
and we are done. □

**Proposition 3.27.** *For $i = 1, 2$ let $E_i$ be an $\alpha_i-$twisted holomorphic vector bundle on $X$. If $E_i$ is approximate $g-$Hermite-Einstein, then $E_1 \otimes E_2$ is approximate $g-$Hermite-Einstein.*

**Proposition 3.28.** *Let $E_1$ and $E_2$ be two $\alpha-$twisted holomorphic vector bundles on $X$ of respective ranks $r_1$ and $r_2$. If they are approximate $g-$Hermite-Einstein and we have*
$$\frac{\deg_g(E_1)}{r_1} = \frac{\deg_g(E_2)}{r_2},$$
*then $E_1 \oplus E_2$ is approximate $g-$Hermite-Einstein.*

**Proposition 3.29.** *Let $E$ be an $\alpha-$twisted holomorphic vector bundle on $X$. If $E$ is approximate $g-$Hermite-Einstein, then*
  (1) *for every $p, q \in \mathbb{N}$, $E^{p,q}$ is approximate $g-$Hermite-Einstein;*
  (2) *for every $p \in \mathbb{N}$, $\wedge^p E$ is approximate $g-$Hermite-Einstein.*

**Proposition 3.30.** *Let $f : X \longrightarrow Y$ be an étale covering, and choose a Kähler matric $g$ on $Y$ and the metric $f^*g$ on $X$.*
  (1) *Let $E$ be an $\alpha-$twisted holomorphic vector bundle on $Y$. If it is approximate $g-$Hermite-Einstein, then $f^*E$ is approximate $f^g-$Hermite-Einstein.*

(2) Let $F$ be an $f^*\alpha-$twisted holomorphic vector bundle on $X$. If it is approximate $f^*g-$Hermite-Einstein, then $f_*F$ is approximate $g-$Hermite-Einstein.

The last property we will need is the following:

**Proposition 3.31.** *Let $E_1$ and $E_2$ be two $\alpha-$twisted holomorphic vector bundles on $X$ of respective ranks $r_1$ and $r_2$, and suppose that such that*
$$\frac{\deg_g(E_1)}{r_1} > \frac{\deg_g(E_2)}{r_2}.$$
*If $E_1$ and $E_2$ are approximate $g-$Hermite-Einstein, then every morphism $f \in \mathscr{H}om(E_1, E_2)$ is zero.*

*Proof.* We know that $f$ is a global section of the untwisted vector bundle $E_1^* \otimes E_2$. As $E_1$ and $E_2$ are approximate $g-$Hermite-Einstein, by Propositions 3.26 and 3.27 we have that $E_1^* \otimes E_2$ is approximate $g-$Hermite-Einstein. Now, notice that
$$\deg_g(E_1^* \otimes E_2) = r_1 \deg_g(E_2) - r_2 \deg_g(E_1) < 0.$$
By Proposition 5.6 in Chapter IV of [17] we then know that $E_1^* \otimes E_2$ has no non-zero global sections, and we are done. □

## 4. Semistability for twisted vector bundles

If $A$ is a domain and $M$ is an $A-$module of finite type, the **homological dimension** of $M$ is the length $d$ of a minimal free resolution
$$0 \longrightarrow E_d \longrightarrow E_{d_1} \longrightarrow \cdots \longrightarrow E_0 \longrightarrow M \longrightarrow 0.$$
The homological dimension of $M$ is denoted $dh(M)$, and we have that $dh(M) = d$ if and only if $Tor_d^A(\mathbb{C}, M) \neq 0$ and $Tor_{d+1}^A(\mathbb{C}, M) = 0$, or equivalently if and only if $Tor_d^A(\mathbb{C}, M) \neq 0$ and $Tor_k^A(\mathbb{C}, M) = 0$ for every $k > d$.

If $A = \mathscr{O}_{\mathbb{C}^n, 0}$, then for every $A-$module $M$ we have $dh(M) \leq n$. As a consequence, the same holds for $A = \mathscr{O}_{X,x}$ where $X$ is a complex manifold of dimension $n$ and $x \in X$.

4.1. **Singularities of twisted sheaves.** Let $\mathscr{E} = \{\mathscr{E}_i, \phi_{ij}\}$ be an $\alpha-$twisted coherent sheaf on $X$. If $x \in U_i$, then we may consider the stalk $\mathscr{E}_{i,x}$ of $\mathscr{E}_i$ at $x$, which is a $\mathscr{O}_{X,x}-$module of finite type. If $x \in U_{ij}$, the isomorphism $\phi_{ij} : \mathscr{E}_i \longrightarrow \mathscr{E}_j$ of $\mathscr{O}_{U_{ij}}-$modules induces an isomorphism $\phi_{ij,x} : \mathscr{E}_{i,x} \longrightarrow \mathscr{E}_{j,x}$ of $\mathscr{O}_{X,x}-$modules. The following definition therefore makes sense:

**Definition 4.1.** *We call **stalk of $\mathscr{E}$ at** $x$ the isomorphism class $\mathscr{E}_x$ of the $\mathscr{O}_{X,x}-$module $\mathscr{E}_{i,x}$, where $i \in I$ is such that $x \in U_i$.*

As the homological dimension of a module is invariant under isomorphism, the homological dimension $dh(\mathscr{E}_x)$ makes sense.

4.1.1. *Singularity sets of twisted sheaves.* For $m \in \mathbb{N}$ let
$$S_m(\mathscr{E}) := \{x \in X \,|\, dh(\mathscr{E}_x) \geq n - m\},$$
which is called $m-$**th singularity set of** $\mathscr{E}$, and we clearly have
$$S_0(\mathscr{E}) \subseteq S_1(\mathscr{E}) \subseteq \cdots \subseteq S_n(\mathscr{E}) = X.$$
The subset $S_{n-1}(\mathscr{E})$ is called **singular locus** of $\mathscr{E}$, and we have
$$S_{n-1}(\mathscr{E}) = \{x \in X \,|\, \mathscr{E}_x \text{ is not a free } \mathscr{O}_{X,x} - \text{module}\}.$$

As the function
$$dh(\mathscr{E}_i) : U_i \longrightarrow \mathbb{Z}, \qquad dh(\mathscr{E}_i)(x) := dh(\mathscr{E}_{i,x})$$
is upper semicontinuous, the same holds for
$$dh(\mathscr{E}) : X \longrightarrow \mathbb{Z}, \qquad dh(\mathscr{E})(x) := dh(\mathscr{E}_x),$$
so that $S_m(\mathscr{E})$ is a closed subset of $X$ for every $m$. Scheja's Theorem tells us that the subset
$$S_m(\mathscr{E}_i) := \{x \in U_i \,|\, dh(\mathscr{E}_{i,x}) \geq n - m\}$$
is a closed analytic subset of $U_i$ of dimension at most $m$ for every $i \in I$ and every $m$. As
$$S_m(\mathscr{E}) \cap U_i = S_m(\mathscr{E}_i),$$
it follows that $S_m(\mathscr{E})$ is a closed analytic subset of $X$ of dimension at most $m$ for every $m$. As a corollary we get

**Lemma 4.2.** *Let $\mathscr{E}$ be an $\alpha-$twisted coherent sheaf on $X$ such that for every $x \in X$ there is an open neighborhood $U \subseteq X$ of $x$ and an exact sequence of $\alpha-$twisted coherent sheaves*
$$0 \longrightarrow \mathscr{E}_{|U} \longrightarrow \mathscr{F}_1 \longrightarrow \cdots \mathscr{F}_k \longrightarrow 0,$$
*where $\mathscr{F}_j$ is a locally free $\alpha_{|U}-$twisted sheaf on $U$. Then*
$$\dim(S_m(\mathscr{E})) \leq m - k.$$

*Proof.* If $x \in U_i$, up to restricting $U$ we may suppose $U \subseteq U_i$, so the exact sequence in the statement becomes
$$0 \longrightarrow \mathscr{E}_{i|U} \longrightarrow \mathscr{F}_1 \longrightarrow \cdots \mathscr{F}_k \longrightarrow 0,$$
where $\mathscr{F}_j$ is a locally free sheaf of $\mathscr{O}_{U_i}-$modules. But then we know that
$$\dim(S_m(\mathscr{F}_i)) \leq m - k,$$
and as this holds for every $i \in I$, we conclude. $\square$

By definition of $S_{n-1}(\mathscr{E})$ we see that $\mathscr{E}_{|X \setminus S_{n-1}(\mathscr{E})}$ is a locally free $\alpha-$twisted coherent sheaf, so the rank of $\mathscr{E}_x$ as a $\mathscr{O}_{X,x}-$module is invariant over $X \setminus S_{n-1}(\mathscr{E})$. We will call this number the **rank** of $\mathscr{E}$.

4.1.2. *Torsion of twisted sheaves.* We recall that if $M$ is an $A-$module of finite type, an element $m \in M$ is a torsion element if there is $a \in A$ such that $am = 0$ in $M$. The set

$$T(M) := \{m \in M \,|\, m \text{ is a torsion element}\}$$

is an $A-$submodule of $M$, the torsion submodule. If $T(M) = 0$, we say that $M$ is torsion-free.

We let $M^* := Hom_A(M, A)$. We then have a natural morphism of $A-$modules

$$\sigma_M : M \longrightarrow M^{**}$$

whose kernel is $T(M)$, so that $M$ is torsion-free if and only if $\sigma_M$ is injective. We say that $M$ is reflexive if $\sigma_M$ is an isomorphism of $A-$modules.

We recall that
- every free $A-$module is reflexive, and every reflexive $A-$module is torsion free;
- every torsion-free $A-$module is a submodule of a free $A-$module of the same rank;
- for every $A-$module $M$ we have that $M^*$ is torsion-free and that $M^{**}$ is reflexive;
- an $A-$module $M$ is torsion-free if and only if $dh(M) \leq n - 1$, and it is reflexive if and only if $dh(M) \leq n - 2$, where $n$ is the Krull dimension of $A$.

If $\mathscr{E}$ is an $\alpha-$twisted coherent sheaf, we let $T(\mathscr{E})$ be the torsion subsheaf of $\mathscr{E}$, i. e. $T(\mathscr{E}) = \{T(\mathscr{E}_i), \phi_{ij|T(\mathscr{E}_i)}\}$, where $T(\mathscr{E}_i)$ is the torsion subsheaf of $\mathscr{E}_i$. It is easy to see that $T(\mathscr{E})$ is an $\alpha-$twisted coherent subsheaf of $\mathscr{E}$, called **torsion of $\mathscr{E}$**.

**Definition 4.3.** *An $\alpha-$twisted coherent sheaf $\mathscr{E}$ is **torsion-free** if $\mathscr{E}_x$ is a torsion-free $\mathcal{O}_{X,x}-$module for every $x \in X$.*

This is equivalent to $T(\mathscr{E}) = 0$, or even to $\mathscr{E}_i$ torsion free for every $i \in I$. As a consequence of the previous properties of $A-$modules, it follows that:
- a locally free $\alpha-$twisted coherent sheaf is torsion-free,
- every $\alpha-$twisted coherent subsheaf of a torsion-free $\alpha-$twisted coherent sheaf if torsion-free,
- every torsion-free $\alpha-$twisted coherent sheaf is an $\alpha-$twisted subsheaf of a locally free $\alpha-$twisted coherent sheaf of the same rank,
- if $\mathscr{E}$ is a torsion-free $\alpha-$twisted coherent sheaf, then for every $x \in X$ we have $dh(\mathscr{E}_x) \leq n-1$, so by Lemma 4.2 we get that $\dim(S_m(\mathscr{E})) \leq m - 1$. In particular, we see that the singular locus of a torsion-free $\alpha-$twisted coherent sheaf has codimension at least 2,
- The natural morphism $\sigma_{\mathscr{E}} : \mathscr{E} \longrightarrow \mathscr{E}^{**}$ has kernel equal to $T(\mathscr{E})$, hence $\mathscr{E}$ is torsion-free if and only if $\sigma_{\mathscr{E}}$ is injective.

**Definition 4.4.** *An $\alpha-$twisted coherent sheaf $\mathscr{E}$ is **reflexive** if $\sigma_{\mathscr{E}}$ is an isomorphism.*

This is equivalent to $\mathscr{E}_i$ reflexive for every $i \in I$. Again, as before we have that:

- all locally-free $\alpha-$twisted coherent sheaves are reflexive,
- all reflexive $\alpha-$twisted coherent sheaves are torsion-free,
- for every $\alpha-$twisted coherent sheaf $\mathscr{E}$, we have that $\mathscr{E}^*$ is reflexive and $\mathscr{E}^{**}$ is torsion-free,
- if $\mathscr{F}$ is reflexive, then $\dim(S_m(\mathscr{F})) \leq m - 2$ for every $m$. In particular, the singular locus of a reflexive $\alpha-$twisted coherent sheaf has codimension at least 3 in $X$,
- if $\mathscr{E}$ and $\mathscr{F}$ are two $\alpha-$twisted coherent sheaves and $\mathscr{F}$ is reflexive, then the coherent sheaf $\mathscr{H}om(\mathscr{E}, \mathscr{F})$ is reflexive.

We end this section with the following:

**Definition 4.5.** *An $\alpha-$twisted coherent sheaf $\mathscr{E} = \{\mathscr{E}_i, \phi_{ij}\}$ is **normal** if $\mathscr{E}_i$ is normal for every $i \in I$, i. e. if for every $i \in I$, every open $U \subseteq U_i$ and every closed analytic subset $A \subseteq U$ of codimension at least 2, the restriction morphism $\mathscr{E}_i(U) \longrightarrow \mathscr{E}_i(U \setminus A)$ is an isomorphism.*

Normal $\alpha-$twisted coherent sheaves are useful for the following property, which is an immediate consequence of the untwisted analogue:

**Lemma 4.6.** *An $\alpha-$twisted coherent sheaf is reflexive if and only if it is torsion-free and normal. Moreover, if*

$$0 \longrightarrow \mathscr{E} \longrightarrow \mathscr{F} \longrightarrow \mathscr{G} \longrightarrow 0$$

*is an exact sequence of $\alpha-$twisted coherent sheaves, if $\mathscr{F}$ is reflexive and $\mathscr{G}$ is torsion-free, then $\mathscr{E}$ is normal (and hence reflexive).*

**4.2. Semistability for twisted sheaves.** Let $\mathscr{E}$ be an $\alpha-$twisted coherent sheaf on a complex manifold $X$. We know that $\mathscr{E}$ admits a finite resolution

$$0 \longrightarrow \mathscr{E}_m \longrightarrow \cdots \longrightarrow \mathscr{E}_0 \longrightarrow \mathscr{F} \longrightarrow 0$$

where $\mathscr{E}_j$ is a locally free $\alpha-$twisted sheaf of rank $r_j$ (see [3]). If $E_j$ is $\alpha-$twisted holomorphic vector bundle associated to $\mathscr{E}_j$, then $\det(E_j)$ is an $\alpha^{r_j}-$twisted holomorphic line bundle, whose associated locally free $\alpha^{r_j}-$twisted coherent sheaf is denoted $\det(\mathscr{E}_j)$ and called determinant of $\mathscr{E}_j$.

We then let

$$\det(\mathscr{E}) := \otimes_{j=0}^m \det(\mathscr{E}_j)^{(-1)^j},$$

which is a locally free sheaf of rank 1 twisted by $\alpha^r$, where $r$ is the rank of $\mathscr{E}$. It is called **determinant** of $\mathscr{E}$, and it does not depend on the chosen resolution.

If $\mathscr{E}$ is a torsion-free $\alpha-$twisted coherent sheaf, then $\det(\mathscr{E})$ is canonically isomorphic to $(\wedge^r \mathscr{E})^{**}$, and we have $\det(\mathscr{E})^* \simeq \det(\mathscr{E}^*)$. Moreover, if $f : \mathscr{E} \longrightarrow \mathscr{F}$ is a morphism of $\alpha-$twisted coherent sheaves, where $\mathscr{E}$ and $\mathscr{F}$ are both torsion free and have the same rank, then $f$ induces a morphism

$$\det(f) : \det(\mathscr{E}) \longrightarrow \det(\mathscr{F}).$$

Let now $\mathscr{E}$ be an $\alpha$-twisted coherent sheaf of rank $r$ and suppose that $X$ is a compact Kähler manifold of dimension $n$. As $\det(\mathscr{E})$ is a locally free $\alpha^r$-twisted sheaf of rank 1, we may associate to it an $\alpha^r$-twisted holomorphic line bundle $L$. Choose a Hermitian metric $h$ on $L$: we will let

$$\gamma_1(\mathscr{E},h) := \gamma_1(L,h),$$

and call it the **first Chern form of** $(\mathscr{E},h)$.

As $g$ is Kähler metric (or more generally if $\sigma_g^{n-1}$ is $\partial\bar{\partial}$-closed), the cohomology class of $\gamma_1(\mathscr{E},h)$ does not depend on $h$, so we write it as

$$c_1(\mathscr{E}) \in H^{1,1}(X) \cap H^2(X,\mathbb{Z}),$$

and call it **first Chern class of** $\mathscr{E}$. The $g$-**degree of** $\mathscr{E}$ is

$$\deg_g(\mathscr{E}) := \int_X \gamma_1(\mathscr{E},h) \wedge \sigma_g^{n-1},$$

which again does not depend on $h$.

If $r > 0$ is the rank of $\mathscr{E}$, the **slope of** $\mathscr{E}$ **with respect to** $g$ is

$$\mu_g(\mathscr{E}) = \frac{\deg_g(\mathscr{E})}{r}.$$

**Definition 4.7.** *A torsion-free $\alpha$-twisted coherent sheaf $\mathscr{E}$ is said to be $g$-**semistable** if for every $\alpha$-twisted coherent subsheaf $\mathscr{F}$ of $\mathscr{E}$ whose rank $r'$ is such that $0 < r' < r$, we have that $\mu_g(\mathscr{F}) \leq \mu_g(\mathscr{E})$. If for every $\alpha$-twisted subsheaf the inaquality is strict, we say that $\mathscr{E}$ is $g$-**stable**. An $\alpha$-twisted holomorphic vector bundle $E$ is $g$-**semistable** (resp. $g$-**stable**) if the associated locally free $\alpha$-twisted coherent sheaf is $g$-semistable (resp. $g$-stable).*

Let us now collect some properties of semistable twisted sheaves that will be useful in what follows.

**Lemma 4.8.** *Let*
$$0 \longrightarrow \mathscr{E} \longrightarrow \mathscr{F} \longrightarrow \mathscr{G} \longrightarrow 0$$
*be an exact sequence of $\alpha$-twisted coherent sheaves on a compact complex manifold $X$ with a Kähler metric $g$. Then*
$$r'(\mu_g(\mathscr{F}) - \mu_g(\mathscr{E})) + r''(\mu_g(\mathscr{F}) - \mu_g(\mathscr{G})) = 0,$$
*where $r'$ is the rank of $\mathscr{E}$ and $r''$ is the rank of $\mathscr{G}$.*

*Proof.* If $r$ is the rank of $\mathscr{F}$, we have that $\det(\mathscr{F})$ is $\alpha^r$-twisted, $\det(\mathscr{E})$ is $\alpha^{r'}$-twisted and $\det(\mathscr{G})$ is $\alpha^{r''}$-twisted, so that $\det(\mathscr{E}) \otimes \det(\mathscr{G})$ is $\alpha^r$-twisted (since $r' + r'' = r$). As the sequence is exact, we have an isomorphism $\det(\mathscr{E}) \simeq \det(\mathscr{E}) \otimes \det(\mathscr{G})$, and hence $c_1(\mathscr{F}) = c_1(\mathscr{E}) + c_1(\mathscr{G})$. But then $\deg_g(\mathscr{F}) = \deg_g(\mathscr{E}) + \deg_g(\mathscr{G})$, hence
$$(r' + r'')\mu_g(\mathscr{F}) = r\mu_g(\mathscr{F}) = r'\mu_g(\mathscr{E}) + r''\mu_g(\mathscr{G}),$$
and the statement follows. □

As an immediate consequence we have the following:

**Proposition 4.9.** *Let $\mathscr{E}$ be an torsion-free $\alpha$−twisted coherent sheaf of rank $r$ on a complex manifold $X$ with a Kähler metric $g$.*

(1) *The sheaf $\mathscr{E}$ is $g$−semistable if and only if for every $\alpha$−twisted quotient $\mathscr{G}$ of $\mathscr{F}$ of rank $0 < r''$ we have $\mu_g(\mathscr{F}) \leq \mu_g(\mathscr{G})$.*
(2) *The sheaf $\mathscr{E}$ is $g$−stable if and only if for every $\alpha$−twisted quotient $\mathscr{G}$ of $\mathscr{F}$ of rank $0 < r'' < r$ we have $\mu_g(\mathscr{F}) < \mu_g(\mathscr{G})$.*

We now have the following property of torsion twisted sheaves.

**Lemma 4.10.** *If $\mathscr{E}$ is a torsion $\alpha$−twisted coherent sheaf on a compact complex manifold $X$ with a Kähler metric $g$, then $\deg_g(\mathscr{E}) \geq 0$.*

*Proof.* We first show that if $\mathscr{E}$ and $\mathscr{E}'$ are two $\alpha$−twisted torsion-free coherent sheaves of the same rank $r$ and $f : \mathscr{E} \longrightarrow \mathscr{E}'$ is injective, then the induced morphism $\det(f) : \det(\mathscr{E}) \longrightarrow \det(\mathscr{E}')$ is injective. Indeed, let $A = S_{n-1}(\mathscr{E})$ and $A' = S_{n-1}(\mathscr{E}')$: over $X \setminus (A \cup A')$ both $\mathscr{E}$ and $\mathscr{E}'$ are locally free, and $f$ is an injection between two locally free $\alpha$−twisted sheaves of the same rank.

The induced morphism $\det(f)$ is then an isomorphism over $X \setminus (A \cup A')$, so $\ker(\det(f))$ is supported on $A \cup A'$, and hence it is a torsion $\alpha^r$−twisted sheaf. But moreover $\ker(\det(f))$ is an $\alpha$−twisted coherent subsheaf of $\det(\mathscr{E})$, which is torsion-free, so $\ker(f) = 0$ and $f$ is injective.

Now, as by hypothesis $\mathscr{E}$ is torsion, over $X \setminus A$ it is locally free and torsion, so it is trivial. It follows that $\mathscr{E}$ is supported on $A$, and it has rank 0. Let

$$0 \longrightarrow \mathscr{E}_m \xrightarrow{f_m} \cdots \xrightarrow{f_1} \mathscr{E}_0 \xrightarrow{f_0} \mathscr{E} \longrightarrow 0$$

be a locally free resolution of $\mathscr{E}$. We let $S_1 := \mathscr{E}_1 / \ker(f_1)$, which is a torsion-free $\alpha$−twisted coherent sheaf, and we have an exact sequence

$$0 \longrightarrow S_1 \xrightarrow{f_1} \mathscr{E}_0 \xrightarrow{f_0} \mathscr{E} \longrightarrow 0.$$

As the rank of $\mathscr{E}$ is 0, $\det(\mathscr{E})$ is an untwisted locally free coherent sheaf of rank 1. Moreover $S_1$ and $\mathscr{E}_0$ are both torsion-free $\alpha$−twisted coherent sheaves of the same rank $r$. As $f_1$ is injective, the previous proof gives that $\det(f_1) : \det(S_1) \longrightarrow \det(\mathscr{E}_0)$ is injective, i. e. it is a non-trivial morphism of $\alpha^r$−twisted coherent sheaves. Since

$$\det(\mathscr{E}) \simeq \det(\mathscr{E}_0) \otimes \det(S_1)^* \simeq \mathscr{H}om(\det(S_1), \det(\mathscr{E}_0)),$$

we then have that

$$\Gamma(X, \det(\mathscr{E})) \simeq Hom(\det(S_1), \det(\mathscr{E}_0)).$$

But then $\det(f_1)$ corresponds to a non-trivial holomorphic section of $\det(\mathscr{E})$.

Let $V$ be the zero locus of this section, which is a divisor on $X$ whose cohomology class is $c_1(\det(\mathscr{E})) = c_1(\mathscr{E})$. It then follows that

$$\deg_g(\mathscr{E}) = \int_X c_1(\mathscr{E}) \wedge \sigma_g^{n-1} = \int_V \sigma_g^{n-1} \geq 0,$$

and we are done. $\square$

**Remark 4.11.** The previous proof works even without the Kähler assumption on $g$, see Proposition 1.3.5 of [21].

As a consequence of this we get the following:

**Proposition 4.12.** *Let $\mathcal{E}$ be an torsion-free $\alpha-$twisted coherent sheaf of rank $r$ on a complex manifold $X$ with a Kähler metric $g$.*
  (1) *The sheaf $\mathcal{E}$ is $g-$semistable if and only if one of the two following conditions is verified:*
    (a) *for every $\alpha-$twisted subsheaf $\mathcal{F}$ of $\mathcal{E}$ such that $\mathcal{E}/\mathcal{F}$ is torsion-free, we have $\mu_g(\mathcal{F}) \leq \mu_g(\mathcal{E})$.*
    (b) *for every torsion-free $\alpha-$twisted quotient $\mathcal{G}$ of $\mathcal{E}$ we have $\mu_g(\mathcal{E}) \leq \mu_g(\mathcal{G})$.*
  (2) *The sheaf $\mathcal{E}$ is $g-$stable if and only if one of the two following conditions is verified:*
    (a) *for every $\alpha-$twisted subsheaf $\mathcal{F}$ of $\mathcal{E}$ with $\mathcal{E}/\mathcal{F}$ torsion-free and such that $\mathcal{F}$ has rank $0 < r' < r$, we have $\mu_g(\mathcal{F}) \leq \mu_g(\mathcal{E})$.*
    (b) *for every torsion-free $\alpha-$twisted quotient $\mathcal{G}$ of $\mathcal{E}$ of rank $0 < r'' < r$ we have $\mu_g(\mathcal{E}) \leq \mu_g(\mathcal{G})$.*

*Proof.* The proof is identical to that of Proposition 7.6 in Chapter V of [17]. □

The following will be useful in what follows:

**Proposition 4.13.** *Let $\mathcal{E}$ be an torsion-free $\alpha-$twisted coherent sheaf of rank $r$ on a complex manifold $X$ with a Kähler metric $g$.*
  (1) *If $r = 1$, then $\mathcal{E}$ is $g-$stable.*
  (2) *If $\mathcal{L}$ is a locally free $\alpha'-$twisted sheaf of rank 1, then the $\alpha-$twisted sheaf $\mathcal{E}$ is $g-$(semi)stable if and only if $\mathcal{F} \otimes \mathcal{L}$ is $g-$(semi)stable.*
  (3) *$\mathcal{E}$ is $g-$(semi)stable if and only if $\mathcal{E}^*$ is $g-$(semi)stable.*

*Proof.* The proof is identical to that of Proposition 7.7 in Chapter V of [17]. □

The last results we need, whose proofs are identical to those of the corresponding results in the untwisted cases (see Propositions 5.7.9 and 5.7.11, and Corollaries 5.7.12 and 5.7.14), are the following:

**Proposition 4.14.** *Let $\mathcal{E}$ and $\mathcal{F}$ be two torsion-free $\alpha-$twisted coherent sheaves on a complex manifold $X$ with a Kähler metric $g$. Then $\mathcal{E} \oplus \mathcal{F}$ is $g-$semistable if and only if $\mathcal{E}$ and $\mathcal{F}$ are $g-$semistable and $\mu_g(\mathcal{E}) = \mu_g(\mathcal{F})$.*

**Proposition 4.15.** *Let $\mathcal{E}$ and $\mathcal{F}$ be two torsion-free $\alpha-$twisted coherent sheaves on a complex manifold $X$ with a Kähler metric $g$, and let $f : \mathcal{E} \longrightarrow \mathcal{F}$ be a morphism of $\alpha-$twisted coherent sheaves.*
  (1) *If $\mu_g(\mathcal{E}) > \mu_g(\mathcal{F})$, then $f = 0$,*
  (2) *If $\mu_g(\mathcal{E}) = \mu_g(\mathcal{F})$ and $\mathcal{E}$ is $g-$stable, then either $f = 0$ or $\mathcal{E}$ and $\mathcal{F}$ have the same rank and $f$ is injective.*

(3) If $\mu_g(\mathscr{E}) = \mu_g(\mathscr{F})$ and $\mathscr{F}$ is $g-$stable, then either $f = 0$ or $\mathscr{E}$ and $\mathscr{F}$ have the same rank and $f$ is generically surjective.

This Proposition has two important corollaries:

**Corollary 4.16.** *Let $\mathscr{E}$ and $\mathscr{F}$ be two torsion-free $\alpha-$twisted coherent sheaves on a complex manifold $X$ with a Kähler metric $g$. If they have the same rank and same degree with respect to $g$, and one of them is $g-$stable, then every morphism of $\alpha-$twisted coherent sheaves between them is either 0 or an isomorphism.*

**Corollary 4.17.** *Let $E$ be an $\alpha-$twisted holomorphic vector bundle on a compact complex manifold $X$ with a Kähler metric $g$. If $E$ is $g-$stable, then $\Gamma(X, \underline{End}(E)) \simeq \mathbb{C} \cdot id_E$.*

4.3. **Hermite-Einstein implies polystable.** Let $X$ be a complex manifold of dimension $n$ and $g$ a Kähler metric on $X$. Let $\mathscr{E}$ be an $\alpha-$twisted coherent sheaf on $X$ of rank $r$. Let $L$ be the $\alpha^r-$twisted holomorphic line bundle associated to $\det(\mathscr{E})$, and choose a Hermitian metric $h$ on $L$.

We call $g-$**degree form of** $(\mathscr{E}, h)$ the $d-$closed real $2n-$form on $X$ defined as
$$d_g(\mathscr{E}, h) := \gamma_1(\mathscr{E}, h) \wedge \sigma_g^{n-1}.$$
If $X$ is compact, we then have
$$\deg_g(\mathscr{E}) = \int_X d_g(\mathscr{E}, h).$$

The first result we need to prove is the following, relating the properties of the $g-$degree form and the $g-$Hermite-Einstein condition.

**Lemma 4.18.** *Let $E$ be an $\alpha-$twisted holomorphic vector bundle on a compact complex manifold $X$ with Kähler metric $g$. Let $h$ be a Hermitian metric on $E$ and suppose that $(E, h)$ verifies the $g-$Hermite-Einstein condition with Einstein factor $c$. Let*
$$0 \longrightarrow E' \longrightarrow E \longrightarrow E'' \longrightarrow 0$$
*be an exact sequence of $\alpha-$twisted holomorphic vector bundles.*

(1) *If $r'$ is the rank of $E'$, $r$ is the rank of $E$ and $h' := h^{E'}$ is the Hermitian metric induced by $h$ on $E'$, the real $d-$closed $2n-$form*
$$\frac{d_g(E', h')}{r'} - \frac{d_g(E, h)}{r}$$
*is everywhere negative on $X$.*

(2) *If the previous $2n-$form is 0, then the exact sequence above splits, and if we let $h''$ be the natural Hermitian metric induced by $h$ on $E''$, then $(E', h')$ and $(E'', h'')$ verify the $g-$Hermite-Einstein condition with Einstein factor $c$.*

*Proof.* The proof is identical to that of Proposition 8.2 in Chapter V of [17], using the results in section 2.4. The same proof works even only if $\sigma_g^{n-1}$ is $\partial\bar{\partial}-$closed (see Proposition 2.3.1 of [21]). □

We are now in the position to prove one of the main results of this section, namely:

**Theorem 4.19.** *Let $E$ be an $\alpha-$twisted holomorphic vector bundle on a compact complex manifold $X$ with Kähler metric $g$. Let $h$ be a Hermitian metric on $E$ and suppose that $(E, h)$ verifies the $g-$Hermite-Einstein condition with Einstein factor $c$. Then $E$ is $g-$semistable, and we have $E = E_1 \oplus \cdots \oplus E_k$ where*

(1) *$E_1, \cdots, E_k$ are $g-$stable $\alpha-$twisted holomorphic subbundles of $E$,*
(2) *if $h^{E_j}$ is the Hermitian metric induced by $h$ on $E_j$, then $(E_j, h^{E_j})$ verifies the $g-$Hermite-Einstein condition with Einstein factor $c$.*

*Proof.* Let us first prove that $E$ is $g-$semistable. To do so, let us consider the locally free $\alpha-$twisted coherent sheaf $\mathscr{E}$ associated to $E$, and let $r$ be its rank: we need to show that $\mathscr{E}$ is $g-$semistable, so let us consider an $\alpha-$twisted coherent subsheaf $\mathscr{F}$ of $\mathscr{E}$ and let $p$ be its rank. By Proposition 4.12, in order to prove that $\mathscr{E}$ is $g-$semistable we may suppose that $0 < p < r$ and that $\mathscr{E}/\mathscr{F}$ is torsion-free.

Let $\iota : \mathscr{F} \longrightarrow \mathscr{E}$ be the inclusion, so taking the determinant we get

$$\det(\iota) : \det(\mathscr{F}) \longrightarrow (\wedge^p \mathscr{E})^{**},$$

since $\det(\mathscr{F}) \simeq (\wedge^p \mathscr{F})^{**}$. As $\mathscr{E}$ is locally free, we have that $\wedge^p \mathscr{E}$ is locally free, and hence that $(\wedge^p \mathscr{E})^{**} \simeq \wedge^p \mathscr{E}$.

The first part of the proof of Lemma 4.10 shows that $\det(\iota)$ is injective. Moreover, notice that $(\wedge^p \mathscr{E}) \otimes \det(\mathscr{F})^*$ is an untwisted sheaf, so tensoring with $\det(\mathscr{F})^*$ we get that $\det(\iota)$ defines a global section

$$f : \mathcal{O}_X \longrightarrow (\wedge^p \mathscr{E}) \otimes \det(\mathscr{F})^*,$$

which is not trivial since $\iota$ is injective. This corresponds to finding a holomorphic section of the holomorphic vector bundle $\wedge^p E \otimes L^*$, where $L$ is the $\alpha^p-$twisted holomorphic line bundle corresponding to $\det(\mathscr{F})$.

As $(E, h)$ verifies the $g-$Hermite-Einstein condition with Einstein factor $c$, by Remark 3.22 we get that $c = c_g(E)$, i. e.

$$c = \frac{2\pi n \mu_g(E)}{\int_X \sigma_g^n}.$$

By Proposition 3.12 we know that $(\wedge^p E, \wedge^p h)$ verifies the $g-$Hermite-Einstein condition with Einstein factor $pc$.

Let $h'$ be a Hermitian metric on $L$. By Propositions 3.8 and 3.19 up to apply a conformal change to $h'$ we know that $(L, h')$ verifies the $g-$Hermite-Einstein condition with Einstein factor $c'$. Again by Remark 3.22 we have

$$c' = \frac{2\pi n \mu_g(L)}{\int_X \sigma_g^n}.$$

Notice that
$$\mu_g(L) = \mu_g(\det(\mathscr{F})) = \deg_g(\det(\mathscr{F})) = \deg_g(\mathscr{F}) = p\mu_g(\mathscr{F}),$$
so
$$c' = \frac{2\pi n p \mu_g(\mathscr{F})}{\int_X \sigma_g^n}.$$

By Proposition 3.9 we then see that $(L^*, (h')^*)$ verifies the $g$−Hermite-Einstein condition with Einstein factor $-c'$, and hence by Proposition 3.10 we see that $(\wedge^p E \otimes L^*, \wedge^p h \otimes (h')^*)$ verifies the $g$−Hermite-Einstein condition with Einstein factor $pc - c'$.

But $\wedge^p E \otimes L^*$ is a holomorphic vector bundle on a compact Kähler manifold which has a non-trivial holomorphic section. Proposition 3.15 gives then $pc - c' \geq 0$, i. e.
$$\frac{2\pi n p \mu_g(E)}{\int_X \sigma_g^n} \geq \frac{2\pi n p \mu_g(L)}{\int_X \sigma_g^n},$$
so that $\mu_g(\mathscr{E}) \geq \mu_g(\mathscr{F})$, showing that $\mathscr{E}$ is $g$−semistable.

Let us now prove the second part of the statement. First, if $E$ is $g$−stable, then we let $k = 1$ and $E = E_1$, and we are done. We will then suppose that $E$ is $g$−semistable but not $g$−stable.

This implies that there is an $\alpha$−twisted coherent subsheaf $\mathscr{F}$ of $\mathscr{E}$ of rank $0 < p < r$ with $\mathscr{G} := \mathscr{E}/\mathscr{F}$ torsion-free and $\mu_g(\mathscr{F}) = \mu_g(\mathscr{E})$. Letting $L$ be the $\alpha^p$−twisted holomorphic vector bundle associated to $\det(\mathscr{F})$ and $h'$ a Hermitian metric on $L$, by the previous part of the proof $(\wedge^p E \otimes L^*, \wedge^p h \otimes (h')^*)$ verifies the $g$−Hermite-Einstein condition with Einstein factor 0.

By Proposition 3.15 the mean curvature of $(\wedge^p E \otimes L^*, \wedge^p h \otimes (h')^*)$ is negative semidefinite and every holomorphic section of $\wedge^p E \otimes L^*$ is parallel with respect to the Chern connection $D$. In particular the morphism $f$ defined before has to be parallel with respect to $D$, so $L$ is an $\alpha^p$−twisted holomorphic line subbundle of $\wedge^p E$ which is parallel with respect to the Chern connection of $(\wedge^p E, \wedge^p h)$.

Now, let $X' := X \setminus S_{n-1}(\mathscr{F})$ be the open subset of $X$ over which $\mathscr{F}$ is locally free, and let $F'$ be the $\alpha$−twisted holomorphic vector bundle associated to $\mathscr{F}_{|X'}$. The inclusion $\iota : \mathscr{F} \longrightarrow \mathscr{E}$ induces an inclusion $\iota' : \mathscr{F}_{|X'} \longrightarrow \mathscr{E}_{|X'}$, and hence an inclusion $j' : F' \longrightarrow E_{|X'}$ of $\alpha_{|X'}$−twisted holomorphic vector bundles. The previous discussion shows that $F'$ is a parallel twisted holomorphic subbundle of $E$ with respect to the Chern connection of $(E, h)$.

By Lemma 2.42 we then get an $\alpha_{|X'}$−holomorphic subbundle $G'$ of $E_{|X'}$ such that $E_{|X'} = F' \oplus G'$. Notice that $\mathscr{G}_{|X'}$ is a locally free $\alpha_{|X'}$−twisted sheaf on $X'$ whose associated vector bundle is $G'$. As the exact sequence
$$0 \longrightarrow F' \longrightarrow E_{|X'} \longrightarrow G' \longrightarrow 0$$
is splitting, the same holds for the exact sequence
$$0 \longrightarrow \mathscr{F}_{|X'} \longrightarrow \mathscr{E}_{|X'} \longrightarrow \mathscr{G}_{|X'} \longrightarrow 0.$$

As $\mathscr{E}$ is locally free, we have that $\mathscr{F}$ is reflexive (as a subsheaf of a reflexive sheaf), and we know that $\mathscr{G}$ is torsion-free by assumption. By Lemma 4.6 we then get that $\mathscr{H}om(\mathscr{E},\mathscr{F})$ and $\mathscr{H}om(\mathscr{F},\mathscr{F})$ are reflexive, and hence normal. This implies that
$$\Gamma(X, \mathscr{H}om(\mathscr{E},\mathscr{F})) = \Gamma(X', \mathscr{H}om(\mathscr{E},\mathscr{F})),$$
and
$$\Gamma(X, \mathscr{H}om(\mathscr{F},\mathscr{F})) = \Gamma(X', \mathscr{H}om(\mathscr{F},\mathscr{F})),$$
since $S_{n-1}(\mathscr{F})$ is a closed analytic subset of codimension at least 2 (since $\mathscr{F}$ is reflexive).

As the exact sequence
$$0 \longrightarrow \mathscr{F}_{|X'} \longrightarrow \mathscr{E}_{|X'} \longrightarrow \mathscr{G}_{|X'} \longrightarrow 0$$
is splitting, there is a splitting morphism $p' : \mathscr{E}_{|X'} \longrightarrow \mathscr{F}_{|X'}$, i. e. such that $p' \circ \iota' = id_{\mathscr{F}_{|X'}}$. The previous property then implies that there is a unique morphism $p : \mathscr{E} \longrightarrow \mathscr{F}$ such that $p \circ \iota = id_{\mathscr{F}}$, i. e. the exact sequence
$$0 \longrightarrow \mathscr{F} \longrightarrow \mathscr{E} \longrightarrow \mathscr{G} \longrightarrow 0$$
splits too, and we have $\mathscr{E} = \mathscr{F} \oplus \mathscr{G}$. But as $\mathscr{E}$ is locally free, it follows that $\mathscr{F}$ and $\mathscr{G}$ are locally free, and hence we have $E = F \oplus G$ where $F$ and $G$ are the $\alpha-$twisted holomorphic vector bundles corresponding to $\mathscr{F}$ and $\mathscr{G}$. By Lemma 2.42 we then may proceed by induction on the rank, getting the statement. □

Let us now present the following definition:

**Definition 4.20.** *An $\alpha-$twisted holomorphic vector bundle $E$ is called $g-$**polystable** if it is $g-$semistable and we have $E = E_1 \oplus \cdots \oplus E_k$ where $E_1, \cdots, E_k$ are $g-$stable $\alpha-$twisted holomorphic subbundles of $E$ such that $\mu_g(E_j) = \mu_g(E)$ for every $j = 1, \cdots k$.*

This definition allows us to rephrase Theorem 4.19 as follows, proving one direction of the first part of Theorem 1.1:

**Theorem 4.21.** *Let $X$ be a compact Kähler manifold with a Kähler metric $g$ and $E$ an $\alpha-$twisted holomorphic vector bundle on $X$. If $E$ is $g-$Hermite-Einstein, then $E$ is $g-$polystable.*

**Remark 4.22.** The proof of Theorem 4.19 works more generally if we suppose $g$ to be such that $\sigma_g^{n-1}$ is $\partial\bar{\partial}-$closed (see Theorem 2.3.2 of [21]).

4.4. **Approximate Hermite-Einstein implies semistable.** An adaptation of the proof of Theorem 4.19 allows us to prove the following, which is one direction of the second item of Theorem 1.1:

**Theorem 4.23.** *Let $E$ be an $\alpha-$twisted holomorphic vector bundle on a compact complex manifold $X$ with Kähler metric $g$. If $E$ is approximate $g-$Hermite-Einstein, then $E$ is $g-$semistable.*

*Proof.* Consider an $\alpha$–twisted coherent subsheaf $\mathscr{F}$ of the locally free $\alpha$–twisted coherent sheaf $\mathscr{E}$ associated to $E$, such that $\mathscr{E}/\mathscr{F}$ is torsion-free and such that $\mu_g(\mathscr{F}) = \mu_g(\mathscr{E})$. We let $r$ be the rank of $\mathscr{E}$ and $p$ be the rank of $\mathscr{F}$, and we suppose that $0 < p < r$. Moreover we let $L$ be the $\alpha^p$–twisted holomorphic vector bundle associated to $\det(\mathscr{F})$.

Recall that as $E$ is approximate $g$–Hermite-Einstein, by Proposition 3.29 then $\wedge^p E$ is approximate $g$–Hermite-Einstein. The $\alpha^p$–twisted holomorphic line bundle $L$ is approximate $g$–Hermite-Einstein by Propositions 3.8 and 3.25, so by Proposition 3.26 we see that $L^*$ is approximate $g$–Hermite-Einstein. By Proposition 3.27 we then conclude that $\wedge^p E \otimes L^*$ is approximate $g$–Hermite-Einstein.

As in the proof of Theorem 4.19, we have a nontrivial holomorphic section of $\wedge^p E \otimes L^*$, hence $\deg_g(\wedge^p E \otimes L^*) \geq 0$ by Proposition 3.31. Notice that

$$\deg_g(\wedge^p E \otimes L^*) = \int_X c_1(\wedge^p E \otimes L^*) \wedge \sigma_g^{n-1} =$$

$$= \int_X (c_1(\wedge^p E) - r_p c_1(L)) \wedge \sigma_g^{n-1},$$

where $r_p$ is the rank of $\wedge^p E$, so that

$$0 \leq \deg_g(\wedge^p E \otimes L^*) = \deg_g(\wedge^p E) - r_p \deg_g(L) = \deg_g(\wedge^p E) - r_p \deg_g(\mathscr{F}).$$

But then we get

$$0 \leq \mu_g(\wedge^p E \otimes L^*) = \mu_g(\wedge^p E) - p\mu_g(\mathscr{F}) = p(\mu_g(E) - \mu_g(\mathscr{F})),$$

so that $\mu_g(E) \geq \mu_g(\mathscr{F})$, proving that $E$ is $g$–semistable. $\square$

## 5. The Kobayashi-Hitchin correspondence

The aim of this section is to complete the proof of point 1 of Theorem 1.1, namely:

**Theorem 5.1.** *Let $X$ be a compact Kähler manifold with a Kähler metric $g$, and $E$ an $\alpha$–twisted holomorpic vector bundle on $X$. Then $E$ is $g$–polystable if and only if it is $g$–Hermite-Einstein.*

By Theorem 4.21 we know that if $E$ is $g$–Hermite-Einstein, then $E$ is $g$–polystable. We are left with the proof of the opposite direction, and to do so we follow closely section 3 of [21].

Consider two Hermitian metrics $h_0$ and $h$ on $E$, and fix a Kähler metric $g$ on $X$. We let $D_0$ be the Chern connection of $(E, h_0)$ and $D$ the Chern connection of $(E, h)$. For every $i \in I$ we have

$$D_i^{1,0} = D_{0,i}^{1,0} + (f_i^{h_0,h})^{-1} \circ D_{0,i}^{1,0}(f_i^{h_0,h}),$$

where $D_{0,i}^{1,0}$ denotes the $(1,0)$–part of the connection induced by $D_{0,i}$ on $End(E_i)$ (see as instance section (1.9) of [26]).

If $R_i$ (resp. $R_{0,i}$) is the Chern curvature of $(E_i, h_i)$ (resp. of $(E_i, h_{0,i})$), we then have
$$R_i = \overline{\partial}_i \circ D_i^{1,0} = \overline{\partial}_i \circ D_{0,i}^{1,0} + \overline{\partial}_i((f_i^{h_0,h})^{-1} \circ D_{0,i}^{1,0}(f_i^{h_0,h})) =$$
$$= R_{0,i} + \overline{\partial}_i((f_i^{h_0,h})^{-1} \circ D_{0,i}^{1,0}(f_i^{h_0,h})).$$

If we now let $R$ (resp. $R_0$) be the Chern curvature of $(E, h)$ (resp. of $(E, h_0)$) we then get
$$R_{|U_i} = R_i - B_i \cdot id_{E_i} = R_{0,i} + \overline{\partial}_i((f_i^{h_0,h})^{-1} \circ D_{0,i}^{1,0}(f_i^{h_0,h})) - B_i id_{E_i} =$$
$$= R_{0|U_i} + \overline{\partial}_i((f_i^{h_0,h})^{-1} \circ D_{0,i}^{1,0}(f_i^{h_0,h})).$$

It then follows that $R = R_0 + \overline{\partial}((f^{h_0,h})^{-1} \circ D_0^{1,0}(f^{h_0,h}))$, and hence that
$$K_g(E, h) = K_g(E, h_0) + i\Lambda_g(\overline{\partial}((f^{h_0,h})^{-1} \circ D_0^{1,0}(f^{h_0,h}))).$$

Letting
$$K_g^0(E, h_0) := K_g(E, h_0) - c_g(E) \cdot id_E,$$
we then see that
$$K_g(E, h) - c_g(E) \cdot id_E = K_g^0(E, h_0) + i\Lambda_g(\overline{\partial}((f^{h_0,h})^{-1} \circ D_0^{1,0}(f^{h_0,h}))).$$

In conclusion, we see that $h$ is $g$−Hermite-Einstein if and only if
$$K_g^0(E, h_0) + i\Lambda_g(\overline{\partial}((f^{h_0,h})^{-1} \circ D_0^{1,0}(f^{h_0,h}))) = 0.$$

By Remark 2.58, it follows that $E$ is $g$−Hermite-Einstein if and only if there is $h_0 \in Herm^+(E)$ and $f \in End_{h_0}^+(E)$ such that
$$K_g^0(E, h_0) + i\Lambda_g(\overline{\partial}(f^{-1} \circ D_0^{1,0}(f))) = 0.$$

This will be called **Hermite-Einstein equation**.

Now, recall that if $f \in End_{h_0}^+(E)$ then it makes sense to consider $\log(f)$, which is a $h_0$−Hermitian endomorphism of $E$. For every $\epsilon > 0$ we now let
$$L_\epsilon^{h_0} : End_{h_0}^+(E) \longrightarrow End_{h_0}(E),$$
where
$$L_\epsilon^{h_0}(f) := K_g^0(E, h_0) + i\Lambda_g(\overline{\partial}(f^{-1} \circ D_0^{1,0}(f))) + \epsilon \log(f).$$

We will call $L_\epsilon^{h_0} = 0$ the **perturbed equation**.

5.1. **Solution of the equation $L_1 = 0$.** The first step in the proof of the Kobayashi-Hitchin correspondence is that there is $f \in End_{h_0}^+(E)$ such that $L_1(f) = 0$, i. e. the perturbed equation for $\epsilon = 1$ has a solution.

**Lemma 5.2.** *Let $E$ be an $\alpha$−twisted holomorphic vector bundle on $X$.*

(1) *If $h, h_0 \in Herm^+(E)$, $\epsilon > 0$ and $f \in End_h^+(E)$, then*
$$Tr(L_\epsilon^{h_0}(f)) = Tr(K_g^0(E, h_0)) + i\Lambda_g \overline{\partial}\partial(Tr(\log(f))) + \epsilon Tr(\log(f)).$$

(2) *There is $h_0 \in Herm^+(E)$ such that $Tr(K_g^0(E, h_0)) = 0$ and such that the perturbed equation $L_1^{h_0} = 0$ has a solution.*

(3) If $Tr(K_g^0(E,h)) = 0$ and $L_\epsilon^h(f) = 0$, then $\det(f) = 1$.

*Proof.* It is easy to see that $\exp(Tr(f)) = \det(\exp(f))$ for every endomorphism $f$ of $E$. In particular, if $f \in End_h^+(E)$ it follows that
$$\det(f) = \exp(Tr(\log(f))),$$
and that
$$\partial(\det(f)) = \det(f) \cdot Tr(f^{-1} \circ D_h^{1,0}(f)).$$
It follows that
$$\partial(Tr(\log(f))) = \partial(\log(\det(f))) = Tr(f^{-1} \circ D_h^{1,0}(f)).$$

Now, as in the proof of Lemma (3.2.1) [21] a local calculation shows that
$$Tr(f^{-1} \circ D_0^{1,0}(f)) = Tr(f^{-1} \circ D_h^{1,0}(f)) = \partial(Tr(\log(f))),$$
so that
$$Tr(i\Lambda_g(\overline{\partial}(f^{-1} \circ D_0^{1,0}(f)))) = i\Lambda_g \overline{\partial}\partial(Tr(\log(f))).$$
But then
$$Tr(L_\epsilon^{h_0}(f)) = Tr(K_g^0(E,h_0)) + Tr(i\Lambda_g(\overline{\partial}(f^{-1} \circ D_0^{1,0}(f)))) + Tr(\epsilon \log(f)) =$$
$$= Tr(K_g^0(E,h_0)) + i\Lambda_g \overline{\partial}\partial(Tr(\log(f))) + \epsilon Tr(\log(f)).$$
This completes the proof of the first point of the statement.

For the second point, take any Hermitian metric $h$ on $E$, and notice that if we let $K_g^0(E,h) = K_g(E,h) - c_g(E) \cdot id_E$, then we have
$$\int_X Tr(K_g^0(E,h))\sigma_g^n = 0.$$
It follows that there is a smooth function $\psi$ such that
$$i\Lambda_g \overline{\partial}\partial \psi = K_g^0(E,h),$$
and hence there is a function $\phi$ such that
$$i\Lambda_g \overline{\partial}\partial \phi = -\frac{1}{r}K_g^0(E,h),$$
where $r$ is the rank of $E$.

We let $h_1 := h^{\exp(\phi)}$, which is then a Hermitian metric on $E$ such that
$$K_g(E,h_1) = K_g(E,h) + i\Lambda_g \overline{\partial}\partial(\phi) \cdot id_E$$
(see the proof of Lemma 3.18). We then get
$$Tr(K_g^0(E,h_1)) = Tr(K_g^0(E,h)) + ri\Lambda_g \overline{\partial}\partial(\phi) = 0.$$

By Lemma 3.5 we know that $K_g(E,h_1)$ is a $h_1$-Hermitian endomorphism, so that $K_g^0(E,h_1)$ is too. By Lemma 2.10 we then see that $\exp(K_g(E,h_1)) \in End_{h_1}^+(E)$, hence
$$h_0 := \exp(\widehat{K_g(E,h_1)})_{h_1} \in Herm^+(E).$$
It is easy to see that $Tr(K_g^0(E,h_0)) = 0$. Moreover, let
$$f_1 := \exp(-K_g^0(E,h_1)) \in End_{h_0}^+(E),$$

so that an easy calculation shows that $L_1^{h_0}(f) = 0$.

The final point of the statement is proved exactly as in point (iii) of Lemma (3.2.1) of [21]. □

Now, the same of Lemma (3.2.3) of [21] shows that if $h_0$ is as in point 2 of Lemma 5.2, and if for every $\epsilon > 0$ and $f \in End_{h_0}^+(E)$ we let

$$\widehat{L}^{h_0}(\epsilon, f) := f \circ L_\epsilon^{h_0}(f),$$

then $\widehat{L}^{h_0}(\epsilon, f) \in End_{h_0}(E)$.

By the identification of $Herm(E)$ with $End_{h_0}(E)$, the norm $L_k^p$ on $Herm(E)$ induced a $L_k^p$−norm on $End_{h_0}(E)$. We will write $L_k^p End_{h_0}(E)$ for the Banach space completion of $End_{h_0}(E)$, and $L_k^p End_{h_0}^+(E)$ for the interior of the closure of $End_{h_0}^+(E)$ in $L_k^p End_{h_0}(E)$.

Notice that if $\epsilon > 0$, then $\widehat{L}^{h_0}(\epsilon, \cdot)$ is a second order differential operator. It follows that

$$\widehat{L}^{h_0} : (0,1] \times L_k^p End_{h_0}^+(E) \longrightarrow L_{k-2}^p End_{h_0}(E),$$

and by the Sobolev Multiplication Theorem, the Left Composition Lemma and the fact that a multilinear continuous map between Banach spaces is smooth, it follows that the map $\widehat{L}^{h_0}$ is differentiable. Moreover, as in section 3.2 of [21] for every $\epsilon \in (0,1]$ and every $f \in End_{h_0}^+(E)$ we have that $d\widehat{L}^{h_0}(\epsilon, f)$ is a linear, second order differential operator extending to

$$d\widehat{L}^{h_0}(\epsilon, f) : L_k^p End_{h_0}(E) \longrightarrow L_{k-2}^p End_{h_0}(E).$$

The same proof of Lemma (3.2.4) shows that this differential operator is elliptic, and it is an isomorphism if and only if it is injective or surjective.

Now, let $f \in End_{h_0}^+(E)$, so that $f^{\frac{1}{2}} \in End_{h_0}^+(E)$ (see Lemma 2.10). For every $\psi \in End(E)$, write $\psi = \{\psi_i\}_{i \in I}$, and if $f = \{f_i\}_{i \in I}$ we have $f^{\frac{1}{2}} = \{f_i^{\frac{1}{2}}\}_{i \in I}$. We define

$$Ad(f_i^{\frac{1}{2}})(\psi_i) := f_i^{\frac{1}{2}} \circ \psi_i \circ f_i^{-\frac{1}{2}} \in End(E_i).$$

Notice that

$$Ad(f_i^{\frac{1}{2}})(\psi_i) = f_i^{\frac{1}{2}} \circ \psi_i \circ f_i^{-\frac{1}{2}} =$$
$$= (\phi_{ij}^{-1} \circ f_j^{\frac{1}{2}} \circ \phi_{ij}^{-1}) \circ (\phi_{ij}^{-1} \circ \psi_j \circ \phi_{ij}) \circ (\phi_{ij}^{-1} \circ f_j^{-\frac{1}{2}} \circ \phi_{ij}) =$$
$$= \phi_{ij}^{-1} \circ (f_j^{\frac{1}{2}} \circ \psi_j \circ f_j^{-\frac{1}{2}}) \circ \phi_{ij} = \phi_{ij}^{-1} \circ Ad(f_i^{\frac{1}{2}})(\psi_j) \circ \phi_{ij}.$$

It follows that

$$Ad(f^{\frac{1}{2}})(\psi) \in End(E),$$

hence we have

$$Ad(f^{\frac{1}{2}}) : End(E) \longrightarrow End(E).$$

In a similar way one defines

$$Ad(f^{-\frac{1}{2}}) : End(E) \longrightarrow End(E).$$

Moreover, we define a new connection $D^f$ on $End(E)$ as follows: recall that if $D_0$ is the Chern connection of $(E, h_0)$, it induces the Chern connection of $(End(E), End(h_0))$, that we still write $D_0$ (which is a connection on an untwisted bundle). We let

$$D^f := Ad(f^{\frac{1}{2}}) \circ D_0 \circ Ad(f^{-\frac{1}{2}}),$$

which is again a connection on $End(E)$ compatible with $End(h_0)$.

Finally, for $f \in End_{h_0}^+(E)$ and $\varphi \in End_{h_0}(E)$ we let

$$\eta_f(\varphi) := f^{-\frac{1}{2}} \circ \varphi \circ f^{-\frac{1}{2}}.$$

The following requires the same proof of Propositions (3.2.5) and (3.2.6) in [21]:

**Lemma 5.3.** *Let $h_0 \in Herm^+(E)$, $(\epsilon, f) \in (0,1] \times L_k^p End_{h_0}^+(E)$ and $\varphi \in End_{h_0}(E)$. Suppose that $\widehat{L}^{h_0}(\epsilon, f) = 0$.*

  (1) *If there is $\alpha \in \mathbb{R}$ such that $d\widehat{L}^{h_0}(\epsilon, f)(\varphi) + \alpha f \circ \log(f) = 0$, then*

$$i\Lambda_g \overline{\partial} \partial (|\eta_\varphi(f)|^2) + 2\epsilon |\eta_\varphi(f)|^2 + |D^f(\eta_\varphi(f))|^2 \leq -2\alpha h_0(\log(f), \eta_\varphi(f)),$$

  *where $h_0$ denotes the Hermitian metric induced by $h_0$ on $End(E)$.*
  (2) *We have that $d\widehat{L}^{h_0}(\epsilon, f) : L_k^p End_{h_0}(E) \longrightarrow L_{k-2}^p End_{h_0}(E)$ is an isomorphism.*

This allows us to prove the following, as done in Corollary (3.2.7) of [21]:

**Proposition 5.4.** *Let $h_0 \in Herm^+(E)$, and suppose that $\epsilon_0 \in (0, 1]$ and $f_0 \in End_{h_0}^+(E)$ are such that $L_{\epsilon_0}^{h_0}(f_0) = 0$. Then there are $\delta > 0$ and a unique differentiable map*

$$f : (0,1] \cap (\epsilon_0 - \delta, \epsilon_0 + \delta) \longrightarrow End_{h_0}^+(E)$$

*such that $f(\epsilon_0) = f_0$ and $L_\epsilon^{h_0}(f(\epsilon)) = 0$ for every $\epsilon \in (0,1) \cap (\epsilon_0 - \delta, \epsilon_0 + \delta)$.*

Now, let $f_1 \in End_{h_0}^+(E)$ be a solution of $L_1^{h_0} = 0$, and define $J(f_1)$ to be the subset of $(0, 1]$ given by all those $\epsilon \in (0, 1]$ for which there is a differentiable map $f : [\epsilon, 1] \longrightarrow End_{h_0}^+(E)$ such that $f(1) = f_1$ and $L_{\epsilon'}^{h_0}(f(\epsilon')) = 0$ for every $\epsilon' \in [\epsilon, 1]$.

By Lemma 5.2 we know that $1 \in J(f_1)$, so that $J(f_1) \neq \emptyset$, and Proposition 5.4 it follows that $J(f_1)$ is an open subset of $(0, 1]$. It then follows that there is a maximal open interval $(\epsilon_0, 1]$ over which the perturbed equation has a unique differentiable solution.

5.2. **Closure of $J(f_1)$.** Let now $\epsilon_0 \in [0, 1)$ be such that there is a unique differentiable map

$$f : (\epsilon_0, 1] \longrightarrow End_{h_0}^+(E)$$

such that $f(1) = f_1$ and $L_\epsilon^{h_0}(f(\epsilon)) = 0$ for every $\epsilon \in (\epsilon_0, 1]$. We will write $f_\epsilon := f(\epsilon)$, and by Lemma 5.2 we may and will suppose that $\det(f_\epsilon) = 1$.

The aim of this section is to show that if $\epsilon_0 > 0$, then $f_\epsilon$ converges to $f_{\epsilon_0} \in End^+_{h_0}(E)$ as $\epsilon$ converges to $\epsilon_0$, and that $L^{h_0}_{\epsilon_0}(f_{\epsilon_0}) = 0$. Proposition 5.4 allows us then to extend $f$ to an interval of the form $(\epsilon_1, 1]$ for $\epsilon_1 < \epsilon_0$, contradicting the maximality of $\epsilon_0$. It will then follow that $J = (0, 1]$.

We will let
$$m_\epsilon := \max_{x \in X} |\log(f_\epsilon)|(x), \quad \varphi_\epsilon := f'(\epsilon), \quad \eta_\epsilon := f_\epsilon^{\frac{1}{2}} \circ \varphi_\epsilon \circ f_\epsilon^{-\frac{1}{2}}.$$

Moreover, and this will be essential in the whole section, we will suppose that $E$ is simple.

**Lemma 5.5.** *Let $E$ be a simple $\alpha$–twisted holomorphic vector bundle on $X$. For every $\epsilon > \epsilon_0$ we have that $Tr(\eta_\epsilon) = 0$, and that there is a positive real number $C(m_\epsilon)$ depending only on $m_\epsilon$ such that*
$$||D^{f_\epsilon}(\eta_\epsilon)||^2_{L^2} \geq C(m_\epsilon)||\eta_\epsilon||^2_{L^2}.$$

*Proof.* Since we supposed that $\det(f_\epsilon) = 1$ for every $\epsilon$, we have
$$\frac{d}{d\epsilon} \det(f_\epsilon) = 0.$$
But as $det(f_\epsilon) = \exp(Tr(\log(f_\epsilon)))$ (see the proof of Lemma 5.2) we get
$$0 = \frac{d}{d\epsilon} \exp(Tr(\log(f_\epsilon))) = \det(f_\epsilon) \cdot Tr\left(f_\epsilon^{-1} \circ \frac{d}{d\epsilon} f_\epsilon\right) =$$
$$= Tr(f_\epsilon^{-1} \circ \varphi_\epsilon) = Tr(f_\epsilon^{-\frac{1}{2}} \circ \varphi_\epsilon \circ f_\epsilon^{-\frac{1}{2}}) = Tr(\eta_\epsilon),$$
and the first point of the statement is proved.

For the second point, let $\psi_\epsilon := Ad(f_\epsilon^{\frac{1}{2}})(\eta_\epsilon)$. By definition of $D^{f_\epsilon}$ and $Ad(f_\epsilon^{\frac{1}{2}})$ we have
$$|D^{f_\epsilon}(\eta_\epsilon)|^2 \geq |(D^{f_\epsilon})^{0,1}(\eta_\epsilon)|^2 = |f_\epsilon^{\frac{1}{2}} \circ \overline{\partial}(\psi_\epsilon) \circ f_\epsilon^{-\frac{1}{2}}|^2 \geq C(m_\epsilon)|\overline{\partial}(\psi_\epsilon)|^2.$$
It follows that
$$||D^{f_\epsilon}(\eta_\epsilon)||^2_{L^2} = \int_X |D^{f_\epsilon}(\eta_\epsilon)|^2 \sigma_g^n \geq C(m_\epsilon) \int_X |\overline{\partial}(\psi_\epsilon)|^2 \sigma_g^n = C(m_\epsilon)||\overline{\partial}(\psi_\epsilon)||^2_{L^2}.$$

By definition we have
$$||\overline{\partial}(\psi_\epsilon)||^2_{L^2} = (\Delta_{\overline{\partial}}(\psi_\epsilon), \psi_\epsilon)_{L^2}.$$
As $Tr(\eta_\epsilon) = 0$ it follows that $Tr(\psi_\epsilon) = 0$, so $\psi_\epsilon$ is $L^2$–orthogonal to the identity, and hence since $E$ is simple to $\ker(\Delta_{\overline{\partial}})$ (see Remark 7.2.2 of [21]).

As $\Delta_{\overline{\partial}}$ is self-adjoint and elliptic, all its eigenvalues are non-negative. If $\lambda_1$ is the smallest positive eigenvalue of $\Delta_{\overline{\partial}}$, it follows that
$$(\Delta_{\overline{\partial}}(\psi_\epsilon), \psi_\epsilon)_{L^2} \geq \lambda_1 ||\psi_\epsilon||^2_{L^2}.$$
It then follows that
$$||D^{f_\epsilon}(\eta_\epsilon)||^2_{L^2} \geq C'(m_\epsilon)||\psi_\epsilon||^2_{L^2} \geq C'(m_\epsilon)||\eta_\epsilon||^2_{L^2},$$
and we are done. □

Another important property is the following:

**Lemma 5.6.** *Let $E$ be a simple $\alpha-$twisted holomorphic vector bundle. For every $\epsilon$ we have*
$$\max_{x\in X}|\varphi_\epsilon|(x)\leq C(m_\epsilon)$$
*for some positive real number $C(m_\epsilon)$ depending only on $m_\epsilon$.*

*Proof.* The proof is as the one of Proposition (3.3.3) of [21], when one has to replace Proposition (3.2.5) with Lemma 5.3 and Lemma (3.3.1) with Lemma 5.5. □

We moreover need the following two Lemmas:

**Lemma 5.7.** *Let $E$ be a simple $\alpha-$twisted holomorphic vector bundle on $X$, $h_0$ a Hermitian metric on $E$ and $f\in End_{h_0}^+(E)$ a solution of $L_\epsilon^{h_0}=0$ for some $\epsilon>0$.*

  (1) *We have*
  $$\frac{1}{2}i\Lambda_g\overline{\partial}\partial(|\log(f)|^2)+\epsilon|\log(f)|^2\leq |K_g^0(E,h_0)|\cdot|\log(f)|.$$

  (2) *If $m:=\max_{x\in X}|\log(f)|(x)$, then we have*
  $$m\leq\frac{1}{\epsilon}\max_{x\in X}|K_g^0(E,h_0)|(x).$$

  (3) *There is real number $C$ (depending only on $g$ and $h_0$) such that*
  $$m\leq C(||\log(f)||_{L^2}+\max_{x\in X}|K_g^0(E,h_0)|(x))^2).$$

*Proof.* The proof is identical to that of Lemma (3.3.4) of [21]. □

**Lemma 5.8.** *Let $E$ be a simple $\alpha-$twisted holomorphic vector bundle on $X$. Suppose that there is $m\in\mathbb{R}$ such that $m_\epsilon\leq m$ for every $\epsilon\in(\epsilon_0,1]$. Then for every $p>1$ and every $\epsilon\in(\epsilon_0,1]$ we have that*
$$||\varphi_\epsilon||_{L_2^p}^2\leq C(m)\cdot(1+||f_\epsilon||_{L_2^p}^2),\quad ||f_\epsilon||_{L_2^p}^2\leq e^{C(m)(1-\epsilon)}(1+||f_1||_{L_2^p}^2).$$

*Proof.* The proof is identical to that of Proposition (3.3.5) of [21], where one has to replace Proposition (3.3.3) by Lemma 5.6. □

All these results together allow us to prove the following:

**Proposition 5.9.** *Let $E$ be a simple $\alpha-$twisted holomorphic vector bundle on $X$, $h_0$ a Hermitian metric on $E$ and $f_1$ a solution of $L_1^{h_0}=0$. Then $J(f_1)=(0,1]$.*

*Proof.* The proof is identical to that of Proposition (3.3.6), point i), where one has to replace Lemma (3.3.4) with Lemma 5.7 and Proposition (3.3.5) with Lemma 5.8. □

5.3. **Existence of Hermite-Einstein metrics.** The first result we prove is the following:

**Proposition 5.10.** *Let $E$ be a simple $\alpha-$twisted holomorphic vector bundle on $X$, $h_0$ a Hermitian metric on $E$ and $f_1 \in End^+_{h_0}(E)$ a solution of the perturbed equation $L_1^{h_0} = 0$. Let $f : (0,1] \longrightarrow End^+_{h_0}(E)$ be the unique differentiable solution of the perturbed equation, and suppose that there is a real number $C$ such that $||f_\epsilon||_{L^2} \leq C$ for every $\epsilon \in (0,1]$. Then the equation $L_0^{h_0} = 0$ has a solution.*

*Proof.* As by hypothesis we have a uniform bound for $||f_\epsilon||_{L^2}$, then point 3 of Lemma 5.7 provides a uniform bound for $m_\epsilon \leq m$, and hence Lemma 5.8 provides a uniform bound for $||f_\epsilon||^2_{L^p_2}$.

This implies the existence of a sequence $\epsilon_k$ converging to 0 such that $f_{\epsilon_k}$ converges weakly to some $f_0 \in L^p_2 End(E)$. The fact that $m_\epsilon \leq m$ implies that the eigenvalues of $f_0$ take values in $[e^{-m}, e^m]$, which implies that $f_0 \in L^p_2 End^+(E)$ (see the proof of Lemma 7.3.10 of [21]).

Now we know that $L^p_2$ has a compact embedding in $L^p_1$, hence we may suppose that $f_{\epsilon_k}$ converges strongly to $f_0$ in $L^p_1 End^+(E)$.

Suppose that that $L_0^{h_0}(f_0) = 0$. This implies that $i\Lambda_g \overline{\partial} D_0^{1,0}(f_0)$ is a multilinear algebraic expression in $f$, $\log(f)$, $D_0^{1,0}(f)$ and $\overline{\partial}(f)$. But since $f_0 \in L^p_2 End(E)$, it follows that $i\Lambda_g \overline{\partial} D_0^{1,0}(f) \in L^p_1 End(E)$. The Elliptic Regularity Theorem then implies that $f_0 \in L^p_3 End(E)$. Repeating this process, by Rellich's Theorem we then see that $f_0 \in A^0(\underline{End}(E))$, i. e. $f_0$ is a solution of $L_0^{h_0} = 0$.

We are then left to prove that $L_0^{h_0}(f_0) = 0$. To do so, we just need to prove that if $\zeta$ is any smooth endomorphism of $E$, we have $(L_0^{h_0}(f_0), \zeta)_{L^2} = 0$. Recall that $L_{\epsilon_k}^{h_0}(f_{\epsilon_k}) = 0$ for every $k$, hence we get

$$(L_0^{h_0}(f_0), \zeta)_{L^2} = (L_0^{h_0}(f_0) - L_{\epsilon_k}^{h_0}(f_{\epsilon_k}), \zeta)_{L^2}.$$

If we let $\Psi_k := f_0^{-1} \circ D_0^{1,0}(f_0) f_{\epsilon_k}^{-1} \circ D^{1,0}(f_{\epsilon_k})$, we then see that

$$(L_0^{h_0}(f_0), \zeta)_{L^2} = (i\Lambda_g \overline{\partial}(\Psi_k), \zeta)_{L^2} + (\epsilon_k \log(f_{\epsilon_k}), \zeta)_{L^2}.$$

Now, recall that if $p \gg n$, then the map

$$L^p \longrightarrow L^2, \quad f \mapsto \log(f)$$

is continuous, hence we see that $(\epsilon_k \log(f_{\epsilon_k}), \zeta)_{L^2}$ converges to 0. Moreover we have

$$(\Lambda_g \overline{\partial}(\Psi_k), \zeta)_{L^2} = \int_X Tr(\Lambda_g \overline{\partial}(\Psi_k) \cdot *\zeta) = \frac{1}{n!} \int_X \Lambda_g(Tr(\overline{\partial}(\Psi_k) \cdot \zeta))\sigma_g^n =$$

$$= \frac{1}{(n-1)!} \int_X Tr(\overline{\partial}(\Psi_k) \cdot \zeta)\sigma_g^{n-1} = \frac{1}{(n-1)!} \int_X Tr(\Psi_k \wedge \overline{\partial}(\zeta \cdot \sigma_g^{n-1})) = (\Psi_k, \beta)_{L^2},$$

where $\beta \in A^{1,0}(\underline{End}(E))$ does not depend on $k$.

Remark that the map
$$L_1^p \longrightarrow L^2, \quad f \mapsto f^{-1} \circ D_0^{1,0}(f)$$
is continuous (as $p \gg n$). As $f_{\epsilon_k}$ converges strongly to $f_0$ in $L_1^p$, we have that $\Psi_k$ converges to 0, so $(\Psi_k, \beta)_{L^2}$ converges to 0. But this implies that
$$(L_0^{h_0}(f_0), \zeta)_{L^2} = (i\Lambda_g \overline{\partial}(\Psi_k), \zeta)_{L^2} + (\epsilon_k \log(f_{\epsilon_k}), \zeta)_{L^2}$$
converges to 0, i. .e $(L_0^{h_0}(f_0), \zeta)_{L^2} = 0$ for every $\zeta$, so $L_0^{h_0}(f_0) = 0$. □

As a consequence, if $||f_\epsilon||_{L^2}$ is uniformly bounded, then the Hermitian metric $\widehat{(f_0)}_{h_0}$ is a $g$–Hermite-Einstein metric, and hence $E$ is $g$–Hermite-Einstein. In the remaining part of this section we prove that if $||f_\epsilon||_{L^2}$ are not uniformly bounded, then $E$ is not $g$–stable. These two results together will prove that if $E$ is a $g$–stable $\alpha$–twisted holomorphic vector bundle, then $E$ is $g$–Hermite-Einstein, proving Theorem 5.1.

The main definition we will need is the following:

**Definition 5.11.** *If $E$ is an $\alpha$–twisted holomorphic vector bundle on $X$, an element $\pi \in L_1^2 End(E)$ is a **weakly holomorphic $\alpha$–twisted subbundle** of $E$ if $\pi^* = \pi = \pi^2$ and $(id_E - \pi) \circ \overline{\partial}(\pi) = 0$ almost everywhere on $X$.*

The reason for the name is the following:

**Lemma 5.12.** *Let $E$ be an $\alpha$–twisted holomorphic vector bundle on $X$ whose associated locally free $\alpha$–twisted sheaf is $\mathscr{E}$. If $\pi \in L_1^2 End(E)$ is a weakly holomorphic $\alpha$–twisted subbundle of $E$, there is a coherent $\alpha$–twisted subsheaf $\mathscr{F}$ of $\mathscr{E}$ and an analytic subset $S$ of $X$ such that:*

(1) *$S$ has codimension at least 2 in $X$,*
(2) *$\pi_{|X\setminus S} \in A^0(E_{|X\setminus S})$ and we have $\pi_{|X\setminus S}^* = \pi_{|X\setminus S} = \pi_{|X\setminus S}^2$ and $(id_{E|X\setminus S} - \pi_{|X\setminus S}) \circ \overline{\partial}(\pi_{|X\setminus S}) = 0$,*
(3) *$\mathscr{F}_{|X\setminus S}$ is the image of $\pi_{|X\setminus S}$ and an $\alpha$–twisted holomorphic subbundle of $E$.*

*Proof.* We let $\pi = \{\pi_i\}_{i\in I}$, and notice that $\pi_i \in L_1^2 End(E_i)$ is such that $\pi_i^* = \pi_i = \pi_i^2$ and $(id_{E_i} - \pi_i) \circ \overline{\partial}(\pi_i) = 0$ almost everywhere on $U_i$. As $E_i$ as a holomorphic vector bundle on $U_i$, by [27] (see Theorem 3.4.3 of [21]) there is a coherent subsheaf $\mathscr{F}_i$ of $\mathscr{E}_i$ and an analytic subset $S_i$ of $U_i$ such that

(1) $S_i$ has codimension at least 2 in $U_i$,
(2) $\pi_{i|U_i\setminus S_i} \in A^0(E_{i|U_i\setminus S_i})$ and

$$\pi_{i|U_i\setminus S_i}^* = \pi_{i|U_i\setminus S_i} = \pi_{i|U_i\setminus S_i}^2, \quad (id_{E_i|U_i\setminus S_i} - \pi_{i|U_i\setminus S_i}) \circ \overline{\partial}(\pi_{i|U_i\setminus S_i}) = 0,$$

(3) $\mathscr{F}_{i|U_i\setminus S_i}$ is the image of $\pi_{i|U_i\setminus S_i}$, and a holomorphic subbundle of $E_i$.

Now, the fact that $\pi$ is an endomorphism of $E$ imply that $\mathscr{F}$ is an $\alpha$–twisted coherent sheaf and that $S_i \cap U_{ij} = S_j \cap U_{ij}$ (since $S_i$ is the locus of $U_i$ where $\mathscr{F}_i$ is not free). Hence the $S_i$'s glue together to give an analytic subset $S$ of $X$ of codimension at least 2, and we are done. □

Before going on we need two Lemmas, and we recall that if $f \in End_{h_0}^+(E)$, then for every $\sigma \in (0,1]$ we may define $f^\sigma \in End_{h_0}^+(E)$ (see section 2.2.1 and Lemma 2.53).

**Lemma 5.13.** *Let $E$ be an $\alpha$−twisted holomorphic vector bundle and $h_0 \in Herm^+(E)$. Consider $f \in End_{h_0}^+(E)$ and $\sigma \in (0,1]$.*

(1) *We have*
$$i\Lambda_g(h_0(f^{-1} \circ D_0^{1,0}(f), D_0^{1,0}(f^\sigma))) \geq |f^{-\frac{\sigma}{2}} \circ D_0^{1,0}(f^\sigma)|^2.$$

(2) *If there is $\epsilon > 0$ such that $L_\epsilon^{h_0}(f) = 0$, then*
$$\frac{1}{\sigma}i\Lambda_g\overline{\partial}\partial(Tr(f^\sigma)) + \epsilon h_0(\log(f), f^\sigma) + |f^{-\frac{\sigma}{2}} \circ D_0^{1,0}(f^\sigma)|^2 \leq -h_0(K_g^0(E, h_0), f^\sigma).$$

*Proof.* The proof is as that of Lemma (3.4.4) of [21]. $\square$

For $\epsilon > 0$ and $x \in X$ let $\lambda(\epsilon, x)$ be the largest eigenvalue of $\log(f_{\epsilon,x})$ (which is well-defined by Remark 2.9). We let
$$M_\epsilon := \max_{x \in X} \lambda(\epsilon, x), \qquad \rho(\epsilon) := e^{-M_\epsilon}.$$
As $Tr(\log(f_\epsilon)) = 0$, it follows that $M_\epsilon$ grows as $m_\epsilon$ and $\rho(\epsilon) \leq 1$. The second Lemma we need is the following:

**Lemma 5.14.** *Let $E$ be an $\alpha$−twisted holomorphic vector bundle on $X$, $h_0 \in Herm^+(E)$, $f_1$ a solution of $L_1^{h_0} = 0$ and $f : (0,1] \longrightarrow End_{h_0}^+(E)$ the unique differentiable map such that $f(1) = f_1$ and $L_\epsilon^{h_0}(f_\epsilon) = 0$ for every $\epsilon \in (0,1]$. If $\limsup_{\epsilon \to 0} ||\log(f_\epsilon)||_{L^2} = +\infty$, then:*

(1) *for every $x \in X$, if $\lambda$ is an eigenvalue of $\rho(\epsilon)f_{\epsilon,x}$, then $\lambda \leq 1$.*
(2) *For every $x \in X$ there is an eigenvalue $\lambda$ of $\rho(\epsilon)f_{\epsilon,x}$ such that $\lambda \leq \rho(\epsilon)$.*
(3) *We have $\max_{x \in X}(\rho(\epsilon)|f_\epsilon|(x)) \geq 1$.*
(4) *If $\epsilon_k$ is a sequence converging to 0, then $\rho(\epsilon_k)$ converges to 0.*

*Proof.* The proof is as that of Lemma (3.4.5) of [21]. $\square$

The two previous Lemmas 5.13 and 5.14 allows us to produce a weakly holomorphic $\alpha$−twisted subbundle of $E$.

**Proposition 5.15.** *Let $E$ be an $\alpha$−twisted holomorphic vector bundle on $X$, $h_0 \in Herm^+(E)$, $f_1$ a solution of $L_1^{h_0} = 0$ and $f : (0,1] \longrightarrow End_{h_0}^+(E)$ the unique differentiable map such that $f(1) = f_1$ and $L_\epsilon^{h_0}(f_\epsilon) = 0$ for every $\epsilon \in (0,1]$. If $\limsup_{\epsilon \to 0} ||\log(f_\epsilon)||_{L^2} = +\infty$, then:*

(1) *for $k \to +\infty$ the endomorphisms $\rho(\epsilon_k)f_{\epsilon_k}$ converge weakly in $L_1^2$ to an element $f_\infty \neq 0$,*
(2) *there is a sequence $\{\sigma_l\}_{l \in \mathbb{N}}$ converging to 0 such that $f_\infty^{\sigma_l}$ converge weakly in $L_1^2$ to an element $f_\infty^0$,*
(3) *letting $\pi := id_E - f_\infty^0$, then $\pi$ is a weakly holomorphic $\alpha$−twisted subbundle of $E$.*

*Proof.* The proof is exactly as the one of Proposition (3.4.6) of [21], where one has to replace Lemma (3.4.4) by Lemma 5.13, Lemma (3.3.4) with Lemma 5.7, and Lemma (3.4.5) with Lemma 5.14. □

Now, by Lemma 5.12 we see that the weakly holomorphic $\alpha-$twisted subbundle of $E$ defines an $\alpha-$twisted coherent subsheaf $\mathscr{F}$ of $\mathscr{E}$, and with the same proof of Corollary (3.4.7) we see that $0 < rk(\mathscr{F}) < rk(E)$. We conclude this section with the following:

**Proposition 5.16.** *Let $E$ be an $\alpha-$twisted holomorphic vector bundle on $X$, $h_0 \in Herm^+(E)$, $f_1$ a solution of $L_1^{h_0} = 0$ and $f : (0,1] \longrightarrow End_{h_0}^+(E)$ the unique differentiable map such that $f(1) = f_1$ and $L_\epsilon^{h_0}(f_\epsilon) = 0$ for every $\epsilon \in (0,1]$. If $\limsup_{\epsilon \to 0} ||\log(f_\epsilon)||_{L^2} = +\infty$, then $E$ is not $g-$stable.*

*Proof.* The proof is identical to that of Proposition (3.4.8) of [21]. □

As a consequence we finally conclude with the proof of Theorem 5.1

*Proof.* By Theorem 4.19 we have that if $E$ is $g-$Hermite-Einstein, then $E$ is $g-$polystable.

Conversely, if $E$ be $g-$stable, then Corollary 4.17 implies that $E$ is simple, hence by Lemma 5.2 and Propositions 5.4 and 5.9 there is a Hermitian metric $h_0$ on $E$ for which there is $f_1 \in End_{h_0}^+(E)$ which is solution of the perturbed equation $L_1^{h_0} = 0$, and we let $f : (0,1] \longrightarrow End_{h_0}^+(E)$ be the unique differentiable solution of the perturbed equation.

If there is no real number $C$ such that $||f_\epsilon||_{L^2} \leq C$ for every $\epsilon \in (0,1]$, then
$$\limsup_{\epsilon \to 0} ||\log(f_\epsilon)||_{L^2} = +\infty.$$
By Proposition 5.16 then $E$ is not $g-$stable, which is impossible. As a consequence there must be a constant $C$ such that $||f_\epsilon||_{L^2} \leq C$ for every $\epsilon \in (0,1]$. By Proposition 5.10 then $E$ is $g-$Hermite-Einstein.

If now $E$ is $g-$polystable, then $E = E_1 \oplus \cdots \oplus E_k$ where $E_1, \cdots, E_k$ are all $g-$stable of the same $g-$slope. Hence $E_1, \cdots, E_k$ are all $g-$Hermite-Einstein with the same Einstein factor, so by Proposition 3.11 $E$ is $g-$Hermite-Einstein. □

The definition of mean curvature, of Chern classes, of degree, of $g-$stability and of $g-$Hermite-Einstein metrics do not depend on the fact that $g$ is a Kähler metric on $X$, but only on the fact that $\sigma_g^{n-1}$ is $\partial\bar{\partial}-$closed, i. e. on the fact that $g$ is a Gauduchon metric on $X$. By [9] we know that on every compact complex manifold there is a Gauduchon metric.

By Remark 4.22, and by the fact that all the results of this section go through if we suppose that $g$ is a Gauduchon metric on $X$ (see [21]), we conclude the following, providing a generalization of [29].

**Theorem 5.17.** *Let $X$ be a compact complex manifold and $g$ a Gauduchon metric on $X$. An $\alpha-$twisted holomorphic vector bundle is $g-$Hermite-Einstein if and only if $E$ is $g-$polystable.*

## 6. Approximate Kobayashi-Hitchin

The aim of this section is to complete the proof of the approximate Kobayashi-Hitchin correspondence for twisted vector bundles, i. e. the following

**Theorem 6.1.** *Let $E$ be an $\alpha-$twisted holomorphic vector bundle over a compact Kähler manifold with Kähler metric $g$. Then $E$ is $g-$semistable if and only if it is approximate $g-$Hermite-Einstein.*

Theorem 4.23 tells us that if $E$ is approximate $g-$Hermite-Einstein, then it is $g-$semistable. This section is devoted to prove the converse, and for this we follow closely [17] and [13].

First we introduce the Donaldson Lagrangian, and show that the associated evolution equation has a unique smooth solution on $\mathbb{R}_+$ once the starting Hermitian metric is fixed. As a consequence of this, we show that if the Donaldson Lagrangian of $E$ is bounded from below, then $E$ is approximate $g-$Hermite-Einstein.

In order to prove Theorem 6.1 we will then just need to prove that the Donaldson Lagrangian for $E$ is bounded from below. To do so we will adapt to twisted sheaves the regularization process described in [13].

### 6.1. Donaldson's Lagrangian.
Let $h, k \in Herm^+(E)$, and consider the $k-$Hermitian endomorphism $f^{k,h} \in End(E)$. The determinant $\det(f^{k,h})$ of $f^{k,h}$ is a smooth function on $X$, and since by Remark 2.47 we know that $f^{k,h}$ is invertible we see that $\det(f^{k,h})$ is never zero. If we let

$$Q_1(h,k) := \log(\det(f^{k,h})),$$

it follows that this is a smooth function on $X$.

**Lemma 6.2.** *For every $h, k, l \in Herm^+(E)$ we have*

$$Q_1(h,k) = -Q_1(k,h), \quad Q_1(h,l) = Q_1(h,k) + Q_1(k,l).$$

*Proof.* By definition, for every $h \in Herm^+(E)$ we have $f^{h,h} = id_E$, so that $Q_1(h,h) = 0$. Moreover, by Remark 2.47 we have $f^{h,k} = (f^{k,h})^{-1}$, hence $Q_1(h,k) = -Q_1(k,h)$. By Remark 2.48 we moreover have $f^{l,k} \circ f^{k,h} = f^{l,h}$ for every $h, k, l \in Herm^+(E)$, hence we have

$$Q_1(h,l) = \log(\det(f^{l,h})) = \log(\det(f^{l,k} \circ f^{k,h})) = \log(\det(f^{l,k})\det(f^{k,h})) =$$

$$= \log(\det(f^{l,k})) + \log(\det(f^{k,h})) = Q_1(k,l) + Q_1(h,k),$$

and we are done. □

Consider now $\hbar \in \Omega^A_{h,k}(E)$, where $A$ is an interval in $\mathbb{R}$. For every $t \in A$ we let $h_t := \hbar(t)$, which is a Hermitian metric on $E$, and we let $D_t$ be the Chern connection of $(E, h_t)$ and $R_t$ its curvature. We then have a function

$$\mathscr{R} : A \longrightarrow A^{1,1}(\underline{End}(E)), \quad \mathscr{R}(t) := R_t.$$

If $\hbar$ is differentiable, we have
$$f^\hbar : A \longrightarrow End(E), \quad f^\hbar(t) := f^{h_t, h'_t},$$
where $h'_t = \hbar'(t)$. We know that $\hbar$ is a geodetic in $Herm^+(E)$ from $h$ to $k$ if and only if $f^\hbar$ is constant, i. e. if and only if $\partial_t f^\hbar = 0$.

Since $f^{h_t, h'_t} \in End(E)$ and $R_t \in A^{1,1}(\underline{End}(E))$ for every $t \in A$, we see that $f^{h_t, h'_t} \cdot R_t \in A^{1,1}(\underline{End}(E))$. Its trace is then a smooth $(1,1)$–form on $X$ depending on $t$, and we define
$$Q_2^\hbar(h, k) := i \int_A Tr(f^{h_t, h'_t} \cdot R_t) dt,$$
a $(1,1)$–form on $X$ depending on $\hbar$.

We now introduce the following notation: if $\hbar \in \Omega_{h,k}^{a,b}(E)$, we let
$$\hbar^{-1} : [a, b] \longrightarrow Herm^+(E), \quad \hbar^{-1}(t) := \hbar(b - t + a).$$
If $\hbar \in \Omega_{a,b}^{h,k}(E)$ and $\hbarbar \in \Omega_{k,l}^{b,c}(E)$, we define
$$\hbar * \hbarbar : [a, c] \longrightarrow Herm^+(E), \quad \hbar * \hbarbar(t) := \begin{cases} \hbar(t), & t \in [a, b] \\ \hbarbar(t), & t \in [b, c] \end{cases}$$
which is then an element of $\Omega_{h,l}^{a,c}(E)$.

We have the following:

**Lemma 6.3.** *For every $h, k, l \in Herm^+(E)$ and for every $\hbar \in \Omega_{h,k}^{a,b}(E)$ and $\hbarbar \in \Omega_{k,l}^{b,c}$ we have*
$$Q_2^{\hbar^{-1}}(k, h) = -Q_2^\hbar(h, k), \quad Q_2^{\hbar * \hbarbar}(h, l) = Q_2^\hbar(h, k) + Q_2^{\hbarbar}(k, l).$$

*Proof.* Notice that $\hbar^{-1} \in \Omega_{k,h}^{a,b}(E)$ and we have
$$(h^{-1})'_t = (\hbar^{-1})'(t) = \frac{d}{dt} \hbar(b - t + a) = -\hbar'(b - t + a) = -h'_{b-t+a}.$$
By Remark 2.49 it follows that
$$f^{\hbar^{-1}}(t) = f^{h_t^{-1}, (h^{-1})'_t} = f^{h_{b-t+a}, -h'_{b-t+a}} = -f^\hbar(b - t + a).$$
Moreover, as $R(\hbar^{-1}(t)) = R_{b+a-t}$, we get
$$Q_2^{\hbar^{-1}}(k, h) = i \int_a^b Tr(f^{\hbar^{-1}}(t) \cdot R(\hbar^{-1}(t))) dt =$$
$$= -\int_a^b Tr(f^{h_s, h'_s} \cdot R_s) ds = -Q_2^\hbar(h, k).$$
The second equality in the statement is trivial. □

Finally, for every $h, k \in Herm^+(E)$ and every $\hbar \in \Omega_{h,k}^A(E)$ we let
$$L_g^\hbar(h, k) := \int_X \left( Q_2^\hbar(h, k) - \frac{c_g(E)}{n} Q_1(h, k) \sigma_g \right) \wedge \sigma_g^{n-1},$$

where $n$ is the dimension of $X$. The function $L_g^\hbar(h,k)$ will be called **Donaldson Lagrangian of $E$ between $h$ and $k$ along $\hbar$**.

The following is an immediate consequence of Lemmas 6.2 and 6.3:

**Lemma 6.4.** *For every $h, k, l \in Herm^+(E)$ and for every $\hbar \in \Omega_{h,k}^{a,b}(E)$ and $\hbar \in \Omega_{k,l}^{b,c}$ we have*

$$L_g^{\hbar^{-1}}(k,h) = -L_g^\hbar(h,k), \quad L_g^{\hbar * \hbar}(h,l) = L_g^\hbar(h,k) + L_g^\hbar(k,l).$$

6.1.1. *The Donaldson Lagrangian is path-independent.* In this section we prove that $L_g^\hbar(h,k)$ does not depend on $\hbar$, but only on $h$ and $k$. To do so, we first prove the following:

**Lemma 6.5.** *Let $h$ be a Hermitian metric on $E$ and $\hbar \in \Omega_{h,h}^{a,b}(E)$. Then*

$$Q_2^\hbar(h,h) \in \partial A^{0,1}(X) + \overline{\partial} A^{1,0}(X).$$

*Proof.* Let $a_0, \cdots, a_m \in [a,b]$ such that $a_0 = a$, $a_m = b$ and $a_j < a_{j+1}$ for every $0 \leq j \leq m-1$ be the points where $\hbar$ is not differentiable. Fix now $k \in Herm^+(E)$, and for every $0 \leq j \leq m$ consider the segment

$$\hbar_j : [0,1] \longrightarrow Herm^+(E), \quad \hbar_j(t) := (1-t)k + th_{a_j}$$

joining $k$ to $h_{a_j}$. Let $\widetilde{\hbar}_j := \hbar_j * \hbar_{|[a_j, a_{j+1}]} * \hbar_{j+1}^{-1}$, which is a piecewise differentiable closed curve based at $k$. By Lemma 6.3 we have

$$Q_2^{\widetilde{\hbar}_j}(k,k) = Q_2^{\hbar_j}(k, h_{a_j}) + Q_2^{\hbar_{|[a_j, a_{j+1}]}}(h_{a_j}, h_{a_{j+1}}) - Q_2^{\hbar_{j+1}}(k, h_{a_{j+1}})$$

for every $j = 0, \cdots, m-1$, so that

$$Q_2^{\widetilde{\hbar}_0 * \cdots * \widetilde{\hbar}_{m-1}}(k,k) = \sum_{j=0}^{m-1} Q_2^{\hbar_{|[a_j, a_{j+1}]}}(h_{a_j}, h_{a_{j+1}}) = Q_2^\hbar(h,h).$$

As a consequence, if we know that $Q_2^{\widetilde{\hbar}_j}(k,k) \in \partial A^{0,1}(X) + \overline{\partial} A^{1,0}(X)$, it follows that the same holds for $Q_2^\hbar(h,h)$.

We then just need to prove the statement for the curve $\widetilde{\hbar}_j$. We let $\Delta_j := [a_j, a_{j+1}] \times [0,1]$ and

$$H_j : \Delta_j \longrightarrow Herm^+(E), \quad H_j(t,s) = sh_t + (1-s)k,$$

which is a smooth function such that

- $H_j(t,0) = k$ for every $t \in [a_j, a_{j+1}]$,
- $H_j(t,1) = h_t$ for every $t \in [a_j, a_{j+1}]$
- $H_j(a_j, s) = \hbar_j(s) = (1-s)k + sh_{a_j}$ for every $s \in [0,1]$,
- $H_j(a_{j+1}, s) = \hbar_{j+1}(s) = (1-s)k + sh_{a_{j+1}}$ for every $s \in [0,1]$.

We now let

$$\partial_s H_j(t) := \left. \frac{d}{ds} H_j(t,s) \right|_{s=0}, \quad \partial_t H_j(s) := \left. \frac{d}{dt} H_j(t,s) \right|_{t=0},$$

which are Hermitian forms on $E$, and let

$$u_j(t,s) := f^{H_j(t,s), \partial_s H_j(t)}, \qquad v_j(t,s) := f^{H_j(t,s), \partial_t H_j(s)}$$

$$R_j(t,s) := \overline{\partial}(H_j(t,s)^{-1} \partial H_j(s,t)),$$

i. e. $R_j(t,s)$ is the Chern curvature of $(E, H_j(t,s))$. We then have three functions

$$u_j, v_j : \Delta_j \longrightarrow End(E), \quad R_j : \Delta_j \longrightarrow A^{1,1}(\underline{End}(E)).$$

We now define a $(1,0)$−form $\widetilde{d}H_j$ on $\Delta_j$ with coefficients in $E^* \otimes \overline{E}^*$ as follows:

$$\widetilde{d}H_j := \partial_s H_j(t) \cdot ds + \partial_t H_j(s) \cdot dt$$

(notice that both $\partial_s H_j(t)$ and $\partial_t H_j(s)$ are Hermitian forms on $E$, i. e. smooth global sections of the untwisted vector bundle $\overline{E}^* \otimes E^*$). We moreover let

$$H_j^{-1} \widetilde{d}H_j := u_j ds + v_j dt,$$

which is a $(1,0)$−form on $\Delta_j$ with coefficients in $\underline{End}(E)$. We then let

$$H_j^{-1} \widetilde{d}H_j \cdot R_j := u_j \cdot R_j ds + v_j \cdot R_j dt,$$

which is then a $(1,0)$ form on $\Delta_j$ whose coefficients are $(1,1)$−forms on $X$ with coefficients in $\underline{End}(E)$. Finally, let

$$\phi_j := iTr(H_j^{-1} \widetilde{d}H_j \cdot R_j) = iTr(u_j \cdot R_j)ds + iTr(v_j \cdot R_j)dt,$$

which is a $(1,0)$−form on $\Delta_j$ whose coefficients are $(1,1)$−forms on $X$.

Now, we have

$$\int_{\partial \Delta_j} \phi_j = \int_a^b (\phi_j(t,1) - \phi_j(t,0))dt + \int_0^1 (\phi_j(a_j,s) - \phi_j(a_{j+1},s))ds.$$

As $v_j(t,0) = 0$, we get

$$\phi_j(t,0)dt = iTr(v_j(t,0) \cdot R_j(t,0))dt = 0.$$

Moreover we have

$$\phi_j(a_j, s)ds = iTr(u_j(a_j, s) \cdot R_j(a_j, s))ds = iTr(f^{\hbar_j(s), \hbar'_j(s)} \cdot R_j(\hbar_j(s)))ds,$$

$$\phi_j(a_{j+1}, s)ds = iTr(u_j(a_{j+1}, s) \cdot R_j(a_{j+1}, s))ds =$$
$$= iTr(f^{\hbar_{j+1}(s), \hbar'_{j+1}(s)} \cdot R_j(\hbar_{j+1}(s)))ds,$$

$$\phi_j(t,1)dt = iTr(v_j(t,1) \cdot R_j(t,1))dt = iTr(f^{\hbar(t), \hbar'(t)} \cdot R_j(\hbar(t)))dt.$$

It then follows that

$$\int_{\partial \Delta_j} \phi_j = Q_2^{\hbar_j}(h_{a_j}, k) + Q_2^{\hbar_{|[a_j, a_{j+1}]}}(h_{a_j}, h_{a_{j+1}}) - Q_2^{\hbar_{j+1}}(h_{a_{j+1}}, k).$$

We then have

$$Q_2^{\widetilde{\hbar_j}}(k,k) = \int_{\partial \Delta_j} \phi_j = \int_{\Delta_j} \widetilde{d}\phi_j,$$

and we just need to prove that
$$\int_{\Delta_j} \widetilde{d}\phi_j \in \partial A^{0,1}(X) + \overline{\partial} A^{1,0}(X).$$
To do so, we start by letting
$$\alpha_j : \Delta_j \longrightarrow A^{0,1}(X), \quad \alpha_j(s,t) := iTr(v_j(s,t) \cdot \overline{\partial} u_j(s,t)).$$
We notice that for every $s, t \in \Delta_j$ we have $v(s,t), u(s,t) \in End(E)$, so that $v(s,t) \cdot \overline{\partial} u(s,t) \in A^{0,1}(End(E))$, and its trace is then a $(0,1)$−form on $X$.

For every $(s,t) \in \Delta_j$ let $D_j(s,t)$ be the Chern connection of $(E, H_j(s,t))$, whose curvature is $R_j(s,t)$. We then have $D_j^{0,1}(s,t) = \overline{\partial}$, so
$$\alpha_j(s,t) = iTr(v_j(s,t) \cdot D_j^{0,1}(s,t) u_j(s,t)).$$
It follows that
$$\partial \alpha_j = iTr(D_j^{1,0} v_j \wedge D_j^{0,1} u_j + v_j \cdot D_j^{1,0} D_j^{0,1} u_j).$$
Similarly, we have
$$\overline{\alpha}_j : \Delta_j \longrightarrow A^{1,0}(X),$$
and
$$\overline{\partial} \overline{\alpha}_j = -iTr(D_j^{0,1} v_j \wedge D_j^{1,0} u_j + v_j \cdot D_j^{0,1} D_j^{1,0} u_j).$$
As a consequence we see that
$$\partial \alpha_j + \overline{\partial} \overline{\alpha}_j = iTr(-D_j^{0,1} D_j^{1,0}(v_j \cdot u_j) + (D_j^{0,1} D_j^{1,0} v_j) \cdot u_j + v_j \cdot (D_j^{1,0} D_j^{0,1} u_j)).$$

Or aim is to prove that
$$\widetilde{d}\phi_j = -(\partial \alpha_j + \overline{\partial} \overline{\alpha}_j) ds \wedge dt + i\overline{\partial}\partial Tr(v_j \cdot u_j) ds \wedge dt,$$
which will then conclude the proof. In order to do this we need the following formulas:
$$\partial_s R_j = D_j^{0,1} D_j^{1,0} u_j, \quad \partial_t R_j = D_j^{0,1} D_j^{1,0} v_j,$$
$$\partial_t u_j = -v_j \cdot u_j + f^{H_j, \partial_t \partial_s H_j}, \quad \partial_s v_j = -u_j \cdot v_j + f^{H_j, \partial_s \partial_t H_j}$$
(the last two formulas are immediate from the definition of $u_j$ and $v_j$, while the first two formulas may be proved locally as in the proof of Lemma 3.6 in Chapter VI of [17], and then by gluing the local formulas since the $B-$field $B$ does not depend neither on $t$ nor on $s$).

Now recall that
$$\phi_j(s,t) = iTr(u_j \cdot R_j) ds + iTr(v_j \cdot R_j) dt, \quad \widetilde{d} = \partial_s ds + \partial_t dt,$$
hence we have
$$\widetilde{d}\phi_j = iTr(-\partial_t u_j \cdot R_j - u_j \cdot \partial_t R_j + \partial_s v_j \cdot R_j + v_j \cdot \partial_s R_j) ds \wedge dt =$$
$$= iTr(-v_j(D_j^{1,0} D_j^{0,1} + D_j^{0,1} D_j^{1,0}) u_j + v_j D_j^{0,1} D_j^{1,0} u_j - u_j D_j^{0,1} D_j^{1,0} v_j) ds \wedge dt.$$
The formula we are looking for then follows immediately. □

The previous result has an important corollary, which tells us that the Donaldson Lagrangian between two Hermitian metrics does not depend on the chosen path.

**Corollary 6.6.** *Let $h, k \in Herm^+(E)$. For $i = 1, 2$ and $a_i < b_i$ two real numbers, let $\hbar_i \in \Omega_{h,k}^{a_i,b_i}(E)$. Then $L_g^{\hbar_1}(h, k) = L_g^{\hbar_2}(h, k)$.*

*Proof.* Consider $\hbar := \hbar_1 * \hbar_2^{-1} \in \Omega_{h,h}^{a,b}(E)$. We then have
$$L_g^{\hbar}(h, h) = \int_X \left( Q_2^{\hbar}(h, h) - \frac{c_g(E)}{n} Q_1(h, h) \sigma_g \right) \wedge \sigma_g^{n-1}.$$
Lemma 6.2 we get $Q_1(h, h) = 0$, and hence
$$L_g^{\hbar}(h, h) = i \int_X \left( \int_a^b Tr(f^{h_t, h_t'} \cdot R_t) dt \right) \wedge \sigma_g^{n-1}.$$
Now, by Lemma 6.5 there are a $(0, 1)-$form $\phi$ and a $(1, 0)-$form $\psi$ such that
$$i \int_a^b Tr(f^{h_t, h_t'} \cdot R_t) dt = \partial \phi + \overline{\partial} \psi,$$
hence
$$L_g^{\hbar}(h, h) = \int_X (\partial \phi + \overline{\partial} \psi) \wedge \sigma_g^{n-1} = 0.$$
But now notice that by Lemma 6.4 we have
$$L_g^{\hbar}(h, h) = L_g^{\hbar_1}(h, k) - L_g^{\hbar_2}(h, k),$$
which concludes the proof. $\square$

As a consequence on the previous Corollary we will be allowed to use the notation $L_g(h, k)$ instead of $L_g^{\hbar}(h, k)$, that will simply be called **Donaldson Lagrangian between** $h$ **and** $k$. In particular, for every $k \in Herm^+(E)$ we get a function
$$L_{g,k} : Herm^+(E) \longrightarrow \mathbb{R}, \quad L_{g,k}(h) := L_g(h, k),$$
called **Donaldson Lagrangian at** $k$.

6.1.2. *Critical points of the Donaldson Lagrangian.* We now want to relate the Donaldson Lagrangian and Hermite-Einstein metrics. We first need the following:

**Lemma 6.7.** *Let $\hbar$ be a differentiable curve in $Herm^+(E)$ and $k \in Herm^+(E)$. For every $t$ consider the segment $\hbar_t$ in $Herm^+(E)$ connecting $h_t$ and $k$. Then we have*
$$\partial_t Q_1(h_t, k) = Tr(f^{h_t, h_t'}),$$
$$\partial_t Q_2^{\hbar_t}(h_t, k) = i Tr(f^{h_t, h_t'} \cdot R_t) \bmod \partial A^{0,1}(X) + \overline{\partial} A^{1,0}(X).$$

*Proof.* Along the proof of Lemma 6.5 we proved that if $\hbar : [a, b] \longrightarrow Herm^+(E)$ is a piecewise differentiable curve and $k \in Herm^+(E)$, we have
$$i \int_a^b Tr(f^{h_s, h_s'} \cdot R_s) ds + Q_2^{\hbar_a}(h_a, k) - Q_2^{\hbar_b}(h_b, k) \in \partial A^{0,1}(X) + \overline{\partial} A^{1,0}(X).$$
This holds for every $t \in [a, b]$, and if we derive this with respect to $t$ we get the statement. $\square$

Using this, we are able to prove the following:

**Lemma 6.8.** *Let $\hbar$ be a differentiable curve in $Herm^+(E)$, and fix $k \in Herm^+(E)$. Then we have*

$$\partial_t L_{g,k}(h_t) = \frac{1}{n!}\int_X Tr((K_g(E,h_t) - c_g(E)\cdot id_E)\circ f^{h_t,h'_t})\sigma_g^n.$$

*Proof.* By definition of $L_{g,k}(h_t)$ and Lemma 6.7 we have

$$\partial_t L_{g,k}(h_t) = \int_X \left(iTr(f^{h_t,h'_t}\cdot R_t) - \frac{c_g(E)}{n}Tr(f^{h_t,h'_t})\sigma_g\right)\frac{\sigma_g^{n-1}}{(n-1)!}.$$

Since, by definition, we have

$$iR_t \wedge \sigma_g^{n-1} = \frac{1}{n}K_g(E,h_t)\sigma_g^n,$$

we then get

$$\partial_t L_{g,k}(h_t) = \frac{1}{n!}\int_X (Tr(f^{h_t,h'_t}\circ K_g(E,h_t)) - c_g(E)Tr(f^{h_t,h'_t}))\sigma_g^n,$$

which proves the statement. $\square$

Let $\hbar$ be a differentiable curve in $Herm^+(E)$ and let $\widehat{K}_g(E,h_t)$ be the mean curvature Hermitian form of $(E,h_t)$: we notice that $\widehat{K}_g(E,h_t) - c_g(E)h_t$ is a Hermitian form on $E$ as well. Recall that for every $t$ we have that $h'_t$ is a Hermitian form on $E$, hence $h'_t$ and $\widehat{K}_g(E,h_t) - c_g(E)h_t$ are both tangent vectors of $Herm^+(E)$ at $h_t$.

Recall that $K_g(E,h_t)$ is a $h_t$–Hermitian endomorphism by Lemma 3.5: by Example 2.54 we then get

$$K_g(E,h_t) = f^{h_t,\widehat{K}_g(E,h_t)},$$

so

$$K_g(E,h_t) - c_g(E)id_E = f^{h_t,\widehat{K}_g(E,h_t)} - c_g(E)f^{h_t,h_t} = f^{h_t,\widehat{K}_g(E,h_t) - c_g(E)h_t}.$$

But this implies that

$$\int_X Tr((K_g(E,h_t) - c_g(E)\cdot id_E)\circ f^{h_t,h'_t})\sigma_g^n =$$

$$= \int_X Tr(f^{h_t,\widehat{K}_g(E,h_t) - c_g(E)h_t}\circ f^{h_t,h'_t})\sigma_g^n = (\widehat{K}_g(E,h_t) - c_g(E)h_t, h'_t)_{h_t},$$

where the last is the Riemannian metric at $h_t$. Lemma 6.8 then gives

$$\partial_t L_{g,k}(h_t) = (\widehat{K}_g(E,h_t) - c_g(E)h_t, h'_t)_{h_t}.$$

This allows us to conclude the following:

**Proposition 6.9.** *Let $k \in Herm^+(E)$. An element $h \in Herm^+(E)$ is a critical point for $L_{g,k}$ if and only if $(E,h)$ is $g$–Hermite-Einstein.*

*Proof.* The point $h$ is critical for $L_{g,k}$ if and only if for every differentiable curve $\hbar$ from $h$ we have that
$$\left.\frac{d}{dt}L_{g,k}(h_t)\right|_{t=0} = 0.$$
The previous discussion shows that $h$ is critical for $L_{g,k}$ if and only if
$$(\widehat{K}_g(E,h) - c_g(E)h, h'_0)_h = 0.$$
As this has to be verified for every differentiable curve $\hbar$ from $h$, this happens if and only if $\widehat{K}_g(E,h) - c_g(E)h = 0$, i. e. if and only if $(E,h)$ is $g-$Hermite-Einstein. $\square$

6.2. **The evolution equation.** Let now $k \in Herm^+(E)$, and consider $\hbar : [a,b] \longrightarrow Herm^+(E)$ a curve of Hermitian metrics on $E$. Consider the function
$$L_{g,k} : Herm^+(E) \longrightarrow \mathbb{R},$$
and let us look at its differential $dL_{g,k}$: for a given $h \in Herm^+(E)$, we have a linear map
$$dL_{g,k,h} : T_h Herm^+(E) \longrightarrow \mathbb{R}.$$
In particular, for every $t \in [a,b]$ we have $h'_t \in T_{h_t} Herm^+(E)$, and we may calculate the value of $dL_{g,k,h_t}$ on $h'_t$. By definition of the differential, and by the discussion in the previous section, we have
$$dL_{g,k,h_t}(h'_t) = \partial_t L_{g,k}(h_t) = (\widehat{K}_g(E,h_t) - c_g(E)h_t, h'_t)_{h_t}.$$
Let now $\mathrm{grad} L_{g,k}$ be the gradient of $L_{g,k}$, i. e. the vector field on $Herm^+(E)$ which is dual to the differential form $dL_{g,k}$: the previous formula gives then
$$\mathrm{grad} L_{g,k}(h_t) = \widehat{K}_g(E,h_t) - c_g(E)h_t.$$
We now consider the **evolution equation**
$$\partial_t \hbar = -\mathrm{grad} L_{g,k}(\hbar) = -(\widehat{K}_g(E,\hbar) - c_g(E)\hbar)$$
in terms of Hermitian forms, or equivalently
$$f^{\hbar} = -(K_g(E,\hbar) - c_g(E)id_E)$$
in terms of endomorphisms. Our present aim is to study the solutions of the evolution equation for a given $k \in Herm^+(E)$: we will show that there is always a unique smooth solution $\hbar : [0,+\infty) \longrightarrow Herm^+(E)$ with $\hbar(0) = k$.

6.2.1. *Uniqueness and convergence of solutions.* The first result is about uniqueness. More precisely, we prove the following:

**Proposition 6.10.** *Let $k \in Herm^+(E)$.*
  (1) *If $\hbar_1, \hbar_2 : [0,a) \longrightarrow Herm^+(E)$ are two smooth solutions of the evolution equation (for the given $k$) such that $h_{1,0} = h_{2,0}$, then $h_{1,t} = h_{2,t}$ for every $t \in [0,a)$.*

(2) *If $\hbar : [0, a) \longrightarrow Herm^+(E)$ is a smooth solution of the evolution equation (for the given k), then $h_t$ converges uniformly in the $C^0$-topology to a Hermitian metric $h_a$ as t converges to a.*

*Proof.* Fix $b < a$ and let $\Delta := [0, b] \times [0, 1]$. We let $t$ be the variable of $[0, b]$ and $s$ the variable on $[0, 1]$, and we define

$$H : \Delta \longrightarrow Herm^+(E), \quad H(t, s) := (1 - s)h_{1,t} + sh_{2,t},$$

which is a smooth function. Notice that $H(t, 0) = h_{1,t}$ and $H(t, 1) = h_{2,t}$.

As in the proof of Lemma 6.5, we let

$$u(t, s) := f^{H(t,s), \partial_s H(t,s)}, \quad v(t, s) := f^{H(t,s), \partial_t H(t,s)},$$

so that we get

$$u, v : \Delta \longrightarrow A^0(\underline{End}(E)),$$

i. e. $u(t, s)$ and $v(t, s)$ are endomorphisms of $E$ for every $(t, s) \in \Delta$. We moreover let

$$R : \Delta \longrightarrow A^{1,1}(\underline{End}(E))$$

be the function mapping $(t, s) \in \Delta$ to the Chern curvature $R(t, s)$ of $(E, H(t, s))$.

As already remarked in the proof of Lemma 6.5, we have

$$\partial_s R = D^{0,1} D^{1,0} u, \quad \partial_t R = D^{0,1} D^{1,0} v,$$

$$\partial_t u = -v \cdot u + f^{H(t,s), \partial_t \partial_s H(t,s)}, \quad \partial_s v = -u \cdot v + f^{H(t,s), \partial_s \partial_t H(t,s)}$$

where $D = D^{1,0} + D^{0,1}$ is the family of Chern connections defined by $H$ on $E$. We notice that $D^{0,1} = \overline{\partial}$, i. e. $D^{0,1}$ does not depend neither on $t$ nor on $s$, while $D^{1,0}$ varies with $(t, s)$.

If we fix $t_0 \in [0, b]$, for every $s \in [0, 1]$ we have $u(t_0, s) \in End(E)$. Hence we have $u(t_0, s)^2 = u(t_0, s) \circ u(t_0, s) \in End(E)$, thus $Tr(u(t_0, s)^2)$ is a smooth function on $X$ for every $s$. We then let

$$e_{t_0} := \frac{1}{2} \int_0^1 Tr(u(t_0, s)^2) ds,$$

which is then a smooth function on $X$: for every $x \in X$, the value $e_{t_0}(x)$ is the energy (with respect to the Riemannian metric) of the path $H(t_0, s)_x$ running in $Herm^+(E_x)$ from $h_{1,t_0,x}$ to $h_{2,t_0,x}$. We then notice that

$$e_{t_0} : X \longrightarrow \mathbb{R}_+,$$

and we have $e_{t_0}(x) = 0$ if and only it $h_{1,t_0,x} = h_{2,t_0,x}$. Finally, we let

$$e : [0, b] \longrightarrow C^\infty(X), \quad e(t) := e_t.$$

Suppose now that for every $t_0 \in [0, b]$ and every $x \in X$ the path $H(t_0, s)_x$ is a geodetic in $Herm^+(E_x)$. As shown in the proof of Lemma 6.5, we have that $u(t, s) = u(t)$ does not depend on $s$, i. e. $\partial_s u = 0$. Moreover, we have

$$\partial_s v = -u \cdot v + f^{H(t,s), \partial_s \partial_t H(t,s)}.$$

As a consequence we have
$$\partial_t e = \frac{1}{2}\int_0^1 Tr(\partial_t u \cdot u + u \cdot \partial_t u)ds = \int_0^1 Tr(u \cdot \partial_t u)ds.$$
Since
$$\partial_t u = -v \cdot u + f^{H(t,s),\partial_t\partial_s H(t,s)},$$
we get
$$\partial_t e = \int_0^1 Tr(-u \cdot v \cdot u + u \cdot f^{H(t,s),\partial_t\partial_s H(t,s)})ds = \int_0^1 Tr(u\partial_s v)ds.$$
Integrating by parts and using $\partial_s u = 0$ we then get
$$\partial_t e = Tr(u(t) \cdot (v(t,1) - v(t,0))) = Tr(u(t) \cdot (f^{h_{2,t},h'_{2,t}} - f^{h_{1,t},h'_{1,t}})).$$

As we are supposing that $\hbar_1$ and $\hbar_2$ are two solutions of the evolution equation, we have that
$$f^{\hbar_2} = K_g(E,h_{2,t}) - c_g(E)id_E, \quad f^{\hbar_1} = K_g(E,h_t) - c_g(E)id_E.$$
We then get
$$\partial_t e = Tr(u \cdot (K_g(E,h_{2,t}) - K_g(E,h_{1,t}))).$$

Similar calculations give that
$$\overline{\partial}\partial e = -\int_0^1 Tr(D^{1,0}u \wedge D^{0,1}u)ds + Tr(u \cdot (R(h_{2,t}) - R(h_{1,t}))),$$
where $R(h)$ denotes the Chern curvature of $(E,h)$. Notice that $\overline{\partial}\partial e$ is then a $(1,1)$-form on $X$, so we may consider $i\Lambda_g(\overline{\partial}\partial e)$, which is a smooth function on $X$: we will write $\square_g e := i\Lambda_g(\overline{\partial}\partial e)$, and we have
$$\square_g e = -\int_0^1 |D^{1,0}u|^2 ds + Tr(u \cdot (K_g(E,h_{2,t}) - K_g(E,h_{2,t}))).$$

As a consequence we see that
$$\partial_t e + \square_g e = -\int_0^1 |D^{1,0}u|^2 ds \leq 0.$$
Let us now consider the function
$$m : [0,b] \longrightarrow \mathbb{R}, \quad m(t) := \max_{x \in X}(e_t(x)).$$
By the Maximum Principle for Parabolic Equations (see Lemma 4.1 in Chapter VI of [17]) we know that the function $m$ is monotone decreasing in $t$.

Now, consider $m(0) = \max_{x \in X} e_0(x)$: recall that $e_0(x)$ is the energy of the path $H(0,s)_x$ in $Herm^+(E_x)$ connecting $h_{1,0,x}$ and $h_{2,0,x}$. Since $h_{1,0} = h_{2,0}$ by hypothesis, we then get $e_0 = 0$, so that $m(0) = 0$. Since $m$ is monotone decreasing, it follows that $m(t) \leq 0$ for every $t \in [0,b]$, so that $e_t(x) \leq 0$ for every $t \in [0,b]$ and every $x \in X$. But since $e_t(x) \geq 0$, we then get $e_t(x) = 0$ for every $x \in X$ and every $t \in [0,b]$, i. e. $h_{1,t} = h_{2,t}$ for every $t \in [0,b]$. Since this holds for every $b < a$, we then see that $\hbar_1 = \hbar_2$, completing the proof of the first point of the statement.

For the second point, let $x \in X$ and consider $t, t' \in [0, a)$. We let $\rho_x(t, t')$ be the distance between $h_{t,x}$ and $h_{t',x}$. Moreover, we let $e_x(t, t')$ be the energy of the (unique) geodesic path in $Herm^+(E_x)$ connecting $h_{t,x}$ to $h_{t',x}$. We have
$$e_x(t, t') = \frac{1}{2}\rho_x(t, t')^2.$$
Since the Riemannian metric is complete, in order to prove the statement we just need to prove that the family $\{h_t\}_{t \in [0,a)}$ is uniformly Cauchy, i. e. that for every $\epsilon > 0$ there is $\delta > 0$ such that for every $t, t' \in [0, a)$ such that $|t - t'| < \delta$ we have
$$\max_{x \in X} e_x(t, t') < \epsilon.$$
Fix then $\epsilon > 0$. Recall that the family $\{h_t\}_{t \in [0,a)}$ is continuous at 0, hence there is $\delta > 0$ such that for every $\tau \in [0, \delta)$ we have
$$\max_{x \in X} e_x(0, \tau) < \epsilon.$$
Take now any two $t, t' \in [0, a)$ such that $t < t'$ and $t' - t < \delta$. Write $\tau := t' - t \in [0, \delta)$. Let
$$\hbar_1 : [0, a) \longrightarrow Herm^+(E), \quad \hbar_1(s) := \hbar(s + \tau),$$
so that in particular $h_{1,t} = h_{t'}$. Consider the function
$$f : X \times [0, a) \longrightarrow \mathbb{R}, \quad f(x, t) := e_x(t, t').$$
If we let
$$H : \Delta \longrightarrow Herm^+(E), \quad H(t, s) := (1 - s)h_t + sh_{1,t},$$
consider the energy function $e : [0, a) \longrightarrow C^\infty(X)$ associated to $H$ (that we used in the previous part of the proof), i. e. the function mapping $t \in [0, a)$ to the energy function $e_t : X \longrightarrow \mathbb{R}_+$, defined by letting $e_t(x)$ be the energy of the path $H(t, s)_x$ in $Herm^+(E_x)$ connecting $h_{t,x}$ to $h_{1,t,x} = h_{t+\tau,x} = h_{t',x}$. We then see that for every $x \in X$ and every $t \in [0, a)$ we have $e_t(x) = f(x, t)$.

But now we know from the previous part of the proof that the function
$$m : [0, a) \longrightarrow \mathbb{R}, \quad m(t) := \max_{x \in X} f(x, t) = \max_{x \in X} e_t(x)$$
is monotone decreasing in $t$. This implies that
$$\max_{x \in X} e_x(t, t') = \max_{x \in X} f(x, t) \leq \max_{x \in X} f(x, 0) = \max_{x \in X} e_x(0, \tau) < \epsilon,$$
and we are done. □

Proposition 6.10 tells us that if the evolution equation has a smooth solution $\hbar$ over an interval $[0, a)$, then this solution is unique and may be extended to a continuous family $\overline{\hbar}$ defined over $[0, a]$.

6.2.2. *Short-time solution.* We now want to show that given two Hermitian metrics $h, k$ on $E$, then there are $\delta > 0$ and $\hbar : [0, \delta] \longrightarrow Herm^+(E)$ which solves the evolution equation (with respect to $k$) and which is such that $h_0 = h$.

If $\hbar : [0, a] \longrightarrow Herm^+(E)$ is a smooth curve, we let
$$R_\hbar : [0, a] \longrightarrow A^{1,1}(\underline{End}(E)), \quad R_\hbar(t) := R_{k_t},$$
were $R_{k_t}$ is the Chern curvature of $(E, k_t)$, and
$$\widehat{K}_\hbar : [0, a] \longrightarrow Herm(E) \quad \widehat{K}_\hbar(t) := \widehat{K}_g(E, k_t).$$
Moreover, we let
$$P_\hbar : [0, a] \longrightarrow Herm(E), \quad P_\hbar(t) := k_t' + \widehat{K}_g(E, k_t) - c_g(E)k_t.$$

Consider now a smooth curve $\hbar : [0, a] \longrightarrow Herm^+(E)$, and let $v : [0, a] \longrightarrow Herm(E)$. For $s \ll 1$ we have that $h_t + sv_t \in Herm^+(E)$ for every $t \in [0, a]$, i. e. we have
$$\hbar + sv : [0, a] \longrightarrow Herm^+(E).$$
We view $\hbar, v, \hbar + sv$ as smooth global sections of $\pi^* H_E$. We will let
$$P_\hbar(v) := P_{\hbar+v}$$
for every $v$ such that $\hbar + v$ is a family of Hermitian metrics.

For every global section $v$ of $\pi^* H_E$ we let
$$dP_\hbar(v) := \lim_{s \to 0} \frac{P_\hbar(v) - P_\hbar(0)}{s}.$$
It is easy to see that
$$dP_\hbar(v) = v' + \frac{d}{ds} \widehat{K}_{\hbar+sv} \bigg|_{s=0} - c_g(E)v.$$

Notice that $f^{\hbar,v}$ is a smooth family of endomorphisms of $E$. Then $i\Lambda_g \overline{\partial}\partial f^{\hbar,v}$ is a family of endomorphisms of $E$, and the family of Hermitian forms associated to it and to $\hbar$ will be denoted $\square_\hbar v$.

By Lemma 3.6 we have that
$$\frac{d}{ds} \widehat{K}_{\hbar+sv} \bigg|_{s=0} = \square_\hbar v,$$
so that
$$dP_\hbar(v) = v' + \square_\hbar v - c_g(E)v,$$
i. e.
$$dP_\hbar = \partial_t + \square_\hbar - c_g(E) \cdot id_{Herm(E)}.$$

If $v \in L_q^p(X/0, a, H_E)$, then $dP_\hbar(v) \in L_{q-2}^p(X/0, a, H_E)$, i. e. we may view $dP_\hbar$ as a linear operator
$$dP_\hbar : L_q^p(X/0, a, H_E) \longrightarrow L_{q-2}^p(X/0, a, H_E).$$

The proof of the following is as that of Lemma 6.5 in Chapter VI of [17].

**Lemma 6.11.** *Let $p \geq 2n+2$ and $q \geq 2$. Then*
$$dP_\hbar : L^p_q(X/0, a, H_E) \longrightarrow L^p_{q-2}(X/0, a, H_E)$$
*is an isomorphism.*

As a consequence we get the following:

**Proposition 6.12.** *Let $E$ be an $\alpha$-twisted holomorphic vector bundle on a compact Kähler manifold $X$. There is $\delta > 0$ such that the evolution equation has a smooth solution $\hbar : [0, \delta] \longrightarrow Herm^+(E)$.*

*Proof.* Let $a > 0$ and choose a Hermitian metric $h$ on $E$. Let
$$\hbar : [0, a] \longrightarrow Herm^+(E), \quad \hbar(t) := h$$
be the constant function of value $h$. Consider an integer $p > 2n+2$, where $2n$ is the real dimension of $X$. We know by Lemma 6.11 that
$$dP_\hbar : L^p_2(X/0, a, H_E) \longrightarrow L^p_0(X/0, a, H_E)$$
is an isomorphism.

The Implicit Function Theorem implies that $P_\hbar$ maps a neighborhood $U$ of $0 \in L^p_2(X/0, a, H_E)$ onto a neighborhood $U'$ of $P_\hbar(0)$ in $L^p_0(X/0, a, H_E)$. Let now
$$w : [0, a] \longrightarrow Herm(E), \quad w(t) := \begin{cases} 0, & t \in [0, \delta] \\ P_\hbar(0), & t \in (\delta, a] \end{cases}$$
for some $\delta > 0$ such that $w \in U'$ (which is possible if $\delta \ll 1$).

But then there is $v \in U$ such that $P_\hbar(v) = w$. Define then
$$\hbar := \hbar + v : [0, \delta] \longrightarrow Herm^+(E).$$
Notice that as $w_{|[0,\delta]} = 0$, we get
$$0 = w_{|[0,\delta]} = P_\hbar(v)_{|[0,\delta]} = P(\hbar) = \hbar' + \widehat{K}\hbar - c_g(E)\hbar,$$
so that
$$\hbar' = -(\widehat{K}_\hbar - c_g(E)\hbar),$$
i. e. $\hbar$ is a solution of the evolution equation. The fact that $\hbar$ is smooth can be proved exactly as in Theorem 7.1 in Chapter VI of [17]. □

We then see that the evolution equation has always a unique smooth solution for short time intervals $[0, \delta]$ (for $\delta \ll 1$).

6.2.3. *All-time solution.* The aim of this section is the show that the unique smooth solution of the evolution equation for a short time interval found in the previous section can be extended to a unique smooth solution on $[0, +\infty)$.

Let now $\hbar : [a, b] \longrightarrow Herm^+(E)$ be a family of Hermitian metrics on an $\alpha$-twisted holomorphic vector bundle $E$, and suppose it is a solution of the evolution equation, i. e.
$$\hbar' = -(\widehat{K}_g(E, \hbar) - c_g(E)\hbar).$$

For every $t \in [a,b]$ let $h_t = \hbar(t)$, and write $h_t = \{h_{t,i}\}$ where $h_{t,i}$ is a Hermitian metric on $E_i$ and we have
$$h_{t,i} = {}^T\phi_{ij} h_{t,j} \overline{\phi}_{ij}.$$

Let $D_t$ be the Chern connection of $(E, h_t)$, so that $D_t = \{D_{t,i}\}$ where $D_{t,i}$ is the Chern connection of $(E_i, h_{t,i})$. For every $i \in I$ consider $\partial_t D_{t,i}$: if $\Gamma_{t,i}$ is the connection form of $D_{t,i}$ with respect to a given local frame, we have
$$\Gamma_{t,i} = \phi_{ij}^{-1} \Gamma_{t,j} \phi_{ij} + \phi_{ij}^{-1} d\phi_{ij} + \omega_{ij} id_{E_i},$$
we see that
$$\partial_t \Gamma_{t,i} = \phi_{ij}^{-1} \partial_t \Gamma_{t,j} \phi_{ij},$$
so the $\partial_t D_{t,i}$'s glue together to give $\partial_t D_t \in A^{1,0}(\underline{End}(E))$.

Now, for every $i \in I$ we have
$$\partial_t D_{t,i} = \partial_t(h_{t,i}^{-1} \cdot \partial h_{t,i}) = -h_{t,i}^{-1} \partial_t h_{t,i} \cdot \partial h_{t,i} + h_{t,i}^{-1} \partial \partial_t h_{t,i}.$$

As $\hbar$ is a solution of the evolution equation, we have $\partial_t h_t = -(\widehat{K}_t - c_g h_t)$, where we let $\widehat{K}_t := \widehat{K}_g(E, h_t)$ and $c_g = c_g(E)$. Restricting to $E_i$ we get
$$\partial_t h_{t,i} = -(\widehat{K}_{t|E_i} - c_g h_{t,i}).$$

We then get
$$\partial_t D_{t,i} = h_{t,i}^{-1}(\widehat{K}_{t|E_i} - c_g(E) h_{t,i}) \partial h_{t,i} - h_{t,i}^{-1} \partial(\widehat{K}_{t|E_i} - c_g h_{t,i}) =$$
$$= h_{t,i}^{-1}(\widehat{K}_{t|E_i} - c_g h_{t,i}) \cdot h_{t,i}^{-1} \partial h_{t,i} - \partial(h_{t,i}^{-1}(\widehat{K}_{t|E_i} - c_g h_{t,i})) - h_{t,i}^{-1} \partial h_{t,i} \cdot h_{t,i}^{-1}(\widehat{K}_{t|E_i} - c_g(E) h_{t,i}).$$

Notice now that
$$h_{t,i}^{-1}(\widehat{K}_{t|E_i} - c_g h_{t,i}) = K_{t|E_i} - c_g(E) id_{E_i},$$
where $K_t = K_g(E, h_t)$, and that $h_{t,i}^{-1} \partial h_{t,i} = D_{t,i}$, hence we get
$$\partial_t D_{t,i} = (K_{t|E_i} - c_g id_{E_i}) \cdot D_{t,i} - \partial(K_{t|E_i} - c_g id_{E_i}) - D_{t,i} \cdot (K_{t|E_i} - c_g(E) id_{E_i}) =$$
$$= -D_{t,i}^{1,0}(K_{t|E_i} - c_g(E) id_{E_i}) = -D_{t,i}^{1,0}(K_{t|E_i}),$$
where $D_{t,i}^{1,0}$ here denotes the $(1,0)$-part of the Chern connection induced by $D_{t,i}$ on $End(E_i)$.

As the $K_{t|E_i}$'s glue together to give $K_t$ of $E$, and as $D_{t,i}^{1,0}$ is the restriction to $E_i$ of the $(1,0)$-part of the Chern connection of $(End(E), End(h))$, we then see that
$$\partial_t D_t = -D_t^{1,0} K_t.$$

Now, if $V$ is any holomorphic vector bundle, $h$ is a Hermitian metric on it and $D$ is the Chern connection of $(V, h)$, we let
$$\delta_h^{0,1} : A^{p,q}(V) \longrightarrow A^{p,q-1}(V), \qquad \delta_h^{0,1}(\sigma) := -*(D^{0,1}(*\sigma)).$$

Recall that (see section 2 in Chapter III of [17])
$$\Lambda_g D^{1,0} - D^{1,0} \Lambda_g = i\delta_h^{0,1}.$$

**Lemma 6.13.** *We have $D_t^{1,0} K_t = \delta_{h_t}^{0,1} R_t$.*

*Proof.* In the statement $R_t$ is the Chern curvature of $(E, h_t)$. The Hermitian metric $h_t$ induces a Hermitian metric $End(h_t)$ on $\underline{End}(E)$, whose Chern connection is the connection induced by the Chern connection $D_t$ of $(E, h_t)$.

The Bianchi identity (see Lemma 3.1) gives $D_t R_t = 0$. We then get

$$D_t^{1,0} K_t = D^{1,0}(i\Lambda_g R_t) = iD_t^{1,0}\Lambda_g R_t = i(-i\delta_{h_t}^{0,1} + \Lambda_g D^{1,0})R_t =$$
$$= \delta_{h_t}^{0,1} R_t + i\Lambda_g(D_t^{1,0} R_t) = \delta_{h_t}^{0,1} R_t,$$

and we are done. □

Since $\partial_t D_t = -D_t^{1,0} K_t$, Lemma 6.13 gives $\partial_t D_t = -\delta_{h_t}^{0,1} R_t$. Now, let

$$\Box_h := \overline{\partial}\delta_h^{0,1} - \delta_h^{0,1}\overline{\partial}.$$

**Corollary 6.14.** *We have $\partial_t R_t = -\Box_{h_t} R_t$ and $\partial_t K_t = -\Box_{h_t} K_t$. Moreover, we have*

$$(\partial_t + \Box_{h_t})Tr(R_t) = 0.$$

*Proof.* To prove the statement, recall that $R_{t|E_i} = R_{t,i} - B_i \cdot id_{E_i}$, so that $\partial_t R_{t|E_i} = \partial_t R_{t,i}$. Hence

$$(\partial_t R_t)_{|E_i} = \partial_t R_{t,i} = \partial_t \overline{\partial} D_{t,i} = \overline{\partial}\partial_t D_{t,i} = \overline{\partial}(\partial_t D_t)_{|E_i} = -\overline{\partial}(\delta_{h_t}^{0,1} R_t)_{|E_i}.$$

Hence we get

$$\partial_t R_t = -\overline{\partial}\delta_{h_t}^{0,1} R_t = -\Box_{h_t} R_t + \delta_{h_t}^{0,1}\overline{\partial} R_t = -\Box_{h_t} R_t.$$

Applying $i\Lambda_g$ we get $\partial_t K_t = -\Box_{h_t} K_t$, and taking the trace we get the last part of the statement. □

Now, let $\hbar : [0, a) \longrightarrow Herm^+(E)$ be a solution of the evolution equation for some $a > 0$, whose existence is known since Proposition 6.12. Let

$$f : X \times [0, a] \longrightarrow \mathbb{R}, \quad f(x, t) := |R_t|^2(x),$$

and for $k \geq 0$ we let

$$f_k : X \times [0, a] \longrightarrow \mathbb{R}, \quad f_k(x, t) := |D_t^k R_t|^2(x),$$

where $D_t^k$ is the composition of $k$ copies of $D_t$. The following may be proved exactly as in Lemma 8.7 in Chapter VI of [17], using Corollary 6.14 in place of equation (6.8.5), and Lemma 6.13 in place of equation (6.8.4).

**Lemma 6.15.** *There are real numbers $c, c_k \in \mathbb{R}$ depending only on $g$ such that*

$$(\partial_t + \Box_{h_t})f \leq c(f^{\frac{3}{2}} + f),$$
$$(\partial_t + \Box_{h_t})f_k \leq c_k f_k^{\frac{1}{2}}\left(\sum_{i+j=k} f_i^{\frac{1}{2}}(f_j^{\frac{1}{2}} + 1)\right).$$

*Moreover, we have $(\partial_t + \Box_{h_t})|K_t|^2 \leq 0$.*

As a consequence, we get the following:

**Corollary 6.16.** *Let $E$ be an $\alpha-$twisted holomorphic vector bundle and consider a smooth solution $\hbar : [0, a) \longrightarrow Herm^+(E)$ of the evolution equation.*

(1) *The functions*
$$s_R : [0, a) \longrightarrow \mathbb{R}, \quad s_R(t) := \sup_{x \in X} |Tr(R_t)|(x),$$
*and*
$$s_K : [0, a) \longrightarrow \mathbb{R}, \quad s_K(t) := \sup_{x \in X} |K_t|(x)$$
*are both monotone decreasing, and in particular bounded.*

(2) *If there is $\beta \in \mathbb{R}$ such that $|R_t|(x) \leq \beta$ for every $x \in X$ and $t \in [0, a)$, then for each $k \in \mathbb{N}$ there is $\beta_k \in \mathbb{R}$ such that $|D_t^k R_t|(x) \leq \beta_k$ for every $x \in X$ and $t \in [0, a)$.*

*Proof.* By Corollary 6.14 we know that $(\partial_t + \square_{h_t})|Tr(R_t)| = 0$, and by Lemma 6.15 we know that that $(\partial_t + \square_{h_t})|K_t|^2 \leq 0$. It then follows from the Maximum Principle for Parabolic Equations (see Lemma 4.1 in Chapter VI of [17]) that both $s_R$ and $s_K$ are monotone decreasing functions.

For the remaining part, we proceed by induction on $k$. If $k = 1$, this is the hypotesis. Consider now $k \in \mathbb{N}$, and suppose that for every $j < k$ there is a constant $\beta_j$ such that $|D_t^j R_t|(x) \leq \beta_j$ for every $x \in X$ and $t \in [0, a)$, i. e. $f_j^{\frac{1}{2}}(x, t) \leq \beta_j$ for every $(x, t) \in X \times [0, a)$. By Lemma 6.15 there is then a constant $A_k$ such that
$$(\partial_t + \square_{h_t}) f_k \leq A_k(1 + f_k).$$

Now, consider the following Cauchy problem
$$\begin{cases} (\partial_t + \square_{h_t})(u) = A_k(1 + u) \\ u(0) = f_k(0) \end{cases}$$
and notice that $(\partial_t + \square_{h_t})(u) = (\partial_t + \square_{h_t})(1 + u)$. The differential equation is linear in $1 + u$, hence the Cauchy problem has a unique smooth solution $\widetilde{u}$ defined for every $t \geq 0$.

An easy calculation then shows
$$(\partial_t + \square_{h_t})((f_k - \widetilde{u})e^{-A_k t}) = e^{-A_k t}((\partial_t + \square_{h_t})(f_k) - A_k(f_k + 1)) \leq 0.$$

But then the Maximum Principle for Parabolic Equation (see Lemma 4.1 in Chapter VI of [17]) gives us that the function
$$g : [0, a) \longrightarrow \mathbb{R}, \quad g(t) := \sup_{x \in X}((f_k(t) - \widetilde{u}(t))e^{-A_k t})$$
is monotone decreasing. As $f_k(0) = \widetilde{u}(0)$ we get that $g(0) = 0$, so that $g(t) \leq 0$ for every $t \in [0, a)$. Since $e^{-A_k t} > 0$ for every $t$, it follows that $f_k(t) \leq \widetilde{u}(t)$ for every $t$, and hence the statement holds. $\square$

Another useful consequence of Lemma 6.15 is the following, whose proof is identical to that of Lemma 8.16 in Chapter VI of [17]:

**Lemma 6.17.** *Let $E$ be an $\alpha-$twisted holomorphic vector bundle and consider a smooth solution $\hbar : [0, a) \longrightarrow Herm^+(E)$ of the evolution equation. If there is $q > 3n$ for which $|R_t|$ is uniformly bounded in $L^q(X)$ (i. e. independently of $t$), then $|R_t|$ is uniformly bounded in $L^\infty(X)$.*

This allows us to prove the following:

**Lemma 6.18.** *Let $E$ be an $\alpha-$twisted holomorphic vector bundle and consider a smooth solution $\hbar : [0, a) \longrightarrow Herm^+(E)$ of the evolution equation. We let*
$$\mathscr{R} : [0, a) \longrightarrow C^\infty(X), \quad \mathscr{R}(t) := |R_t|.$$
*Suppose that:*
   (1) *$h_t$ converges in the $C^0-$topology to a Hermitian metric $h_a$ for $t$ converging to $a$;*
   (2) *the function $\sup_{x \in X} |K_t|(x)$ is uniformly bounded on $[0, a)$.*

*Then for $p < +\infty$, we have that $\hbar$ is bounded in $C^1(X, a, H_E)$ and in $L_2^p(X, a, H_E)$, and that $\mathscr{R}$ is bounded in $L^p(X, a)$.*

*Proof.* The proof is almost identical to that of Lemma 8.22 in Chapter VI of [17]. Suppose that $\hbar$ is not bounded in $C^1(X, a, H_E)$. This implies that there is a sequence $\{t_k\}_{k \in \mathbb{N}}$ of points in $[0, a)$ with the two following properties:
   (1) $\lim_{k \to +\infty} t_k = a$ and
   (2) if we let $M_k := \sup_{x \in X} |\partial h_{t_k}|$, then $\lim_{k \to +\infty} M_k = +\infty$.

We let $x_k \in X$ be a point where $|\partial h_{t_k}|$ attains its maximum $M_k$. Taking a subsequence we may suppose that $\lim_{k \to +\infty} x_k = x_0 \in X$. Let $i_0 \in I$ be such that $x_0 \in U_{i_0}$, and choose an open neighborhood $U$ of $x_0$ contained in $U_{i_0}$. Fix local holomorphic coordinates $z_1, \cdots, z_n$ on $U$.

Choose a local frame $s$ of $E_{i_0}$ over $U$, and represent $K_{t|U}$ by a matrix of smooth functions whose entries are denoted $K_{t,rs}$, and $h_{t,i_0}$ by a matrix $H_{t,i_0}$ of smooth functions whose entries are denoted $h_{t,i_0,pq}$. We let $h_{t,i_0}^{pq}$ be the entries of $H_{t,i_0}^{-1}$, and represent the Kähler metric $g$ by a matrix $G$ on $U$, whose entries are denoted $g_{\alpha\beta}$. We let $g^{\alpha\beta}$ be the entries of $G^{-1}$.

The proof now follows the same lines as that of Lemma 8.22 in Chapter VI of [17], where one has to replace the formula (8.23) (expressing $\widehat{K}_t$ in terms of $h_t$) with the following formula

(1) $\widehat{K}_{t,rs} = \Delta(h_{t,i_0,rs}) + \sum_{\alpha,\beta=1}^{n} \sum_{p,q=1}^{r} g^{\alpha\beta} h_{t,i_0}^{pq} \partial_\alpha h_{t,i_0,rq} \partial_{\overline{\beta}} h_{t,i_0,ps} - i\Lambda_g B_{i_0} h_{t,i_0,rs},$

where
$$\Delta = -\sum_{\alpha,\beta} g^{\alpha\beta} \frac{\partial^2}{\partial z_\alpha \partial \overline{z}_\beta},$$
and
$$\partial_\alpha = \frac{\partial}{\partial z_\alpha}, \quad \partial_{\overline{\beta}} = \frac{\partial}{\partial \overline{z}_\beta}.$$

The proof works since $B_{i_0}$ does not depend on $t$, so $i\Lambda_g B_{i_0}$ is bounded on $[0, a)$. □

We are now in the position to prove the main result of this section:

**Proposition 6.19.** *Let $E$ be an $\alpha-$twisted holomorphic vector bundle over a compact Kähler manifold and let $h_0$ be a Hermitian metric on $E$. Then the evolution equation has a unique smooth solution $\hbar : [0, +\infty) \longrightarrow Herm^+(E)$ such that $\hbar(0) = h_0$.*

*Proof.* Proposition 6.10 gives unicity, and Proposition 6.12 gives a smooth solution $\hbar : [0, a) \longrightarrow Herm^+(E)$ of the evolution equation for some $a \in \mathbb{R}$. Suppose now that $[0, a)$ is the largest interval on which a solution exists.

By point 2 of Proposition 6.10 we know that $h_t$ converges in the $C^0$ topology to a Hermitian metric $h_a$ as $t$ converges to $a$. By point 1 of Corollary 6.16 we know that $\sup_{x \in X} |K_t|(x)$ is bounded on $[0, a)$. We then may apply Lemma 6.18, getting that $\hbar$ is bounded in $C^1(X, a, H_E)$, and that $\mathscr{R}$ is bounded in $L^p(X, a)$ for every $p < +\infty$.

The boundedness of $\mathscr{R}$ implies by Lemma 6.17 that $|R_t|$ is uniformly bounded in $L^\infty$ on $[0, a)$, so Point 2 of Corollary 6.16 implies that for every $k$ the functions $|D_t^k R_t|$ are uniformly bounded for $t \in [0, a)$.

The boundedness of the family $\hbar$ allows us to show that $\hbar$ is bounded in $C^k(X, a, H_E)$ for every $k \in \mathbb{N}$. The proof goes by induction on $k$. Suppose that $\hbar$ is bounded in $C^{k-1}(X, a, H_E)$. Then the families of the first order partial derivatives of $h_t$ are bounded in $C^{k-2}(X, a, H_E)$. Since $\widehat{K}_t$ is uniformly bounded in $C^l(X, a, H_E)$ for all $l \in \mathbb{N}$, and hence in particular in $C^{k-2}(X, a, H_E)$, equation (1) and Elliptic Regularity imply that $\hbar$ is bounded in $C^k(X, a, H_E)$.

It follows that $\hbar$ is bounded in $C^\infty(X, a, H_E)$. But since $h_t$ converges to $h_a$ in the $C^0-$topology, it follows that the convergence is in the $C^\infty-$topology, i. e. we can extend $\hbar : [0, a) \longrightarrow Herm^+(E)$ to a smooth solution over $[0, a]$.

By Proposition 6.12, starting with $h_a$ we may extend $\hbar$ to a unique smooth solution of the evolution equation defined over $[a, a')$ for some $a' > a$, and hence we extend $\hbar$ to a smooth solution on $[0, a')$. But this contradicts the fact that $[0, a)$ is the largest interval over which $\hbar$ exists, concluding the proof. □

6.2.4. *Properties of the solution.* The main property we show in this section is the following:

**Proposition 6.20.** *Let $E$ be an $\alpha-$twisted holomorphic vector bundle and $\hbar : [0, +\infty) \longrightarrow Herm^+(E)$ be a smooth solution of the evolution equation.*

(1) *For every $k \in Herm^+(E)$ the function*

$$L_{g,k,\hbar} : [0, +\infty) \longrightarrow \mathbb{R} \quad L_{g,k,\hbar}(t) := L_{g,k}(h_t)$$

*is monotone decreasing.*

(2) *The function*
$$m_K : [0, +\infty) \longrightarrow \mathbb{R}, \quad m_K(t) := \max_{x \in X} |K_g(E, h_t) - c_g(E)id_E|^2$$
*is monotone decreasing.*

(3) *If there is $A \in \mathbb{R}$ such that $L_{g,k,\hbar} \geq A$, then*
$$\lim_{t \to +\infty} m_K(t) = 0.$$

*Proof.* As seen in section 5.1.3 we have
$$\partial_t L_{g,k,\hbar}(t) = \partial_t L_{g,k}(h_t) = (\widehat{K}_g(E, h_t) - c_g(E)h_t, h_t')_{h_t}.$$

Since $\hbar$ is a solution of the evolution equation, we have
$$h_t' = -(\widehat{K}_g(E, h_t) - c_g(E)h_t),$$
hence
$$\partial_t L_{g,k,\hbar} = -||\widehat{K}_t - c_g(E)h_t||_{h_t}^2 \leq 0,$$
so that $L_{g,k,\hbar}$ is monotone decreasing.

For the second point, by Corollary 6.14 we know that $\partial_t K_t = -\square_{h_t} K_t$. Let us now calculate $\square ||K_t - c_g(E)id_E||^2$. To do so, we first have
$$D_t^{1,0} D_t^{0,1} ||K_t - c_g(E)id_E||^2 = D_t^{1,0} D_t^{0,1} (Tr((K_t - c_g(E)id_E)^2)) =$$
$$= 2Tr((K_t - c_g(E)id_E) \circ (D_t^{1,0} D_t^{0,1} K_t)) + 2Tr(D_t^{1,0} K_t \circ D_t^{0,1} K_t).$$

Then we get
$$\square ||K_t - c_g(E)id_E||^2 = i\Lambda_g D_t^{1,0} D_t^{0,1} ||K_t - c_g(E)id_E||^2 =$$
$$= 2Tr((K_t - c_g(E)id_E) \circ \square_{h_t} K_t) - 2|D_t^{1,0} K_t|^2 =$$
$$= -2Tr((K_t - c_g(E)id_E)\partial_t K_t) - 2|D_t^{1,0} K_t|^2 =$$
$$= -2\partial_t ||K_t - c_g(E)id_E||^2 - 2|D_t^{1,0} K_t|^2.$$

It then follows that
$$(\partial_t + \square)||K_t - c_g(E)id_E||^2 = -2|D_t^{1,0} K_t|^2 \leq 0.$$

By the Maximum Principle for Parabolic Equations (see Lemma 4.1 in Chapter VI of [17]) we then get the statement.

The remaining part of the proof is exactly identical to the proof of point (iii) of Proposition 9.1 in Chapter VI of [17]. □

An immediate corollary of point 3 of Proposition 6.20 is:

**Proposition 6.21.** *Let $E$ be an $\alpha$−twisted holomorphic vector bundle on a compact Kähler manifold $X$ with Kähler metric $g$. If for every $k \in Herm^+(E)$ there is a real number $A_k$ such that $L_{g,k}(h) \geq A_k$ for every $h \in Herm^+(E)$, then $E$ is approximate $g$−Hermite-Einstein.*

6.3. **Regularization of twisted sheaves.** We start by describing a regularization process that was first used by Buchdahl in [2] for holomorphic vector bundles on surfaces, and then by Jacob in [13] for holomorphic vector bundles on compact Kähler manifolds. We adapt it here for holomorphic twisted vector bundles. As we will see, this construction will allows us to prove the boundedness of the Donaldson Lagrangian of a $g-$semistable twisted vector bundle.

6.3.1. *Blow-ups and regularization of subsheaves.* Let $X$ be a compact Kähler manifold with Kähler metric $g$. The starting point of the construction is an exact sequence
$$0 \longrightarrow S \xrightarrow{f} E \xrightarrow{p} Q \longrightarrow 0$$
of $\alpha-$twisted coherent sheaves, where $E$ is locally free of rank $r$ (we will consider it as an $\alpha-$twisted holomorphic vector bundle), $S$ is torsion-free of rank $s$ and $Q$ is torsion-free of rank $q = r - s$.

We let $Z$ be the singular set of $Q$: on $X \setminus Z$ the sheaves $S$, $E$ and $Q$ are all $\alpha-$twisted holomorphic vector bundles, and the morphisms $f$ and $p$ are morphisms of $\alpha-$twisted holomorphic vector bundles. As the locus of $X$ where $Q$ is locally free is where the rank of $f$ is maximal, i. e. equal to $s$, we get
$$Z = \{x \in X \,|\, rk_x(f) \leq s - 1\}.$$
For every $k \in \mathbb{N}$ we then let
$$Z_k := \{x \in X \,|\, rk_x(f) \leq k\}.$$
We notice that $Z_k \subseteq Z_{k+1}$ for every $k \in \mathbb{N}$, that $Z_0 = \emptyset$, that $Z_{s-1} = Z$ and that $Z_k = X$ for $k \geq s$.

Let us now choose $k_0 \in \mathbb{N}$ such that $Z_{k_0}$ is the smallest non-empty set among the $Z_k$'s, and choose $x \in Z_{k_0}$, so that $rk_x(f) = k_0$. Let us write $S = \{S_i, \psi_{ij}\}_{i,j \in I}$ and $E = \{E_i, \phi_{ij}\}_{i,j \in I}$.

If $x \in U_i$, we may choose local frames $s_i$ of $S_i$ and $e_i$ of $E_i$ so that the matrix $F_i$ representing $f_i$ with respect to these local frames is
$$F_i = \begin{bmatrix} I_{k_0} & 0 \\ 0 & G_i \end{bmatrix}$$
where $G_i$ is a $(s - k_0) \times (r - k_0)-$matrix whose entries are holomorphic functions vanishing identically on $Z_{k_0}$.

Suppose that $x \in U_{ij}$, and represent $\phi_{ij}$ (resp. $\psi_{ij}$) by a $r \times r-$matrix $a_{ij}$ (resp. a $s \times s-$matrix $b_{ij}$) with respect to $e_i$ and $e_j$ (resp. $s_i$ and $s_j$). As $f$ is a morphism of twisted sheaves, we have $a_{ij}F_i = F_j b_{ij}$. From this relation and the form of $F_i$ and $F_j$, we get an invertible $(r - k_0) \times (r - k_0)-$matrix $a'_{ij}$ and an invertible $(s - k_0) \times (s - k_0)-$matrix $b'_{ij}$ such that $a'_{ij} G_i = G_j b'_{ij}$.

Let us now consider the blow up $\pi : \widetilde{X} \longrightarrow X$ of $X$ along $Z_{k_0}$ (with reduced structure). We let $\widetilde{U}_i := \pi^{-1}(U_i)$, and $\pi_i := \pi_{|\widetilde{U}_i}$: we then have that $\pi_i : \widetilde{U}_i \longrightarrow U_i$ is the blow-up of $U_i$ along $Z_{k_0} \cap U_i$ (with reduced structure).

Let $U$ be an open subset of $\widetilde{X}$ over which the exceptional divisor of $\pi$ is given by equation $w = 0$.

Over $U \cap \widetilde{U}_i$ we then have that $\pi_i^* f_i$ is represented by the matrix

$$\pi_i^* F_i = \begin{bmatrix} I_{k_0} & 0 \\ 0 & \pi_i^* G_i \end{bmatrix},$$

where $\pi_i^* G_i$ is a $(s - k_0) \times (r - k_0)$−matrix whose entries are holomorphic functions which are all multiples of $w$. We let $m_i$ be the largest power of $w$ one can pull out of $\pi_i^* G_i$. Since $a'_{ij} G_i = G_j b'_{ij}$ and $a'_{ij}$ and $b'_{ij}$ are invertible, we get $m_i = m_j$. We then let $m_{k_0}(f) := m_i$, and call it the **vanishing multiplicity of** $f$ **along** $Z_{k_0}$. We will use the notation $m(f)$ for $m_{k_0}(f)$ if no confusion is possible.

Notice that there is a $(r - k_0) \times (s - k_0)$−matrix $\widetilde{G}_{ij}$ which does not vanish identically on the exceptional divisor of $\pi_i$, such that

$$\pi_i^* F_i = \begin{bmatrix} I_{k_0} & 0 \\ 0 & \widetilde{G}_i \end{bmatrix} \cdot \begin{bmatrix} I_{k_0} & 0 \\ 0 & w^a I_{s-k_0} \end{bmatrix}.$$

We let

$$\widetilde{F}_i := \begin{bmatrix} I_{k_0} & 0 \\ 0 & \widetilde{G}_i \end{bmatrix}, \quad T_i := \begin{bmatrix} I_{k_0} & 0 \\ 0 & w^a I_{s-k_0} \end{bmatrix}$$

so that $\pi_i^* F_i = \widetilde{F}_i T_i$.

The matrix $T_i$ gives rise to a morphism $t_i : S_i \longrightarrow S_i$ of $\mathcal{O}_{U_i}$−modules, and the matrix $\widetilde{F}_i$ gives rise to a morphism $\widetilde{f}_i : S_i \longrightarrow E_i$ of $\mathcal{O}_{U_i}$−modules. Let us define

$$\widetilde{S}_i := t_i(S_i)$$

which is a coherent subsheaf of $\mathcal{O}_{U_i}$−modules of $S_i$.

**Lemma 6.22.** *For every $i, j \in I$ there is an isomorphism $\widetilde{\psi}_{ij} : \widetilde{S}_i \longrightarrow \widetilde{S}_j$ such that $\widetilde{S} = \{\widetilde{S}_i, \widetilde{\phi}_{ij}\}_{i,j \in I}$ is an $\alpha$−twisted coherent sheaf, and such that*

$$\widetilde{f} = \{\widetilde{f}_{i|\widetilde{S}_i}\}_{i \in I} : \widetilde{S} \longrightarrow E_i$$

*is an injective morphism of $\alpha$−twisted coherent sheaves.*

*Proof.* First, we define $\psi_{ii} := id_{\widetilde{S}_i}$. Now, let us consider $i \neq j \in I$. We represent $\psi_{ij}$ on a local frame by a $s \times s$−matrix $b_{ij} = [\beta_{ij,pq}]$. We let $\widetilde{b}_{ij} = [\widetilde{\beta}_{ij,pq}]$ be defined as

$$\widetilde{\beta}_{ij,pq} = \begin{cases} \beta_{ij,pq}, & p, q \leq k_0, \ p, k > k_0 \\ w^{m(f)} \beta_{ij,pq}, & p > k_0, \ q \leq k_0 \\ w^{-m(f)} \beta_{ij,pq}, & p \leq k_0, \ q > k_0 \end{cases}$$

To this matrix we associate a morphism $\widetilde{\psi}_{ij} : \widetilde{S}_i \longrightarrow \widetilde{S}_j$, and we prove that it is an isomorphism whose inverse is $\widetilde{\psi}_{ji}$. To do so, let us calculate $\widetilde{b}_{ij} \cdot \widetilde{b}_{ji}$, i. e. we calculate the entry $\gamma_{p,q}$ in position $(p, q)$ of this product.

- If $p, q \leq k_0$ we have

$$\gamma_{p,q} = \sum_{l=1}^{s} \widetilde{\beta}_{ij,pl}\widetilde{\beta}_{ji,lq} = \sum_{l=1}^{k_0} \beta_{ij,pl}\beta_{ji,lq} + \sum_{l=k_0+1}^{s} w^{-m(f)}\beta_{ij,pl}w^{m(f)}\beta_{ji,lq} =$$

$$= \sum_{l=1}^{s} \beta_{ij,pl}\beta_{ji,lq} = \delta_{pq},$$

where the last equality comes from the fact that the last sum is the entry in position $(p,q)$ of $b_{ij} \cdot b_{ji}$, which is a matrix representing $\phi_{ij} \circ \phi_{ji} = id$.

- If $p > k_0$ and $q \leq k_0$ we have

$$\gamma_{p,q} = \sum_{l=1}^{s} \widetilde{\beta}_{ij,pl}\widetilde{\beta}_{ji,lq} = \sum_{l=1}^{k_0} w^{m(f)}\beta_{ij,pl}\beta_{ji,lq} + \sum_{l=k_0+1}^{s} \beta_{ij,pl}w^{m(f)}\beta_{ji,lq} =$$

$$= w^{m(f)}\sum_{l=1}^{s} \beta_{ij,pl}\beta_{ji,lq} = 0,$$

where the last equality is again as before (since $p \neq q$).

- If $p \leq k_0$ and $q > k_0$ we have

$$\gamma_{p,q} = \sum_{l=1}^{s} \widetilde{\beta}_{ij,pl}\widetilde{\beta}_{ji,lq} = \sum_{l=1}^{k_0} \beta_{ij,pl}w^{-m(f)}\beta_{ji,lq} + \sum_{l=k_0+1}^{s} w^{-m(f)}\beta_{ij,pl}\beta_{ji,lq} =$$

$$= w^{-m(f)}\sum_{l=1}^{s} \beta_{ij,pl}\beta_{ji,lq} = 0.$$

- If $p, q > k_0$ we have

$$\gamma_{p,q} = \sum_{l=1}^{s} \widetilde{\beta}_{ij,pl}\widetilde{\beta}_{ji,lq} = \sum_{l=1}^{k_0} w^{m(f)}\beta_{ij,pl}w^{-m(f)}\beta_{ji,lq} + \sum_{l=k_0+1}^{s} \beta_{ij,pl}\beta_{ji,lq} =$$

$$= \sum_{l=1}^{s} \beta_{ij,pl}\beta_{ji,lq} = \delta_{pq}.$$

In conclusion $\widetilde{b}_{ij} \cdot \widetilde{b}_{ji} = I_s$, so $\widetilde{\psi}_{ij}$ is an isomorphism whose inverse is $\widetilde{\psi}_{ji}$.

We are left to show that $\widetilde{\psi}_{ij} \circ \widetilde{\psi}_{jk} \circ \widetilde{\psi}_{ki} = \alpha_{ijk} \cdot id_{E_i}$. A calculation similar to the previous one shows that $\widetilde{b}_{ij}\widetilde{b}_{jk}\widetilde{b}_{ki} = \alpha_{ijk}I$, which completes the proof of the fact that $\widetilde{S}$ is an $\alpha$−twisted coherent sheaf.

The fact that $\widetilde{f}$ is a morphism of $\alpha$−twisted coherent sheaves is as follows: represent $f_i$ by a matrix $F_i$, $\widetilde{f}_i$ by a matrix $\widetilde{F}_i$ and $t_i$ by a matrix $T_i$. Moreover, represent $\phi_{ij}$ by $a_{ij}$: we then have

$$\pi^*a_{ij} \cdot \widetilde{F}_i \cdot T_i = \pi^*(a_{ij} \cdot F_i) = \pi^*(F_j \cdot b_{ij}) = \widetilde{F}_j \cdot T_j \cdot \pi^*b_{ij}.$$

A calculation similar to the previous one for the product of $\widetilde{b}_{ij}$'s shows that $T_j \cdot \pi^*b_{ij} = \widetilde{b}_{ij}$. This shows that $\pi^*\phi_{ij} \circ \widetilde{f}_i = \widetilde{f}_j \circ \widetilde{\psi}_{ij}$, so that $\widetilde{f}$ is a morphism of $\alpha$−twisted coherent sheaves. As $f_i$ and $t_i$ are both injective, it follows that $\widetilde{f}$ is injective. □

As a consequence, we now have a new exact sequence

$$0 \longrightarrow \widetilde{S} \xrightarrow{\widetilde{f}} \pi^*E \xrightarrow{\widetilde{p}} \widetilde{Q} \longrightarrow 0$$

of $\pi^*\alpha$−twisted coherent sheaves on $\widetilde{X}$, where $\widetilde{S}$ is again torsion-free of rank $s$ and $\widetilde{Q}$ is torsion-free of rank $q$. The matrix $\widetilde{F}_i$ representing $\widetilde{f}_i$ is

$$\widetilde{F}_i = \begin{bmatrix} I_{k_0} & 0 \\ 0 & \widetilde{G}_i \end{bmatrix}$$

hence for every $x \in \widetilde{X}$ we have $rk_x(\widetilde{f}) \geq k_0$. We have two possible cases:
  i) either $rk_x(\widetilde{f}) > k_0$ for every $x \in \pi^{-1}(Z_k)$,
  ii) or if $x \in \pi^{-1}(U_i)$, $\mathfrak{m}_x$ is the maximal ideal of $x$, $p$ is the smallest integer such that $\mathfrak{m}_x^p \subseteq V(\widetilde{G}_i)$ (the ideal generated by the entries of $\widetilde{G}_i$), and $q$ is the smallest integer such that $\mathfrak{m}_{\pi(x)}^q \subseteq V(G_i)$, then $p < q$.

We are now in the position to prove the following:

**Proposition 6.23.** *Let $X$ be a compact manifold, $E$ and $\alpha$−twisted holomorphic vector bundle and $S$ a torsion-free coherent subsheaf of $E$ with torsion-free quotient. Then there exists a finite number of blow-ups*

$$\widetilde{X}_N \xrightarrow{\pi_N} \cdots \xrightarrow{\pi_2} \widetilde{X}_1 \xrightarrow{\pi_1} X,$$

*and for each $1 \leq k \leq N$ there is a $\pi_k^* \circ \cdots \circ \pi_1^*(\alpha)$−twisted torsion-free coherent sheaf $\widetilde{S}_k$ with an injective morphism $f_k : \widetilde{S}_k \longrightarrow \pi_k^* \circ \cdots \circ \pi_1^*(E)$ verifying the two following properties:*

(1) *for every $i \in I$ there is a morphism $t_{k,i} : \pi_k^*S_{k-1} \longrightarrow \pi_k^*S_{k-1}$ such that*
  (a) $\pi_k^*f_{k-1,i} = f_{k,i}t_{k,i}$ *and*
  (b) *for every $x \in X_k$, if the exceptional divisor of $\pi_k$ around $x$ has equation $w = 0$, then $t_{k,i}$ is represented (with respect to a local frame) by a diagonal matrix whose entries are monomials in $w$.*
(2) *The rank of $f_N$ is constant, so $\widetilde{S}_N$ is a $\pi_N^* \circ \cdots \circ \pi_1^*(\alpha)$−twisted holomorphic subbundle of $\pi_N^* \circ \cdots \circ \pi_1^*E$, and the corresponding quotient is a $\pi_N^* \circ \cdots \circ \pi_1^*(\alpha)$−twisted holomorphic bundle.*

*Proof.* The construction provided above shows us that after the blow-up of $Z_{k_0}$ in $X$ we get $\widetilde{f} : \widetilde{S} \longrightarrow \pi^*E$. Now, for every $k \in \mathbb{N}$ let

$$\widetilde{Z}_k := \{x \in \widetilde{X} \mid rk_x(\widetilde{f}) \leq k\},$$

and let $\widetilde{k}_0$ be the smallest $k$ such that $\widetilde{Z}_k \neq \emptyset$. Then i) and ii) give
  i') either $\widetilde{k}_0 > k_0$,
  ii') or $\widetilde{k}_0 = k_0$, $\widetilde{Z}_{k_0} \subseteq \pi^{-1}(Z_{k_0})$, and $m_{k_0}(\widetilde{f}) < m_{k_0}(f)$.

If this last case happens, with a finite number of blow-ups we reduce to the case where $\widetilde{k}_0 > k_0$. We now repeat the construction, and the statement is proved. □

The sequence of blow-ups described in the statement of Proposition 6.23 is called **regularization** of the exact sequence

$$0 \longrightarrow S \xrightarrow{f} E \xrightarrow{p} Q \longrightarrow 0,$$

and each blow-up in the process is a **regularization step**.

6.3.2. *Metrics, curvatures and regularization.* As seen in section 2.7, if

$$0 \longrightarrow S \xrightarrow{f} E \xrightarrow{p} Q \longrightarrow 0$$

is an exact sequence of $\alpha-$twisted holomorphic vector bundles on $X$ and $h$ is a Hermitian metric on $E$, then $h$ induces Hermitian metrics $h^S$ on $S$ and $h^Q$ on $Q$. Moreover, we have a splitting morphism $\varphi : Q \longrightarrow E$ of $\alpha-$twisted $C^\infty$ vector bundles.

If we choose a local frame for $E_i$ and represent $h_i$ by a matrix $H_i$, $f_i$ by a matrix $F_i$ and $\varphi_i$ by a matrix $\Phi_i$, let $H_i^S$ be the matrix representing $h_i^S$ and $H_i^Q$ be the matrix representing $h_i^Q$: we then have

$$H_i^S = {}^T F_i H_i \overline{F}_i, \quad H_i^Q = {}^T \Phi_i H_i \overline{\Phi}_i.$$

Using the same notation of the previous section, let $\pi : \widetilde{X} \longrightarrow X$ be the blow-up of $X$ along $Z_{k_0}$ with reduced structure, and let

$$0 \longrightarrow \widetilde{S} \xrightarrow{\widetilde{f}} \pi^* E \xrightarrow{\widetilde{p}} \widetilde{Q} \longrightarrow 0$$

be the induced exact sequence. Represent $h_i^{\widetilde{S}}$ by a matrix $H_i^{\widetilde{S}} = [h_{i,pq}^{\widetilde{S}}]$, and $h_i^{\widetilde{Q}}$ by a matrix $H_i^{\widetilde{Q}} = [h_{i,pq}^{\widetilde{Q}}]$ with respect to a local frame of $E_i$.

The following is an immediate consequence of Proposition 4 of [13]:

**Lemma 6.24.** *If $w$ is the local equation of the exceptional divisor of $\pi$, we then have*

$$\pi_i^* h_{i,pq}^S = \begin{cases} h_{i,pq}^{\widetilde{S}}, & p,q \leq k_0 \\ w^{-m(f)} h_{i,pq}^{\widetilde{S}}, & p > k_0, \ q \leq k_0 \\ \overline{w^{-m(f)}} h_{i,pq}^{\widetilde{S}}, & p \leq k_0, \ q > k_0 \\ w^{m(f)} \overline{w^{m(f)}} h_{i,pq}^{\widetilde{S}}, & p,q > k_0 \end{cases}$$

*and*

$$\pi_i^* h_{i,pq}^Q = \begin{cases} h_{i,pq}^{\widetilde{Q}}, & p,q \leq k_0 \\ w^{m(f)} h_{i,pq}^{\widetilde{Q}}, & p > k_0, \ q \leq k_0 \\ \overline{w^{m(f)}} h_{i,pq}^{\widetilde{Q}}, & p \leq k_0, \ q > k_0 \\ w^{-m(f)} \overline{w^{-m(f)}} h_{i,pq}^{\widetilde{Q}}, & p,q > k_0 \end{cases}$$

As a consequence of the previous Lemma we get:

**Corollary 6.25.** *Let $R$ be the Chern curvature of $(Q, h^Q)$, and $\widetilde{R}$ the Chern curvature of $(\widetilde{Q}, h^{\widetilde{Q}})$. Let $w$ be a local equation of the exceptional divisor of $\pi$, and let*

$$m_p := \begin{cases} 0, & p \leq k_0 \\ m_{k_0}(f), & p > k_0 \end{cases}$$

*Then we have*

$$\pi^* Tr(R) = \sum_{p=1}^{q} m_p \partial \overline{\partial} \log |w|^2 + Tr(\widetilde{R}).$$

*Proof.* Consider the curvature $R_i$ of $(Q_i, h_i^Q)$. By Lemma 1 of [13] we have

$$\pi_i^* R_i = \sum_{p=1}^{q} m_p \partial \overline{\partial} \log |w|^2 + Tr(\widetilde{R}_i),$$

where $\widetilde{R}_i$ is the curvature of $(\widetilde{Q}_i, h_i^{\widetilde{Q}})$. As the Chern curvature of $(Q, h)$ is obtained locally as $R_i - B_i id_{E_i}$ and the Chern curvature of $(\widetilde{Q}, h^{\widetilde{Q}})$ is obtained locally as $\widetilde{R}_i - \pi_i^* B_i id_{\pi_i^* E_i}$, the statement follows. $\square$

Let now $g$ be a Kähler metric on $X$ and $\sigma_g$ its Kähler form. The metric $\pi^* g$, whose associated $(1,1)$−form is $\pi^* \sigma_g$, is not a Kähler metric on $\widetilde{X}$ since it is degenerate on the exceptional divisor of $\pi$. Anyway, if $E$ is an $\alpha$−twisted holomorphic vector bundle, and $h$ is a Hermitian metric on it, then we may define the degree of $E$ with respect to $\pi^* g$ as follows:

**Definition 6.26.** *The $\pi^* g$−**degree** of $\pi^* E$ is*

$$\deg_{\pi^* g}(\pi^* E) := \int_{\widetilde{X}} \pi^* \gamma_1(E, h) \wedge \pi^* \sigma_g^{n-1}.$$

Again, as $\sigma_g$ is closed and the exceptional divisor of $\pi$ is contracted by $\pi$, the previous definition does not depend on $h$. Using the fact that by definition we have

$$\gamma_1(E, h) = \frac{i}{2\pi} Tr(R)$$

where $R$ is the Chern curvature of $(E, h)$, we get

$$\deg_{\pi^* g}(\pi^* E) = \frac{i}{2\pi} \int_{\widetilde{X}} \pi^* Tr(R) \wedge \pi^* \sigma_g^{n-1}.$$

The same definition makes sense for every torsion-free $\alpha$−twisted coherent sheaf whose Chern curvature is $L^1$ on the locally-free locus. We show that this is the case for $S$ and $Q$ as before. To do so, recall that given

$$0 \longrightarrow S \xrightarrow{f} E \xrightarrow{p} Q \longrightarrow 0$$

as before, on $X \setminus Z$ there is an element $C \in A^{0,1}(\underline{Hom}(Q, S))$, the second fundamental form of the induced metric. We then have $C \wedge C^* \in A^{1,1}(\underline{End}(Q))$, so that we let

$$\gamma_g(Q) := i\Lambda_g(C \wedge C^*) \in End(Q).$$

Notice that $Tr(\gamma_g(Q))$ is a smooth function on $X \setminus Z$, and we have:

**Proposition 6.27.** *The second fundamental form $C$ is in $L^2$, and*
$$\int_{X \setminus Z} Tr(\gamma_g(Q))\sigma_g^n$$
*is finite. Moreover, the Chern curvatures of $(S, h^S)$ and $(Q, h^Q)$ are in $L^1$.*

*Proof.* The proof of the first part is identical to that of Proposition 1 of [13] (where one replaces Lemma 1 by Corollary 6.25). The second part is a consequence of the first part, of the Gauss-Codazzi equations for the curvatures on the locally free part, and of the fact that the Chern curvature of an $\alpha-$twisted locally free sheaf is smooth. $\square$

As a consequence, the definition of degree of $E$ with respect to $g$ (resp. of degree of $\pi^*E$ with respect to $\pi^*g$) extends to $S$ and $Q$ (resp. to $\widetilde{S}$ and $\widetilde{Q}$). The following is then an immediate consequence of this definition and of Corollary 6.25 (see Lemma 2 of [13]).

**Proposition 6.28.** *Let*
$$0 \longrightarrow S \xrightarrow{f} E \xrightarrow{p} Q \longrightarrow 0$$
*be an exact sequence of $\alpha-$twisted coherent sheaves, where $E$ is locally free of rank $r$ and $S$ and $Q$ are torsion-free of respective rank $s$ and $q$. Let $\pi : \widetilde{X} \longrightarrow X$ be a regularization, and let*
$$0 \longrightarrow \widetilde{S} \xrightarrow{\widetilde{f}} \pi^*E \xrightarrow{\widetilde{p}} \widetilde{Q} \longrightarrow 0$$
*be the regularized sequence. If $g$ be a Kähler metric on $X$, then*
$$\deg_{\pi^*g}(\widetilde{S}) = \deg_g(S), \quad \deg_{\pi^*g}(\widetilde{Q}) = \deg_g(Q).$$

6.4. **Regularization and Donaldson's Lagrangian.** We now define the Donaldson Lagrangian for the subsheaves $S$ and $Q$ in the exact sequence

(2) $$0 \longrightarrow S \xrightarrow{f} E \xrightarrow{p} Q \longrightarrow 0$$

where $E$ is an $\alpha-$twisted locally free coherent sheaf and $S$ and $Q$ are torsion-free.

Recall that if $h, k \in Herm^+(E)$ and $\hbar \in \Omega^{0,1}_{h,k}(E)$, then we defined
$$Q_1(h, k) := \log(\det(f^{k,h})), \quad Q_2^\hbar(h, k) := i\int_0^1 Tr(f^{h_t, h'_t} \cdot R_t)dt,$$

where $h_t = \hbar(t)$ and $R_t$ is the Chern curvature of $(E, h_t)$. For a Kähler metric $g$ on $X$ whose associated Kähler form is $\sigma_g$, we defined
$$L_g(h, k) = \int_X \left(Q_2^\hbar(h, k) - \frac{c_g(E)}{n}Q_1(h, k)\sigma_g\right) \wedge \sigma_g^{n-1},$$

which is independent of the chosen path $\hbar$. We recall that

$$c_g(E) = \frac{2\pi n \mu_g(E)}{\int_X \sigma_g^n}.$$

6.4.1. *Definition of the Lagrangian for torsion-free subsheaves.* Let now $\pi : \widetilde{X} \longrightarrow X$ be a blow-up. For a $\pi^*\alpha$–twisted holomorphic vector bundle $\widetilde{E}$ on $\widetilde{X}$, one may define the Donaldson functional using $\pi^*g$: using the same notation as above, but for $h, k \in Herm^+(\pi^*E)$, we let

$$L_{\pi^*g}(h,k) := \int_X \left( Q_2^\hbar(h,k) - \frac{c_g(E)}{n} Q_1(h,k)\pi^*\sigma_g \right) \wedge \pi^*\sigma_g^{n-1},$$

which is again independent of the chosen path $\hbar$ since $\pi^*\sigma_g$ is closed.

If $\pi : \widetilde{X} \longrightarrow X$ is a regularization of the exact sequence (2), let the regularized exact sequence be

$$(3) \qquad 0 \longrightarrow \widetilde{S} \xrightarrow{\widetilde{f}} \pi^*E \xrightarrow{\widetilde{p}} \widetilde{Q} \longrightarrow 0$$

where $\widetilde{S}$ and $\widetilde{Q}$ are both locally free.

Let $h, k \in Herm^+(E)$, so that by Lemma 2.36 we have that $\pi^*h, \pi^*k \in Herm^+(\pi^*E)$. The exact sequence (3) implies by section 2.7 that $\pi^*h$ and $\pi^*k$ induce Hermitian metrics $h^{\widetilde{S}} := (\pi^*h)^{\widetilde{S}}$ and $k^{\widetilde{S}} := (\pi^*k)^{\widetilde{S}}$ on $\widetilde{S}$, and $h^{\widetilde{Q}} := (\pi^*h)^{\widetilde{Q}}$ and $k^{\widetilde{Q}} := (\pi^*k)^{\widetilde{Q}}$ on $\widetilde{Q}$. We then let

$$L_g^S : Herm^+(E) \times Herm^+(E) \longrightarrow \mathbb{R}, \quad L_g^S(h,k) := L_{\pi^*g}^{\widetilde{S}}(h^{\widetilde{S}}, k^{\widetilde{S}}),$$

$$L_g^Q : Herm^+(E) \times Herm^+(E) \longrightarrow \mathbb{R}, \quad L_g^Q(h,k) := L_{\pi^*g}^{\widetilde{Q}}(h^{\widetilde{Q}}, k^{\widetilde{Q}}).$$

**Remark 6.29.** We used the notation $L_{\pi^*g}^{\widetilde{S}}$ instead of $L_{\pi^*g}$ for the Donaldson Lagrangian of $\widetilde{S}$ in order to avoid confusion between $\widetilde{S}$ and $\widetilde{Q}$. Moreover, we notice that $L_g^S$ and $L_g^Q$ are both defined on Hermitian metrics on $E$.

**Remark 6.30.** Even if all the Donaldson Lagrangian involved in the definition of $L_g^S$ and $L_g^Q$ do not depend on the chosen path, in order to calculate them we need to fix one. A natural choice would be $\hbar \in \Omega_{h,k}^{0,1}(E)$: if we let

$$\pi^*\hbar : [0,1] \longrightarrow Herm^+(\pi^*E), \quad \pi^*\hbar(t) := \pi^*h_t,$$

then $\pi^*\hbar \in \Omega_{\pi^*h,\pi^*k}^{0,1}(\pi^*E)$, so by section 2.7 we then let

$$\hbar^{\widetilde{S}} : [0,1] \longrightarrow Herm^+(\widetilde{S}), \quad \hbar^{\widetilde{S}}(t) = h_t^{\widetilde{S}} := (\pi^*h_t)^{\widetilde{S}}$$

and

$$\hbar^{\widetilde{Q}} : [0,1] \longrightarrow Herm^+(\widetilde{Q}) \quad \hbar^{\widetilde{Q}}(t) = h_t^{\widetilde{Q}} = (\pi^*h_t)^{\widetilde{Q}},$$

which are the curves of Hermitian metrics induced by $\pi^*\hbar$ on $\widetilde{S}$ and $\widetilde{Q}$. We then may use $\hbar^{\widetilde{S}}$ in order to calculate $L_g^S(h,k)$ and $\hbar^{\widetilde{Q}}$ in order to calculate $L_g^Q(h,k)$. Anyway, we may use any path $\hbar \in \Omega_{h^{\widetilde{S}},k^{\widetilde{S}}}^{0,1}(\widetilde{S})$ to calculate $L_g^S(h,k)$, and any path $\hbar \in \Omega_{h^{\widetilde{Q}},k^{\widetilde{Q}}}^{0,1}(\widetilde{Q})$ to calculate $L_g^Q(h,k)$.

Notice that the previous definition may however depend on the choice of the regularization, which is not unique. The following tells us that it is not the case:

**Proposition 6.31.** *If $\pi : \widetilde{X} \longrightarrow X$ and $\pi' : \widetilde{X}' \longrightarrow X$ are two blow-ups producing regularizations of $S$ and $Q$, then $L^S_{\pi^*g} = L^S_{(\pi')^*g}$ and $L^Q_{\pi^*g} = L^Q_{(\pi')^*g}$.*

*Proof.* The proof is identical to that of Proposition 5 of [13]. $\square$

This allows us to simplify the notation, and let $L^S_g = L^S_{\pi^*g}$ and $L^Q_g = L_{\pi^*g}$.

6.4.2. *Relations between the Lagrangians.* We now describe the relation between $L_g$, $L^S_g$ and $L^Q_g$. We will use the following notation: for every $h \in Herm^+(E)$ we let $C_h \in A^{0,1}(\underline{Hom}(Q, S))$ be the second fundamental form of the metric $h^Q$ induced by $h$ on $Q$.

If $\pi : \widetilde{X} \longrightarrow X$ is a regularization of the exact sequence (2), and if $C_{\pi^*h} \in A^{0,1}(\underline{Hom}(\widetilde{Q}, \widetilde{S}))$ is the second fundamental form of the metric $h^{\widetilde{Q}}$ induced by $\pi^*h$ on $\widetilde{Q}$, then the same proof of Proposition 1 of [13] shows that

$$\int_{\widetilde{X}} Tr(C_{\pi^*h} \wedge C^*_{\pi^*h}) \wedge \pi^*\sigma_g^{n-1}$$

is a real number which does not depend on the chosen regularization.

We let

$$||C_h||^2_{L^2} := \int_{\widetilde{X}} Tr(C_{\pi^*h} \wedge C^*_{\pi^*h}) \wedge \pi^*\sigma_g^{n-1},$$

and call it the $L^2$−norm of $C_h$ (we know that $C_h$ is in $L^2$ by Proposition 6.27). Using this notation we state now the following, which is the relation between the various Donaldson Lagrangians.

**Proposition 6.32.** *In the exact sequence (2) suppose that $\mu_g(S) = \mu_g(E)$. Then for every $h, k \in Herm^+(E)$ we have*

$$L_g(h, k) = L^S_g(h, k) + L^Q_g(h, k) + ||C_h||^2_{L^2} - ||C_k||^2_{L^2}.$$

*Proof.* Let us first suppose that $S$ and $Q$ are both locally free. We let

$$Q^S_1(h, k) := \log(\det(f^{k^S, h^S})), \quad Q^Q_1(h, k) := \log(\det(f^{k^Q, h^Q})).$$

For every $i \in I$, choose a local frame of $S_i$ and a local frame of $Q_i$, whose images in $E_i$ under $f_i$ and $\varphi_i$ give a local frame of $E_i$, where $\varphi = \{\varphi_i\}$ is the splitting morphism $\varphi : Q \longrightarrow E$.

With respect to these local frames we represent $h_i$ by a matrix $H_i$, $k_i$ by a matrix $K_i$, $f_i$ by a matrix $F_i$, $\varphi_i$ by a matrix $P_i$, $h^S_i$ by a matrix $H^S_i$, $k^S_i$ by a matrix $K^S_i$, $h^Q_i$ by a matrix $H^Q_i$ and $k^Q_i$ by a matrix $K^Q_i$.

Over the open subset where the local frame is defined we then have

$$Q_1(h, k) = \log(\det(K_i^{-1} H_i)), \quad Q^S_1(h, k) = \log(\det((K^S_i)^{-1} H^S_i)),$$
$$Q^Q_1(h, k) = \log(\det((K^Q_i)^{-1} H^Q_i)).$$

By the definition of $f^{k_i,h_i}$, $f^{k_i^S,h_i^S}$ and $f^{k_i^Q,h_i^Q}$ we have that

$$f^{k_i,h_i} \circ f_i = f_i \circ f^{k_i^S,h_i^S}, \qquad f^{k_i,h_i} \circ \varphi_i = \varphi_i \circ f^{k_i^Q,h_i^Q},$$

hence

$$K_i^{-1} H_i = \begin{bmatrix} (K_i^S)^{-1} H_i^S & 0 \\ 0 & (K_i^Q)^{-1} H_i^Q \end{bmatrix}.$$

This implies that

$$\det(K_i^{-1} H_i) = \det((K_i^S)^{-1} H_i^S) \det((K_i^Q)^{-1} H_i^Q),$$

and hence $Q_1(h,k) = Q_1^S(h,k) + Q_1^Q(h,k)$.

We are now interested in $Q_2$. To do so, let $\pi: E \longrightarrow S$ be the splitting morphism, and remark that if $h$ is a Hermitian metric on $E$, we have that $C_h = \pi \circ \overline{\partial} \circ \varphi$ (see section 2.4.2). Let now $\hbar \in \Omega_{h,k}^{0,1}(E)$ and $h_t := \hbar(t)$. Notice that for every $t \in [0,1]$ the Hermitian metric $h_t$ induces a different splitting of the exact sequence (2), namely

$$0 \longrightarrow Q \xrightarrow{\varphi_t} E \xrightarrow{\pi_t} S \longrightarrow 0,$$

and we have $C_{h_t} = \pi_t \circ \overline{\partial} \circ \varphi_t$.

For every $t \in [0,1]$ we let $s_t : Q \longrightarrow S$ be such that $f \circ s_t = \varphi_t - \varphi_0$ (or equivalently $s_t \circ p = \pi_0 - \pi_t$). The same argument in the proof of Proposition 10.2 in Chapter VI of [17] gives

$$\partial_t C_{h_t} = \overline{\partial}(\partial_t s_t)$$

and

$$f^{h_t,h_t'} = \begin{bmatrix} f^{h_t^S,(h^S)_t'} & -\partial_t s_t \\ -\partial_t s_t^* & f^{h_t^Q,(h^Q)_t'} \end{bmatrix}.$$

The Gauss-Codazzi equations (see section 2.7.2) give

$$R_t = \begin{bmatrix} R_t^S - C_{h_t} \wedge C_{h_t}^* & -D_t^{1,0} C_{h_t} \\ -D_t^{0,1} C_{h_t}^* & R_t^Q - C_{h_t}^* \wedge C_{h_t} \end{bmatrix},$$

where $R_t$, $R_t^S$ and $R_t^Q$ are the Chern curvatures of $(E, h_t)$, $(S, h_t^S)$ and $(Q, h_t^Q)$ respectively.

As a consequence we get

$$Tr(f^{h_t,h_t'} \cdot R_t) - Tr(f^{h_t^S,(h^S)_t'} \cdot R_t^S) - Tr(f^{h_t^Q,(h^Q)_t} \cdot R_t^Q) =$$
$$= -Tr(f^{h_t^S,(h^S)_t'} \cdot C_{h_t} \wedge C_{h_t}^*) - Tr(f^{h_t^Q,(h^Q)_t'} \cdot C_{h_t}^* \wedge C_{h_t}) +$$
$$+ Tr(\partial_t s_t \circ D_t^{0,1} C_{h_t}^*) - Tr(\partial_t s_t^* \circ D_t^{1,0} C_{h_t}).$$

Since

$$Tr(\partial_t s_t \circ D_t^{0,1} C_{h_t}^*) = \overline{\partial}(Tr(\partial_t s_t \circ C_{h_t}^*)) - Tr(\partial_t C_{h_t} \wedge C_{h_t}^*),$$
$$Tr(\partial_t s_t^* \circ D_t^{1,0} C_{h_t}) = \partial(Tr(\partial_t s_t^* \circ C_{h_t})) - Tr(\partial_t C_{h_t}^* \wedge C_{h_t}),$$

we finally get

$$Tr(f^{h_t,h_t'} \cdot R_t) - Tr(f^{h_t^S,(h^S)_t'} \cdot R_t^S) - Tr(f^{h_t^Q,(h^Q)_t'} \cdot R_t^Q) =$$

$$= -Tr(f^{h_t^S,(h^S)'_t} \cdot C_{h_t} \wedge C_{h_t}^*) - Tr(f^{h_t^Q,(h^Q)'_t} \cdot C_{h_t}^* \wedge C_{h_t}) +$$
$$-Tr(\partial_t C_{h_t} \wedge C_{h_t}^*) + Tr(\partial_t C_{h_t}^* \wedge C_{h_t}) + \delta_t,$$

where

$$\delta_t = \overline{\partial}(Tr(\partial_t s_t \circ C_{h_t}^*)) - \partial(Tr(\partial_t s_t^* \circ C_{h_t})) \in \partial A^{0,1}(X) + \overline{\partial} A^{1,0}(X).$$

The same proof of Proposition 10.2 in Chapter VI of [17] gives that

$$Tr(f^{h_t,h'_t} \cdot R_t) - Tr(f^{h_t^S,(h^S)'_t} \cdot R_t^S) - Tr(f^{h_t^Q,(h^Q)'_t} \cdot R_t^Q) = -\partial_t(Tr(C_{h_t} \wedge C_{h_t}^*)) + \delta,$$

hence

$$Tr(f^{h_t,h'_t} \cdot R_t) = Tr(f^{h_t^S,(h^S)'_t} \cdot R_t^S) + Tr(f^{h_t^Q,(h^Q)'_t} \cdot R_t^Q) - \partial_t(Tr(C_{h_t} \wedge C_{h_t}^*)) - \delta.$$

Finally, we have

$$Q_2^{\hbar}(h,k) = i \int_0^1 Tr(f^{h_t,h'_t} \cdot R_t) \wedge \sigma_g^{n-1} =$$

$$= i \int_0^1 Tr(f^{h_t^S,(h^S)'_t} \cdot R_t^S) \wedge \sigma_g^{n-1} + i \int_0^1 Tr(f^{h_t^Q,(h^Q)'_t} \cdot R_t^Q) \wedge \sigma_g^{n-1} +$$
$$-iTr(C_h \wedge C_h^*) + iTr(C_k \wedge C_k^*) + \delta,$$

where

$$\delta = \int_0^1 \delta \wedge \sigma_g^{n-1} \in \partial A^{0,1}(X) + \overline{\partial} A^{1,0}(X).$$

But now notice that

$$Q_2^{\hbar,S}(h,k) = i \int_0^1 Tr(f^{h_t^S,(h^S)'_t} \cdot R_t^S) \wedge \sigma_g^{n-1},$$

$$Q_2^{\hbar,Q}(h,k) = i \int_0^1 Tr(f^{h_t^Q,(h^Q)'_t} \cdot R_t^Q) \wedge \sigma_g^{n-1},$$

hence

$$Q_2^{\hbar}(h,k) = Q_2^{\hbar,S}(h,k) + Q_2^{\hbar,Q}(h,k) - iTr(C_h \wedge C_h^*) + iTr(C_k \wedge C_k^*) + \delta.$$

Now, by hypothesis we have $\mu_g(E) = \mu_g(S)$, and hence $\mu_g(Q) = \mu_g(E)$. It follows that $c_g(E) = c_g(E) = c_g(Q)$, so that the previous relations prove the statement in the case $S$ and $Q$ are both locally free.

If $S$ and $Q$ are not locally free, we may regularize them by taking a blow-up $\pi : \widetilde{X} \longrightarrow X$, getting an exact sequence

$$0 \longrightarrow \widetilde{S} \xrightarrow{\widetilde{f}} \pi^* E \xrightarrow{\widetilde{p}} \widetilde{Q} \longrightarrow 0$$

of $\alpha$−twisted locally free sheaves. As $\mu_g(E) = \mu_g(S) = \mu_g(Q)$, by Proposition 6.28 we have $\mu_{\pi^*g}(\pi^* E) = \mu_{\pi^*g}(\widetilde{S}) = \mu_{\pi^*g}(\widetilde{Q})$. By the previous part of the proof we then have

$$L_{\pi^*g}(h,k) = L^{\widetilde{S}}_{\pi^*g}(h,k) + L^{\widetilde{Q}}_{\pi^*g}(h,k) + ||C_{\pi^*k}||^2_{L^2} - ||C_{\pi^*h}||^2_{L^2}.$$

But by Proposition 6.31 we have

$$L^{\widetilde{S}}_{\pi^*g}(h,k) = L^S_g(h,k), \quad L^{\widetilde{Q}}_{\pi^*g}(h,k) = L^Q_g(h,k),$$

and as in the proof of Proposition 1 of [13] we have
$$||C_{\pi^*k}||^2_{L^2} = ||C_k||^2_{L^2}, \quad ||C_{\pi^*h}||^2_{L^2} = ||C_h||^2_{L^2},$$
and we are done. $\square$

6.5. **A lower bound.** In this section we prove the following result:

**Proposition 6.33.** *Let $X$ be a compact Kähler manifold with Kähler metric $g$, and consider the exact sequence (2), where we suppose that $\mu_g(S) = \mu_g(E)$ and that $S$ has minimal rank among the torsion-free $\alpha$−twisted coherent subsheaves of $E$ with this property. If $\pi : \widetilde{X} \longrightarrow X$ is regularization of the exact sequence (2) and $\widetilde{S}$ the twisted sheaf induced by $S$, then there is $M \in \mathbb{R}$ such that $L^{\widetilde{S}}_{\pi^*g}(h,k) \geq M$ for every $h, k \in Herm^+(\widetilde{S})$.*

The proof will be as follows: we first assume that $g$ is a Kähler metric on $X$ with volume 1, and suppose that $\pi$ consists of a single blow-up. If we take $h, k \in Herm^+(E)$ and let $h^{\widetilde{S}}$ and $k^{\widetilde{S}}$ be the induced Hermitian metrics on $\widetilde{S}$, by definition we have
$$L^S_g(h,k) = L^{\widetilde{S}}_{\pi^*g}(h^{\widetilde{S}}, k^{\widetilde{S}}).$$

By Remark 6.30 we may assume that $L^{\widetilde{S}}_{\pi^*g}(\widetilde{h}^S, \widetilde{k}^S)$ is calculated using a path $\hbar^{\widetilde{S}} \in \Omega^{0,1}_{h^{\widetilde{S}}, k^{\widetilde{S}}}(\widetilde{S})$. We will let $\widetilde{h}_t := \hbar^{\widetilde{S}}(t)$ for every $t \in [0,1]$: notice that this Hermitian metric is not necessarily induced by a Hermitian metric on $\pi^*E$.

Let $\widetilde{R}_t$ be the Chern curvature of $(\widetilde{S}, \widetilde{h}_t)$, and let $\widetilde{K}^S_t = i\Lambda_{\pi^*g}\widetilde{R}^S_t$: we remark that the metric that we are using to define the mean curvature is not a Kähler metric, since it blows-up along the exceptional divisor of $\pi$.

Now, the evolution equation has the form
$$f^{\widetilde{h}_t, \widetilde{h}'_t} = -(\widetilde{K}_t - c_{\pi^*g}(\widetilde{S})id_{\widetilde{S}}).$$

By Proposition 6.28 we know that $c_g(S) = c_{\pi^*g}(\widetilde{S})$, hence the evolution equation is

(4) $$f^{\widetilde{h}_t, \widetilde{h}'_t} = -(\widetilde{K}_t - c_g(S)id_{\widetilde{S}}).$$

Suppose now that $\hbar^{\widetilde{S}}$ is a solution of the evolution equation, and let
$$L(t) : [0,1] \longrightarrow \mathbb{R}, \quad L(t) := L^{\widetilde{S}}_{\pi^*g}(\widetilde{h}_t, k^{\widetilde{S}}).$$
Then
$$\partial_t L = \frac{1}{n!}\int_{\widetilde{X}} Tr((\widetilde{K}_t - c_g(S)id_{\widetilde{S}}) \circ f^{\widetilde{h}_t, \widetilde{h}'_t})\pi^*\sigma^n_g = -||\widetilde{K}_t - c_g(S)id_{\widetilde{S}}||^2 \leq 0.$$

But then $L$ is decreasing along any solution of the evolution equation, so that $L^{\widetilde{S}}_{\pi^*g}$ is bounded from below along a solution of the evolution equation for every initial Hermitian metric on $\widetilde{S}$ coming from a Hermitian metric on $E$, and then it is bounded from below in general.

The previous strategy works if we are able to show that the evolution equation
$$f^{\widetilde{h}_t,\widetilde{h}'_t} = -(\widetilde{K}_t - c_g(S)id_{\widetilde{S}})$$
has a solution starting from any Hermitian metric on $\widetilde{S}$ coming from a Hermitian metric on $E$. Even if $\widetilde{S}$ is locally free, this cannot be concluded from 6.19 since $\widetilde{K}_t$ is defined using $\pi^*g$, which is not a Kähler metric.

6.5.1. *Existence of a solution of the evolution equation.* We now prove that the evolution equation has a solution for all times:

**Proposition 6.34.** *The evolution equation (4) has a smooth solution*
$$\hbar : [0,+\infty) \longrightarrow Herm^+(\widetilde{S}).$$

*Proof.* Let us first suppose that the regularization $\pi : \widetilde{X} \longrightarrow X$ of the exact sequence (2) is a single blow-up. Fix a Fubini-Study metric $\sigma$ on the exceptional divisor of $\pi$, and take $\epsilon > 0$ a small real number such that $g_\epsilon := \pi^*g + \epsilon\sigma$ is a Kähler metric on $\widetilde{X}$.

For every Hermitian metrics $\widetilde{h}, \widetilde{k} \in Herm^+(\widetilde{S})$ we let
$$L_\epsilon(\widetilde{h}, \widetilde{k}) := L_{g_\epsilon}^{\widetilde{S}}(\widetilde{h}, \widetilde{k}),$$
which is a Donaldson Lagrangian of a $\pi^*\alpha$-twisted holomorphic vector bundle with respect to a Kähler metric. By Proposition 6.19, the evolution equation
$$(5) \qquad f^{\widetilde{h}_t,\widetilde{h}'_t} = -(K_{g_\epsilon}(\widetilde{S},h_t) - c_{g_\epsilon}(S)id_{\widetilde{S}})$$
has then a unique smooth solution
$$\widetilde{\hbar}_\epsilon : [0,+\infty) \longrightarrow Herm^+(\widetilde{S})$$
for a given initial Hermitian metric. We show that there is a sequence $\epsilon_m$ converging to 0 such that $\widetilde{\hbar}_{\epsilon_m}$ converges to a solution of the evolution equation (4).

To do so, let $\widetilde{K}_\epsilon$ the mean curvature of $(\widetilde{S}, \widetilde{h})$ with respect to $g_\epsilon$, i. e. $\widetilde{K}_\epsilon = i\Lambda_{g_\epsilon}\widetilde{R}$ where $\widetilde{R}$ is the Chern curvature of $(\widetilde{S}, \widetilde{h})$. By definition we have
$$\int_{\widetilde{X}} Tr(-\widetilde{K}_\epsilon + c_{g_\epsilon}(\widetilde{S})id_{\widetilde{S}})\pi^*\sigma_{g_\epsilon}^n = 0,$$
hence it follows that the equation
$$\Delta f = Tr(-i\widetilde{K}_\epsilon + c_{g_\epsilon}(\widetilde{S})id_{\widetilde{S}})$$
has a smooth solution, denoted $\varphi_\epsilon$.

We now let $\widetilde{h}_\epsilon := e^{\varphi_\epsilon}\widetilde{h}$, which is a Hermitian metric on $\widetilde{S}$ (since it is a conformal change of $\widetilde{h}$. The evolution equation (5) has then a unique smooth solution
$$\widetilde{\hbar}_\epsilon : [0,+\infty) \longrightarrow Herm^+(\widetilde{S})$$

such that $\widetilde{\hbar}_\epsilon(0) = \widetilde{h}_\epsilon$. We let $\widetilde{R}_{\epsilon,t}$ be the Chern curvature of $(\widetilde{S}, \widetilde{h}_{\epsilon,t})$, where $\widetilde{h}_{\epsilon,t} = \widetilde{\hbar}_\epsilon(t)$, and $\widetilde{K}_{\epsilon,t}$ will be its mean curvature. We now make the following:

**Claim**: *for every $0 < t_1 \leq t_2 < +\infty$ there is a constant $N \in \mathbb{R}$ such that $|\widetilde{R}_{\epsilon,t}| \leq N$ for every $\epsilon$ and every $t \in [t_1, t_2]$.*

Let us first show that if this claim holds, then we are done. Indeed, as in Corollary 6.16 the claim implies that for every $k$ there is then a constant $N_k \in \mathbb{R}$ such that $|\widetilde{D}_{\epsilon,t}^k \widetilde{R}_{\epsilon,t}| \leq N_k$, where $\widetilde{D}_{\epsilon,t}$ is the Chern connection of $(\widetilde{S}, \widetilde{h}_{\epsilon,t})$.

This gives a a uniform bound for the $C^k$−norm of $\widetilde{R}_{\epsilon,t}$. By compactness we then see that there is a subsequence $\epsilon_m$ converging to 0 such that the solution $\widetilde{h}_{\epsilon_m,t}$ converges to a Hermitian metric $\widetilde{h}_{0,t}$ for every $t \in [t_1, t_2]$, and hence we get a solution

$$\widetilde{\hbar} : [t_1, t_2] \longrightarrow Herm^+(\widetilde{S}),$$

of the evolution equation (4) defined on $[t_1, t_2]$.

The subsequence $\epsilon_m$ depends on $t_1$ and $t_2$, so this does not yet give the desired solution. Anyway, if we take a sequence $t_n$ going to $+\infty$, for each $n$ we will find a subsequence $\epsilon_{n,m}$ of $\epsilon_{n-1,m}$ converging to 0 and such that $\widetilde{h}_{\epsilon_{n,m},t}$ converges to $h_{0,t}$ for every $t \in [t_1, t_n]$. Letting $t_n$ going to $+\infty$ we get a solution

$$\widetilde{\hbar} : [t_1, +\infty) \longrightarrow Herm^+(\widetilde{S})$$

of the evolution equation (4). Up to rescaling, we are then done.

In conclusion, we are left to prove the claim, i. e. that for every $0 < t_1 \leq t_2 < +\infty$ there is a constant $N \in \mathbb{R}$ such that $|\widetilde{R}_{\epsilon,t}| \leq N$ for every $\epsilon$ and every $t \in [t_1, t_2]$. This requires some steps.

*Step 1*: there is a uniform bound for the $||\widetilde{K}_{\epsilon,t}||_{L^1}$ (i. e. independent of $\epsilon > 0$ and $t \geq 0$). To prove this, recall that

$$(\partial_t + \Box_{\widetilde{h}_{\epsilon,t}})|\widetilde{K}_{\epsilon,t}|^2 \leq 0$$

by Lemma 6.15. The same argument of the proof of Proposition 6 of [13] shows that

$$(\partial_t - \Box_{\widetilde{h}_{\epsilon,t}})|\widetilde{K}_{\epsilon,t}| \leq 0,$$

so

$$\partial_t ||\widetilde{K}_{\epsilon,t}||_{L^1} = \partial_t \int_{\widetilde{X}} |\widetilde{K}_{\epsilon,t}| \sigma_{g_\epsilon}^n \leq 0.$$

It follows that the function mapping $t$ to $||\widetilde{K}_{\epsilon,t}||_{L^1}$ is decreasing, hence one just needs to prove a uniform bound for $\widetilde{K}_{\epsilon,0}$ (i. e. independent of $\epsilon > 0$).

Notice that

$$\widetilde{K}_{\epsilon,0} = K_{g_\epsilon}(\widetilde{S}, \widetilde{h}_\epsilon) = K_{g_\epsilon}(\widetilde{S}, \widetilde{h}) + \Delta_\epsilon \varphi_\epsilon \cdot id_{\widetilde{S}} = \widetilde{K}_\epsilon + \Delta_\epsilon \varphi_\epsilon \cdot id_{\widetilde{S}},$$

so
$$||\widetilde{K}_{\epsilon,0}||_{L^1} = \int_{\widetilde{X}} |\widetilde{K}_{\epsilon,0}|\sigma_{g_\epsilon}^n \leq \int_{\widetilde{X}} |\Delta_\epsilon \varphi_\epsilon|\sigma_{g_\epsilon}^n + ||\widetilde{K}_\epsilon||_{L^1} =$$
$$= \int_{\widetilde{X}} |Tr(-\widetilde{K}_\epsilon + c_{g_\epsilon}(\widetilde{S})id_{\widetilde{S}})\sigma_{g_\epsilon}^n + ||\widetilde{K}_\epsilon||_{L^1} \leq 2||\widetilde{K}_\epsilon||_{L^1} + \Gamma,$$

for some constant number $\Gamma$.

The same proof of Proposition 6 of [13] (replacing Proposition 2 with Proposition 6.27) shows that $||\widetilde{K}_\epsilon||_{L^1}$ is uniformly bounded, so we get a uniform bound on $||\widetilde{K}_{\epsilon,t}||_{L^1}$.

*Step 2*: the same proof of Proposition 7 of [13] shows that the existence of a uniform bound for $||\widetilde{K}_{\epsilon,t}||_{L^1}$ implies the existence of a uniform bound for $||\widetilde{K}_{\epsilon,t}||_{L^\infty}$. As in Proposition 8 of [13], this implies a uniform bound on $Tr(f^{\widetilde{h}_{\epsilon,t},\widetilde{h}'_{\epsilon,t}})$ (i. e. independent of $\epsilon > 0$ and $t \in [t_1, t_2]$ for every $0 < t_1 \leq t_2 < +\infty$). As in Lemma 6.18 this proves the uniform bound of the Chern curvatures we are looking for.

This concludes the proof of the statement in the case of a single blow-up. Suppose now that $\pi = \pi_k \circ \pi_{k-1} \circ \cdots \circ \pi_1$, where $\pi_j : X_j \longrightarrow X_{j-1}$ is a single blow-up (here we let $X_0 := X$, so that $X_k = \widetilde{X}$). For every $j$ we let $\sigma_j$ be a Fubini-Study metric on the exceptional divisor of $\pi_j \circ \cdots \circ \pi_1$, and for every choice of $\epsilon_1, \cdots, \epsilon_k > 0$ we define a Kähler metric on $X_j$ in a recursive way as $g_j := \pi_j^* g_{j-1} + \epsilon_j \sigma_j$, where we let $g_0 = g$.

If we let $\epsilon_k$ go to 0, by the previous part of the proof we find a smooth solution of the evolution equation for the metric $\pi_k^* g_{k-1}$, defined on some interval $[t_{k-1}, +\infty)$. Repeating this process we will then find a smooth solution of the evolution equation with respect to $\pi^* g$, defined over some interval $[t_0, +\infty)$. Up to rescaling, we are done. □

6.5.2. *Bounds and subsheaves.* Before giving the proof of Proposition 6.33 we need to recall some results of [25] and [26]. The first result is the following:

**Lemma 6.35.** *Let $E$ be an $\alpha-$twisted holomorphic vector bundle on a compact Kähler manifold $X$, and let $h, k \in Herm^+(E)$. Then*
$$\Delta(\log(Tr(f^{h,k}))) \leq 2(|K_g(E,h)| + |K_g(E,k)|).$$

*Proof.* The proof is essentially the same of that of point (d) of Lemma 3.1 in [25]. See even the proof of equation (1.9.2) of [26]. □

The second result, which is the content of Propositions 2.1 and 2.2 of [25], is the following:

**Lemma 6.36.** *Suppose that $X$ is either a compact Kähler manifold with a Kähler metric $g$, or a Zariski open subset of a compact Kähler manifold with a metric $g$ which is the restriction of a smooth metric. Let $n$ be the dimension of $X$. Then the following properties hold.*

(1) *$X$ has finite volume with respect to $g$.*

(2) $X$ has a positive exhaustion function $\phi$ such that $\Delta \phi$ is bounded.
(3) There is an increasing function $a: \mathbb{R}_+ \longrightarrow \mathbb{R}_+$ such that
   (a) $a(0) = 0$ and $a(x) = x$ if $x > 1$;
   (b) if $f : X \longrightarrow \mathbb{R}$ is bounded, positive and $\Delta(f) \leq M$ for some real number $M$, then there is a constant $C(M)$ depending only on $M$ such that
$$\sup_{x \in X} |f|(x) \leq C(M) \cdot a\left( \int_X |f| \sigma_g^n \right);$$
   (c) if $f$ is as in the previous point and $M = 0$, then $\Delta(f) = 0$.

A third useful ingredient is a different formula for the Donaldson Lagrangian that can be found in [25] and [26]. Let us take an $\alpha-$twisted holomorphic vector bundle $E$ on $X$, and suppose that $h_0, h \in Herm^+(E)$. We let $f^{h_0,h} \in End_{h_0}^+(E)$ be the endomorphism associated to $h_0$ and $h$, and we will write $s^{h_0,h} := \log(f^{h_0,h})$. We know by Lemma 2.56 that $s^{h_0,h} \in End_{h_0}(E)$, and we notice that $f^{h_0,h} = \exp(s^{h_0,h})$, and that we have
$$h = \widehat{f^{h_0,h}}_{h_0} = \widehat{\exp(s^{h_0,h})}_{h_0}.$$

Now, for $t \in [0,1]$ we have that $t \cdot s^{h_0,h} \in End_{h_0}(E)$, so we have
$$\mathcal{s} : [0,1] \longrightarrow End_{h_0}(E), \quad \mathcal{s}(t) = s_t := t \cdot s^{h_0,h}.$$
We then see that $\exp(t \cdot s^{h_0,h}) \in End_{h_0}^+(E)$, hence we have
$$\mathcal{f} : [0,1] \longrightarrow End_k^+(E), \quad \mathcal{f}(t) = f_t := \exp(t \cdot s^{h_0,h}).$$
Finally, we let
$$h_t := \widehat{\exp(t \cdot s^{h_0,h})}_{h_0} = \widehat{f_t}_{h_0},$$
which is a Hermitian metric on $E$, and we have that
$$\hbar : [0,1] \longrightarrow Herm^+(E), \quad \hbar(t) := h_t$$
is a smooth family of Hermitian metrics such that $\hbar(0) = h_0$ and $h_1 = h$. Notice that $f_t = f^{h_0,h_t}$ (see Example 2.54).

Now, recall from Example 2.52 that $f^{h_0,h}$ is both $h_0-$Hermitian and $h-$Hermitian. Moreover, we have
$$\partial_t \mathcal{f} = s^{h_0,h} \circ \mathcal{f}(t) = s^{h_0,h} \circ f^{h_0,h_t},$$
hence
$$f^{h_t,h_t'} = f_t' \circ (f_t)^{-1} = \mathcal{f}'(t) \circ \mathcal{f}(t)^{-1} = \partial_t \mathcal{f}(t) \circ \mathcal{f}(t)^{-1} = s^{h_0,h}.$$

Recall that by Lemma 6.8 that
$$\partial_t L_g(h_0, h_t) = \frac{1}{n!} \int_X Tr((K_g(E, h_t) - c_g(E) \cdot id_E) \circ f^{h_t,h_t'}) \sigma_g^n,$$
hence it follows that
$$\partial_t L_g(h_0, h_t) = \frac{1}{n!} \int_X Tr((K_g(E, h_t) - c_g(E) \cdot id_E) \circ s^{h_0,h}) \sigma_g^n.$$

From this we get

$$\frac{d^2}{dt^2}L_g(h_0,h_t) = \frac{1}{n!}\int_X Tr(\partial_t K_g(E,h_t)\circ s^{h_0,h})\sigma_g^n =$$

$$= \frac{1}{n!}\int_X Tr(\partial_t(i\Lambda_g R_t)\circ s^{h_0,h})\sigma_g^n = \frac{1}{n!}\int_X iTr(\partial_t R_t\circ s^{h_0,s})\sigma_g^n.$$

Now, we know that

$$\partial_t R_t = \overline{\partial}(D_h^{1,0}(f'(t)\circ f(t)^{-1})) = \overline{\partial}(D_h^{1,0}(s^{h_0,h})),$$

hence we get

$$\frac{d^2}{dt^2}L_g(h_0,h) = \int_X iTr(\overline{\partial}(D_h^{1,0}(s^{h_0,h}))\circ s^{h_0,s})\sigma_g^n.$$

An easy calculation (see section (5.4) in [26]) shows that

$$\int_X iTr(\overline{\partial}(D_h^{1,0}(s^{h_0,h}))\circ s^{h_0,h})\sigma_g^n = n!\int_X (\overline{\partial}s^{h_0,h},\overline{\partial}s^{h_0,h})_h.$$

Letting $S_i^{h_0,h} = [s_{a,b,i}]$ be the matrix representing $s_i^{h_0,h}$ with respect to a local frame, and if $\lambda_1,\cdots,\lambda_r$ are the eigenvalues of $s^{h_0,h}$ (and hence of $s_i^{h_0,h}$), we then get (see again section (5.4) of [26]) that

$$L_g(h_0,h) = i\int_X Tr(K_g(E,h_0)\circ s^{h_0,h})\sigma_g^n + \int_X \sum_{a,b=1}^r |\overline{\partial}s_{a,b,i}|^2 \frac{e^{\lambda_b-\lambda_a}-(\lambda_b-\lambda_a)-1}{(\lambda_b-\lambda_a)^2}.$$

We are now in the position to prove the following:

**Lemma 6.37.** *Let $E$ be an $\alpha$-twisted holomorphic vector bundle on a compact Kähler manifold $X$ with Kähler metric $g$, and consider $h_0 \in Herm^+(E)$. Suppose that $E$ is $g$-stable and that there is $M \in \mathbb{R}$ such that $\sup|K_g(E,h_0)|(x) \leq M$. Then there are constants $C_1, C_2 \in \mathbb{R}$ such that for every $s \in End_{h_0}(E)$ such that $Tr(s) = 0$, $\sup_{x\in X}|s|(x) < +\infty$ and $\sup_{x\in X}|K_g(S,h)|(x) \leq M$ (for $h = \widehat{\exp(s)}_{h_0}$), we have*

$$\sup_{x\in X}|s|(x) \leq C_1 + C_2 L_g(h,h_0).$$

*Proof.* The proof follows the same strategy as that of Proposition 5.3 of [25]. First, using the fact that

$$\sup|K_g(E,h_0)|(x), \sup_{x\in X}|K_g(S,h)|(x) \leq M,$$

by Lemmas 6.35 and 6.36 it follows that

$$\sup_{x\in X}|s|(x) \leq C_1 + C_2||s||_{L^1}$$

(see section (5.4) in [26], in particular equation (5.4.2)).

Suppose now that the statement is false. We can then find $s \in End_{h_0}(E)$ as in the statement, which contradicts the inequality we want to prove, and such that $||s||_{L^1}$ is arbitrarily large. This means that there is a sequence of real numbers $C_m$ going to $+\infty$, and a sequence of $s_m \in End_{h_0}(E)$ such that $Tr(s_m) = 0$, $||s_m||_{L^1}$ going to infinity, and for which

$$||s_m||_{L^1} \geq C_m \cdot L_g(h_m, h_0),$$

where we let

$$h_m := \widehat{\exp(s_m)}_{h_0}.$$

We now let

$$u_m := \frac{s_m}{||s_m||_{L^1}},$$

so that $||u_m||_{L^1} = 1$ and $\sup_{x \in X} |u_m|(x)$ is bounded. The same proof of Lemmas 5.4 and 5.5 of [25] shows that the $u_m$'s converge weakly in $L_1^2$ to a non-trivial $u_\infty \in L_1^2 End(E)$ which verifies the following properties:

(1) If $\Phi : \mathbb{R} \times \mathbb{R} \longrightarrow \mathbb{R}_+$ is a smooth function such that $\Phi(x,y) \leq \frac{1}{x-y}$ for every $x > y$, and if $U_{\infty,i} = [u_{a,b,i}]$ is the matrix representing $u_{\infty,i}$ with respect to a local frame of $E_i$, we have

$$i \int_X Tr(K_g(E, h_0) \circ u_\infty)\sigma_g^n + \int_X \sum_{a,b=1}^r |\overline{\partial} u_{a,b,i}|^2 \Phi(\lambda_a, \lambda_b) \leq 0,$$

where $\lambda_1, \cdots, \lambda_r$ are the eigenvalues of $u_\infty$.

(2) The eigenvalues $\lambda_1, \cdots, \lambda_r$ of $u_\infty$ are constant and not all equal.

Now, the same argument in the proof of Lemma 5.3 of [25] provides a weakly holomorphic subbundle $S$ of $E$ such that $\mu_g(S) \geq \mu_g(E)$ (see the proof of Lemma 5.7 of [25]), getting a contradiction. $\square$

6.5.3. *Proof of Proposition 6.33.* We are now ready to present a proof of Proposition 6.33

*Proof.* Consider a solution $\widetilde{\hbar} : [0, +\infty) \longrightarrow Herm^+(\widetilde{S})$ of the evolution equation, whose existence is proved in Proposition 6.34. We let $\widetilde{h}_0 := \widetilde{\hbar}(0)$ and $\widetilde{h}_t := \widetilde{\hbar}(t)$ for every $t > 0$. We let $\widetilde{f}_t := f^{\widetilde{h}_0, \widetilde{h}_t}$ and $\widetilde{s}_t := \log(\widetilde{f}_t)$. We may moreover suppose that $\det(\widetilde{f}_t) = 1$, so that $Tr(\widetilde{s}_t) = 0$.

As we saw in the previous section, we have

$$L_{\pi^*g, \widetilde{h}_0}^{\widetilde{S}}(\widetilde{h}_t) = \int_{\widetilde{X}} iTr(K_{\pi^*g}(\widetilde{S}, \widetilde{h}_0) \circ \widetilde{s}_t)\sigma_{\pi^*g}^n + \int_{\widetilde{X}} \sum_{a,b=1}^r |\overline{\partial}\widetilde{s}_{a,b,i}|^2 \frac{e^{\lambda_b - \lambda_a} - (\lambda_b - \lambda_a) - 1}{(\lambda_b - \lambda_a)^2}.$$

Consider now the push-forward $\pi_*\widetilde{S}$: it is an $\alpha$-twisted coherent sheaf on $X$ which is locally free outside the closed subset $Z$ we blow-up to obtain $\pi : \widetilde{X} \longrightarrow X$. As $\pi^{-1}(Z)$ has measure zero and $X \setminus Z$ is isomorphic, via $\pi$, to $\widetilde{X} \setminus \pi^{-1}(Z)$, the same argument presented in section 6.5.2 shows that

$$L_{\pi^*g, \widetilde{h}_0}^{\widetilde{S}}(\widetilde{h}_t) = \int_{X \setminus Z} iTr(K_g(\widetilde{S}, \widetilde{h}_0) \circ \widetilde{s}_t)\sigma_g^n + \int_{X \setminus Z} \sum_{a,b=1}^r |\overline{\partial}\widetilde{s}_{a,b,i}|^2 \frac{e^{\lambda_b - \lambda_a} - (\lambda_b - \lambda_a) - 1}{(\lambda_b - \lambda_a)^2}.$$

We now prove that there are two constants $C_1, C_2 > 0$ such that

$$||\widetilde{s}_t||_{L^1} \leq C_1 + C_2 L^{\widetilde{S}}_{\pi^*g, \widetilde{h}_0}(\widetilde{h}_t),$$

at least for $t$ varying in a sequence going to $+\infty$.

Indeed, suppose that there are not such constants. Notice that $Tr(\widetilde{s}_t) = 0$ and, by the proof of Proposition 6.34, we have that $\sup_{x \in X} |K_g(\widetilde{S}, \widetilde{h}_t)|(x)$ is uniformly bounded. We then may apply the same proof of Lemma 6.37 to get an $\alpha-$twisted coherent subsheaf (over $X \setminus Z$) $\mathscr{F}$ of $\pi_*\widetilde{S}$ such that $0 < rk(\mathscr{F}) < rk(\pi_*\widetilde{S})$.

As $\widetilde{S}$ is a $\pi^*\alpha-$twisted subbundle of $\pi^*E$, it follows that $\pi_*\widetilde{S}$ is an $\alpha-$twisted coherent subsheaf of $E$ over $X \setminus Z$. It follows that over $X \setminus Z$ the sheaf $\mathscr{F}$ is an $\alpha-$twisted coherent subsheaf of $E$. One may then apply the same construction of section 7 of [27] (locally and then gluing, like in the proof of Lemma 5.12) to prove that $\mathscr{F}$ extends to an $\alpha-$twisted coherent subsheaf of $E$ with $0 < rk(\mathscr{F}) < rk(E)$, and such that $\mu_g(\mathscr{F}) \geq \mu_g(E)$.

As $E$ is $g-$semistable we then get $\mu_g(\mathscr{F}) = \mu_g(E)$. But since $rk(\mathscr{F}) < rk(\pi_*\widetilde{S}) = rk(S)$, we get a contradiction by the minimality of $S$. In conclusion, it follows that there are two constants $C_1, C_2 > 0$ such that

$$||\widetilde{s}_t||_{L^1} \leq C_1 + C_2 L^{\widetilde{S}}_{\pi^*g, \widetilde{h}_0}(\widetilde{h}_t),$$

at least for $t$ varying in a sequence going to $+\infty$, getting the statement.

But this implies that $L^{\widetilde{S}}_{\pi^*g, \widetilde{h}_0} \geq -\frac{C_1}{C_2}$, at least for $t$ varying in a sequence going to $\infty$. Hence $L^{\widetilde{S}}_{\pi^*g, \widetilde{h}_0}$ is bounded from below along any solution of the evolution equation, and hence we conclude the proof. $\square$

Another result that we will use in the proof of the approximate Kobayashi-Hitchin correspondence is the following:

**Proposition 6.38.** *Let $X$ be a compact Kähler manifold with Kähler metric $g$ and $\pi : \widetilde{X} \longrightarrow X$ a blow-up. Let $E$ be a $\pi^*g-$semistable $\pi^*\alpha-$twisted holomorphic vector bundle on $\widetilde{X}$, and let $S$ be a $\pi^*\alpha-$twisted coherent subsheaf of maximal rank among the $\pi^*\alpha-$twisted coherent subsheaves of $\pi^*E$ such that $\mu_{\pi^*g}(\pi^*E) = \mu_{\pi^*g}(S)$. Then there is a constant $C$ such that for every $\widetilde{h}, \widetilde{k} \in Herm^+(\pi^*E)$ we have $L^S_{\pi^*g}(\widetilde{h}, \widetilde{k}) \geq C$.*

*Proof.* We construct a regularization $\pi_1 : \widetilde{X}_1 \longrightarrow \widetilde{X}$ of the exact sequence

$$0 \longrightarrow S \longrightarrow \pi^*E \longrightarrow Q \longrightarrow 0.$$

Using the same proof of Proposition 6.33 we provide a proper $\alpha-$twisted torsion-free coherent subsheaf $\mathscr{F}$ of $\pi_*\pi_{1*}\widetilde{S}$ such that $\mu_g(\mathscr{F}) \geq \mu_g(\pi_*\pi_{1*}\widetilde{S})$.

We let $Z$ be the closed subset of $X$ we blow-up to get $\widetilde{X}$. Recall that since $\widetilde{S}$ is a $\pi_1^*\pi^*\alpha-$twisted holomorphic subbundle of $\pi_1^*E$ on $X \setminus Z$, the $\alpha-$twisted sheaf $\mathscr{F}$ is an $\alpha-$twisted subsheaf of $\pi_*E$.

This induces a map $\mathscr{F} \longrightarrow \pi_*E$ of $\alpha$–twisted sheaves, and composing with the natural map $\pi^*\pi_*E \longrightarrow E$ we then get a morphism $j : \pi^*\mathscr{F} \longrightarrow E$ of $\pi^*\alpha$–twisted sheaves, which is injective over $\pi^{-1}(Z)$. Hence $\ker(j)$ is a $\pi^*\alpha$–twisted sheaf whose support is contained in $\pi^{-1}(Z)$.

Consider now the exact sequence
$$0 \longrightarrow \ker(j) \longrightarrow \pi^*\mathscr{F} \longrightarrow \pi^*\mathscr{F}/\ker(j) \longrightarrow 0,$$
and recall that by Lemma 4.10 and its proof, we have
$$\deg_{\pi^*g}(\ker(j)) = \int_V \pi^*\sigma_g^n = 0,$$
since $V$ is the closure of the locus where a holomorphic section of $\ker(j)$ vanishes, which is contained in $Z$. It follows that
$$\deg_{\pi^*g}(\pi^*\mathscr{F}) = \deg_{\pi^*g}(\pi^*\mathscr{F}/\ker(j)),$$
and as these two sheaves have the same rank we get
$$\mu_{\pi^*g}(\pi^*\mathscr{F}) = \mu_{\pi^*g}(\pi^*\mathscr{F}/\ker(j)).$$

But now notice that $\pi^*\mathscr{F}/\ker(j)$ is a $\pi^*\alpha$–twisted coherent subsheaf of $E$, and as $E$ is $\pi^*g$–semistable we get finally that $\mu_{\pi^*g}(\pi_g^*\mathscr{F}/\ker(j)) = \mu_{\pi^*g}(E)$. Since the rank of $\pi^*\mathscr{F}/\ker(j)$ equals the rank of $\pi^*\mathscr{F}$, and this is smaller than the rank of $S$, we get a contradiction with respect to the minimality of $S$. We may now apply the same strategy as in the proof of Proposition 6.33 to get the lower boundedness of $L_{\pi^*g}^S$. $\square$

6.6. **Conclusion of the proof.** We are now in the position to prove Theorem 6.1, namely that an $\alpha$–twisted holomorphic vector bundle $E$ on a compact Kähler manifold with Kähler metric $g$ is $g$–semistable if and only if it is approximate $g$–Hermite-Einstein.

By Theorem 4.23 we know that if $E$ is approximate $g$–Hermite-Einstein, then it is $g$–semistable. We need to prove the converse, and by Proposition 6.21 it is sufficient to prove that the Donaldson Lagrangian $L_{g,k}$ is bounded below for every $k \in Herm^+(E)$. We will need the following:

**Lemma 6.39.** *Let $X$ be a compact Kähler manifold with Kähler metric $g$, $\pi : \widetilde{X} \longrightarrow X$ a blow-up.*
  (1) *If $E$ is a $g$–semistable $\alpha$–twisted holomorphic vector bundle on $X$, then $\pi^*E$ is a $\pi^*g$–semistable $\pi^*\alpha$–twisted holomorphic vector bundle on $\widetilde{X}$.*
  (2) *If we moreover have that $\widetilde{Q}$ is a torsion-free quotient of $\pi^*E$ such that $\mu_{\pi^*g}(\pi^*E) = \mu_{\pi^*g}(\widetilde{Q})$, then $\widetilde{Q}$ is $\pi^*g$–semistable.*

*Proof.* As degree and slope are defined even with respect to $\pi^*g$, the notion of $\pi^*g$–semistability is available as well, and has the same definition as the one with respect to a Kähler metric.

Suppose first that $\pi^*E$ is not $\pi^*g$–semistable. There is then a proper $\pi^*\alpha$–twisted coherent subsheaf $\mathscr{F}$ of $\pi^*E$ such that $0 < rk(\mathscr{F}) < rk(\pi^*E)$ and $\mu_{\pi^*g}(\mathscr{F}) > \mu_{\pi^*g}(\pi^*E)$. As $\pi$ is an isomorphism between $X \setminus Z$ and $\widetilde{X} \setminus \pi^{-1}(Z)$, it follows that $\mu_g(\pi_*\mathscr{F}) > \mu_g(E)$.

On $X \setminus Z$ we have that $\pi_*\mathscr{F}$ is a proper $\alpha$−twisted coherent subsheaf of $E$. Since the codimension of $Z$ is at least two, as already explained (see the proof of Proposition 6.38) we then have that $\pi_*\mathscr{F}$ is a proper $\alpha$−twisted coherent subsheaf of $E$ on $X$, getting a contradiction, and hence showing the first point of the statement.

Suppose now that $\widetilde{Q}$ is not $\pi^*g$−semistable, and let $\mathscr{G}$ be a torsion-free $\pi^*\alpha$−twisted coherent subsheaf of $\widetilde{Q}$ with $0 < rk(\mathscr{G}) < rk(\widetilde{Q})$ and such that $\mu_{\pi^*g}(\mathscr{G}) > \mu_{\pi^*g}(\widetilde{Q})$. Then
$$\mu_{\pi^*g}(\widetilde{Q}/\mathscr{G}) < \mu_{\pi^*g}(\widetilde{Q}) = \mu_{\pi^*g}(\pi^*E).$$
Consider now the morphism
$$f : \pi^*E \longrightarrow \widetilde{Q} \longrightarrow \widetilde{Q}/\mathscr{G}.$$
Then $\ker(f)$ is a torsion-free $\pi^*\alpha$−twisted coherent subsheaf of $\pi^*E$ such that $\mu_{\pi^*g}(\ker(f)) > \mu_{\pi^*g}(\pi^*E)$.

Notice that if $rk(\ker(f)) = 0$, then $rk(\pi^*E) = rk(\widetilde{Q}) - rk(\mathscr{G}) < rk(\widetilde{Q})$, which is impossible since $\widetilde{Q}$ is a quotient of $\pi^*E$. Similarly, if $rk(\ker(f)) = rk(\pi^*E)$, then we get that $rk(\widetilde{Q}) = rk(\mathscr{G})$, which is not possible. It follows that $\ker(f)$ is a torsion-free $\pi^*\alpha$−twisted coherent subsheaf of $\pi^*E$ such that $0 < rk(\ker(f)) < rk(\pi^*E)$ and $\mu_{\pi^*g}(\ker(f)) > \mu_{\pi^*g}(\pi^*E)$. As $\pi^*E$ is $\pi^*g$−semistable, this is impossible, and we are done. □

We are now in the position to prove Theorem 6.1:

*Proof.* By Theorem 4.23 we just need to show that if $E$ is $g$−semistable, then it is approximate $g$−Hermite-Einstein. By Proposition 6.21 we just need to prove that the Donaldson Lagrangian of $E$ with respect to $g$ is bounded below, i. e. we need to prove that for every $k \in Herm^+(E)$ there is a constant $M_k \in \mathbb{R}$ such that $L_g(h, k) \geq M_k$ for every $h \in Herm^+(E)$. Fix then $k \in Herm^+(E)$.

If $E$ is $g$−stable, by Theorem 5.1 we know that $E$ is $g$−Hermite-Einstein, and hence approximate $g$−Hermite-Einstein by Proposition 3.25, and we are done. Suppose then that $E$ is not $g$−stable, and let $S$ be a torsion-free $\alpha$−twisted coherent subsheaf of $E$ such that $\mu_g(S) = \mu_g(E)$ and such that $Q = E/S$ is torsion-free. Suppose moreover that $S$ has minimal rank among all the subsheaves with these properties.

Consider the exact sequence

(6) $$0 \longrightarrow S \longrightarrow E \longrightarrow Q \longrightarrow 0.$$

By minimality of $S$, we see that $S$ is $g$−stable, hence by Proposition 6.21 the Donaldson Lagrangian $L_g^S$ is bounded below, i. e. for every Hermitian metric $k$ on $E$ there is a constant $B_k \in \mathbb{R}$ such that $L_{g,k}^S(h) \geq B_k$ for every $h \in Herm^+(E)$.

By Proposition 6.32 we have
$$L_g(h,k) = L_g^S(h,k) + L_g^Q(h,k) + ||C_h||_{L^2}^2 - ||C_k||_{L^2}^2.$$

Notice that $||C_h||_{L^2}^2$ is positive and $||C_k||_{L^2}^2$ is fixed once $k$ is fixed. We then have
$$L_g(h,k) \geq B_k + ||C_k||_{L^2}^2 + L_g^Q(h,k),$$

so we just need to prove that there is a constant $N_k \in \mathbb{R}$ such that $L_g^Q(h,k) \geq N_k$ for every $h \in Herm^+(E)$.

Take now a regularization $\pi : \widetilde{X} \longrightarrow X$ of the exact sequence (6), so that we get an exact sequence
$$0 \longrightarrow \widetilde{S} \longrightarrow \pi^*E \longrightarrow \widetilde{Q} \longrightarrow 0.$$

We know that $L_g^Q(h,k) = L_{\pi^*g}^{\widetilde{Q}}(\widetilde{h}^Q, \widetilde{k}^Q)$.

As $E$ is $g$-semistable we know by Lemma 6.39 that the pull-back $\pi^*E$ is $\pi^*g$-semistable. Moreover, since $\mu_g(S) = \mu_g(E)$, we get that
$$\mu_{\pi^*g}(\widetilde{Q}) = \mu_g(Q) = \mu_g(S) = \mu_{\pi^*g}(\pi^*E).$$

By Lemma 6.39 we conclude $\widetilde{Q}$ is $\pi^*g$-semistable. If it is $\pi^*g$-stable, we are done.

Suppose then that $\widetilde{Q}$ is not $\pi^*g$-semistable. Let $S_1$ be a torsion-free $\pi^*\alpha$-twisted coherent subsheaf of $\widetilde{Q}$ such that $0 < rk(S_1) < rk(\widetilde{Q})$, $\mu_{\pi^*g}(S_1) = \mu_{\pi^*g}(\widetilde{Q})$, and the quotient $\widetilde{Q}/S_1$ is torsion-free. Choose $S_1$ to have minimal rank among all the subsheaves with these properties. We then see that $S_1$ is $\pi^*g$-stable, and hence from Proposition 6.38 we know that there is a constant $P_{\widetilde{k}}$ such that $L_{\pi^*g}^{S_1}(h,\widetilde{k}) \geq P_{\widetilde{k}}$ for every $h \in Herm^+(\widetilde{Q})$.

Again we have
$$L_{\pi^*g}^{\widetilde{Q}}(h,\widetilde{k}) = L_{\pi^*g}^{S_1}(h,\widetilde{k}) + L_{\pi^*g}^{Q_1}(h,\widetilde{k}) + ||C_h||_{L^2}^2 - ||C_{\widetilde{k}}||_{L^2}^2,$$

where $Q_1 = \widetilde{Q}/S_1$. As before, the problem is reduced to prove that there is a constant $W_{\widetilde{k}}$ such that $L_{\pi^*g}^{Q_1}(h,\widetilde{k}) \geq W_{\widetilde{k}}$ for every $h \in Herm^+(\widetilde{Q})$.

To do so we blow up again $\pi_1 : X_1 \longrightarrow \widetilde{X}$ in order to provide a regularization of $\widetilde{Q}_1$. By Lemma 6.39 we know that $\widetilde{Q}_1$ is $\pi_1^*\pi^*g$-semistable, and we have $rk(\widetilde{Q}_1) < rk(\widetilde{Q})$. After a finite number of steps of this type we then have to stop, concluding the proof. $\square$

### 6.7. Corollaries.
We conclude the present paper with some easy corollaries of the (approximate) Kobayashi-Hitchin correspondence which may prove to be useful for some applications. The first two corollaries are immediate applications of Theorem 1.1 and Propositions 3.12, 3.29, 3.10 and 3.27

**Corollary 6.40.** *Let $E$ be an $\alpha$-twisted holomorphic vector bundle.*
   (1) *If $E$ is $g$-stable, then $\wedge^p E$ and $Sym^p E$ are $g$-polystable.*
   (2) *If $E$ is $g$-semistable, then $\wedge^p E$ and $Sym^p E$ are $g$-semistable.*

**Corollary 6.41.** *Let $E$ be an $\alpha-$twisted holomorphic vector bundle and $F$ a $\beta-$twisted holomorphic vector bundle.*
  (1) *If $E$ and $F$ are $g-$stable, then $E \otimes F$ are $g-$polystable.*
  (2) *If $E$ and $F$ are $g-$semistable, then $E \otimes F$ are $g-$semistable.*

The following is known as Bogomolov's inequality.

**Corollary 6.42.** *Let $E$ be a $g-$semistable $\alpha-$twisted holomorphic vector bundle of rank $r$. Then*
$$\int_X ((r-1)c_1(E)^2 - 2rc_2(E)) \wedge \sigma_g^{n-2} \leq 0.$$

*Proof.* The proof is identical to that of Theorems 4.7 and 5.7 in Chapter IV of [17]. Se even Theorem 2.2.3 and Corollary 2.2.4 of [21]. □

Arvid Perego, Dipartimento di Matematica, Università di Genova, via Dodecaneso 35, 16146 Genova (Italy). E-mail: perego@dima.unige.it